\documentclass{article}
\usepackage[english]{babel}
\usepackage{amsthm,amsfonts,amsbsy, amssymb,amsmath,graphicx}
\usepackage{mathrsfs}
\usepackage{graphics}
\usepackage{afterpage}
\usepackage[all]{xy}
\usepackage{xcolor}
\usepackage[hyphens]{url}
\usepackage{hyperref}
\usepackage{bm}

\DeclareMathAlphabet{\mathdutchcal}{U}{dutchcal}{m}{n}

\newtheorem{theorem}{Theorem}
\newtheorem*{theorem*}{Theorem}
\newtheorem{lemma}{Lemma}
\newtheorem{corollary}{Corollary}

\newtheorem{proposition}{Proposition}

\theoremstyle{definition}
\newtheorem{definition}{Definition}

\theoremstyle{remark}
\newtheorem{remark}{Remark}
\newtheorem{example}{Example}

\def\Z{{\mathbb Z}}
\def\N{{\mathbb N}}

\def\R{{\mathbb R}}

\newcommand{\difd}{\mathrm{Diff}_\partial^0}

\DeclareMathOperator{\uast}{\underline{\ast}}
\DeclareMathOperator{\oast}{\overline{\ast}}

\newcommand{\sgn}{\mathrm{sgn}}
\newcommand{\mor}{\mathrm{Mor}}

\date{}

\title{The crossing and the arc from the topological viewpoint}
\author{Igor Nikonov\footnote{nikonov at mech.math.msu.su}}

\begin{document}

\maketitle

\begin{abstract}
The combinatorial approach to knot theory treats knots as diagrams modulo Reidemeister moves. Many constructions of knot invariants (e.g., index polynomials, quandle colorings, etc.) use elements of diagrams such as arcs and crossings by assigning invariant labels to them. The universal invariant labels, which carry the most information, can be thought of as equivalence classes of arcs and crossings modulo the relation, which identifies corresponding elements of diagrams connected by a Reidemeister move. We can call these equivalence classes the arcs and crossings of the knot. In the paper, we give a topological description of sets of these classes as the isotopy classes of probes of diagram elements.

In the second part of the paper, we discuss homotopy classes of diagram elements. We demonstrate that the sets of these classes are fundamental for algebraic objects that are responsible for coloring diagrams of tangles on a given surface. For arcs, these algebraic objects are quandles; for regions, they are partial ternary quasigroups; for semiarcs, they are biquandloids; and for crossings, they are crossoids. The definitions of the last three algebraic structures are given in the paper.

Additionally, we introduce the multicrossing complex of a tangle and define the crossing homology class. In a sense, the multicrossing complex unifies tribracket, biquandle and crossoid homologies; and the tribracket, biquandle and crossoid cycle invariants are actually the result of pairing a tribracket (biquangle, crossoid) cocycle with the crossing homology class. 
\end{abstract}


\section{Introduction}

Knot theory can be approached from two different sides. The topological approach defines a knot as an embedding of the circle in $\R^3$. On this way one gets decomposition of knots into prime ones~\cite{Schubert}; Thurston's trichotomy of torus, satellite and hyperbolic knots~\cite{Thurston}; and knot invariants like genus~\cite{Seifert}, knot group~\cite{Dehn}, hyperbolic volume and Alexander polynomials~\cite{Alexander}.

\begin{figure}[h]
\centering
  \includegraphics[width=0.15\textwidth]{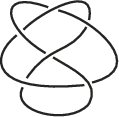}
  \caption{A knot diagram}\label{pic:knot_diagram}
\end{figure}

\begin{figure}[h]
\centering
  \includegraphics[width=0.5\textwidth]{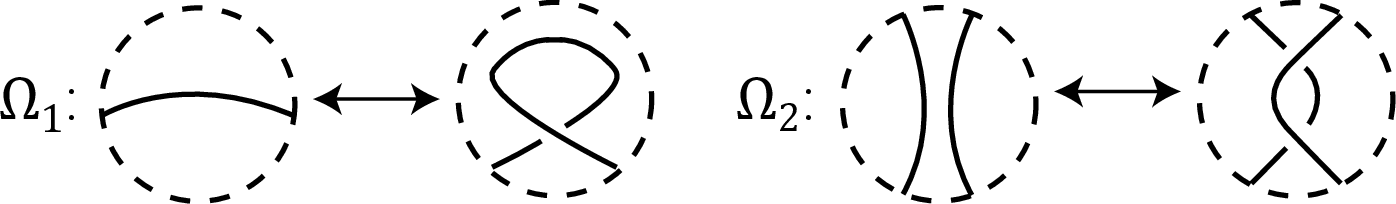}\\ \vbox{\phantom{1em}}
  \includegraphics[width=0.3\textwidth]{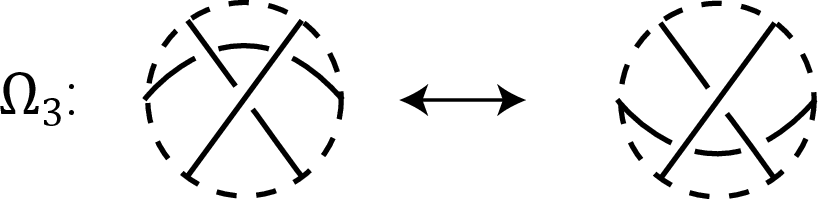}
  \caption{Reidemeister moves}\label{pic:reidmove}
\end{figure}

Combinatorially, a knot can be defined as an equivalence class of diagrams (Fig.~\ref{pic:knot_diagram}) modulo Reidemeister moves (Fig.~\ref{pic:reidmove}). This approach leads to Tait conjectures, skein polynomials, (bi)quandle cocycle invariants, and Khovanov homology. The growing popularity of diagram methods allows one to talk about the combinatorial revolution in knot theory~\cite{Nel}.

\begin{figure}[h]
\centering
  \includegraphics[width=0.4\textwidth]{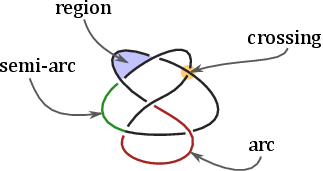}
  \caption{Elements of a knot diagram}\label{pic:diagram_parts}
\end{figure}

A knot diagram as an embedded $4$-valent graph determines sets of its elements such as arcs, semiarcs, crossings, and regions (Fig.~\ref{pic:diagram_parts}).

A Reidemeister move between two diagrams induces a correspondence between their elements. This correspondence is a bijection on the elements of the diagrams that are not involved in the move. On the other hand, the move (e.g. a second Reidemeister move) can split, merge, eliminate, or create diagram elements.

Constructions of many knot invariants use labels assigned to diagram elements. For example, the well-known formula for the linking number of a two-component link $L=K_1\cup K_2$
$$
lk(K_1,K_2)=\frac 12\sum_{c\in K_1\cap K_2} sgn(c)
$$
is a sum of quantitative labels (signs) of the crossings distinguished by another label (mixed component type) of the crossings.
Other examples exploit labelings of diagram arcs (quandle colorings), semiarcs (biquandle colorings), crossings (parity brackets and index polynomials) and regions (shadow quandle cocycles).

The labels used are assumed to be invariant under Reidemeister moves. This means that the labels of any correspondent elements (under some Reidemeister move) must coincide. In other words, invariant labels are maps to some coefficient set from the set of equivalence classes of diagram elements modulo the correspondences induced by the Reidemeister moves. Depending on the type of diagram elements, we call these equivalence classes \emph{(semi)arcs, crossings} and \emph{regions of the knot} because they are no longer linked to a concrete diagram. The aim of this paper is to provide a topological description of these classes.

\begin{figure}[h]
\centering
  \includegraphics[width=0.4\textwidth]{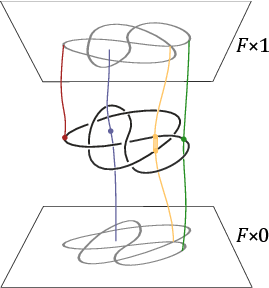}
  \caption{Topological interpretation of diagram elements}\label{pic:diagram_parts_top}
\end{figure}

The geometric definition of knot quandle provides a clue to the topological description of knot elements. Conceptually, the answer is as follows (Fig.~\ref{pic:diagram_parts_top}). Let $K$ be a knot in the thickening $F\times I$ (where $I=[0,1]$) of a compact oriented surface $F$. Then
\begin{itemize}
    \item the arcs of the knot $K$ are the isotopy classes of paths from $K$ to $F\times\{1\}$;
    \item the semiarcs of the knot are the isotopy classes of paths from $F\times\{0\}$ to $F\times\{1\}$ which intersect $K$ once;
    \item the crossings of the knot are the isotopy classes of paths from $F\times\{0\}$ to $F\times\{1\}$ which intersect $K$ twice; the segment between the intersection points is framed;
    \item the regions of the knot are the isotopy classes of paths from $F\times\{0\}$ to $F\times\{1\}$ which do not intersect $K$.
\end{itemize}

We will refer to the paths above as arc probes, semiarc probes, crossing probes, and region probes, respectively.

To give a more precise formulation of the statement, it is necessary to pay attention to the fact that knot invariants can be divided into two types. Invariants of the first type produce the same result regardless of the chosen diagram of the knot. This type includes numerical and polynomial invariants, such as the genus of the knot, the number of crossings, the Jones polynomial, and the Alexander polynomial.

The second type includes invariants (which we will refer to as coinvariants) whose values calculated on different knot diagrams are formally different but are isomorphic. Examples of these invariants are the knot group, the set of biquandle colorings, and the Khovanov homology.

Coinvariants form part of a categorical structure and, therefore, possess the property of functoriality and define a monodromy, which is the group of isomorphisms of the values of the coinvariant. The process of transitioning from invariants to coinvariants is known as categorification. A natural number, which represents the value of an invariant, can then be viewed as a dimension of some space, that is, an isomorphism class of finite-dimensional spaces. This interpretation is used in the construction of Khovanov homology.
To reverse this process, we can use decategorification, which involves the transition to orbits of coinvariant under the action of the monodromy group. 

Another refinement involves considering homotopy of sequences of Reidemeister moves. These are sequences of moves that transform one diagram to another, but they can be altered locally without affecting the final result, such as contracting forward and reverse Reidemeister moves. Coinvariants that remain unchanged under these homotopic transformations are called homotopy coinvariants or $h$-coinvariants.

Now, we can give a more accurate formulation of the result described by Fig.~\ref{pic:diagram_parts_top}: \emph{the set of isotopic classes of arc (crossing, region) probes of a tangle is the universal $h$-coinvariant of arcs (crossings, regions)}. Similarly, the set of orbits of the isotopic classes of diagram elements, under the action of the motion group of the knot, is the universal invariant of diagram elements.

In summary, the content of this article is as follows. Section~\ref{subsect:probe_spaces} describes the probe spaces that are used later when exploring the probes of the diagram elements.

In Section~\ref{subsect:diagram_categories}, we define categories of tangle diagrams. There are two types of diagram category: the rigid category and the homotopy category. In the rigid category, morphisms are sequences of Reidemeister moves, while in the homotopy category, the sequences of moves are considered up to homotopy.
Each category has two ways of describing it: a combinatorial one, where diagrams are objects and Reidemeister moves are morphisms, or topological one, where objects are tangles and morphisms are isotopies. The rigid category is useful for describing diagram elements as functors, but the homotopy category seems better suited for studying tangle invariants.

We introduce the notation for diagram elements in Section~\ref{subsect:diagram_elements}, and (combinatorial) strong and weak equivalences of diagram elements in Section~\ref{subsect:diagram_element_equivalence}.

Section~\ref{subsect:category_invariants} provides a formal definition of invariants and coinvariants for functors into the category of relations, specifically the functors associated with diagram elements.  We use strong and weak equivalences on sets of diagram elements to define a universal (co)invariant for the diagram element functors.

Sections~\ref{sect:arc}--\ref{sect:crossing} are dedicated to various elements of diagrams and follow a similar structure. Within these chapters, we introduce topological strong and weak equivalences and demonstrate that these equivalences are equivalent to the combinatorial ones. Consequently, the set of isotopy classes of diagram elements becomes a universal homotopical coinvariant of the corresponding diagram elements.

Additionally, Section~\ref{subsect:arc_motion_group} introduces the notion of tangle motion group, which can be understood as the tangle symmetry group. By factoring the set of isotopy classes of diagram elements under the action of the motion group, we obtain a universal invariant of the diagram elements.
Section~\ref{sect:crossing_R2_equivalence} discusses two types of equivalences on the set of crossings: the twisting equivalence and the equivalence generated by the second Reidemeister move. We show that these two types of equivalences coincide and describe the corresponding universal $h$-invariant.

Section~\ref{sect:midcrossing} discusses the midcrossings, i.e., the classes of crossings that arise when one admits the upper and lower third Reidemeister moves to be applied to the crossing. Midcrossings induce local transformations of tangles that can be used to define skein modules and skein invariants of tangles. Two examples of midcrossings are nugatory midcrossings that appear in the prime decomposition of a knot, and ribbon midcrossing that emerge as ribbon singularities of a ribbon knot. 

Section~\ref{sect:trait} considers the trait functor that is a modification of the crossing functor where a crossing survives after a third Reidemeister move. We give a description of the universal $h$-(co)invariant of the trait functor in terms of the fundamental group of the surface $F$ and thus reproduce the results of~\cite{Ntribe}.

Section~\ref{sect:elements} is devoted to relations between the diagram elements. In Section~\ref{sect:element_incidence} we consider the incidence relations between elements, and in Section~\ref{sect:element_compatibility_relations} we give a topological formulation for compatibility relations of diagram elements in Reidemeister moves. Section~\ref{subsect:element_representing_graph} shows how the elements of knots can be presented by probe diagrams. 

We define elements of tangles as isotopy classes of probes, and the second part of this paper is devoted to the study of homotopy classes. We find that the sets of homotopy classes of diagram elements is closely connected to colorings of diagrams using elements of a set with a given algebraic structure such as quandle. Moreover, the sets of homotopy classes of diagram elements appear to be the fundamental object with the given algebraic structure.

We introduce the homotopy classes of diagram elements in Section~\ref{sect:colorings}. In Section~\ref{subsect:quandle} we consider colorings of arcs of diagrams. We define the topological quandle of a tangle in a thickened surface and prove that it is the fundamental quandle. The topological quandle turns out to be an invariant of virtual tangles.

In Section~\ref{subsect:partial_tribracket} we consider colorings of regions of diagrams and the corresponding algebraic structure, that is a tribracket in partial ternary quasigroups.  We define the topological partial ternary quasigroup and show that it is fundamental. We also extend the definition of tribracket homology and tribracket cycle invariant to the case of partial ternary quasigroups.

In Section~\ref{subsect:biquandloid} we introduce a modification of biquandle that we call a biquandloid. The biquandloid structure is more suited for describing colorings of semiarc of diagrams in a fixed surface than that of biquandle. We show that the topological biquandloid is fundamental, and extend the constuction of biquandle cycle invariant to biquandloids.

In Section~\ref{subsect:crossoid} we define the crossoid structure that can be considered as a generalization of parities on knots. Once again, we define the topological crossoid and show that it is fundamental. We define the crossoid homology and the crossoid cycle invariant. The crossoid cycle invariant generalizes both the biquandle cycle invariant and the index polynomial. 

In Section~\ref{subsect:multicrossing_complex} we introduce the multicrossing complex and define the crossing homology class. In a sense, the multicrossing complex unifies tribracket, biquandle and crossoid homologies, and the tribracket, biquandle and crossoid cycle invariants are actually the result of pairing a tribracket (biquangle, crossoid) cocycle with the crossing homology class. 

The paper ends with two speculative sections devoted to invariants which use diagram elements in their construction, and further research directions.

{\color{red}
}

\section{Diagram}\label{sect:diagrams}

\subsection{Tangles and tangle diagrams}

\begin{figure}[h]
\centering
  \includegraphics[width=0.4\textwidth]{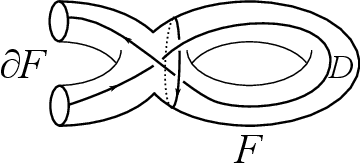}
  \caption{A tangle diagram}\label{pic:tangle}
\end{figure}

\begin{definition}\label{def:tangle}
Let $F$ be an oriented compact connected surface and $Y\subset \partial F\times (0,1)$ a finite set.
A \emph{tangle} is an embedding $T\colon M\to F\times (0,1)$ of an oriented compact $1$-manifold $M$ into the thickening of the surface $F$ such that $T(M)\cap\partial F\times I=T(\partial M)=Y$ (Fig.~\ref{pic:tangle}). We assume that the embedding is transversal to the boundary.
\end{definition}

Let $\Sigma=\Sigma(F,Y)$ be the space of tangles equipped with the Whitney topology. The space $\Sigma$ is stratified with respect to the natural projection $p\colon F\times I\to F$: $\Sigma=\bigcup_{k\in\N\cup\{\infty\}}\Sigma_k$. We need only the first three strats.

The subset $\Sigma_0$ consists of tangles $T$ for whose the projection $p_T\colon T\to p(T)$ is a bijection except a finite number of double points. Any tangle $T\in\Sigma_0$ determines a tangle diagram.

\begin{figure}[h]
\centering
  \includegraphics[width=0.3\textwidth]{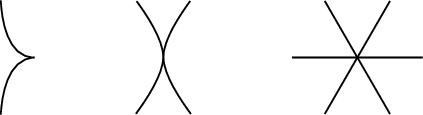}
  \caption{Codimension $1$ singularities: an ordinary cusp, a simple tangency and a triple point}\label{pic:singularity_codim1}
\end{figure}

The subset $\Sigma_1$ consists of tangles $T$ whose projection includes one singular point of codimension $1$ (Fig.~\ref{pic:singularity_codim1}) besides a finite number of double point.

\begin{figure}[h]
\centering
  \includegraphics[width=0.5\textwidth]{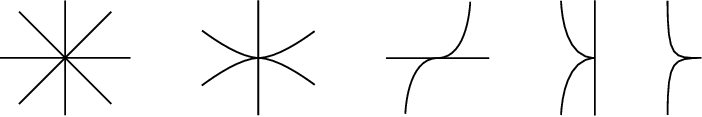}
  \caption{Codimension $2$ singularities: a quadruple point, a tangent triple point, a cubic tangency, an intersected cusp and a ramphoidal cusp~\cite{FiedKur}}\label{pic:singularity_codim2}
\end{figure}

The subset $\Sigma_2$ consists of tangles whose projection contains one singular point of codimension $2$ (Fig.~\ref{pic:singularity_codim2}) or two points of codimension $1$.

\begin{definition}\label{def:sigma_transversility}
    A smooth family of tangles $\tau_t$, $t\in[0,1]$, is \emph{$\Sigma_1$-transversal} if
\begin{enumerate}
    \item $\tau_t\in\Sigma_0\cup\Sigma_1$ for all $t$;
    \item $\tau_0,\tau_1\in\Sigma_0$;
    \item $\tau_t$ is transversal to $\Sigma_1$.
\end{enumerate}

A homotopy $\Gamma=(\tau^s)_{s\in[0,1]}$ of paths with fixed ends is \emph{$\Sigma_2$-transversal} if
\begin{enumerate}
    \item $\tau^0$ and $\tau^1$ are $\Sigma_1$-transversal paths;
    \item for all $s\in [0,1]$ except a finite set the path $\tau^s$ is $\Sigma_1$-transversal;
    \item for any exceptional $s$ there exist a unique $t^*\in (0,1)$ such that
\begin{itemize}
    \item $\tau^s$ satisfies the conditions of $\Sigma_1$-transversality on $[0,1]\setminus\{t^*\}$;
    \item either $\tau^s_{t^*}\in\Sigma_2$ or $\tau^s_{t^*}$ is a non-transversal intersection of $\tau^s$ with $\Sigma_1$.
\end{itemize}
\end{enumerate}
\end{definition}

The transversality theorem implies the following. 
\begin{proposition}\label{prop:sigma_transversality}
\begin{enumerate}
    \item the set $\Sigma_0$ is dense in $\Sigma$;
    \item any path $\tau_t, t\in[0,1]$, in $\Sigma$ such that $\tau_0,\tau_1\in\Sigma_0$ can be approximated by a $\Sigma_1$-transversal path $\tilde\tau_t$ such that $\tilde\tau_0=\tau_0$ and $\tilde\tau_1=\tau_1$;
    \item any homotopy $\Gamma=(\tau^s)_{s\in[0,1]}$ between $\Sigma_1$-transversal paths $\tau^0$ and $\tau^1$ can be approximated by a $\Sigma_2$-transversal homotopy $\tilde\Gamma=(\tilde\tau^s)_{s\in[0,1]}$.
\end{enumerate}
\end{proposition}

\begin{definition}\label{def:tangle_diagram}
Let $F$ be an oriented compact connected surface and $X\subset\partial F$ a finite set. A \emph{tangle diagram} is an embedding $D\colon G\to F$ of a finite graph $G$  into $F$ such that
\begin{itemize}
    \item the set of vertices $V(G)=\mathcal C(G)\sqcup \partial G$ splits into the set $\mathcal C(G)$ of vertices of valency $4$ called \emph{crossings}, and the set $\partial G$ of vertices of valency $1$ called \emph{boundary vertices};
    \item $D(G)\cap\partial F=D(\partial G)=X$, and $D$ is transversal to $\partial F$;
    \item each crossing $c\in\mathcal C(G)$ possesses an undercrossing-overcrossing structure: two half-edges incident to $c$ are marked as the \emph{undercrossing}, and the other two are marked as the \emph{overcrossing}; the embedding $D$ maps the undercrossing edges to a pair of opposite edges.
\end{itemize}
\end{definition}

Note that the image $D(G)\subset F$ determines the graph $G$ up to isomorphism. Therefore, below we will often identify the diagram map $D$ and the set $D(G)$ (with the under-overcrossings structure indicated).

Diffeomorphisms of $F$ constant on the boundary act on tangle diagrams by compositions: a tangle diagram $D$ and $f\in \mathrm{Diff}_0(F)$ yield a tangle diagram $f\circ D$.

\begin{definition}\label{def:top_reidemeister_move}
Two tangle diagrams $D_1$ and $D_2$ are \emph{connected by a Reidemeister move} if there is an embedded disc $B\subset F$ such that
\begin{enumerate}
    \item the diagrams $D_1$ and $D_2$ coincide outside $B$;
    \item there is a diffeomorphism $\phi\colon B\to\mathbb D^2$ to the standard disc in $\mathbb R^2$ such that $\phi(D_1\cap B)$ and $\phi(D_2\cap B)$ are the left- and right-hand side diagrams of one of the three standard Reidemeister moves shown in Fig.~\ref{pic:reidmove}.
\end{enumerate}
The described Reidemeister move will be denoted by $\Omega_i(B,\phi)$, $i=1,2,3$. Analogously, one defines the inverse Reidemeister moves denoted by $\Omega_i^\dag(B,\phi)$, $i=1,2,3$.
\end{definition}

\begin{remark}\label{rem:canonical_lifting}
Note that any tangle $T\in\Sigma_0$ determines a tangle diagram $p(T)$ by projection. On the other hand, given a tangle diagram $D\colon G\to F$, one defines its \emph{canonical lifting} $\lambda(D)$ as follows.

Let $e\colon[0,l]\to F$ be the restriction of $D$ to an edge $e$ of the graph $G$. Consider the map $\tilde e=(e,h)\colon[0,l]\to F\times(0,1)$, where
$$
h(t)=\frac 12+\frac 14\epsilon_0\cdot\alpha\left(\frac{2t}l\right)+ \frac 14\epsilon_1\cdot\alpha\left(2-\frac{2t}l\right),
$$
$\epsilon_0=-1$ (resp. $\epsilon_0=1$) if the beginning half-edge of $e$ is undercrossing (resp. overcrossing), $\epsilon_1=-1$ (resp. $\epsilon_1=1$) if the ending half-edge of $e$ is undercrossing (resp. overcrossing), and $\alpha\colon\mathbb R\to\mathbb R$ is a fixed smooth descending function such that $\alpha(x)=1$ for all $x<0$ and $\alpha(x)=0$ for all $x>1$.

Concatenation of the lifts $\tilde e$ for all the edges of $G$ forms a tangle $\lambda(D)$ in the thickening $F\times I$.
\end{remark}

Let $\tau_t$, $t\in[0,1]$, be a smooth family of tangles. By isotopy extension theorem, there exists an isotopy $H_t\in Diff(F\times I)$ such that $H_0=id$ and $\tau_t=H_t\circ\tau_0$, $t\in[0,1]$. We call $H_t$ an \emph{extension of $\tau_t$}. From the definition we have the following statement.

\begin{proposition}\label{prop:isotopy_lifting_ambiguity}
Let $H_t$ and $H'_t$ be extensions of a tangle path $\tau_t$. Then $H'_1=f\circ H_1$ for some $f$ in the stabilizer
\[
\difd(F\times[0,1], \tau_1)=\{f\in\difd(F\times I)\mid f|_{\tau_1}=id\}.
\]
\end{proposition}

\subsection{The probe spaces}\label{subsect:probe_spaces}

For a connected compact oriented surface $F$ denote:
\begin{itemize}
    \item the group of diffeomorphisms of $F\times I$ rel $\partial(F\times I)=\partial F\times I\cup F\times\{0,1\}$ by  $Diff(F\times I)$;
    \item the group of pseudoisotopies, i.e. diffeomorphisms of $F\times I$ rel $\partial F\times I\cup F\times 0$ by  $P(F)$;
    \item the group of diffeomorphisms $\psi$ of $F\times I$ rel $\partial F\times I$ such that $\psi(F\times i)=F\times i$, $i=0,1$, by $\widetilde{Diff}(F\times I)$.
\end{itemize}

The homotopy type of these groups can be described as follows~\cite{Budney}.

\begin{proposition}\label{prop:diffeomorphism_groups_type}
\begin{enumerate}
    \item $Diff(F\times I)\sim \Omega(Diff_0(F))$;
    \item $P(F)$ is contractible;
    \item $\widetilde{Diff}(F\times I)\sim Diff(F)$.
\end{enumerate}
Here $Diff(F)$ is the group of diffeomorphisms of $F$ rel $\partial F$ and $Diff_0(F)$ is the component of the identity.
\end{proposition}

\begin{proof}
   1. For $F=S^2$, the equivalence $Diff(F\times I)\sim \Omega(Diff_0(F))$ is considered in~\cite{Hatcher}, and for $F=T^2$, in~\cite{HK}. Let $F\ne S^2, T^2$. We will show that $Diff(F\times I)\sim * \sim \Omega(Diff_0(F))$.

   Let $\gamma$ be a simple non-contractible closed curve in $F$ or a simple arc connecting different components of $\partial F$. Then the annulus $A=\gamma\times I$ is incompressible in $F\times I$. There is a fibration
\[
Diff(F'\times I)\to Diff(F\times I)\to Emb(A, F\times I \ \mathrm{rel}\ \partial A)
\]
where $F'$ is the surface obtained by cutting $F$ along $\gamma$. The embedding space $Emb(A, F\times I\ \mathrm{rel}\ \partial A)$ is contractible by~\cite{Hatcher1}. Hence, $Diff(F\times I)\sim Diff(F'\times I)$. Then the surface $F$ can be reduces to a disjoint union of disks. By Smale conjecture~\cite{Hatcher}, $Diff(\bigsqcup_i D^2\times I)=\prod_i Diff(D^2\times I)\sim *$. Thus,  $Diff(F\times I)\sim *$.

2. The fibration
$ 
Diff(F\times I)\to P(F)\to Diff_0(F),
$ 
defined by the restriction of diffeomorphisms in $P(F)$ to $F\times 1$, induces the fibration
\[
\Omega(Diff_0(F))\to Diff(F\times I)\to P(F).
\]
Since $\Omega(Diff_0(F))\sim Diff(F\times I)$, the space $P(F)$ is contractible.

3. The restriction of a diffeomorphism in $\widetilde{Diff}(F\times I)$ to $F\times 0$ defines the fibration
\[
P(F)\to\widetilde{Diff}(F\times I)\to Diff(F).
\]
Since $P(F)$ is contractible, $\widetilde{Diff}(F\times I)\sim Diff(F)$.
\end{proof}

Note that $Diff(F)$ is a retract of $\widetilde{Diff}(F\times I)$.

\begin{definition}\label{def:probe}
1. A \emph{probe} is an unknotted embedding
\[
    \gamma\colon (I;0;1)\to (Int(F)\times I; Int(F)\times 0; Int(F)\times 1)
\]
transversal to $F\times\partial I$. A \emph{vertical probe} is a probe of the form $\gamma_x=x\times I$, $x\in F$.
Denote the space of probes by ${Pr}(F)$ and the space of vertical probes by ${VPr}(F)$.

2. An \emph{overprobe} is an embedding $\gamma\colon (I;1)\to (Int(F)\times (0,1]; Int(F)\times 1)$ transversal to $F\times 1$. An overprobe $\gamma$ is \emph{vertical} if $\gamma\subset x\times I$ for some $x\in F$. Denote the space of overprobes by ${Pr}^o(F)$ and the space of vertical overprobes by ${VPr}^o(F)$.

3. An \emph{underprobe} is an embedding $\gamma\colon (I;1)\to (Int(F)\times [0,1)]; Int(F)\times 0)$ transversal to $F\times 0$. An underprobe $\gamma$ is \emph{vertical} if $\gamma\subset x\times I$ for some $x\in F$. Denote the space of underprobes by ${Pr}^u(F)$ and the space of vertical underprobes by ${VPr}^u(F)$.

4. A \emph{midprobe} is an embedding $\gamma\colon I\to Int(F)\times (0,1)$. A midprobe $\gamma$ such that $\gamma\subset x\times I$ for some $x\in F$, is called \emph{vertical}. The space of midprobes and vertical midprobes are denoted by ${Pr}^m(F)$ and ${VPr}^m(F)$.
\end{definition}

The goal of the section is to describe the homotopy type of the spaces of probes and vertical probes.

\begin{proposition}\label{prop:basic_probe_spaces}
\begin{enumerate}
    \item $VPr(F)\sim VPr^o(F)\sim VPr^u(F)\sim VPr^m(F)\sim F$;
    \item $Pr^o(F)\sim Pr^u(F)\sim F$;
    \item $Pr^m(F)\sim F\times S^2$.
\end{enumerate}
\end{proposition}

\begin{proof}
    The first statement follow from the definition of vertical probes.

    Given an overprobe, one can contract it to the end in $F\times 1$. This process induces the homotopy eqivalence $Pr^o(F)\sim F\times 1\sim F$. Analogously, we have $Pr^u(F)\sim F$.

    Contracting a midprobe to its beginning, we get a homotopy equivalence between $Pr^m(F)$ and the unit tangent bundle $UT(F\times (0,1))$. Since $F$ is oriented, $F\times (0,1)$ can be embedded to $\mathbb R^3$. Hence,
\[
UT(F\times (0,1))=UT(\mathbb R^3)|_{F\times (0,1)}=(F\times (0,1))\times S^2\sim F\times S^2.
\]
\end{proof}

\begin{remark}\label{rem:midprobe_homotopy_generator}
The equivalence $Pr^m(F)\sim F\times S^2$ implies $\pi_2(Pr^m(F))=\pi_2(F)\times \pi_2(S^2)=\pi_2(F)\times\mathbb Z$. In particular, if the genus of $F$ is positive then $\pi_2(Pr^m(F))=\mathbb Z$. Let us describe the generator of $\pi_2(Pr^m(F))$.

Consider a small neighborhood $U$ of an inner point of $y_0\in F\times I$. We identify $U$ with $\R^3$ so that $y_0$ goes to $0$ and the third coordinate in $\R^3$ corresponds to the vertical direction of $F\times I$. Consider the map $R\colon\mathbb D^2\to SO(3)$ defined by the formula

\begin{multline*}
R(re^{i\phi})=\left(\begin{array}{ccc}
     \cos \phi & -\sin \phi & 0 \\
     \sin \phi & \cos \phi & 0 \\
    0 & 0 & 1
\end{array}\right)
\left(\begin{array}{ccc}
     \cos \pi(1-r) & 0 & -\sin  \pi(1-r) \\
     0 & 1 & 0 \\
     \sin  \pi(1-r) & 0 & \cos  \pi(1-r)
\end{array}\right)\cdot\\
\cdot\left(\begin{array}{ccc}
     \cos \phi & -\sin \phi & 0 \\
     \sin \phi & \cos \phi & 0 \\
    0 & 0 & 1
\end{array}\right).
\end{multline*}

Let $\gamma_0=\{(0,0,t)\mid t\in[0,1]\}$ be a vertical midprobe. Consider the map $f\colon\mathbb D^2\to Pr^m$ defined by the formula $f(z)=R(z)\gamma_0$. Since $f(z)=\gamma_0$ for any $z\in\partial\mathbb D^2$, the map $f$ defines an element $[f]\in\pi_2(Pr^m)$.

Let us extend the map $R$ to a map $\hat R\colon\mathbb D^2\to Diff(\R^3)$.
Consider a monotone function $g\in C^\infty(\R)$ such that $g(t)=0$ for $t<1$ and $g(t)=1$ for $t>2$. Then set $\hat R$ by the formula
\[
\hat R(z)(x)=R((1-z)g(|x|)+z)x,\ z\in\mathbb D^2, x\in\mathbb R^3.
\]
Note that $\hat R(z)(x)=R(z)x$ for $x\in\mathbb D^3$, and $\hat R(z)(x)=x$ when $|x|\ge 2$. For any $z\in\partial\mathbb D^2$, $\hat R(z)(\gamma_0)=\gamma_0$, and the restriction of $\hat R$ to the circle $\partial\mathbb D^2$ illustrates the fact that the process of the rotation by $4\pi$ around an axis is homotopic to the identity in $SO(3)$.

The diffeomorphism $\hat R(z)$ can be thought of as a diffeomorphism of $U$. One can extend it by the identity to a diffeomorphism of $F\times I$. Thus, we get a map $\hat R\colon\mathbb D^2\to Diff(F\times I)$.
\end{remark}

\begin{proposition}\label{prop:probe_space}
    The inclusion ${VPr}(F)\hookrightarrow Pr(F)$ is a homotopy equivalence.
\end{proposition}

\begin{proof}

1. For the sphere $F=S^2$, the homotopy equivalence can be constructed as follows. Consider the thickened sphere as a ball $B$ with its core removed.  Given a probe $\gamma$, consider the isotopy $\varphi_t$ which pulls the core sphere along $\gamma$ in $B$ (Fig.~\ref{pic:sphere_probe_verticalize}). Let $z(t)$ be the center of the core in the moment $t$.  Consider a family $\psi_z$, $z\in B$, of diffeomorphisms of $B$ (for example, M{\"o}bius transformations) which move the core with the center $z$ to the center of the ball. Then the isotopy $\psi_{z(t)}\circ\varphi_t$ verticalizes the probe.

\begin{figure}[ht]
\centering
  \includegraphics[width=0.4\textwidth]{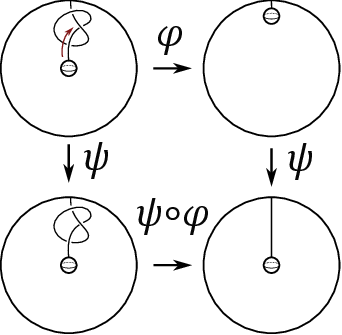}
  \caption{Verticalizing a probe in the thickened sphere}\label{pic:sphere_probe_verticalize}
\end{figure}


2. Let $F\ne S^2$.

\begin{lemma}
    Let $F\ne S^2$ and $Emb_0(I,F\times I)$ denote the space of embeddings $\sigma\colon I\hookrightarrow F\times I$, $\sigma(k)=(x_0,k)$, $k=0,1$, which are isotopic with fixed ends to the vertical probe $x_0\times I$. Then $Emb_0(I,F\times I)$ is contractible.
\end{lemma}

\begin{proof}
    Consider the natural action of $Diff_0(F\times I)$ on the framed probes. We have the exact sequence
\[
1\to Diff((F\setminus D)\times I)\to Diff_0(F\times I)\to \mathcal K^{fr}_0\to 1
\]
where $D\subset F$ is a small disk containing $x_0$ and $\mathcal K^{fr}_0$ is the orbit of the  (framed) vertical probe $x_0\times I$. Since $Diff((F\setminus D)\times I)\sim*$, we have $\mathcal K^{fr}_0\sim Diff_0(F\times I)\sim *$.

By forgetting the framing, we get the covering $f:\mathcal K_0^{fr}\to Emb_0(I,F\times I)$.
We claim that the covering is an equivalence. Let $U^{fr}_0$ be the unknot $U=x_0\times I$ with the initial framing, and $U^{fr}_1\in f^{-1}(U)$. Then there is an isotopy of long framed knots $U^{fr}_t$, $t\in[0,1]$. Consider the covering $\tilde F$ of $F$ which embeds in $\mathbb R^2$. Then the isotopy $U^{fr}_t$ lifts to $\tilde U^{fr}_t$ in $\tilde F\times I\subset\mathbb R^2\times I$. Hence, the framings of $\tilde U^{fr}_0$ and $\tilde U^{fr}_1$ (thus, of $U^{fr}_0$ and $U^{fr}_1$) coincide. Then $f^{-1}(U)=\{U^{fr}_0\}$ and $Emb_0(I,F\times I)\simeq\mathcal K^{fr}_0\sim *$.
\end{proof}

Now, return the proof of Proposition. Consider the fibration
\[
Emb_0(I,F\times I) \to Pr(F) \to \Pi_1(F)
\]
induced by the composition with the projection $p\colon F\times I\to F$. Here $\Pi_1(F)$ is the fundamental groupoid of the surface $F$. Then $Pr(F)\sim\Pi_1(F)$.

The source map of the groupoid  $\Pi_1(F)$ defines a fibration
\[
\tilde F\to \Pi_1(F)\to F
\]
where $\tilde F$ is the universal covering of $F$. Since $\tilde F$ is contractible, we have
\[
Pr(F)\sim \Pi_1(F)\sim F\sim VPr(F).
\]
\end{proof}

The equivalence between $Pr(F)$ and $VPr(F)$ implies the following statement.
\begin{corollary}\label{cor:probes_verticalizing}
    1. Let $\gamma_t$, $t\in I$, is a continuous family of probes such that $\gamma_0$ and $\gamma_1$ are vertical. Then there is a family of vertical probes $\gamma'_t$ which is homotopic to $\gamma_t$ with fixed ends.

    2. Let $\gamma_t$ and $\gamma'_t$, $t\in[0,1]$ be families of vertical probes such that $\gamma_0=\gamma'_0$ and $\gamma_1=\gamma'_1$. If $\gamma_t$ and $\gamma'_t$ are homotopic (with the ends fixed) in the space of probes then they are homotopic in the  space of vertical probes.
\end{corollary}

\begin{proposition}\label{prop:probes_verticalizing}

  1. For any (over,under,mid)probe $\gamma_0$ there exists an isotopy $h\colon I\to \widetilde{Diff}(F\times I)$ such that $h_0=id$ and the (over,under,mid)probe $h_1(\gamma_0)$ is vertical.

  2.  Let $\gamma_0$ and $\gamma_1$ be vertical (over,under,mid)probes and $h_t$, $t\in I$, an isotopy of $F\times I$ such that $h_0=id$ and $h_1(\gamma_0)=\gamma_1$. Then there exists a verticalizing homotopy of the isotopy $h_t$, i.e. a family of diffeomorphisms $H\colon I\times I\to \widetilde{Diff}(F\times I)$ such that $H_{t0}=h_t$, $H_{0s}=id$, $H_{1s}=h_1$ and $H_{t1}(\gamma_0)$ are vertical (over,under,mid)probes for all $t\in[0,1]$ and $s\in I$.
\end{proposition}

\begin{proof}
    1. Assume that $\gamma_0$ is a probe. (The proof for over-, under- and midprobes is analogous.) The inclusion $VPr\hookrightarrow Pr$ induces an isomorphism $\pi_0(VPr)\simeq \pi_0(Pr)$. Hence, there is a family of probes $\gamma_t$, $t\in[0,1]$, such that $\gamma_1\in VPr$. By isotopy extension theorem, there exists $h\colon [0,1]\to \widetilde{Diff}(F\times I)$ such that $h_0=id$ and $h_t(\gamma)=\gamma_t$. In particular, $h_1(\gamma)=\gamma_1\in VPr$.

    2. Let $\gamma_0=x_0\times I$, $\gamma_1=x_1\times I\in VPr$. Consider the path $\gamma=(\gamma_t)_{t\in I}$ in $Pr$ where $\gamma_t=h_t(\gamma_0)$. Since $F$ is connected, there is a path $\alpha$ from $x_0$ to $x_1$. It induces a path $\tilde\alpha$ from $\gamma_0$ to $\gamma_1$ in $VPr$. There is a homotopy $\gamma\sim \gamma\tilde\alpha^{-1}\tilde\alpha$.

    Due to the isomorphism $\pi_1(VPr)\simeq \pi_1(Pr)$, the loop $\gamma\tilde\alpha^{-1}$ is homotopic to a loop $\gamma'$ in $VPr$. Hence, $\gamma\tilde\alpha^{-1}\tilde\alpha\sim\gamma'\tilde\alpha$. Then $\gamma\sim\gamma'\tilde\alpha$, i.e. there is a homotopy $\Gamma\colon I\times I\to Pr$ such that $\Gamma_{t0}=\gamma_t$, $\Gamma_{is}=\gamma_i$, $i=0,1$, and $\Gamma_{t1}\in VPr$ for each $s,t\in I$. By isotopy excision theorem, there exists a isotopy $H\colon I\times I\to \widetilde{Diff}(F\times I)$ such that $H_{t0}=h_t$, $H_{0s}=id$, $H_{1s}=h_1$ and $H_{ts}(\gamma_0)=\Gamma_{ts}$ for all $s,t\in I$. In particular, $H_{t1}(\gamma_0)=\gamma_{t1}\in VPr$ for all $t\in I$.
\end{proof}

\begin{proposition}\label{prop:probe_verticalizing_homotopy}
    Let $\gamma=(\gamma_t)_{t\in I}$ and $\gamma'=(\gamma'_t)_{t\in I}$ be continuous families of vertical (over,under)probes such that $\gamma_0=\gamma'_0$%
, $\gamma_1=\gamma'_1$
    and there is a map $H\colon I\times I\to\widetilde{Diff}(F\times I)$ such that $H_{t0}(\gamma_0)=\gamma_t$, $H_{t1}(\gamma_0)=\gamma'_t$, $H_{0s}(\gamma_0)=\gamma_0$
 and $H_{1s}(\gamma_0)=\gamma_1$
    for all $t,s\in I$. Then there is a isotopy $\tilde H\colon I\times I\times I\to\widetilde{Diff}(F\times I)$ such that $\tilde H_{ts0}=H_{ts}$,
    $\tilde H_{ksu}=H_{ks}$, $\tilde H_{tku}=H_{tk}$, and $\tilde H_{ts1}(\gamma_0)$ is vertical for all $t,s,u\in I$, $k=0,1$.
\end{proposition}

\begin{proof}
    Let $\gamma_t$ and $\gamma'_t$ be probes. (The proof for over- and underprobes is analogous). Denote $\Gamma=(\gamma_{ts})_{t,s\in I}$ where $\gamma_{ts}=H_{ts}(\gamma_0)$.

    Since $\pi_2(Pr,VPr)=0$, there is a family of probes $\tilde\Gamma=(\gamma_{tsu})_{t,s,u\in I}$ such that $\gamma_{ts0}=\gamma_{ts}$, $\gamma_{t0u}=\gamma_t$, $\gamma_{t1u}=\gamma'_t$, $\gamma_{0su}=\gamma_0$, $\gamma_{1su}=\gamma_1$ and $\gamma_{ts1}\in VPr$ for all $t,s,u\in I$. By isotopy excision theorem, there exists a isotopy $\tilde H\colon I^3\to\widetilde{Diff}(F\times I)$ such that $\tilde H_{ts0}=H_{ts}$, $\tilde H_{ksu}=H_{ks}$, $\tilde H_{tku}=H_{tk}$, and $\tilde H_{tsu}(\gamma_0)=\gamma_{tsu}$ for all $t,s,u\in I$. Then $\tilde H_{ts1}(\gamma_0)\in VPr$ for all $t,s\in I$.
\end{proof}

Since $\pi_2(Pr^m,VPr^m)=\mathbb Z$, an analogous proposition for midprobes is more elaborated.

\begin{proposition}\label{prop:midprobe_verticalizing_homotopy}
    Let $\gamma=(\gamma_t)_{t\in I}$ and $\gamma'=(\gamma'_t)_{t\in I}$ be continuous families of vertical midprobes such that $\gamma_0=\gamma'_0$, $\gamma_1=\gamma'_1$ and there is a map $H\colon I\times I\to\widetilde{Diff}(F\times I)$ such that $H_{t0}(\gamma_0)=\gamma_t$, $H_{t1}(\gamma_0)=\gamma'_t$, $H_{0s}(\gamma_0)=\gamma_0$ and $H_{1s}(\gamma_0)=\gamma_1$ for all $t,s\in I$. Then there is a map $H'\colon I\times I\to\widetilde{Diff}(F\times I)$ such that $H'_{t0}=H_{t0}$, $H'_{t1}=H_{t1}$, $H'_{0s}=H_{0s}$, $H'_{1s}(\gamma_0)=\gamma_1$ and $H'_{ts}(\gamma_0)\in VPr^m(F)$ for all $t,s\in I$.
\end{proposition}

\begin{figure}[h]
\centering
  \includegraphics[width=0.5\textwidth]{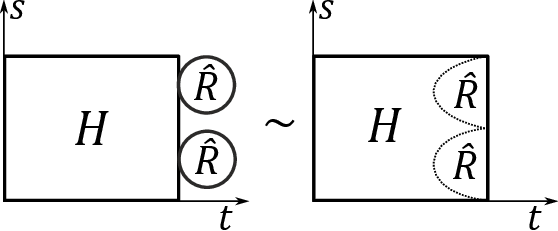}
  \caption{Correction of the isotopy $H$}\label{pic:midprobe_isotopy_correction}
\end{figure}

\begin{proof}
    The map $\Gamma\colon I\times I\to Pr^m(F)$, $\Gamma(t,s)=H_{ts}(\gamma_0)$, defines a class $[\Gamma]\in\pi_2(Pr^m,VPr^m)=\mathbb Z$. If $[\Gamma]=0$ then we can use the reasoning in the proof of Proposition~\ref{prop:probe_verticalizing_homotopy} to get a isotopy $\tilde H_{tsu}$, $t,s,u\in I$. Then the isotopy $H'_{ts}=\tilde H_{ts1}$ satisfies the required conditions.

    If the class $[\Gamma]\ne0$, compensate it by attaching isotopies $\hat R$ from Remark~\ref{rem:midprobe_homotopy_generator} (Fig.~\ref{pic:midprobe_isotopy_correction} left). The composite isotopy is homotopic to an isotopy $H''\colon I\times I\to \widetilde{Diff}(F\times I)$ such that $H''_{t0}=H_{t0}$, $H''_{t1}=H_{t1}$, $H''_{0s}=H_{0s}$, $H''_{1s}(\gamma_0)=\gamma_1$ for all $t,s\in I$ (Fig.~\ref{pic:midprobe_isotopy_correction} right). By definition the family of midprobes $H''_{ts}(\gamma_0)$, $s,t\in I$, defines a trivial class in $\pi_2(Pr^m,VPr^m)$. By applying the reasoning of the previous case to $H''_{ts}$, we get the desired isotopy.
\end{proof}

{
\color{red} 
}

\subsection{Diagram categories}\label{subsect:diagram_categories}

Below we define two categories. The first one (strict category) is a combinatorial category of diagrams and sequences of Reidemeister moves between them. The second (homotopical category) is a topological category of tangle embeddings and isotopies. 



\begin{definition}\label{def:strict_tangle_category}
    Consider the category $\mathfrak{T}_s$ such that
\begin{itemize}
    \item $Ob(\mathfrak T_s)=\Sigma_0$;
    \item for tangles $T_0,T_1\in\Sigma_0$, the morphism set $\mathrm{Mor}_{\mathfrak T_s}(T_0,T_1)$ is the set of equivalence classes of $\Sigma_1$-transversal paths from $T_0$ to $T_1$ modulo homotopies with fixed ends in the class of $\Sigma_1$-transversal paths.
\end{itemize}
The category $\mathfrak{T}_s$ is called the \emph{strict category of tangles}.
\end{definition}

Note that a homotopy though $\Sigma_1$-transversal paths does not change the number of the $\Sigma_1$-points in the path (which correspond to Reidemeister moves) nor the types of the correspondent moves.

\begin{definition}\label{def:strict_diagram_category}
    Consider the category $\mathfrak{D}_s$ such that
\begin{itemize}
    \item $Ob(\mathfrak D_s)$ are tangle diagrams;
    \item for tangle diagrams $D,D'$, the morphism set $\mathrm{Mor}_{\mathfrak D_s}(D_0,D_1)$ is the set of sequences
$$
D=D_0\stackrel{f_1}{\rightarrow} D_1\stackrel{f_2}{\rightarrow}\cdots\stackrel{f_n}{\rightarrow} D_n.
$$
    of diagram isotopies and Reidemeister moves $f_i$ modulo the relations:
    \begin{itemize}
        \item $D\stackrel{f}{\rightarrow} D'\ \sim\ D\stackrel{f'}{\rightarrow} D'$ for isotopies $f$ and $f'$ which are homotopic;
        \item $D_1\stackrel{f}{\rightarrow} D_2\stackrel{g}{\rightarrow} D_3\ \sim\ D_1\stackrel{g\circ f}{\rightarrow} D_3$ for isotopies $f$ and $g$;
        \item $\xymatrix{D_1\ar[r]^f &D_2\ar[r]^{\Omega(B,\phi)} &D_3\ar[r]^{f^{-1}} &D_4}\ \sim\ \xymatrixcolsep{5pc}\xymatrix{D_1\ar[r]^{\Omega(f^{-1}(B),\phi\circ f)} &D_4}$ for an isotopy $f$ and a Reidmeister move $\Omega(B,\phi)$, $\Omega=\Omega_i$ or $\Omega^\dag_i$, $i=1,2,3$.
    \end{itemize}
\end{itemize}
The category $\mathfrak{D}_s$ is called the \emph{strict category of tangles diagrams}.
\end{definition}

\begin{theorem}\label{thm:strict_categories_equivalence}
    The categories $\mathfrak{T}_s$ and $\mathfrak{D}_s$ are equivalent.
\end{theorem}

\begin{proof}
Let us construct functors $\Phi\colon \mathfrak T_s\to\mathfrak D_s$ and  $\Psi\colon \mathfrak T_s\to\mathfrak D_s$.

Given a tangle $T\in\Sigma_0$, its projection $p(T)$ into $F$ has a structure of a tangle diagram which we denote by $\Phi(T)$.

Let $T_0,T_1\in\Sigma_0$ be tangles and $\tau_t$ a $\Sigma_1$-transversal path from $T_0$ to $T_1$. Let $\{t_k\}_{k=1}^n=\tau^{-1}(\Sigma_1)$ be the set of singular points arranged in the ascending order. For a small $\epsilon>0$ denote $D_0=\Phi(\tau_0)$, $D_{2k-1}=\Phi(\tau_{t_k-\epsilon})$,  $D_{2k}=\Phi(\tau_{t_k+\epsilon})$, $k=1,\dots,n$ and $D_{2n+1}=\Phi(\tau_1)$. Then diagrams $D_{2k}$ and $D_{2k+1}$ differ by an isotopy $f_{k+1}$, and diagrams $D_{2k-1}$ and $D_{2k}$ differ by a Reidemeister move $r_k$. Define the morphism $\Phi(\tau_t)\in Mor_{\mathfrak D_s}(\Phi(T_0),\Phi(T_1))$ as the composition
$$
\xymatrix{D_0\ar[r]^{f_1} & D_1\ar[r]^{r_1} & D_2 \ar[r]^{f_2} & \cdots  \ar[r]^{r_n} & D_{2n} \ar[r]^{f_{n+1}} & D_{2n+1}}.
$$
The equivalence class of $\Phi(\tau_t)$ is well defined and $\Phi(\tau_t)=\Phi(\tau'_t)$ for a $\Sigma_1$-transversal path $\tau'_t$  $\Sigma_1$-homotopic to $\tau_t$. Indeed, a homotopy of $\tau_t$ within $\Sigma_1$-transversal paths does not change the sequence of Reidemeister moves.

On the other hand, there is a lifting functor $\Psi\colon \mathfrak D_s\to\mathfrak T_s$.

For a tangle diagram $D$, let $\Psi(D)=\lambda(D)$ be its canonical lifting. For a diagram isotopy $f\colon D\to D'$, define $\Psi(f)=\tau_t$ where $\tau_t=\lambda(f_t(D))$.

Let $f\colon D\to D'$ be a Reidemeister move $\Omega_i(B,\phi)$. For the move $\Omega_i$ fix a standard isotopy $\rho_i$ in $\mathbb D^2\times[0,1]$ that realize the move. Define $\Psi(f)$ as the isotopy which is constant and coincides with $\lambda(D)$ in $(F\setminus B)\times [0,1]$ and coincides with $(\phi\times id)^{-1}\circ\rho_i\circ(\phi\times id)$ in $B\times[0,1]$.

Then $\Phi\circ\Psi=id_{\mathfrak D_s}$. On the other hand, for any tangle $T$ there is a ``vertical'' isotopy $\eta_T$ from $T$ to $\Psi\circ\Phi(T)$ which does not change the projection $p(T)$. For any tangle path $\tau_t$ the composition of isotopies $\eta_{\tau_0}^{-1}\tau_t\eta_{\tau_1}$ is $\Sigma_1$-homotopic to $\Psi\circ\Phi(\tau_t)$. Thus, the isotopies $\eta_T$ establish a natural isomorphism between $id_{\mathfrak T_s}$ and $\Psi\circ\Phi$. Hence, the categories $\mathfrak T_s$ and $\mathfrak D_s$ are equivalent.
\end{proof}

\begin{remark}\label{rem:diagram_dagger_category}
    The categories $\mathfrak T_s$ and $\mathfrak D_s$ are dagger categories~\cite{HV}. In $\mathfrak T_s$, the dagger involution of an isotopy $\tau_t\in\mathrm{Mor}_{\mathfrak T_s}(T_0,T_1)$ is $(\tau_t)^\dag=\tau_{1-t}$. In $\mathfrak D_s$, for a diagram isotopy $f=\{D_t\}$ one has $f^\dag=\{D_{1-t}\}$, and for a Reidemeister move $f=\Omega_i(B,\phi)$ one has $f^\dag=\Omega_i^\dag(B,\phi)$.

    Note that the morphisms $\Omega_i(B,\phi)$ and $\Omega_i^\dag(B,\phi)$ are not inverse in $\mathfrak D_s$. Thus, $\mathfrak T_s$ and $\mathfrak D_s$ are not groupoids.
\end{remark}

Now we define homotopy tangle and diagram categories.

\begin{definition}\label{def:homotopy_tangle_category}
    Consider the category $\mathfrak{T}_h$ such that
\begin{itemize}
    \item $Ob(\mathfrak T_h)=\Sigma_0$;
    \item for tangles $T_0,T_1\in\Sigma_0$, the morphism set $\mathrm{Mor}_{\mathfrak T_h}(T_0,T_1)$ is the set of equivalence classes of isotopies from $T_0$ to $T_1$ modulo homotopy.
\end{itemize}
The category $\mathfrak{T}_h$ is called the \emph{homotopy category of tangles}.
\end{definition}

\begin{remark}\label{rem:homotopy_tangle_category}
    By Proposition~\ref{prop:sigma_transversality}, the morphism set $\mathrm{Mor}_{\mathfrak T_h}(T_0,T_1)$ between tangles $T_0$ and $T_1$ coincides with the set of $\Sigma_1$-transversal paths from $T_0$ to $T_1$ considered up to $\Sigma_2$-transversal homotopies.
\end{remark}

There is a natural projection morphism $\pi\colon\mathfrak T_s\to\mathfrak T_h$ between the strict and the homotopy categories of diagrams.

\begin{figure}[p]
\centering
  \includegraphics[width=0.5\textwidth]{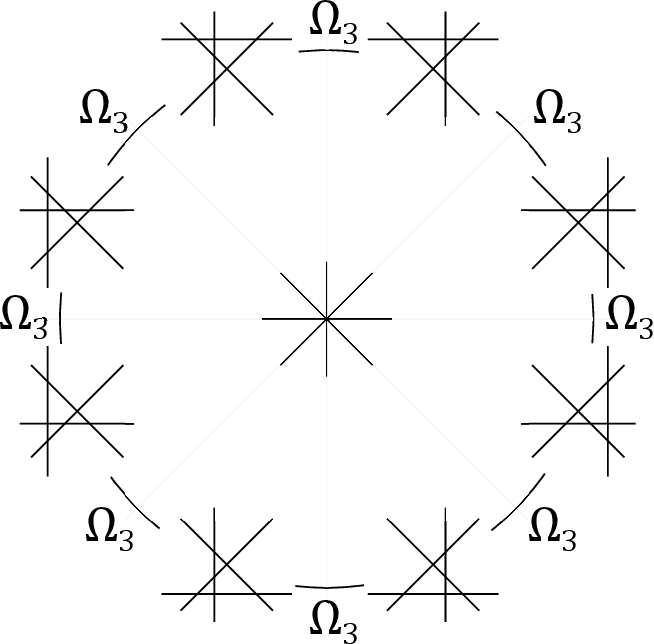}
  \caption{Quadruple point resolution}\label{pic:sing_quadruple_perturbation}
\end{figure}

\begin{figure}[p]
\centering
  \includegraphics[width=0.4\textwidth]{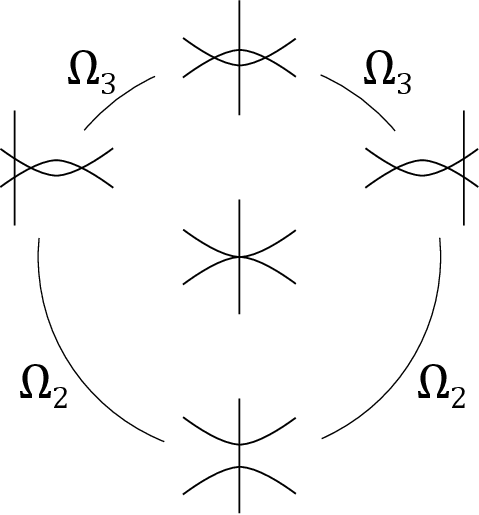}\qquad \includegraphics[width=0.35\textwidth]{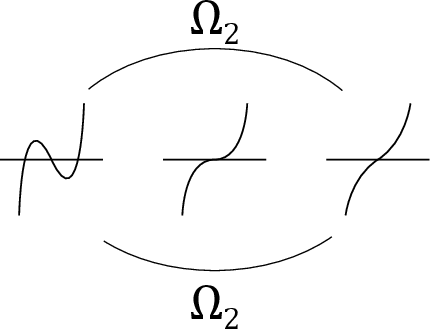}
  \caption{Tangent triple point and cubic tangency resolution}\label{pic:sing_tangent_triple_point_perturbation}
\end{figure}

\begin{figure}[p]
\centering
  \includegraphics[width=0.4\textwidth]{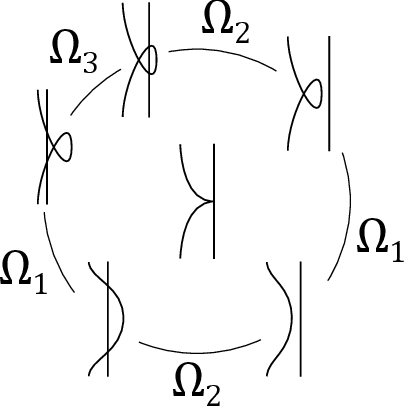}\qquad   \includegraphics[width=0.35\textwidth]{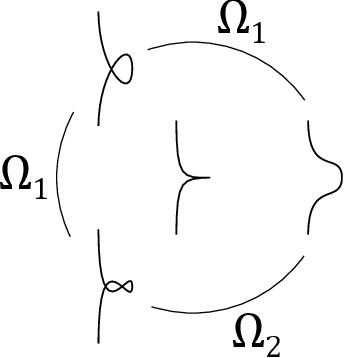}
  \caption{Intersected cusp and ramphoidal cusp resolution}\label{pic:sing_intersected_cusp_perturbation}
\end{figure}

\begin{definition}\label{def:homotopy_diagram_category}
    Consider the category $\mathfrak{D}_h$ such that
\begin{itemize}
    \item $Ob(\mathfrak D_h)$ are tangle diagrams;
    \item for tangle diagrams $D,D'$, the morphism set $\mathrm{Mor}_{\mathfrak D_h}(D,D')$ is the set of equivalence classes of the morphisms $\mathrm{Mor}_{\mathfrak D_s}(D,D')$ modulo the relations:
    \begin{itemize}
        \item contraction of mutually inverse Reidemeister moves
\[
\Omega_i(B,\phi)\Omega_i^\dag(B,\phi)=id,\quad \Omega_i^\dag(B,\phi)\Omega_i(B,\phi)=id,\quad i=1,2,3;
\]
        \item commutativity of distinct Reidemeister moves
\[
\Omega_i(B,\phi)\Omega_j(B',\phi')=\Omega_j(B',\phi')\Omega_i(B,\phi),\quad B\cap B'=
\emptyset;
\]
        \item relations induced by the resolutions of singularities of codimension $2$ (Fig.~\ref{pic:sing_quadruple_perturbation},\ref{pic:sing_tangent_triple_point_perturbation},\ref{pic:sing_intersected_cusp_perturbation}).
    \end{itemize}
\end{itemize}
The category $\mathfrak{D}_h$ is called the \emph{homotopy category of tangles diagrams}.
\end{definition}

\begin{theorem}\label{thm:homotopical_category_equivalence}
    The categories $\mathfrak T_h$ and $\mathfrak D_h$ are equivalent.
\end{theorem}

\begin{proof}
    We need to check that functors $\Phi$ and $\Psi$ induce well-defined functors between $\mathfrak T_h$ and $\mathfrak D_h$.

    Let $\tau_t$ and $\tau'_t$ be homotopical $\Sigma_1$-transversal paths of tangles. By Proposition~\ref{prop:sigma_transversality} there is a $\Sigma_2$-transversal homotopy $\Gamma=(\tau^s)_{s\in[0,1]}$, $\tau^0=\tau$ and $\tau^1=\tau'$. Split $[0,1]$ into intervals $[s_{k-1},s_{k}]$, $0=s_0<s_1<\cdots<s_n=1$, such that the homotopy $\Gamma$ has at most one exceptional point $s$ in $[s_{k-1},s_{k}]$.

    If the interval $[s_{k-1},s_{k}]$ has no exceptional point then $\tau^{s_{k-1}}$ and $\tau^{s_k}$ are homotopic through $\Sigma_1$-transversal paths. Hence, $\Phi(\tau^{s_{k-1}})=\Phi(\tau^{s_k})$.

    If the interval $[s_{k-1},s_{k}]$ contains an exceptional point $s^*$ then $\tau^{s^*}$ has either a non transversal intersection with $\Sigma_1$ or a singular point of codimension $2$. In the first case $\Phi(\tau^{s_{k-1}})$ and $\Phi(\tau^{s_k})$ differ by a relation $\Omega_i(B,\phi)\Omega_i^\dag(B,\phi)=id$ or $\Omega_i^\dag(B,\phi)\Omega_i(B,\phi)=id$. In the second case, $\Phi(\tau^{s_{k-1}})$ and $\Phi(\tau^{s_k})$ differ by a relation induced by resolution of the singular point of codimension $2$. In both cases, $\Phi(\tau^{s_{k-1}})=\Phi(\tau^{s_k})$ as morphisms in $\mathfrak D_h$.

    Thus, $\Phi(\tau_t)=\Phi(\tau'_t)$ in $\mathfrak D_h$, and the functor $\Phi\colon\mathfrak T_h\to\mathfrak D_h$ is well-defined.

    Analogously, the lifting functor $\Psi\colon\mathfrak D_h\to\mathfrak T_h$ is well-defined. The functors $\Phi$ and $\Psi$ establish equivalence between the categories $\mathfrak T_h$ and $\mathfrak D_h$.
\end{proof}

\begin{remark}
    1. The connected components of the categories $\mathfrak D_s$, $\mathfrak D_h$ correspond to isotopy classes of tangles. For a tangle $T$, the connected component of the tangle categories containing $T$ will be denoted $\mathfrak D_s(T)$, $\mathfrak T_s(T)$, $\mathfrak D_h(T)$ and $\mathfrak T_h(T)$. Note that these subcategories are strongly connected, i.e. for any two objects $c, c'$ of a category there exists morphism $f\colon c\to c'$.

    2. The categories $\mathfrak D_h$ and $\mathfrak T_h$ are groupoids, and $\mathfrak T_h$ is a full subcategory of the fundamental groupoid $\Pi_1(\Sigma)$ of the space of tangles $\Sigma$. In some sense, we can consider the groupoid $\Pi_1(\Sigma)$ as a proper category of tangle diagrams.
\end{remark}



\subsection{Diagram elements}\label{subsect:diagram_elements}

Let us define sets of elements of a tangle diagram.

\begin{definition}\label{def:diagram_elements}
    Let $D\colon G\to F$ be a tangle diagram. Then
    \begin{itemize}
        \item the set of crossings $\mathcal C(D)=\mathcal C(G)$ of the graph $G$ is the set of \emph{crossings} of the diagram $D$;
        \item the set of edges $\mathcal{SA}(D)=E(G)$ of $G$ is the set of \emph{semiarcs} of the diagram $D$;
        \item the set of connected components $\mathcal R(D)=\pi_0(F\setminus D)$ is the set of \emph{regions} of $D$;
        \item the set of classes of edges $E(G)$ modulo the equivalence relation which identifies opposite overcrossing edges at every crossing of $G$, is called the set of \emph{arcs} of the diagram $D$ and denoted by $\mathcal A(D)$.
    \end{itemize}
\end{definition}

Examples of diagram elements are presented in Fig.~\ref{pic:diagram_parts}.

Maps $D\mapsto \mathcal C(D)$,\dots, $D\mapsto \mathcal A(D)$ extend to functors $\mathcal C, \mathcal{SA}, \mathcal R, \mathcal A$ from the strict diagram category $\mathfrak D_s$ to the relation category $Rel$. A diagram isotopy $f\colon D\to D'$ induces bijections between the crossings, semiarcs, regions and arcs of the diagram $D$ and $D'$ which are taken for the morphisms $\mathcal C(f)$, $\mathcal{SA}(f)$, $\mathcal R(f)$ and  $\mathcal A(f)$ correspondingly.

Let $f\colon D\to D'$ be a Reidemeister move in a disc $B\subset F$. The morphism $f$ acts on elements of $D$ by the following rule: an element $x$ of the diagram $D$ maps to an element $x'$ of $D'$ if and only if $x$ and $x'$ have common points outside $B$.

This definition implies the following two rules:
\begin{enumerate}
    \item an element $x$ of $D$ which does not participate in the move $f$ maps to the unique corresponding element in $D'$;
    \item inner elements of the move $f$ map to nothing.
\end{enumerate}

\begin{figure}[h]
\centering
  \includegraphics[width=0.4\textwidth]{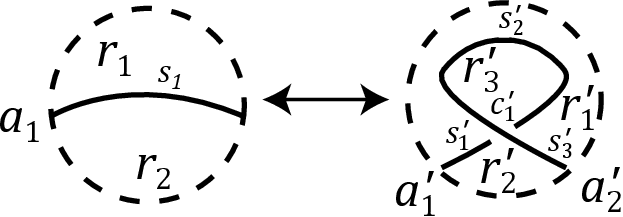}
  \caption{A first Reidemeister move}\label{pic:reidmove1elem}
\end{figure}

For example, for a first Reidemeisted move (Fig.~\ref{pic:reidmove1elem}) we have
\begin{align*}
    \mathcal A(\Omega_1): & a_1\mapsto \{a_1',a_2'\}, & \mathcal A(\Omega_1^\dag): & a_1'\mapsto a_1, a_2'\mapsto a_1,\\
    \mathcal {SA}(\Omega_1): & s_1\mapsto \{s_1',s_3'\}, & \mathcal {SA}(\Omega_1^\dag): & s_1'\mapsto s_1, s_2'\mapsto\emptyset, s_3'\mapsto s_1,\\
    \mathcal R(\Omega_1): & r_1\mapsto r_1', r_2\mapsto r_2', & \mathcal R(\Omega_1^\dag): & r_1'\mapsto r_1, r_2'\mapsto r_2, r_3'\mapsto\emptyset,\\
     & & \mathcal C(\Omega_1^\dag): & c_1'\mapsto\emptyset.
\end{align*}

\begin{figure}[h]
\centering
  \includegraphics[width=0.5\textwidth]{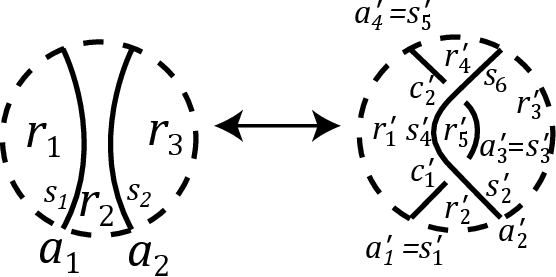}
  \caption{A second Reidemeister move}\label{pic:reidmove2elem}
\end{figure}

For a second Reidemeisted move (Fig.~\ref{pic:reidmove2elem}) we have
\begin{align*}
    \mathcal A(\Omega_2): & a_1\mapsto \{a_1',a_4'\}, a_2\mapsto a_2' & \mathcal A(\Omega_2^\dag): & a_1',a_4'\mapsto a_1, a_2'\mapsto a_2,\\
     & & & a_3'\mapsto\emptyset,\\
    \mathcal {SA}(\Omega_2): & s_1\mapsto \{s_1',s_5'\}, s_2\mapsto \{s_2',s_6'\}, & \mathcal {SA}(\Omega_2^\dag): & s_1', s_5'\mapsto s_1,\ s_2',s_6'\mapsto\emptyset,\\
    & & & s_3',s_4'\mapsto\emptyset,\\
    \mathcal R(\Omega_2): & r_1\mapsto r_1', r_2\mapsto \{r_2',r_4\}, r_3\mapsto r_3' & \mathcal R(\Omega_2^\dag): & r_1'\mapsto r_1, r_2',r_4'\mapsto r_2, r_3'\mapsto\emptyset,\\
     & & & r_5'\mapsto\emptyset,\\
     & & \mathcal C(\Omega_2^\dag): & c_1',c_2'\mapsto\emptyset.
\end{align*}

\begin{figure}[h]
\centering
  \includegraphics[width=0.7\textwidth]{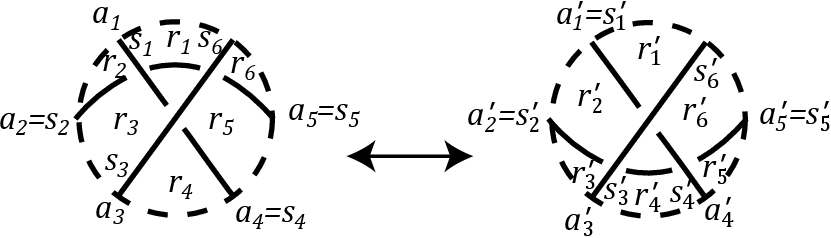}
  \caption{A third Reidemeister move}\label{pic:reidmove3elem}
\end{figure}

For a third Reidemeisted move (Fig.~\ref{pic:reidmove3elem}) we have
\begin{align*}
    \mathcal A(\Omega_3): & a_i\mapsto a_i',\ i=\overline{1,5}, & \mathcal A(\Omega_3^\dag): & a_i'\mapsto a_i,\ i=\overline{1,5},\\
    \mathcal {SA}(\Omega_3): & s_j\mapsto s_j',\ j=\overline{1,6}, & \mathcal {SA}(\Omega_3^\dag): & s_j'\mapsto s_j,\ j=\overline{1,6},\\
    \mathcal R(\Omega_3): & r_j\mapsto r_j',\ j=\overline{1,6}, & \mathcal R(\Omega_3^\dag): & r_j'\mapsto r_j,\ j=\overline{1,6}.
\end{align*}
Unlabelled elements have empty image.

Note that unlike parity or index axioms~\cite{IMN,CFGMX}, there is no correspondence between the crossings in a third Reidemeister move here.

\begin{remark}
    1. The functors $\mathcal A, \mathcal{SA},\mathcal R$ and $\mathcal C$ are dagger functors from the dagger category $\mathfrak D_s$ to the dagger category $Rel$.

    2. The functors $\mathcal A, \mathcal{SA},\mathcal R$ and $\mathcal C$ are not well defined on the homotopical category $\mathfrak D_h$. For example, $\Omega_1\Omega_1^\dag=id$ in  $\mathfrak D_h$ but $\mathcal A(\Omega_1)\mathcal A(\Omega_1^\dag)\ne id$.
\end{remark}


\subsection{The crossing and the arc from the combinatorial viewpoint}\label{subsect:diagram_element_equivalence}

The correspondences of diagram elements induced by Reidemeister moves allow to identify elements of different diagrams. Then one can think about an equivalence class of arcs (semiarcs, regions, crossings) as one and the same arc which bends, splits and merges under isotopies and Reidemeister moves. We will call the sets of equivalence classes of diagram elements the sets of \emph{arcs (semiarcs, regions, crossings) of the tangle}.

We define two types of equivalence relations on diagram elements: a strong and a weak one.

\begin{definition}\label{def:strong_weak_equivalence}
Let $\mathcal F\colon \mathfrak C\to Rel$ be a functor from a category $\mathfrak C$ to the category of relations $Rel$. Two elements $x_1,x_2\in \mathcal F(c)$, $c\in Ob(\mathfrak C)$, are called \emph{affined} if there exists morphism $f\colon c\to c'$ and $y\in \mathcal F(c')$ such that $y\in \mathcal F(f)(x_1)\cap \mathcal F(f)(x_2)$. In this case we write $x_1\approx_{\mathcal F,c} x_2$.

The \emph{strong equivalence relation} associated with the functor $\mathcal F\colon\mathfrak C\to Rel$ is a family of equivalence relations $\sim_{\mathcal F,c}^s$ on the sets $\mathcal F(c)$, $c\in Ob(\mathfrak C)$, generated by the rule:
\begin{itemize}
    \item for any $f\in Mor_{\mathfrak C}(c,c')$, $x_1,x_2\in \mathcal F(c)$ and $y_1\in \mathcal F(f)(x_1)$, $y_2\in \mathcal F(f)(x_2)$  $x_1\sim_{\mathcal F,c}^s x_2$ implies $y_1\sim_{\mathcal F,c'}^s y_2$.
\end{itemize}

The \emph{weak equivalence relation} associated with the functor $\mathcal F\colon\mathfrak C\to Rel$ is the equivalence relation $\sim^w_{\mathcal F}$ on the set $\bigsqcup_{c\in Ob(\mathfrak C)}\mathcal F(c)$ generated by the rule
\begin{itemize}
    \item for any morphism $f\in Mor_{\mathfrak C}(c,c')$,  $x\in \mathcal F(c)$ and  $y\in \mathcal F(f)(x)$ one has $x\sim^w_{\mathcal F} y$.
\end{itemize}

We will use notation $\approx$, $\sim^s$ or $\sim^w$ when the implicit functor $\mathcal F$ is clear.
\end{definition}

\begin{figure}[h]
\centering
  \includegraphics[width=0.6\textwidth]{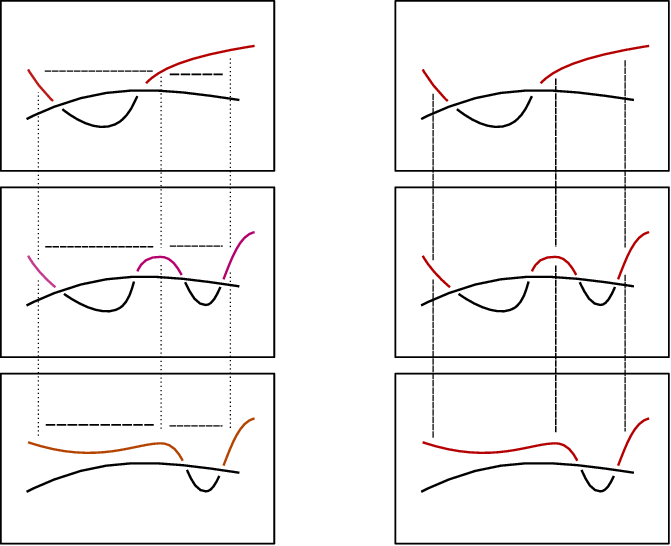}
  \caption{Equivalent arcs in the strong sense (left) and in the weak sense (right). The dashed lines are equivalence relations.}\label{pic:strong_weak_equiv}
\end{figure}

\begin{remark}\label{rem:equivalence_relations_properties}
    \begin{enumerate}
        \item $x_1\approx_{\mathcal F,c} x_2$ implies $x_1\sim^s_{\mathcal F,c} x_2$ implies $x_1\sim^w_{\mathcal F} x_2$;
        \item $\approx_{\mathcal F,c}$ is reflexive and symmetric but not transitive in general. In some sense, the relation $\sim^s_{\mathcal F,c}$ is a transitive closure of the affinity relation.
    \end{enumerate}
\end{remark}

\begin{remark}
    Let $\mathcal F$ be one of the functors $\mathcal A$, $\mathcal SA$ or $\mathcal R$. The relation $x_1\approx_{\mathcal F,D} x_2$ means that the arcs, semiarcs or regions $x_1,x_2\in \mathcal F(D)$ merge after one applies some sequence of Reidemeister moves.
\end{remark}

\begin{example}
    The trefoil diagram in Fig.~\ref{pic:trefoil} has three arcs. They are equivalent in the weak sense because the rotation of the diagram maps one arc to any other. On the other hand, the arcs are not equivalent in the strong sense because they can be distinguished by a quandle coloring.

\begin{figure}[h]
\centering
  \includegraphics[width=0.2\textwidth]{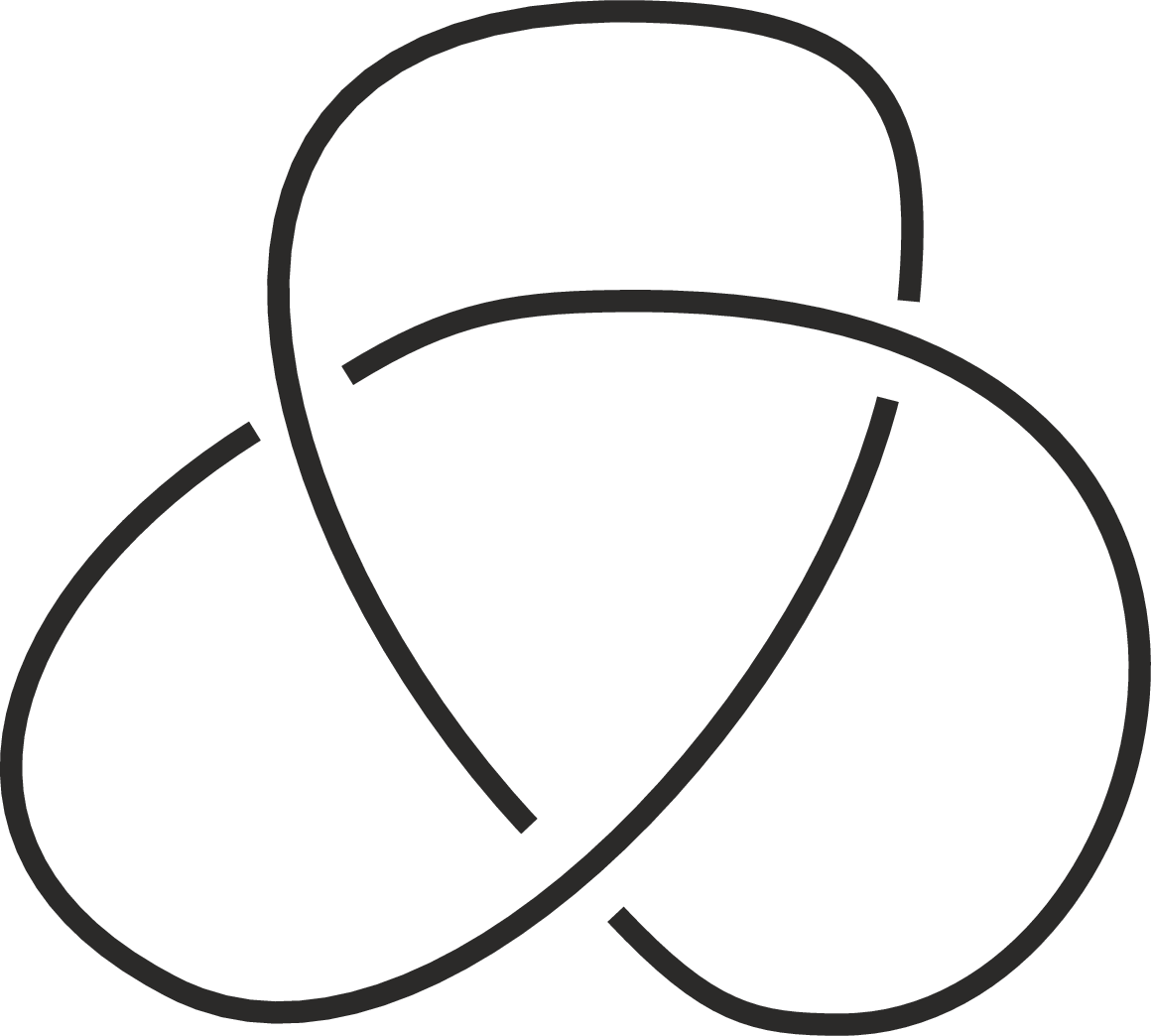}
  \caption{A trefoil diagram}\label{pic:trefoil}
\end{figure}
\end{example}

\begin{definition}\label{def:weak_tangle_elements}
    Let $\mathcal F$ be one of the functors $\mathcal A$, $\mathcal SA$, $\mathcal R$ or $\mathcal C$ considered on the diagrams of a tangle $T$. The set of equivalence classes of the relation $\sim^w_{\mathcal F}$ is called  \emph{the set of arcs (semiarcs, regions, crossings) of the tangle $T$ in the weak sense} and is denoted $\mathcal A^w(T)$ ($\mathcal {SA}^w(T)$, $\mathcal R^w(T)$, $\mathcal C^w(T)$).
\end{definition}

\begin{remark}\label{rem:explicit_implicit_arcs}
    Let $D$ be a diagram of a tangle $T$. The natural projection $\mathcal A(D)\to\mathcal A^w(T)$ from an arc to its weak equivalence class is not a surjection. Then regarding the diagram $D$, the arcs of the tangle $T$ split into a finite set of \emph{explicit arcs} (which correspond to the arcs in $\mathcal A(D)$) and an infinite set of \emph{implicit arcs}. The implicit arcs are not present directly in the diagram $D$ but one can reveal them using the topological description in Section~\ref{sect:arc}.

    The same remark can be applied to the other diagram elements.
\end{remark}

The strong equivalence classes $\bar {\mathcal F}(D)=\mathcal F(D)/\sim^s_{\mathcal F,D}$ describe only explicit arcs (semiarcs, regions, crossings) of the diagram $D$.  In order to reach the set of implicit arcs  (semiarcs, regions, crossings) of the tangle in the strong sense, we need to introduce the notion of a coinvariant of the diagram category.

\subsection{Invariants and coinvariants}\label{subsect:category_invariants}

Let us give a general definition of invariant and coinvariant (cf.~\cite{Armstrong}).

\begin{definition}\label{def:category_invariant}
    Let $\mathfrak C$ be a category. An \emph{invariant} of the category $\mathfrak C$ with values in a category $\mathfrak D$ is a functor $\mathcal F\colon\mathfrak C \to\mathfrak D$ such that for any morphism $f\in Mor_{\mathfrak C}(c,c')$ $\mathcal F(f)=id_{\mathcal F(c)}$.

    A \emph{coinvariant} of the category $\mathfrak C$ with values in a category $\mathfrak D$ is a functor $\mathcal F\colon\mathfrak C \to\mathfrak D$ such that for any morphism $f\in Mor_{\mathfrak C}(c,c')$ the morphism $\mathcal F(f)$ is an isomorphism.

\end{definition}

\begin{remark}
    1. If $\mathfrak C$ is a connected category then an invariant $\mathcal F$ is a constant functor to $\mathfrak D$.

    2. A coinvariant $\mathcal F\colon\mathfrak C\to\mathfrak D$ is a functor from $\mathfrak C$ to the core $Core(\mathfrak D)$ of the category $\mathfrak D$.
\end{remark}

\begin{example}\label{exa:knot_coinvariant}
    1. Let $\mathfrak C$ be a small category. Then any invariant $\mathcal F\colon\mathfrak D_s\to\mathfrak C$ (or $\mathcal F\colon\mathfrak D_h\to\mathfrak C$) is a knot invariant in the usual sense with values in $Ob(\mathfrak C)$.

    2. Let $Q$ be a quandle. Then the correspondence $D\mapsto Col_Q(D)$ which assigns the set of quandle colorings of arcs to a knot diagram, is a coinvariant $Col_Q\colon\mathfrak D_s\to Set$.
\end{example}

Next we define a notion of single-valued (co)invariants for functors into the relation category $Rel$.

\begin{definition}\label{def:single_valued_natural_transformation}
    Let $\mathcal F,\mathcal G\colon\mathfrak C\to Rel$ be functors. A \emph{single-valued natural transformation} $h\colon \mathcal F\Rightarrow \mathcal G$ is a family of functions $h(c)\colon \mathcal F(c)\to \mathcal G(c)$ such that for any morphism $f\colon c\to c'$
\[
h(c')\circ \mathcal F(f)\subset \mathcal G(f)\circ h(c).
\]
\end{definition}

\begin{remark}
We use a weaker condition $h(c')\circ \mathcal F(f)\subset \mathcal G(f)\circ h(c)$ instead of $h(c')\circ \mathcal F(f)= \mathcal G(f)\circ h(c)$ in order to embrace situations when there are several elements which have the same value of an invariant, and one of the elements vanishes after a morphism is applied.

Consider the following example. Let $\mathcal F(c)=\{x_1,x_2\}$, $\mathcal F(c')=\{y\}$, $\mathcal G(c)=\{u\}$, $\mathcal G(c')=\{v\}$. Assume also that $\mathcal F(f)\colon x_1\mapsto y, x_2\mapsto\emptyset$ and $\mathcal G(f)(u)=v$. Then the functions $h(c)(x_i)=u$ and $h(c')(y)=v$ form a single-valued natural transformation. On the other hand, the functions $h(c)$ and $h(c')$ are not a natural transformation in the category $Rel$ in the usual sense, because $h(c')\circ \mathcal F(f)\ne \mathcal G(f)\circ h(c)$.
\end{remark}

%


\begin{definition}\label{def:single_valued_invariant}
    Let $\mathcal F\colon\mathfrak C\to Rel$ be a functor. A \emph{single-valued (co)invariant of the functor} $\mathcal F$ is a pair $(\mathcal G,h)$ where $\mathcal G$ is a (co)invariant of the category $\mathfrak C$ with values in $Rel$ and $h\colon \mathcal F\Rightarrow \mathcal G$ is a single-valued natural transformation.
\end{definition}

\begin{definition}\label{def:universal_invariant}
    A single-valued (co)invariant $(\mathcal G,h)$ of a functor $\mathcal F$ is called \emph{universal} if for any single-valued (co)invariant $(\mathcal G',h')$ there exists a unique single-valued natural transformation $\phi\colon \mathcal G\Rightarrow \mathcal G'$ such that $h'=\phi\circ h$.

    We will use notation $(\hat {\mathcal F}, \hat h)$ and $(\tilde {\mathcal F},\tilde h)$
    for the universal invariant and coinvariant 
    of the functor $F$ correspondingly.
\end{definition}

Recall that a \emph{strongly connected category} is a category $\mathfrak C$ such that for any objects $c,c'\in Ob(\mathfrak C)$ there exists a morphism $f\colon c\to c'$.

\begin{proposition}
    Let $\mathfrak C$ be a strongly connected category, $\mathcal F\colon\mathfrak C\to Rel$ a functor, $(\hat {\mathcal F},\hat h)$ its universal invariant, and $(\tilde {\mathcal F},\tilde h)$ its universal coinvariant. For an object $c\in Ob(\mathfrak C)$ consider the subgroup
$Mon_c(\mathcal F)\subset \mathrm{Aut}_{\mathfrak D}(\tilde {\mathcal F}(c),\tilde {\mathcal F}(c))$ generated by the automorphisms $\tilde{\mathcal F}(f)$, $f\in\mathrm{Mor}_{\mathfrak C}(c,c)$.
Then
\begin{enumerate}
    \item for all $c\in Ob(\mathfrak C)$ the groups $Mon_c(\mathcal F)$ are isomorphic;
    \item for any  $c\in Ob(\mathfrak C)$ $\hat {\mathcal F}(c)=\tilde {\mathcal F}(c)/Mon_c(\mathcal F)$.
\end{enumerate}
\end{proposition}

\begin{proof}
    1. Let $c,c'\in Ob(\mathfrak C)$. By connectivity there are morphisms $f\colon c\to c'$ and $f'\colon c'\to c$. Since $\tilde {\mathcal F}(f), \tilde {\mathcal F}(f')$ are isomorphisms, the map
\[
\tilde {\mathcal F}(g)\mapsto  \tilde {\mathcal F}(f)\tilde {\mathcal F}(g)\tilde {\mathcal F}(f)^{-1}=\tilde{\mathcal F}(fgf')\tilde{\mathcal F}(ff')^{-1}
\]
establishes an isomorphism of the groups $Mon_c(\mathcal F)$ and $Mon_{c'}(\mathcal F)$.

 2. Fix an object $c\in Ob(\mathfrak C)$. Define an invariant $(\mathcal G,h)$ of the functor $\mathcal F$ as follows. For any object $c'$ set $\mathcal G(c')=\tilde {\mathcal F}(c)/Mon_c(\mathcal F)$, and for any morphism $f\colon c'\to c''$ set $\mathcal G(f)=id_{\mathcal G(c)}$. The natural map $h(c')$, $c'\in Ob(\mathfrak C)$, is the composition of the map $\tilde {\mathcal F}(f)\colon \tilde {\mathcal F}(c')\to \tilde {\mathcal F}(c)$ where $f$ is a morphism from $c'$ to $c$, and the projection $\tilde {\mathcal F}(c)\to \mathcal G(c)$. By definition, the map $h(c)$ does not depend on the choice of $f$.

 Let $(\mathcal G',h')$ be a (single-valued) invariant of $\mathcal F$. Then it is a coinvariant of $\mathcal F$. Hence, there exists a unique natural transformation $\phi\colon \tilde{\mathcal F}\Rightarrow \mathcal G'$. For any $f\in\mathrm{Mor}_{\mathfrak C}(c,c)$ we have $\phi(c)=\phi(c)\circ\tilde {\mathcal F}(f)$. Hence, $\phi$ induces a map $\hat\phi$ from $\mathcal G(c)=\tilde {\mathcal F}(c)/Mon_c(\mathcal F)$ to $\mathcal G'(c)$. The map $\hat\phi$ defines a natural transformation from $\mathcal G$ to $\mathcal G'$. Thus, $(\mathcal G,h)$ is the universal invariant of $\mathcal F$.
 \end{proof}

 The group $Mon_c(\mathcal F)$ is called the \emph{monodromy group} of the functor $\mathcal F$. For different $c$ the groups $Mon_c(\mathcal F)$ are isomorphic, so we will use also the notation $Mon(\mathcal F)$.

\subsection{Universal invariants and coinvariants}

Let $\mathcal F\colon\mathfrak C\to Rel$ be a functor into the category of relations. Let us describe a construction for the universal (co)invariant 
of $\mathcal F$.

\subsubsection{Universal invariant}

For a connected category $\mathfrak C$, consider the set $\hat {\mathcal F}(\mathfrak C)=\bigsqcup_{c\in Ob(\mathfrak C)}\mathcal F(c) / \sim^w_{\mathcal F}$ 
and the functor $\hat {\mathcal F}_{\mathfrak C}$ such that $\hat {\mathcal F}_{\mathfrak C}(c)=\hat {\mathcal F}(\mathfrak C)$ for any object $c\in Ob(\mathfrak C)$ and $\hat {\mathcal F}_{\mathfrak C}(f)=id_{\hat {\mathcal F}(\mathfrak C)}$ for any morphism $f\in\mathrm{Mor}(\mathfrak C)$. The map which assigns to an element $x\in \mathcal F(c)$, $c\in Ob(\mathfrak C)$, its equivalence class in $\hat {\mathcal F}(\mathfrak C)$, determines a single-valued natural transformation $\hat h_{\mathfrak C}\colon \mathcal F\Rightarrow\hat {\mathcal F}_{\mathfrak C}$.

In general case, define the functor $\hat {\mathcal F}_{\mathfrak C}$ and the single-valued natural transformation $\hat h_{\mathfrak C}$ by the condition $\hat {\mathcal F}_{\mathfrak C}|_{\mathfrak C'}=\hat {\mathcal F}_{\mathfrak C'}$, $\hat h_{\mathfrak C}|_{\mathfrak C'}=\hat h_{\mathfrak C'}$ for any connected component $\mathfrak C'\in\pi_0(\mathfrak C)$ of the category $\mathfrak C$.

\begin{proposition}\label{prop:single_valued_universal_invariant}
   The pair $(\hat {\mathcal F}_{\mathfrak C}, \hat h_{\mathfrak C})$ is the universal invariant of the functor $\mathcal F$.
\end{proposition}

\begin{proof}
   By definition, $(\hat {\mathcal F}_{\mathfrak C}, \hat h_{\mathfrak C})$ is a single-valued invariant of $\mathcal F$. Let us prove the invariant is universal.

    Let $(\mathcal G,h)$ be a single-valued invariant of $\mathcal F$. Let us show that there exists a unique single-valued natural transformation $\phi\colon\hat {\mathcal F}_{\mathfrak C}\to \mathcal G$ such that $h=\phi\circ\hat h_{\mathfrak C}$.

   Assume first the category $\mathfrak C$ is connected. Let $f\in\mathrm{Mor}_{\mathfrak C}(c,c')$, $x\in \mathcal F(c)$ and $x'\in \mathcal F(f)(x)$. Since $h(c)=h(c')\circ\mathcal  F(f)$, $h(c)(x)=h(c')(\mathcal F(f)(x))\ni h(c')(x')$. Hence, $h(c)(x)=h(c')(x')$ whenever $x'\in \mathcal F(f)(x)$ for some morphism $f$.

    Thus, there is a unique well-defined function $\phi\colon\hat {\mathcal F}(\mathfrak C)\to \mathcal G(\mathfrak C)$ which maps the equivalence class $[x]$ of an element $x\in \mathcal F(c)$, $c\in Ob(\mathfrak C)$, to the element $h(c)(x)\in \mathcal G(c)=\mathcal G(\mathfrak C)$. The map $\phi$ is the required single-valued natural transformation between $\hat {\mathcal F}_{\mathfrak C}$ and the invariant $\mathcal G$.

    In general case, we define the transformation $\phi$ on each connected component of the category $\mathfrak C$ as shown above.
\end{proof}

\subsubsection{Universal local invariant}

In order to define the universal coinvariant we consider a supplement notion.

\begin{definition}\label{def:partial_bijection}
     A relation $f\subset X\times Y$ is called a \emph{partial bijection} if any image $f(x)$, $x\in X$, and preimage $f^\dag(y)$, $y\in Y$, has at most one element. In other words, $f$ is a bijection between the domain and the image.
\end{definition}

\begin{definition}\label{def:local_invariant}
  Let $\mathcal F\colon\mathfrak C\to Rel$ be a functor. A \emph{single-valued local invariant} of the functor $\mathcal F$ is a pair $(\mathcal G,h)$ where $\mathcal G\colon\mathfrak C\to Rel$ is a functor such that any morphism $\mathcal G(f)$, $f\in\mathrm{Mor}(\mathfrak C)$, is a partial bijection, and $h\colon \mathcal F\Rightarrow \mathcal G$ is a single-valued natural transformation.

  A single-valued local invariant $(\mathcal G,h)$ is \emph{universal} if for any single-valued local invariant $(\mathcal G',h')$ there is a unique single-valued natural transformation $\phi\colon \mathcal G\Rightarrow \mathcal G'$  such that $h'=\phi\circ h$.
\end{definition}

A local invariant can be thought of as invariant of the elements of sets $\mathcal F(c)$ with values in local coefficient sets $\mathcal G(c)$.

Let us describe the universal single-valued local invariant of a functor $\mathcal F\colon \mathfrak D_s\to Rel$.

For any $c\in Ob(\mathfrak D_s)$ denote $\bar {\mathcal F}(c)=\mathcal F(c)/ \sim_{\mathcal F,c}^s$. For a morphism $f\colon c\to c'$ define the relation $\bar {\mathcal F}^0(f)$ by the rule $\bar y\in\bar {\mathcal F}(f)(\bar x)$ for elements $\bar x\in\bar {\mathcal F}(c)$, $\bar y\in\bar {\mathcal F}(c')$, if and only if there exist representatives $x\in \mathcal F(c)$ of $\bar x$ and $y\in \mathcal F(c')$ of $\bar y$ such that $y\in \mathcal F(f)(x)$. Then define the completion $\bar {\mathcal F}(f)$ of the map $\bar {\mathcal F}^0(f)$ by the formula
\[
\bar {\mathcal F}(f)=\bigcup_{f_1\circ f_2\circ\cdots\circ f_n=f}\bar {\mathcal F}^0(f_1)\circ\bar {\mathcal F}^0(f_2)\circ\cdots\circ \bar {\mathcal F}^0(f_n)\in\mathrm{Mor}_{Rel}(\bar {\mathcal F}(c),\bar {\mathcal F}(c')).
\]
Denote the projections from $\mathcal F(c)$ to $\bar {\mathcal F}(c)$ by $\bar h(c)$.

The following example shows that $\bar{\mathcal F}$ and $\bar{\mathcal F}^0$ can be different.

\begin{example}
Let $\mathcal F = \mathcal A$. Consider a sequence of decrieasing and increasing second Reidemeister moves in Fig.~\ref{pic:morphism_closure}. Denote $f=f_2\circ f_1$ and $g=g_2\circ g_1$. Then $\mathcal A(f)(x)=\emptyset$, $\mathcal A(f)(y)=y$, $\mathcal A(g)(x)=x$ and $\mathcal A(g)(y)=\emptyset$. Hence, $\mathcal A(g\circ f)(x)=\mathcal A(g\circ f)(y)=\emptyset$.
\begin{figure}[h]
\centering
  \includegraphics[width=0.8\textwidth]{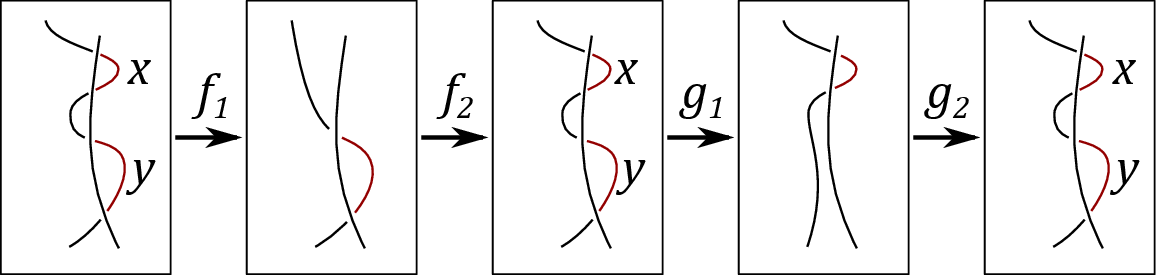}
  \caption{A sequence of second Reidemeister moves}\label{pic:morphism_closure}
\end{figure}

Since $x\sim^s y$, we have $\bar{\mathcal A}^0(g)\circ\bar{\mathcal A}^0(f)(\bar x)=\bar x$, hence, $\bar{\mathcal A}(g\circ f)(\bar x)=\bar x$, but $\bar{\mathcal A}^0(g\circ f)(\bar x)=\emptyset$.
\end{example}


\begin{proposition}\label{prop:universal_local_invariant}
    The pair $(\bar {\mathcal F},\bar h)$ is the universal single-valued local invariant of $\mathcal F$.
\end{proposition}

\begin{proof}
    Let is show first that $\bar {\mathcal F}$ is a well-defined functor. For any $c\in Ob(\mathfrak D_s)$ we have $\bar {\mathcal F}(id_c)=id_{\bar {\mathcal F}(c)}$ by reflexivity of $\sim^s_c$.

    The associativity follows from the fact that any morphism in $\mathfrak D_s$ has a unique presentation as a sequence of Reidemeister moves.

    By definition of $\sim^s$, for an elementary morphism the map $\bar {\mathcal F}(f)$ is injective. Hence, $\mathcal F(f)$ is a partial bijection.  For a general morphism $f$ the map $\bar {\mathcal F}(f)$ is a partial bijection as a composition of partial bijections.

    For an elementary morphism $f\colon c\to c'$ the relation $\bar h(c')\circ \mathcal F(f)\subset \mathcal G(f)\circ \bar h(c)$ follows from the definition of $\bar {\mathcal F}(f)$. For a general morphism $f$ the relation can be proved by induction of the length of a sequence of elementary morphisms representing $f$. Thus, $(\bar {\mathcal F}, \bar h)$ is a local invariant of $\mathcal F$.

    Let $(\mathcal G,h)$ be a local invariant of $\mathcal F$. For any $c\in Ob(\mathfrak D_s)$ and $x,y\in \mathcal F(c)$ denote $x\sim_{\mathcal G,c} y$ if $h(c)(x)=h(c)(y)$. Then $\sim_{\mathcal G,c}$ is an equivalence relation on $\mathcal F(c)$.

    Let $f\in Mor_{\mathfrak C}(c,c')$, $x_1,x_2\in \mathcal F(c)$ and $y_1\in \mathcal F(f)(x_1)$, $y_2\in \mathcal F(f)(x_2)$. Assume that  $x_1\sim_{\mathcal G,c} x_2$, i.e. $h(c)(x_1)=h(c)(x_2)$. Then
\[
h(c')(y_1)=h(c')\mathcal F(f)(x_1)=\mathcal G(f)(h(c)(x_1))=\mathcal G(f)(h(c)(x_2))=h(c')(y_2).
\]
Hence, $y_1\sim_{\mathcal G,c'} y_2$.

Thus, $\sim^s_{\mathcal F}$ is stronger than $\sim_{\mathcal G}$. Then the map $h$ induces well-defined functions $\phi(c)\colon \bar {\mathcal F}(c)\to \mathcal G(c)$. For an elementary morphism $f\colon c\to c'$ we have $\phi(c')\circ\bar {\mathcal F}(f)\subset \mathcal G(f)\circ \phi(c)$ by definition of $\bar {\mathcal F}(f)$. Hence, the relation $\phi(c')\circ\bar {\mathcal F}(f)\subset \mathcal G(f)\circ \phi(c)$ holds for any morphism $f$. Thus, $(\bar {\mathcal F},\hat h)$ is a universal local invariant.
\end{proof}

Since any coinvariant is a local invariant, we have the following statement.

\begin{corollary}
    Let $(\mathcal G,h)$ be a single-valued coinvariant of a functor $\mathcal F\colon\mathfrak C\to Rel$. Then for any $c\in Ob(\mathfrak C)$ and $x,x'\in \mathcal F(c)$  $x\sim^s_{\mathcal F,c}x'$ implies $h_c(x)=h_c(x')$.
\end{corollary}


\subsubsection{Universal coinvariant}\label{subsect:universal_coinvariant}
Let us describe a construction of the universal coinvariant of a functor $\mathcal F\colon\mathfrak C\to Rel$.

Let $\mathfrak C$ be a small category. For any morphism $f\colon c\to c'$ consider its formal opposite arrow by $\bar f\colon c'\to c$. The \emph{zigzag category} of the category $\mathfrak C$ is the category $Z\mathfrak C$ whose objects are $Ob(\mathfrak C)$ and morphisms are compatible sequences of $f,\bar f$, $f\in\mathrm{Mor}(\mathfrak C)$, modulo relations
\[
\overline{id}_c=id_c,\quad fg = f\circ g,\quad \bar g\bar f = \overline{f\circ g},
\]
where $c\in Ob(\mathfrak C)$ and $f\colon c\to c'$, $g\colon c'\to c''$ are morphisms in $\mathfrak C$. The \emph{localization} $\mathfrak C[\mathfrak C^{-1}]$ of the category $\mathfrak C$ is obtained from $Z\mathfrak C$ by imposing relations
\[
f\bar f=id_c,\quad \bar f f=id_{c'}
\]
for any morphism $f\colon c\to c'$ of $\mathfrak C$.

Define a coinvariant $(\tilde {\mathcal F},\tilde h)$ as follows. For $c\in Ob(\mathfrak C)$ consider the set
\[
\tilde {\mathcal F}(c)=\bigsqcup_{\mathrm{Mor}(Z\mathfrak C)\ni\psi\colon c'\to c} \mathcal F_\psi/\sim
\]
where $\mathcal F_\psi=\mathcal F(c')$ and the equivalence $\sim$ is generated by the relations:
\begin{itemize}
    \item for any morphisms $\psi\colon c'\to c$ in $Z\mathfrak C$ and $g\colon c''\to c'$ in $\mathfrak C$ and elements $x\in \mathcal F(c'')$ and $y\in \mathcal F(g)(x)$ one has $x_{\psi g}\sim y_{\psi}$;
    \item for any morphisms $\psi\colon c''\to c$ in $Z\mathfrak C$ and $g\colon c''\to c'$ in $\mathfrak C$ and elements $x\in \mathcal F(c'')$ and $y\in \mathcal F(g)(x)$ one has $x_{\psi }\sim y_{\psi\bar g}$;
    \item for any morphisms $\psi,\psi'\colon c'\to c$ in $Z\mathfrak C$ which coincide in the localization $\mathfrak C[\mathfrak C^{-1}]$, and any $x\in \mathcal F(c')$ one has $x_\psi\sim x_{\psi'}$.
\end{itemize}
Here $x_\psi$ denotes an element $x$ in the set $\mathcal F_\psi$.

For a morphism $f\colon c\to c'$ in $\mathfrak C$ define a map $\tilde {\mathcal F}(f)\colon \tilde {\mathcal F}(c)\to\tilde {\mathcal F}(c')$ by the formula $\tilde {\mathcal F}(f)(x_\psi)=x_{f\psi}$, where $\psi\in\mor_{Z\mathfrak C}(c'',c)$, $c''\in Ob(\mathfrak C)$, and $x\in \mathcal F_\psi$.

Finally, for any object $c\in Ob(\mathfrak C)$ define a function $\tilde h(c)\colon \mathcal F(c)\to\tilde {\mathcal F}(c)$ by the formula $h(c)(x)=x_{id_c}$.

\begin{proposition}
The pair $(\tilde {\mathcal F},\tilde h)$ is the universal coinvariant of $\mathcal F$.
\end{proposition}
\begin{proof}
    A direct check shows that $\tilde {\mathcal F}$ is a functor to $Rel$. For a morphism $f\colon c\to c'$ the map $\tilde {\mathcal F}(f)$ is a bijection: the inverse map is given by the formula $\tilde {\mathcal F}(f)^{-1}(x_\psi)=x_{\bar f\psi}$.

    Let $(\mathcal G,h)$ be a coinvariant. The rule $\bar f\mapsto \mathcal G(f)^\dag$, $f\in\mor(\mathfrak C)$, extends $\mathcal G$ to a functor $\mathcal G\colon Z\mathfrak C\to Rel$ (in fact to a functor $\mathcal G\colon \mathfrak C[\mathfrak C^{-1}]\to Rel$).

    For an object $c\in Ob(\mathfrak C)$ define a function $\phi_c\colon\tilde {\mathcal F}(c)\to \mathcal G(c)$ by the formula
\[
\phi_c(x_\psi)=\mathcal G(\psi)(h(x)),\quad \psi\in\mor(Z\mathfrak C), x\in \mathcal F_\psi.
\]
    The maps $\phi$ form a natural transformation $\tilde {\mathcal F}\Rightarrow \mathcal G$ such that $h=\phi\circ\tilde h$. Hence, $(\tilde {\mathcal F},\tilde h)$ is the universal coinvariant.
\end{proof}


\begin{definition}
    Let $\mathcal F$ be the functor $\mathcal{A}$ ($\mathcal{SA}$, $\mathcal{R}$ or $\mathcal{C}$) from $\mathfrak D_s$ to $Rel$, and $(\tilde F,\tilde h)$ its universal single-valued coinvariant. For a diagram $D\in Ob(\mathfrak D_s)$ of a tangle $T$, the set $\tilde {\mathcal F}(D)$ is called the \emph{set of arcs (semiarcs, regions or crossings) of the tangle $T$ in the strong sense at the diagram $D$} and denoted by $\mathcal A^s_D(T)$ ($\mathcal {SA}^s_D(T)$, $\mathcal R^s_D(T)$ or $\mathcal C^s_D(T)$).
\end{definition}

\subsection{\texorpdfstring{$h$}{h}-coinvariants}

Below we will need a further generalization of the notion of a coinvariant.

\begin{definition}\label{def:relative_functor_coinvariant}
    Let $\mathcal F\colon\mathfrak C\to Rel$ and $\mathcal P\colon \mathfrak C\to \mathfrak C'$ be functors. A \emph{(co)invariant of the functor $\mathcal F$ compatible with $\mathcal P$} is a (co)invariant $(\mathcal G,h)$ of $\mathcal F$ such that $\mathcal G=\mathcal G'\circ \mathcal P$ for some functor $\mathcal G'\colon \mathfrak C'\to Rel$.

    Analogously one defines single-valued and universal single-valued (co)invariants compatible with $\mathcal P$.
\end{definition}

\begin{proposition}\label{prop:universal_relative_functor_coinvariant}
    Let $\mathcal F\colon\mathfrak C\to Rel$ and $\mathcal P\colon \mathfrak C\to \mathfrak C'$ be functors, and $\mathfrak C'$ a groupoid. Let $(\mathcal G^u,h^u)$ be a coinvariant of $\mathcal F$ compatible with $\mathcal P$ such that
\begin{enumerate}
    \item for any $c\in Ob(\mathfrak C)$ and $y\in\mathcal G^u(c)$ there exist $c_1\in Ob(\mathfrak C)$, $f_1\colon c_1\to c$ and $x_1\in F(c_1)$ such that $y=G^u(f_1)(h^u(x_1))$;
    \item for any morphisms $f_i\colon c_1\to c$, $i=1,2$, in $\mathfrak C$ and elements $y\in \mathcal G^u(c)$ and $x_i\in \mathcal F(c_i)$ such that $y=G^u(f_i)(h^u(x_i))$, $i=1,2$, there exists a morphism $f_{12}\colon c_1\to c_2$ such that $\mathcal P(f_{12})=\mathcal P(f_2)^{-1}\circ\mathcal P(f_1)$ and $x_2\in \mathcal F(f_{12})(x_1)$.
\end{enumerate}
Then $(\mathcal G^u,h^u)$ is a universal coinvariant of $\mathcal F$ compatible with $\mathcal P$.
\end{proposition}

\begin{proof}
    Let $(\mathcal G,h)$ be a coinvariant of $\mathcal F$ compatible with $\mathcal P$. We construct a single-valued natural transformation $\phi\colon\mathcal G^u\Rightarrow\mathcal G$.

    Let $c\in Ob(\mathfrak C)$ and $y\in \mathcal G^u(c)$. By the first condition, there is a morphism $f_1\colon c_1\to c$ and an element $x_1\in\mathcal F(c_1)$ such that $y=\mathcal G^u(f_1)(h^u(x_1))$. Then we set $\phi(c)(y)=\mathcal G(f_1)(h(x_1))$.

    Let us show that $\phi(c)(y)$ does not depend on $f_1$ and $x_1$. Let $f_2\colon c_2\to c$ be a morphism and $x_2\in\mathcal F(c_2)$ an element such that $y=\mathcal G^u(f_2)(h^u(x_2))$. Denote $z_i=\mathcal G(f_i)(h(x_i))$, $i=1,2$. We need to check that $z_1=z_2$.

    By the second condition, there exists $f_{12}\colon c_1\to c_2$ such that $\mathcal P(f_{12})=\mathcal P(f_2)^{-1}\circ\mathcal P(f_1)$ and $x_2\in \mathcal F(f_{12})(x_1)$.
    Then
\[
\mathcal G(f_{1})=\mathcal G'(\mathcal P(f_{1}))=\mathcal G'(\mathcal P(f_{2})\circ\mathcal P(f_{12}))=\mathcal G'(\mathcal P(f_{2}))\circ\mathcal G'(\mathcal P(f_{12}))=\mathcal G(f_{2})\circ\mathcal G(f_{12})
\]
and
\begin{multline*}
z_1=\mathcal G(f_1)(h(x_1))=\mathcal G(f_{2})\circ\mathcal G(f_{12})(h(x_1))=\mathcal G(f_{2})(h(\mathcal G(f_{12})(x_1)))\\
=\mathcal G(f_{2})(h(x_2))=z_2.
\end{multline*}
The third and the fourth equalities follows from
\[
h(x_2)\subset h(\mathcal G(f_{12})(x_1))\subset \mathcal G(f_{12})(h(x_1))
\]
and the fact that $\mathcal G(f_{12})(h(x_1))$ is a one-element set.

Thus, the functions $\phi(c)\colon \mathcal G^u(c)\to\mathcal G(c)$, $c\in Ob(\mathfrak C)$, are well defined.

Let $f\colon c\to c'$ be a morphism in the category $\mathfrak C$ and  $z'\in\phi(c')\circ\mathcal G^u(f)\subset \mathcal G(c')$. Then there exists $y\in\mathcal G^u(c)$ such that $z'=\phi(c')(y')$ where $y'=\mathcal G^u(f)(y)$. By the first condition of the statement there exists a morphism $f_1\colon c_1\to c$ and $x_1\in\mathcal F(c_1)$ such that $y=\mathcal G^u(f_1)(h^u(x_1))$. Denote $z=\mathcal G(f_1)(h(x_1))$. Then $\phi(c)(y)=z$ and
\[
\mathcal G(f)(z)=\mathcal G(f)\circ\mathcal G(f_1)(h(x_1))=\mathcal G(f\circ f_1)(h(x_1))=\phi(c')(y')=z'
\]
since $y'=\mathcal G^u(f)(y)=\mathcal G^u(f)\circ\mathcal G^u(f_1)(h^u(x_1))=\mathcal G^u(f\circ f_1)(h^u(x_1))$.

Thus, $\phi(c')\circ\mathcal G^u(f)\subset \mathcal G(f)\circ \phi(c)$, and $\phi$ is a single-valued natural transformation. Hence, $(\mathcal G^u,h^u)$ is a universal coinvariant of the functor $\mathcal F$ compatible with the functor $\mathcal P$.
\end{proof}

We apply Definition~\ref{def:relative_functor_coinvariant} to the projection $\mathcal P\colon \mathfrak D_s\to\mathfrak D_h$ from the strict diagram category to the homotopical diagram category. In this case we refer to (co)invariants of a functor $\mathcal F\colon\mathfrak D_s\to Rel$ compatible with $\mathcal P$ as \emph{h-(co)invariants} of the functor $\mathcal F$.

\begin{definition}\label{def:h_strong_weak_equivalences}
Let $\mathcal F\colon\mathfrak D_s\to Rel$ be a functor. The \emph{homotopical strong (weak) equivalence relation} $\sim^{sh}_{\mathcal F,c}$ ($\sim^{wh}_{\mathcal F}$) associated with $\mathcal F$ is the equivalence obtained from $\sim^s_{\mathcal F,c}$ ($\sim ^w_{\mathcal F}$) by adding the relation: $x\sim y$ if $y\in \mathcal F(f)(x)$ for some morphism $f\colon c\to c$ which is equal to $id_c$ in $\mathfrak D_h$.
\end{definition}

\begin{proposition}\label{prop:h_weak_equivalence}
\begin{enumerate}
\item For any functors $\mathcal F\colon\mathfrak D_s\to Rel$ the relations $\sim^w_{\mathcal F}$ and $\sim^{wh}_{\mathcal F}$ coincide.

\item For any functor $\mathcal F$, any invariant $(\mathcal G,h)$ of $\mathcal F$ is an h-invariant of $\mathcal F$.

\item For any functor $\mathcal F$, the pair $(\hat{\mathcal G},\hat h)$  from Proposition~\ref{prop:single_valued_universal_invariant} is the universal h-invariant of $\mathcal F$.

\end{enumerate}
\end{proposition}
\begin{proof}
    The first statement is due to the fact the new relations of $\sim^{wh}$ follow from $\sim^w$.

    For an invariant $(\mathcal G,h)$ of $\mathcal F$, define a functor $\mathcal G'\colon \mathfrak D_h\to Rel$ by the formulas $\mathcal G'(c)=\mathcal G(c)$, $c\in Ob(\mathfrak D_h)$, and $\mathcal G'(f)=id_{\mathcal G(c)}$, $f\in Mor_{\mathfrak D_h}(c,c')$. Since the projection $\mathcal P\colon\mathfrak D_s\to\mathfrak D_h$ is bijective on the objects, the functor $\mathcal G'$ is well-defined. By definition $\mathcal G=\mathcal G'\circ\mathcal P$, hence, $\mathcal G$ is an $h$-invariant.

    The last statement of the proposition is an immediate consequence of the previous one.
\end{proof}

Now, let us consider the homotopical strong equivalence.

\begin{proposition}\label{prop:h_strong_equivalence}
\begin{enumerate}
\item For the functors $\mathcal A$, $\mathcal {SA}$, $\mathcal R$ the relations $\sim^s$ and $\sim^{sh}$ coincide.
\item For the crossing functor $\mathcal C$ the relation $\sim^s_{\mathcal C}$ is the equality, and $\sim^{sh}_{\mathcal C}$ is generated by the relation $c_1\sim c_2$ in Fig.~\ref{pic:homotopy_crossing_relation}.
\begin{figure}[h]
\centering
  \includegraphics[width=0.2\textwidth]{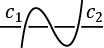}
  \caption{Homotopical strong equivalence relation for crossings}\label{pic:homotopy_crossing_relation}
\end{figure}

\item Let $(\mathcal G,h)$ be an h-coinvariant of a functor $\mathcal F=\mathcal A, \mathcal {SA}, \mathcal C, \mathcal R$. Then $x_1\sim^{sh}_{\mathcal F} x_2$ implies $h(x_1)=h(x_2)$.
\end{enumerate}
\end{proposition}

We postpone the proof until Section~\ref{sect:elements}.





\begin{remark}\label{rem:h_strong_equivalence}
    According to Definition~\ref{def:homotopy_diagram_category}, the additional relations $x\sim y$ of the equivalences $\sim^{sh}$ comes from contraction of inverse Reidemeister moves or resolutions of singularities of codimension $2$. A direct check shows that the new relation $x\sim y$ in these cases implies $x\sim^s y$, except the cubic tangency resolution applied to the crossing functor $\mathcal C$. In the latter case we get the relation in Fig.~\ref{pic:homotopy_crossing_relation}.

    Note that the relation $\sim^s_{\mathcal C}$ is the equality $=$, because crossings cannot split or merge.
\end{remark}


By analogy with the affinity notion of Definition~\ref{def:strong_weak_equivalence}, we give the following.
\begin{definition}\label{def:crossing_homotopy_affinity}
Two crossings $x_1,x_2\in \mathcal C(D')$, $D'\in Ob(\mathfrak D_s)$, are called \emph{homotopically affined} if there exists morphism $f\colon D'\to D$ and $c_i\in \mathcal C(f)(x_i)$, $i=1,2$, such that the crossings $c_1$ and $c_2$ form the configuration shown in Fig.~\ref{pic:homotopy_crossing_relation}. In this case we write $x_1\approx_{\mathcal C,D}^h x_2$.
\end{definition}

\begin{remark}

    Let $\mathcal P\colon \mathfrak C\to \mathfrak C'$ and $\mathcal F\colon\mathfrak C\to Rel$ be functors.  In order to get the construction of the universal coinvariant of $\mathcal F$ compatible with $\mathcal P$, one should add the relation
\[
x_f=x_g \mbox{ whenever } \mathcal P(f)=\mathcal P(g)
\]
 to the list of relations determining the set $\tilde {\mathcal F}(c)$ in Section~\ref{subsect:universal_coinvariant}.

 Let $\mathcal P$ be the projection $\mathcal P\colon \mathfrak D_s\to\mathfrak D_h$ and functor $\mathcal F=\mathcal A$, $\mathcal {SA}$, $\mathcal C$ or $\mathcal R$. Given a tangle diagram, we will denote the universal $h$-coinvariant set of $\mathcal F$ at $D$ by $\mathcal A^{sh}(D)$, $\mathcal {SA}^{sh}(D)$, $\mathcal R^{sh}(D)$, $\mathcal C^{sh}(D)$.
\end{remark}

The aim of the next sections is to describe the sets of diagram elements $\mathcal A^{sh}(D)$, $\mathcal {SA}^{sh}(D)$, $\mathcal R^{sh}(D)$, $\mathcal C^{sh}(D)$ and  $\mathcal A^w(T)$, $\mathcal {SA}^w(T)$, $\mathcal R^w(T)$, $\mathcal C^w(T)$.

\section{Arc}\label{sect:arc}

Let $T\in\Sigma_0$ be a tangle and $D$ its diagram. Consider a tubular neighborhood $N(T)$ of $T$, and denote the complement manifold to the tangle by
\[
M_T=\overline{F\times I\setminus N(T)}\setminus \partial F\times I.
\]

\begin{definition}\label{def:arc_probe}
An \emph{arc probe} is an embedding
\[
\gamma\colon(I; 0; 1)\hookrightarrow(M_T; \partial N(T); F\times\{1\}).
\]
For an arc $a\in\mathcal A(D)$ choose a point $x\in a\subset F$. Let $y=(x,t)\in F\times I$ be the highest point of $(x\times I)\cap N(T)$. Then the embedding
\[
\gamma_a\colon (I;0;1)\hookrightarrow (x\times[t,1]; y; (x,1))
\]
is called a \emph{vertical arc probe} of the arc $a$.
\end{definition}

\begin{figure}[h]
\centering
  \includegraphics[width=0.4\textwidth]{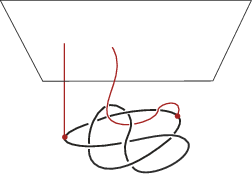}
  \caption{Arc probes}\label{pic:arc_probe}
\end{figure}


\begin{proposition}\label{prop:arc_affinity}
Let $a,b\in\mathcal A(D)$. Then $a\approx_{\mathcal A}b$ if and only if there is an embedded square $\Delta$ between the vertical arc probes $\gamma_a$ and $\gamma_b$
{\small
\[
(I\times I;0\times I;1\times I;I\times 0;I\times 1)\hookrightarrow (M_T;\gamma_a;\gamma_b;\partial N(T);F\times 1).
\]
}
such that the restriction of the projection of $\partial N(T)\to T$ to $\Delta\cap\partial N(T)$ is injective.
\end{proposition}

\begin{proof}
Necessity: if $a\approx b$ then there exists a sequence of moves $f$ which transforms $a$ and $b$ to one arc. We can suppose that  $f$ does not move points $x\in a$ and $y\in b$. Let $\phi$ be the corresponding isotopy of the tangle $T$. Then $\phi$ extends to a spatial isotopy fixed on $x\times [0,1]$ and $y\times [0,1]$. Since $x\in f(a)$ and $y\in f(b)$ belong to one arc then there exists a square $\Delta_0$ between the vertical arc probes $\gamma_{f(a)}$ and $\gamma_{f(b)}$. Then take $\Delta=\phi^{-1}(\Delta_0)$.

\begin{figure}[h]
\centering
  \includegraphics[width=0.7\textwidth]{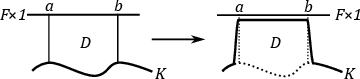}
  \caption{Pulling of an arc}\label{pic:arc_pull}
\end{figure}

Sufficiency: if $\Delta$ is a square then pull the bottom side of $\Delta$ to the top (Fig.~\ref{pic:arc_pull}). The induced isotopy of the knot leads to a diagram where the arcs $a$ and $b$ merge.
\end{proof}

For arc probes $\gamma, \gamma'$,  if there is an embedded square between them, denote $\gamma\approx^t_{\mathcal A, T}\gamma'$ and say that  $\gamma$ and $\gamma'$ are \emph{topologically affined}.

Proposition~\ref{prop:arc_affinity} means two arcs $a$ and $b$ are affined if and only if the corresponding vertical arc probes $\gamma_a$ and $\gamma_b$ are topologically affined.


\begin{example}\label{exa:arc_affinity_nontransitive}
    Consider the unknot in Fig.~\ref{pic:arc_nontransitivity} left. Then $a\approx_{\mathcal A} b$, $b\approx_{\mathcal A} c$ but $a\not\approx_{\mathcal A} c$. Indeed, if there had been an embedded square $\Delta$ between $\gamma_a$ and $\gamma_c$ then the boundary $\partial\Delta$ would have been the unknot. But it is the eight-knot (Fig.~\ref{pic:arc_nontransitivity} right). Thus, the relations $\approx_{\mathcal A}$ and $\approx^t_{\mathcal A}$ are not transitive.

\begin{figure}
\centering
  \includegraphics[width=0.3\textwidth]{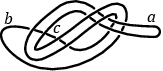}\qquad
  \includegraphics[width=0.3\textwidth]{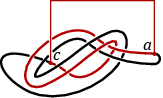}
  \caption{Non-transitivity of arc affinity}\label{pic:arc_nontransitivity}
\end{figure}

\end{example}

Denote the transitive closure of $\approx^t_{\mathcal A, T}$  by $\sim^{ts}_{\mathcal A, T}$. The relation $\sim^{ts}_{\mathcal A, T}$ in the set of arc probes is called  the \emph{topological strong equivalence}.

\begin{proposition}\label{prop:arc_top_affinity_closure}
    For any arc probes $\gamma, \gamma'$  $\gamma\sim^{ts}_{\mathcal A, T}\gamma'$ if and only if $\gamma$ is isotopic to $\gamma'$ (through arc probes).
\end{proposition}

\begin{lemma}\label{lem:embedded_square}
Let an arc probe $\gamma'$ is parallel to an arc probe $\gamma$ (i.e. obtained by shift along a transversal field to $\gamma$) in a manifold $M$. If $\gamma'$ is isotopic to $\gamma''$ in $M\setminus\gamma$ then there is an embedded square $\Delta$ in $M$ between $\gamma$ and $\gamma''$.
\end{lemma}
\begin{proof}
    By construction, there is an embedded square $\Delta_0$ between $\gamma$ and $\gamma'$. The isotopy from $\gamma'$ to $\gamma''$ extends to a spatial isotopy $F$ in $M$ which does not move $\gamma$. Then $\Delta=F(\Delta_0)$ is an embedded square between $\gamma$ and $\gamma'$.
\end{proof}



\begin{proof}[Proof of Proposition~\ref{prop:arc_top_affinity_closure}]
1. Let $\gamma\sim^{ts}_{\mathcal A, T}\gamma'$. Then there exists a sequence of arc probes $\gamma=\gamma_0\approx^t\gamma_1\approx^t\cdots\approx^t\gamma_n=\gamma'$. For any $i$ $\gamma_{i-1}\approx^t\gamma_i$ implies that $\gamma_{i-1}$ and $\gamma_i$ are isotopic. Since isotopy is a transitive relation, $\gamma$ and $\gamma'$ are isotopic.

2. Let $\gamma$ and $\gamma'$ are isotopic. Consider the projections $\bar\gamma=p(\gamma)$ and $\bar\gamma'=p(\gamma')$. The isotopy from $\gamma$ to $\gamma'$ can be described as a sequence of Reidemeister moves $\bar\gamma=\bar\gamma_0\to\bar\gamma_1\to\cdots\to\bar\gamma_n=\bar\gamma'$. Let $\gamma_i$, $i=0,\dots,n$, be a lift of $\bar\gamma_i$ to $F\times(0,1)$. It is sufficient to show that $\gamma_i\sim^{ts}_{\mathcal A, T}\gamma_{i+1}$ for all  $i$.

Note that any two lifts of $\bar\gamma_i$ are topologically strong equivalent.

For a Reidemeister move $\bar\gamma_i\to\bar\gamma_{i+1}$, take a lift $\gamma_i'$ of $\bar\gamma_i$ close to $\gamma_i$. Then $\gamma_i\approx^t_{\mathcal A, T}\gamma_i'$. There is an isotopy from $\gamma_i'$ to a lift $\gamma'_{i+1}$ of $\bar\gamma_{i+1}$ beyond $\gamma_i$.  
Then by Lemma~\ref{lem:embedded_square}
\[
\gamma_i\approx^t_{\mathcal A, T}\gamma_i'\approx^t_{\mathcal A, T}\gamma_{i+1}'\sim^t_{\mathcal A, T}\gamma_{i+1}.
\]
Hence, $\gamma_i\sim^{ts}_{\mathcal A, T}\gamma_{i+1}$.
\end{proof}

\begin{proposition}\label{prop:arc_strong_transitivity}
    For any $a,b\in\mathcal A(D)$  $a\sim^s_{\mathcal A, D} b$ if and only if $\gamma_a\sim^{ts}_{\mathcal A, T}\gamma_b$.
\end{proposition}

\begin{proof}
1. Let $\gamma_a\sim^{ts}_{\mathcal A, T}\gamma_b$. By definition there exists a sequence of arc probes
\[
\gamma_a=\gamma_0\approx^t_{\mathcal A, T}\gamma_1\approx^t_{\mathcal A, T}\cdots\approx^t_{\mathcal A, T}\gamma_n=\gamma_b.
\]
We can suppose that all $\gamma_i$ are distinct. Pull the curves $\gamma_i$, $1\le i\le n-1$, to their ends in $F\times 1$. Extend the transformation of the curves to an isotopy $\phi$ of $F\times I$ which does not move $\gamma_a$ and $\gamma_b$. Denote $T'=\phi(T)$ and $D'=p(T')$. The curves $\gamma'_i=\phi(\gamma_i)$ are vertical probes of some arcs $a'_i\in\mathcal A(D')$ (Fig.~\ref{pic:arc_probe_pull}).

\begin{figure}
\centering
  \includegraphics[width=0.7\textwidth]{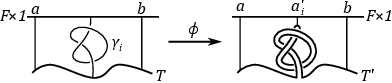}
  \caption{Pulling of an arc probe}\label{pic:arc_probe_pull}
\end{figure}

The condition $\gamma_{i-1}\approx^t_{\mathcal A, T}\gamma_i$ implies $\gamma'_{i-1}\approx^t_{\mathcal A, T'}\gamma'_i$. By Proposition~\ref{prop:arc_affinity}, $a'_{i-1}\approx_{\mathcal A} a'_i$ for all $i$. Hence,
\[
\phi_*(a)=a'_0\sim^s_{\mathcal A,D'} a'_n=\phi_*(b)
\]
where $\phi_*\colon\mathcal A(D)\to \mathcal A(D')$ is the relation induced by the isotopy $\phi$.
By definition of the strong equivalence $a\sim^s_{\mathcal A, D} b$.

2. Let $a\sim_{\mathcal A, D} b$. By Proposition~\ref{prop:arc_top_affinity_closure} it is enough to prove that $\gamma_a$ is isotopic to $\gamma_b$. We prove this by induction.
The case $a=b$ is trivial.

Let $f\colon D\to D'$ is a morphism in $\mathfrak D_s$, $a'\in\mathcal A(f)(a)$ and $b'\in\mathcal A(f)(b)$. Assume that $\gamma_{a'}$ and $\gamma_{b'}$ are isotopic, and $\psi_t$, $t\in[0,1]$, is an isotopy between these arc probes. Let $\phi$ be an isotopy in $F\times I$ which realizes the morphism $f$. Then $\phi^{-1}\circ\psi_t\circ\phi$ is an isotopy between $\gamma_a$ and $\gamma_b$.

Thus, by definition the equivalence $\sim^s_{\mathcal A, T}$ implies $\sim^{ts}_{\mathcal A, T}$.
\end{proof}

\begin{definition}\label{def:strong_arcs_tangle}
   The set of the isotopic classes of arc probes is called the \emph{(topological) set of arcs of the tangle $T$ (in the strong sense)} and denoted by ${\mathscr A}^{ts}(T)$.

   Let $h^{ts}_{\mathcal A,D}$ be the projection $\mathcal A(D)\to\mathscr A^{ts}(T)$ determined by the formula $a\mapsto\gamma_a$.
\end{definition}

\begin{theorem}\label{thm:arc_universal_homotopy_coinvariant}
   The pair $(\mathscr A^{ts}, h^{ts}_{\mathcal A})$ is the universal $h$-coinvariant of the arc functor $\mathcal A$.
\end{theorem}

\begin{proof}
For a morphism $f\colon T\to T'$ in the homotopy category $\mathfrak T_h$, i.e. an isotopy $f={f_t}$, $t\in[0,1]$, between the tangles $T$ and $T'$, the map $\mathscr A^{ts}(f)\colon\mathscr A^{ts}(T)\to\mathscr A^{ts}(T')$, $\gamma\mapsto f_1\circ\gamma$, is a bijection. Since $\mathscr A^{ts}(f)$ depends only on the homotopy class of the isotopy $f$, the functor $\mathscr A^{ts}$ is an $h$-coinvariant of $\mathcal A$.
Let us prove that $\mathscr A^{ts}$ is universal.

Let $(\mathcal G,h)$ be an $h$-coinvariant of $\mathcal A$. Define a single-valued natural transformation $\phi\colon \mathscr A^{ts}\Rightarrow \mathcal G$ as follows.

Let $T$ be a tangle. For an arc probe $\gamma$ of $T$, consider a spatial isotopy $f$ which verticalizes it (for example, by pulling the probe to its end like in the proof of Proposition~\ref{prop:arc_strong_transitivity}). Then $\gamma' = f(\gamma)$ is a vertical arc probe of some arc $a'\in\mathcal A(D')$ where the diagram $D'=p(T')$ is the projection of the tangle $T'=f(T)$ obtained by the pulling the tangle $T$. Define the value $\phi(\gamma)$ by the formula
\[
\phi(\gamma)=\mathcal G(f)^{-1}(h(a'))\in \mathcal G(D).
\]
We show that $\phi$ is well defined. Let $\gamma'$ be an arc probe of $T$ isotopic to $\gamma$ and $H$ an isotopy from $\gamma'$ to $\gamma$ in $M_T$. Denote $f'=f\circ H$. Then
\[
\phi(\gamma_1)=\mathcal G(f')^{-1}(h(a'))=\mathcal G(f)^{-1}(h(a'))=\phi(\gamma)
\]
because $f$ and $f'$ define the same isotopy of the tangle $T$, i.e. the same morphism in $\mathfrak T_h$. Hence, the element $\phi(\gamma)$ does not depend on the choice of a representative in the isotopy class of the arc probe.

Let $f'$ and $f''$ be spatial isotopies which verticalizes $\gamma$. Then $f'(\gamma)=\gamma_{a'}$ and $f''(\gamma)=\gamma_{a''}$ for some arcs $a'\in\mathcal A(D')$ and  $a''\in\mathcal A(D'')$. We can define $\phi(\gamma)$ using either $a'$ or $a''$. Let us show the choice does not affect the result. We need to check that $\mathcal G(f')^{-1}(h(a'))=\mathcal G(f'')^{-1}(h(a''))$, i.e.
\[
h(a'')=\mathcal G(f'')\circ \mathcal G(f')^{-1}(h(a'))=\mathcal G(f''\circ f'^{-1})(h(a')).
\]
Denote the parametrization of the isotopy $f=f''\circ f'^{-1}$ by $h_t$, $t\in[0,1]$. Let $\gamma_t=f_t(\gamma_{a'})$. Consider a spatial isotopy $g_t$, $t\in[0,1)$, which contracts vertically the arc probe $\gamma_{a'}$ to its end in $F\times 1$. Now, modify the isotopy $f$ by contracting the curves $\gamma_t$ to segments short enough to be considered as vertical: $\tilde f_t=f_t\circ g_{\zeta(t)}$ where $\zeta\colon [0,1]\to [0,1)$ is a continuous function such that $\zeta(0)=\zeta(1)=0$ and $\zeta(t)$ is sufficiently close to $1$ for $t\in(0,1)$.

The isotopies $f$ and $\tilde f$ are homotopic. Since $\mathcal G$ is an $h$-coinvariant, $\mathcal G(f)=\mathcal G(\tilde f)$. For any $t\in [0,1]$, the curve $\tilde \gamma_t=\tilde f_t(\gamma)$ is a vertical arc probe $\tilde \gamma_t=\gamma_{a_t}$ of some arc $a_t\in \mathcal A(\tilde f_t(D))$. Hence,
\[
a''=a_1=\mathcal A(\tilde f)(a_0)=\mathcal A(\tilde f)f(a'),
\]
then $h(a'')=\mathcal G(\tilde f)(h(a'))=\mathcal G(f)(h(a'))$.
Thus, the map $\phi$ is well defined.

Let $f\colon D_1\to D_2$ be a morphism in $\mathfrak D_s$, $\gamma_1\in\mathscr A^{ts}(D_1)$ and $\gamma_2=\mathscr A^{ts}(f)(\gamma_1)$. We can consider $f$ as a spatial isotopy such that $\gamma_2=f(\gamma_1)$. Let $g\colon D_2\to D'$ be a spatial isotopy which verticalizes $\gamma_2$, i.e. $g(\gamma_2)=\gamma_{a'}$, $a'\in\mathcal A(D')$. Then $\phi(\gamma_2)=\mathcal G(g)^{-1}(h(a'))$. Since the isotopy $g\circ f$ verticalizes $\gamma_1$,
\[
\phi(\gamma_1)=\mathcal G(g\circ f)^{-1}(h(a'))=\mathcal G(f)^{-1}\mathcal G(g)^{-1}(h(a')).
\]
Then $\phi(\mathscr A^{ts}(f)(\gamma_1))=\phi(\gamma_2)=\mathcal G(f)(\phi(\gamma_1))$. Thus, $\phi$ is a single-valued natural transformation between $\mathscr A^{ts}$ and $\mathcal G$.
\end{proof}

The next statement shows that any set of (implicit) arc of a tangle can be made explicit in some diagram of the tangle.
\begin{proposition}\label{prop:arc_top_realizability}
 For any finite subset of $\{\gamma_i\}_{i=1}^n\subset\mathscr A^{ts}(T)$ there is a morphism $f\colon T\to T'$ such that for any $i=1,\dots, n$, $f(\gamma_i)=\gamma_{a'_i}$ for some $a'_i\in\mathcal A(D')$ where $D'=p(T')$.
\end{proposition}
\begin{proof}
    We can isotope the curves $\gamma_i$ to make them distinct. Consider the isotopies which pull the curves $\gamma_i$ in $F\times I$ along themselves close enough to $F\times 1$. Extend these isotopies to a spatial isotopy $f$ of $F\times I$. Then $f$ is the required morphism.
\end{proof}

\begin{remark}
There is an action of the group of classical long knots $\mathcal L$ on the sets $\mathscr A^{ts}(T)$ given by the formula
$(\lambda,\gamma)\mapsto \gamma\lambda$, $\lambda\in\mathcal L$, $\gamma\in\mathscr A^{ts}(T)$ (Fig.~\ref{pic:arc_long_knot_action}). 

\begin{figure}[t]
\centering
  \includegraphics[width=0.7\textwidth]{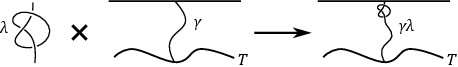}
  \caption{Action of a long knot on an arc of the tangle}\label{pic:arc_long_knot_action}
\end{figure}
\end{remark}


\begin{example}[The arcs of the unknot]
    Let $U\subset D^2\times[0,1]$ be the unknot. Let us show that the set of arcs $\mathscr A^{ts}(U)$ is trivial (consists of one class). Consider an arc probe $\gamma$. The complement $M_U=D^2\times I\setminus N(U)$ is homeomorphic to a full torus with a ball $B$ removed (Fig.~\ref{pic:arc_unknot_homeo}). Consider an isotopy $f_t$, $t\in[0,1]$, which straighten the curve $\gamma$ as shown in Fig.~\ref{pic:arc_unknot_isotopy}. Modify the isotopy $f_t$ to make the ball fixed as follows. Let $(z,\phi)\in D^2\times S^1$ be coordinates in the full torus, where $D^2=\{z\in\mathbb C\mid |z|\le 1\}$. assume the coordinates of the center of the ball $B$ under the isotopy are $(z(t),\phi(t))$ where $z(0)=0$ and $\phi(0)=0$. Then the desired isotopy is $f'_t=g_t\circ f_t$ where $g_t(z,\phi)=\left(\frac{z-z(t)}{1-\overline{z(t)}z}, \phi-\phi(t)\right)$.

\begin{figure}
\centering
  \includegraphics[width=0.7\textwidth]{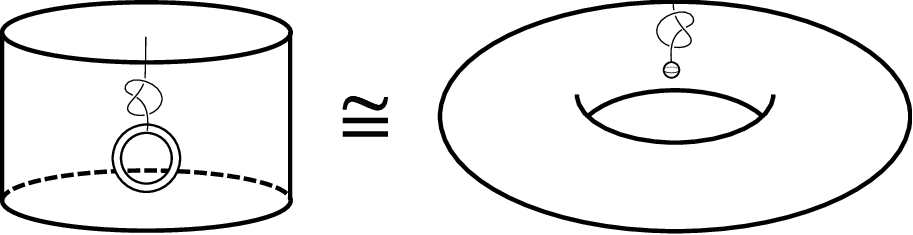}
  \caption{The complement to the unknot}\label{pic:arc_unknot_homeo}
\end{figure}

\begin{figure}
\centering
  \includegraphics[width=0.5\textwidth]{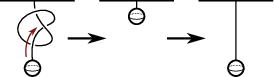}
  \caption{Arc probe straightening}\label{pic:arc_unknot_isotopy}
\end{figure}
    Thus, any arc probe is topologically strong equivalent to the standard arc probe. Hence, $\mathscr A^{ts}(U)$ is trivial.

    Note that the set of arcs of the unknot in another thickened surface can be nontrivial (Fig.~\ref{pic:arc_unknot_torus}).
\begin{figure}
\centering
  \includegraphics[width=0.2\textwidth]{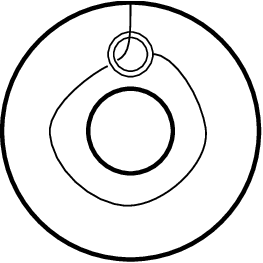}
  \caption{Nontrivial arc probe of the unknot in the thickened torus}\label{pic:arc_unknot_torus}
\end{figure}

    The set of arcs $\mathscr A^{ts}(K)$ for a nontrivial knot $K$ does not consists of one point, because $\mathscr A^{ts}(K)$ projects onto the knot quandle which is nontrivial.
\end{example}

\subsection{Weak equivalence and motion group}\label{subsect:arc_motion_group}

Recall the definition of motion group~\cite{Goldsmith}.

\begin{definition}\label{def:motion_group}
    Let $N$ be a submanifold of a manifold $M$. A \emph{motion} of $N$ in $M$ is a path $f_t$, $t\in[0,1]$, in $Diff(M)$ such that $f_0=id$ and $f_1(N)=N$.

    A motion $f_t$ is \emph{stationary} if $f_t(N)=N$ for all $t\in[0,1]$.

    Two motions $f_t$ and $g_t$ are \emph{equivalent} if the motion $g^{-1}_t\circ f_t$ is homotopic, with endpoints fixed, to a stationary motion. The set $\mathcal M(M,N)$ of the equivalence classes of motions is called the \emph{motion group}.
\end{definition}
\begin{example}\label{exa:motion_stationary_motion}
Consider a classical trefoil $K$. The rotation of the space by $\frac{2\pi}3$ around the symmetry axis of the trefoil is an example of a motion of $K$ (Fig.~\ref{pic:motions} left). A spatial isotopy induced by pulling the trefoil along itself is a stationary motion of $K$ (Fig.~\ref{pic:motions} right).
\begin{figure}
\centering
  \includegraphics[width=0.4\textwidth]{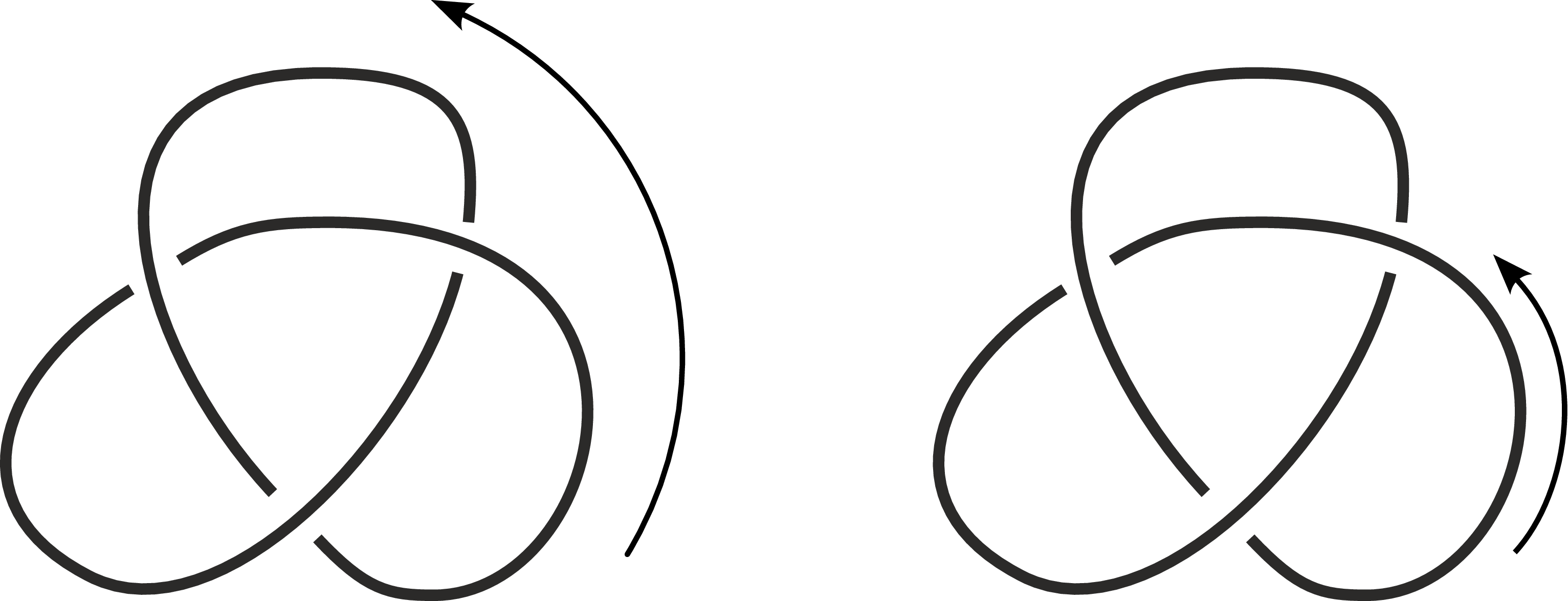}
  \caption{A motion and a stationary motion}\label{pic:motions}
\end{figure}
\end{example}

\begin{definition}\label{def:arc_topological_weak_equivalence}
    An arc probe $\gamma$ of a tangle $T$ is \emph{topologically weak equivalent} to an arc probe $\gamma'$ of a tangle $T'$ ($\gamma\sim^{tw}_{\mathcal A}\gamma'$) if there is an isotopy $\phi$ of $F\times[0,1]$ such that $\phi(T)=T'$ and $\phi(\gamma)=\gamma'$. Denote the set of topologically weak equivalence classes of the arcs of the tangle $T$ by $\mathscr A^{tw}(T)$.
\end{definition}

Note that for isotopic tangles $T$ and $T'$ the sets $\mathscr A^{tw}(T)$ and $\mathscr A^{tw}(T')$ coincide.

For a tangle $T$, its motion group $M(F\times I, T)$ acts on the set $\mathscr A^{ts}(T)$ by composition: $f_t\times\gamma\mapsto f_1(\gamma)$, $f_t\in M(F\times I, T)$, $\gamma\in\mathscr A^{ts}(T)$. By definition of topological weak equivalence, we have the following statement.

\begin{proposition}\label{prop:arc_topological_strong_weak_sets}
$\mathscr A^{tw}(T)= \mathscr A^{ts}(T)/M(F\times I, T)$.
\end{proposition}

\begin{theorem}\label{thm:arc_universal_homotopy_invariant}
The correspondence $T\mapsto\mathscr A^{tw}(T)$ is the universal h-invariant of the arc functor $\mathcal A$.
\end{theorem}

\begin{proof}
A direct check shows that $\mathscr A^{tw}$ is an $h$-invariant. We need to demonstrate the universality property.

Let $(\mathcal G, h)$ be an $h$-invariant of $\mathcal A$. Then it is an $h$-coinvariant, hence, there is a (single-valued) natural transformation $\phi\colon\mathscr A^{ts}\Rightarrow \mathcal G$. We need to show that $\phi(\gamma_1)=\phi(\gamma_2)$ for any $\gamma_1\sim_{\mathcal A}^{tw}\gamma_2$. Then $\phi=\bar\phi\circ\pi$ where $\pi\colon\mathscr A^{ts}(T)\to\mathscr A^{tw}(T)$ is the natural projection, and the single-valued natural map $\bar\phi\colon\mathscr A^{tw}\Rightarrow \mathcal G$ establishes the universality of $\mathscr A^{tw}$.

By definition of $\phi$, there exist isotopies $f_i\colon T_i\to T_i'$, $i=1,2$ which map the arc probes $\gamma_i$ to vertical probes of some arcs $a_i\in\mathcal A(D_i')$ where $D_i'=p(T_i')$, and $\phi(\gamma_i)=G(f_i)^{-1}(h(a_i))=h(a_i)$ because $\mathcal G$ is an invariant. The relation $\gamma_1\sim^{tw}_{\mathcal A}\gamma_2$ implies there is an isotopy $f$ such that $T_2=f(T_2)$ and $\gamma_2=f(\gamma_1)$. Then the composition $f_2\circ f\circ f_1^{-1}$ maps $a_1$ to $a_2$. By proof of Theorem~\ref{thm:arc_universal_homotopy_coinvariant},
\[
h(a_2)=\mathcal G(f_2\circ f\circ f_1^{-1})(h(a_1))=h(a_1),
\]
hence, $\phi(\gamma_2)=\phi(\gamma_1)$.
\end{proof}




\begin{proposition}\label{prop:classical_weak_arc_triviality}
Let $K\subset F\times[0,1]$ be a classical knot (i.e. $F=D^2$ or $S^2$). Then the set of arcs in the weak sense $\mathscr A^{tw}(K)$ of $K$ is trivial.
\end{proposition}
\begin{proof}
Let $\gamma_i$, $i=1,2$, be arc probes of $K$. Contract $\gamma_i$ to its end on $F\times 1$. Extend this isotopy of $\gamma_i$ to a spatial isotopy $f_i$. Then $f_i$ maps $\gamma_i$ to a short vertical segment $\gamma_i'$ and $f_i(K)=K_i$ lies below $\gamma_i'$. We can assume that $\gamma_1'=\gamma_2'$. Denote the length of $\gamma'_1$ by $\varepsilon$ and consider $K_1$ and $K_2$ as long knots in $F\times[0,1-\varepsilon]$. Since $K_1$ and $K_2$ are isotopic to $K$, there is an isotopy $g$ in $F\times[0,1-\varepsilon]$ such that $g(K_1)=K_2$. Extend $g$ by identity to $F\times I$. Then the composition $f=f_2^{-1}\circ g\circ f_1$ maps $K$ to $K$ and $\gamma_1$ to $\gamma_2$. Hence, $\gamma_1\sim^{tw}_{\mathcal A}\gamma_2$. Thus, all arc probes are weak equivalent.
\end{proof}

\begin{corollary}\label{cor:arc_classical_motion_group_transitivity}
   For a classical knot $K$, the motion group $M(F\times I, T)$ acts transitively on the set of arcs $\mathscr A^{ts}(K)$.
\end{corollary}

\begin{corollary}\label{cor:arc_classical_trivial_h-invariant}
    For classical knots, any $h$-invariant of arcs is trivial.
\end{corollary}

\begin{proposition}
Long knot acts trivially on arcs of classical knots.
\end{proposition}
\begin{proof}
    For any knot $K\subset D^2\times I$, arc probe $\gamma$ and a long knot $\lambda$, there is a motion which transforms $\gamma\lambda$ to $\gamma$, see Fig.~\ref{pic:arc_long_knot_triviality}.
\begin{figure}
\centering
  \includegraphics[width=0.6\textwidth]{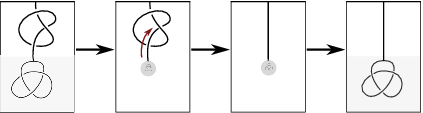}
  \caption{Triviality of the action of long knots}\label{pic:arc_long_knot_triviality}
\end{figure}
\end{proof}

\begin{remark}
    We have defined above arc probe as a curve that connects the tangle with $F\times 1$.  Analogously, one can consider curves connecting the tangle with $F\times 0$. The isotopy classes of such curves corresponds to \emph{under-arcs} of the tangle (which can be identified with the arcs of the mirror tangle).
\end{remark}

{\color{red}
}

\section{Region}\label{sect:region}

Let $T\in\Sigma_0$ be a tangle and $D$ its diagram.

\begin{definition}\label{def:region_probe}
A \emph{region probe} is an embedding
\[
\gamma\colon(I; 0; 1)\hookrightarrow(M_T; F\times 0;  F\times 1)
\]
which presents a trivial long knot in $F\times I$, i.e. isotopic to $x\times I\subset F\times I$.

For a region $r\in\mathcal R(D)$ choose a point $x\in r\subset F$. The embedding $\gamma_r = x\times I$ is called a \emph{vertical region probe} of the region $r$.

\begin{figure}[h]
\centering
  \includegraphics[width=0.4\textwidth]{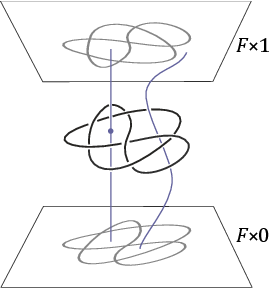}
  \caption{A region probe}\label{pic:region_probe}
\end{figure}
\end{definition}

\begin{proposition}\label{prop:region_affinity}
Let $r,r'\in\mathcal R(D)$. Then $r\approx_{\mathcal R}r'$ if and only if there is an embedded square between the vertical region probes $\gamma_r$ and $\gamma_{r'}$:
{\small
\[
(I\times I; 0\times I; 1\times I; I\times 0; I\times 1)\hookrightarrow (M_T; \gamma_r; \gamma_{r'}; F\times 0; F\times 1).
\]
}
\end{proposition}
\begin{proof}
If $r\approx_{\mathcal R} r'$ then there exists a sequence of moves $f$ which transforms $r$ and $r'$ to one region. We can suppose that  $f$ does not move points $x\in r$ and $y\in r'$. Let $\phi$ be the corresponding isotopy of the tangle $T$. Then $\phi$ extends to a spatial isotopy fixed on $\{x,y\}\times I$. Since $x\in f(r)$ and $y\in f(r')$ belong to one region, there is a curve $\alpha\subset f(r)$ which connects $x$ and $y$. Then $\Delta_0=\alpha\times I$ is a square between the vertical region probes $\gamma_{f(r)}$ and $\gamma_{f(r')}$. Then take $\Delta=\phi^{-1}(\Delta_0)$.

Let $\Delta$ be a square between the vertical region probes $\gamma_r=x\times I$ and $\gamma_{r'}=y\times I$. Isotope the square $\Delta$ rel $\gamma_r\cup\gamma_{r'}$ so that $\Delta\cap F\times k= \alpha\times k$, $k=0,1$,  for some curve $\alpha\subset F$ connecting $x$ and $y$. (Alternatively, we can glue a shrunken mirror image of $\Delta$ below to $\Delta$ as shown in Fig.~\ref{pic:region_mirror_square}.)

\begin{figure}[h]
\centering
  \includegraphics[width=0.7\textwidth]{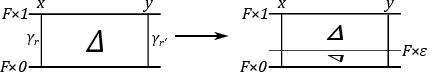}
  \caption{Gluing a mirror copy to a square between $\gamma_r$ and $\gamma_{r'}$}\label{pic:region_mirror_square}
\end{figure}

Let $\Delta_0=\alpha\times I$. One can isotope $\Delta$ so that $\partial\Delta=\partial\Delta_0$ and the interiors of $\Delta$ and $\Delta_0$ don't intersect near the boundary. The disks $\Delta$ and $\Delta_0$ are homotopic rel $\partial\Delta$. By eliminating components of the intersection $\Delta\cap\Delta_0$, one constructs an isotopy $\phi$ from $\Delta$ to $\Delta_0$ rel $\partial\Delta$. Extend $\phi$ to a spatial isotopy of $F\times I$. Then $\phi$ induces a sequence of Reidemeister moves on the tangle $T$, which merges the regions $r$ and $r'$. Thus, $\mathcal R(\phi)(r)\cap\mathcal R(\phi)(r')\ne\emptyset$ and $r\approx_{\mathcal R} r'$.
\end{proof}

For region probes $\gamma, \gamma'$ such that there is an embedded square between them, we denote $\gamma\approx^t_{\mathcal R, T}\gamma'$ and say that  $\gamma$ and $\gamma'$ are \emph{topologically affined}.

\begin{example}\label{ex:region_top_nontransitivity}
    Consider the unknot $K$ in Fig.~\ref{pic:region_nontransitivity} left. Then $a\approx_{\mathcal R} c$, $b\approx_{\mathcal R} c$ but $a\not\approx_{\mathcal R} b$. Indeed, if there had been an embedded square $\Delta$ between $\gamma_a$ and $\gamma_b$ then the boundary $\partial\Delta$ would have not linked with the knot. But the knot and $\partial\Delta$ form the Whitehead link (Fig.~\ref{pic:region_nontransitivity} right).
\begin{figure}
\centering
  \includegraphics[width=0.4\textwidth]{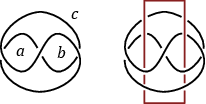}
  \caption{Non-transitivity of region affinity}\label{pic:region_nontransitivity}
\end{figure}
\end{example}

Denote the transitive closure of $\approx^t_{\mathcal R, T}$  by $\sim^{ts}_{\mathcal R, T}$. The relation $\sim^{ts}_{\mathcal R, T}$ in the set of region probes is called  the \emph{topological strong equivalence}.

\begin{proposition}\label{prop:region_top_affinity_closure}
    For any region probes $\gamma, \gamma'$  $\gamma\sim^{ts}_{\mathcal R, T}\gamma'$ if and only if $\gamma$ is isotopic to $\gamma'$ in $M_T$ (through region probes).
\end{proposition}

\begin{proof}
    The necessity is evident.

    Let $\phi=(\gamma_s)$, $s\in[0,1]$, be an isotopy from $\gamma$ to $\gamma'$. There is a sequence $0=s_0<s_1<\cdots<s_n=1$ such that for any $i$ and any $s\in(s_i,s_{i+1})$ we have $\gamma_s\subset U_i$ for some tubular neighbourhood $U_i$ of the region probe $\gamma_i$ in $M_T$. Then $U_i\simeq \mathbb D^2\times I$ where $\gamma_i$ is identified with $0\times I$. We will homotope the isotopy $\phi$ to an isotopy $\tilde\phi=\tilde\gamma_s$ such that $\gamma=\tilde\gamma_0\approx_{\mathcal R}\tilde\gamma_{1/2}\approx_{\mathcal R}\tilde\gamma_1=\gamma'$.

    For simplicity assume $i=0$, $s_0=0$ and $s_1=1$ and identify $U_0$ with $\mathbb D^2\times[0,1]$. Consider the vertical region probe $\hat\gamma=\imath\times I$. Extend the isotopy $\phi=(\gamma_s)$ to a spatial isotopy $(f_s)$, $s\in[0,1]$, of $\mathbb D^2\times I$ such that $f_s|_{\partial\mathbb D^2\times I}=id$.

\begin{figure}
\centering
  \includegraphics[width=0.7\textwidth]{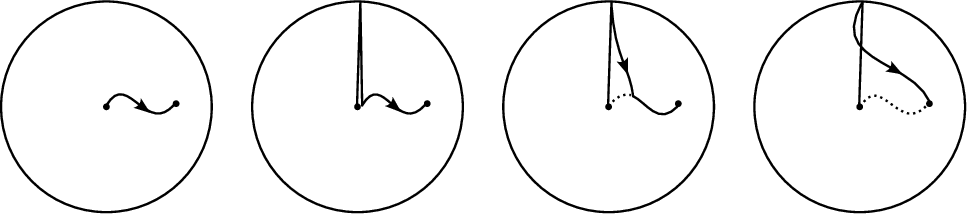}
  \caption{Modification of the isotopy in the tubular neighbourhood}\label{pic:region_probe_small_isotopy}
\end{figure}

    The scheme of the homotopy $\phi_\lambda$ from $\phi$ to $\tilde\phi$ is shown on Fig.~\ref{pic:region_probe_small_isotopy}: we insert a linear isotopy from $\gamma$ to $\hat\gamma$ and back, and then transform the backward linear isotopy to an isotopy from $\hat\gamma$ to $\gamma'$. Denote the linear isotopy from $\gamma$ to $\hat\gamma$ by $\psi=(\hat\gamma_s)$ where $\hat\gamma_s(t)=(s\imath,t)\in\mathbb D^2\times I$. Define the homotopy $\phi_\lambda$ by the formula

\begin{gather*}
\phi_\lambda(s)=\left\{\begin{array}{lc}
    \hat\gamma_{4s}, &  0\le s\le\frac\lambda 2, \\
\hat\gamma_{4(\lambda-s)}   , & \frac\lambda 2\le s\le \lambda,\\
 \gamma_{\frac{s-\lambda}{1-\lambda}}    , & \lambda\le s\le 1,
\end{array}\right.
\qquad \mbox{when }0\le\lambda\le \frac 12.
\end{gather*}

\begin{gather*}
\phi_\lambda(s)=\left\{\begin{array}{lc}
    \hat\gamma_{\frac{2s}{\lambda}}, &  0\le s\le\frac\lambda 2, \\
 f_{2\lambda-1}\left(\hat\gamma_{\frac{2(\lambda-s)}\lambda}\right)    , & \frac\lambda 2\le s\le \lambda,\\
 \gamma_{2s-1}    , & \lambda\le s\le 1,
\end{array}\right.
\qquad \mbox{when }\frac 12\le\lambda\le 1,
\end{gather*}
Denote $\tilde\phi=\phi_1$. Then $\tilde\gamma_0=\gamma$ and $\tilde\gamma_{1/2}=\hat\gamma$ are connected by the embedded square $\Delta_0=[0,\imath]\times I$, and $\tilde\gamma_{1/2}$ and $\tilde\gamma_{1}=\gamma'$ are connected by the embedded square $\Delta_1=f_1(\Delta_0)$.

Thus, we can homotope the isotopy $\phi$ to an isotopy $\tilde\phi$ such that
\[
\gamma=\tilde\gamma_0\approx_{\mathcal R}\tilde\gamma_{\frac 12}\approx_{\mathcal R}\tilde\gamma_1\approx_{\mathcal R}\tilde\gamma_{\frac 32}\approx_{\mathcal R}\cdots\approx_{\mathcal R}\tilde\gamma_{n-\frac{1}2}\approx_{\mathcal R}\tilde\gamma_{n}=\gamma',
\]
hence, $\gamma\sim^{ts}_{\mathcal R, T}\gamma'$.
\end{proof}


\begin{proposition}\label{prop:region_strong_transitivity}
    For any $r,r'\in\mathcal R(D)$  $r\sim^s_{\mathcal R, D} r'$ if and only if $\gamma_{r}\sim^{ts}_{\mathcal R, T}\gamma_{r'}$.
\end{proposition}

\begin{proof}
    Let $\gamma_{r}\sim^{ts}_{\mathcal R, T}\gamma_{r'}$. Then there exist region probes $\gamma_i$, $i=0,\dots,n$, such that
\[
\gamma_r=\gamma_0\approx_{\mathcal R}\gamma_{1}\approx_{\mathcal R}\cdots\approx_{\mathcal R}\gamma_{n}=\gamma_{r'}.
\]
Denote the embedded square between $\gamma_{i-1}$ and $\gamma_i$ by $\Delta_i$, $i=1,\dots,n$. Isotope the squares $\Delta_i$ so that for any $i=1,\dots,n-1$,  $\Delta_i\cup\Delta_{i+1}$ is a smooth surface in the neighbourhood of $\gamma_i$.

We want to show that the region probes $\gamma_i$ can be simultaneously straightened (after some isotopy) in $F\times I\setminus(\gamma_r\cup\gamma_{r'})$. Apply subsequently for $i=1,\dots,n-2$ the following transformation. Consider the intersection points $\Delta_i\cap(\gamma_0\cup\gamma_1\cup\cdots\cup\gamma_{i-2}\cup\gamma_n)$. Redraw the curve $\gamma_i$ in $\Delta_i$ to exclude the intersection points from the square (Fig.~\ref{pic:region_probe_correction}). The new square $\Delta'_i\subset\Delta_i$ between $\gamma_{i-1}$ and the new curve $\gamma'_i$ does not intersect the curves $\gamma_0,\gamma_1,\dots,\gamma_{i-2},\gamma_n$. There is an isotopy from $\gamma_i$ to $\gamma'_i$ in $\Delta_i$ which extends to a spatial isotopy $\phi_i$ in $M_T$ rel $\bigcup_{j=i+1}^{n-1} \gamma_j$. Then replace $\Delta_i$ by $\Delta'_i$, $\gamma_i$ by $\gamma'_i$, and $\Delta_{i+1}$ by $\Delta'_{i+1}=\phi_i(\Delta_{i+1})$.

\begin{figure}[h]
\centering
  \includegraphics[width=0.8\textwidth]{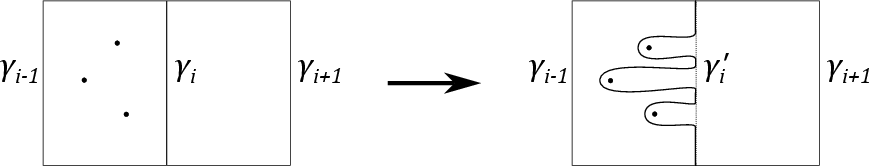}
  \caption{Preliminary isotopy of a region probe}\label{pic:region_probe_correction}
\end{figure}

After modification we have $\Delta_i\cap(\gamma_0\cup\gamma_1\cup\cdots\cup\gamma_{i-2}\cup\gamma_n)=\emptyset$ for any $i=1,\dots,n-2$. Then there is an isotopy $f$ of $F\times I$ rel $\gamma_0\cup\gamma_n$ which straightens consequently the region probes $\gamma_1,\dots,\gamma_{n-2}$ (by pulling $\gamma_i$ along $\Delta_i$ to $\gamma_{i-1}$ which was straightened previously). Denote $\gamma_{r'_i}=\gamma'_i=f(\gamma_i)$, $\Delta'_i=f(\Delta_i)$, $i=1,\dots,n$, and $T'=f(T)$. By Proposition~\ref{prop:region_affinity}
\[
\mathcal R(f)(r)=r'_1\approx_{\mathcal R,D'} r'_2 \approx_{\mathcal R,D'}\cdots \approx_{\mathcal R,D'} r'_{n-2}
\]
where $D'=p(T')$.
Let us show that $r'_{n-2}\sim_{\mathcal R,D'}^s r'_n=\mathcal R(f)(r')$. Then $\mathcal R(f)(r)\sim_{\mathcal R,D'}^s\mathcal R(f)(r')$, hence, $r\sim_{\mathcal R,D}^s r'$.

Isotope $\Delta'_{n-1}$ and $\Delta'_n$ so that $p(\Delta'_{n-1}\cup\Delta'_n)\ne F$, and take $\tilde\gamma=x\times I\in F\times I\setminus (\Delta'_{n-1}\cup\Delta'_n)$. Then there is a spatial isotopy $g_1$ of $F\times I$ rel $\gamma'_{n-2}\cup\tilde\gamma$ which straightens $\gamma'_{n-1}$, and an isotopy $g_2$ of $F\times I$ rel $\gamma'_{n}\cup\tilde\gamma$ which straightens $\gamma'_{n-1}$. Connect the vertical probes $g_1(\tilde\gamma)$ and $g_1(\gamma'_{n-1})$ by an embedded square $\tilde\Delta_0$ in $F\times I\setminus\gamma'_{n-2}$ so that $\tilde\Delta_0$ is transversal to $g_1(\gamma'_n)$ and $g_1(T')$. Denote $\tilde\Delta=g_1^{-1}(\tilde\Delta_0)$.

Isotope the probe $\tilde\gamma$ in $\tilde\Delta$ to eliminate the intersection of the square with $T'$ like in Fig.~\ref{pic:region_probe_correction}. Denote the obtained probe by $\tilde\gamma'$ and the obtained square by $\tilde\Delta'$. Then $\gamma'_{n-1}\approx_{\mathcal R,T}^t\tilde\gamma'$, and we can suppose that $\tilde\gamma'$ is isotopic to $\tilde\gamma$ in $F\times I\setminus(\gamma_{n-2}'\cup\gamma'_{n-1}\cup\gamma'_n)$. Extend this isotopy to a spatial isotopy $f_1$ of $F\times I$ rel $\gamma_{n-2}'\cup\gamma'_{n-1}\cup\gamma'_n$.

Denote $T''=f_1(T')$, $r''_{n-2}=\mathcal R(f_1)(r'_{n-2})$ and $r''_{n}=\mathcal R(f_1)(r'_{n})$. Then $\gamma'_{n-2}=\gamma_{r''_{n-2}}$ and $\gamma'_n=\gamma_{r''_n}$. Since $\tilde\gamma\cap T''=\emptyset$, $\tilde\gamma=\gamma_{\tilde r}$ for some $\tilde r\in\mathcal R(T'')$. Let us show that $r''_{n-2} \sim_{\mathcal R,D''}^s \tilde r \sim_{\mathcal R,D''}^s r''_n$ where $D''=p(T'')$.

Indeed, the isotopy $g_1$ straightens $\gamma'_{n-1}$; the vertical region probes $\gamma'_{n-2}$ and $g_1(\gamma'_{n-1})$ are connected by the embedded square $g_1\circ f_1(\Delta'_{n-1})$, and $g_1(\gamma'_{n-1})$ and $\tilde\gamma$ are connected by the square $g_1\circ f_1(\tilde\Delta')$. Then $\gamma'_{n-2}\approx_{\mathcal R,g_1(T'')}g_1(\gamma'_{n-1})\approx_{\mathcal R,g_1(T'')}\tilde\gamma$, hence, $\mathcal R(g_1)(r''_{n-2})\sim_{\mathcal R, g_1(D'')}^s\mathcal R(g_1)(\tilde r)$ and $r''_{n-2}\sim_{\mathcal R, D''}^s\tilde r$.

Analogously, $r''_{n}\sim_{\mathcal R, D''}^s\tilde r$, hence, $r''_{n-2}\sim_{\mathcal R, D''}^s r''_n$ and $r'_{n-2}\sim_{\mathcal R, D'}^s r'_n$. Thus, $r\sim_{\mathcal R,D}^s r'$.
\end{proof}

\begin{definition}\label{def:region_strong_tangle}
   The set of the isotopic classes of region probes is called the \emph{(topological) set of regions of the tangle $T$ (in the strong sense)} and denoted by ${\mathscr R}^{ts}(T)$.
   Let $h^{ts}_{\mathcal R,D}$ be the projection $\mathcal R(D)\to\mathscr R^{ts}(T)$ determined by the formula $r\mapsto\gamma_r$.
\end{definition}

\begin{theorem}\label{thm:region_universal_homotopy_coinvariant}
   The pair $(\mathscr R^{ts}, h^{ts}_{\mathcal R})$ is the universal $h$-coinvariant of the region functor $\mathcal R$.
\end{theorem}

\begin{proof}
    The proof repeats the arguments of the proof of Theorem~\ref{thm:arc_universal_homotopy_coinvariant}.
    A direct check shows that  $(\mathscr R^{ts}, h^{ts}_{\mathcal R})$ is an $h$-coinvariant of $\mathcal R$.

Let $(\mathcal G,h)$ be an $h$-coinvariant of $\mathcal R$. Define a single-valued natural transformation $\phi\colon \mathscr R^{ts}\Rightarrow \mathcal G$ as follows.

For a region probe $\gamma$ of a tangle $T$, consider a spatial isotopy $f$ which verticalizes it. Then $\gamma' = f(\gamma)$ is a vertical region probe of some region $r'\in\mathcal R(D')$ where $D'=p(T')$ and  $T'=f(T)$. The value $\phi(\gamma)$ is defined as follows
\[
\phi(\gamma)=\mathcal G(f)^{-1}(h(r'))\in \mathcal G(D).
\]
Let $\gamma'$ be a region probe of $T$ isotopic to $\gamma$ and $H$ an isotopy from $\gamma'$ to $\gamma$ in $M_T$. Denote $f'=f\circ H$. Then
\[
\phi(\gamma_1)=\mathcal G(f')^{-1}(h(r'))=\mathcal G(f)^{-1}(h(r'))=\phi(\gamma)
\]
because $f$ and $f'$ define the same isotopy of the tangle $T$, i.e. the same morphism in $\mathfrak T_h$. Hence, the element $\phi(\gamma)$ does not depend on the choice of a representative in the isotopy class of the region probe.

Let $f'$ and $f''$ be spatial isotopies which verticalizes $\gamma$. Then $f'(\gamma)=\gamma_{r'}$ and $f''(\gamma)=\gamma_{r''}$ for some regions $r'\in\mathcal R(D')$ and  $r''\in\mathcal R(D'')$. We need to check that $\mathcal G(f')^{-1}(h(r'))=\mathcal G(f'')^{-1}(h(r''))$, i.e.
\[
h(r'')=\mathcal G(f'')\circ \mathcal G(f')^{-1}(h(r'))=\mathcal G(f''\circ f'^{-1})(h(r')).
\]
Denote the parametrization of the isotopy $f=f''\circ f'^{-1}$ by $f_t$, $t\in[0,1]$. Let $\gamma_t=f_t(\gamma_{r'})$. By Corollary~\ref{cor:probes_verticalizing} there is a 2-parameter isotopy $\gamma_{s,t}$, $s,t\in[0,1]$, of region probes such that $\gamma_{s,0}=\gamma_0$, $\gamma_{s,1}=\gamma_1$, $\gamma_{0,t}=\gamma_t$ and all the region probes $\gamma_{1,t}$ are vertical. Extend the isotopy $\gamma_{s,t}$ to a spatial isotopy $f_{s,t}$ such that $f_{0,t}=f_t$. Denote $\tilde f_t=f_{1,t}$.

The isotopies $f=(f_t)$ and $\tilde f=(\tilde f_t)$ are homotopic. Since $\mathcal G$ is an $h$-coinvariant, $\mathcal G(f)=\mathcal G(\tilde f)$. For any $t\in [0,1]$, $\tilde f_t(\gamma)=\gamma_{r_t}$ is a vertical region probe of some region $r_t\in \mathcal R(\tilde h_t(D))$. Hence,
\[
r''=r_1=\mathcal R(\tilde f)(r_0)=\mathcal R(\tilde f)f(r'),
\]
then $h(r'')=\mathcal G(\tilde f)(h(r'))=\mathcal G(f)(h(r'))$.
Thus, the map $\phi$ is well defined.

Let $f\colon D_1\to D_2$ be a morphism in $\mathfrak D_s$, $\gamma_1\in\mathscr R^{ts}(D_1)$ and $\gamma_2=\mathscr R^{ts}(f)(\gamma_1)$. We can consider $f$ as a spatial isotopy such that $\gamma_2=f(\gamma_1)$. Let $g\colon D_2\to D'$ be a spatial isotopy which verticalizes $\gamma_2$, i.e. $g(\gamma_2)=\gamma_{r'}$, $r'\in\mathcal R(D')$. Then $\phi(\gamma_2)=\mathcal G(g)^{-1}(h(r'))$. Since the isotopy $g\circ f$ verticalizes $\gamma_1$,
\[
\phi(\gamma_1)=\mathcal G(g\circ f)^{-1}(h(r'))=\mathcal G(f)^{-1}\mathcal G(g)^{-1}(h(r')).
\]
Then $\phi(\mathscr R^{ts}(f)(\gamma_1))=\phi(\gamma_2)=\mathcal G(f)(\phi(\gamma_1))$. Thus, $\phi$ is a single-valued natural transformation between $\mathscr R^{ts}$ and $\mathcal G$.
\end{proof}

The following statement is analogous to Proposition~\ref{prop:arc_top_realizability}
\begin{proposition}\label{prop:region_top_realizability}
 For any finite subset of $\{\gamma_i\}_{i=1}^n\subset\mathscr R^{ts}(T)$ there is a morphism $f\colon T\to T'$ such that for any $i=1,\dots, n$, $f(\gamma_i)=\gamma_{r'_i}$ for some $r'_i\in\mathcal R(D')$ where $D'=p(T')$.
\end{proposition}
\begin{proof}
    We can isotope the curves $\gamma_i$ so that they can be verticalized simultaneously. Consider any verticalizing spatial isotopy $f$ for the region probes $\gamma_i$. Then $f$ is the required morphism.
\end{proof}

\begin{definition}\label{def:region_topological_weak_equivalence}
    A region probe $\gamma$ of a tangle $T$ is \emph{topologically weak equivalent} to a region probe $\gamma'$ of a tangle $T'$ ($\gamma\sim^{tw}_{\mathcal R}\gamma'$) if  there is an isotopy $\phi$ of $F\times I$ such that $\phi(T)=T'$ and $\phi(\gamma)=\gamma'$. Denote the set of topologically weak equivalence classes of the regions of the tangle $T$ by $\mathscr R^{tw}(T)$.
\end{definition}

For a tangle $T$, the motion group $M(F\times I, T)$ acts on the set $\mathscr R^{ts}(T)$ by composition: $f_t\times\gamma\mapsto f_1(\gamma)$, $f_t\in M(F\times I, T)$, $\gamma\in\mathscr R^{ts}(T)$.
Then $\mathscr R^{tw}(T)= \mathscr R^{ts}(T)/M(F\times I, T)$.

\begin{theorem}\label{thm:region_universal_homotopy_invariant}
The correspondence $T\mapsto\mathscr R^{tw}(T)$ is the universal h-invariant of the region functor $\mathcal R$.
\end{theorem}
The proof of the theorem is analogous to that of Theorem~\ref{thm:arc_universal_homotopy_invariant}.

\begin{figure}[ht]
\centering
  \includegraphics[width=0.4\textwidth]{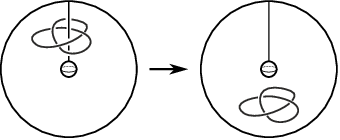}
  \caption{Splitting a spherical knot and a region probe}\label{pic:region_sphere_classical_knot}
\end{figure}

\begin{proposition}\label{prop:region_spherical_trivial}
Let $K\subset S^2\times I$ be a knot in the sphere. Then the set of regions in the weak sense $\mathscr R^{tw}(K)$ is trivial.
\end{proposition}
\begin{proof}
 Indeed, we can split the knot and a region probe (Fig.~\ref{pic:region_sphere_classical_knot}) and isotope them separately to a standard form.
\end{proof}

\begin{remark}\label{rem:region_classical_knots}
1. On the other hand, the set of regions in the strong sense $\mathscr R^{ts}(K)$ is infinite. For region probes $\gamma,\gamma'\in \mathscr R^{ts}(K)$ choose representative which don't intersect and close the curves $\gamma,\gamma'$ to an embedded circle $C$ with arcs connecting the ends of the probes in $S^2\times 0$ and $S^2\times 1$. Denote $d(\gamma,\gamma')=lk(C,K)$. Since the linking coefficient is a link-homotopy invariant, the number $d(\gamma,\gamma')$ depends only on isotopy classes of the probes.

For any region probe $\gamma$ and any $k\in\mathbb Z$ there exists a probe $\gamma_k$ such that $d(\gamma,\gamma_k)=k$. Then $\{\gamma_k\}_{k\in\mathbb Z}\subset \mathscr R^{ts}(K)$ is an infinite sequence of different region probes.

2. For any knot $K\subset \mathbb D^2\times I$ in the disk, the set $\mathscr R^{tw}(K)$ is infinite. For a region probe $\gamma\in\mathscr R^{tw}(K)$ we consider an invariant (the \emph{depth} of the region) $d(\gamma)=d(\gamma, \gamma_\infty)$ where $\gamma_\infty=x_\infty\times I$ for some fixed $x_\infty\in\partial\mathbb D^2$ and $d(\cdot,\cdot)$ is defined as above.

Then for any $k\in\mathbb K$ there exists $\gamma_k\in\mathscr R^{tw}(K)$ such that $d(\gamma_k)=k$. Thus, we have an infinite sequence of different region probes (in the weak sense).
\end{remark}

{\color{red}



}

\section{Semiarc}\label{sect:semiarc}

For a tangle $T$, there is a fibration of $N(T)$ whose fibers are meridian disks. When $T\in\Sigma_0$ we can assume that the meridian disks are orthogonal to the surface $F$. The orientation of $T$ and the right-hand rule determines an orientation of the meridians.

Let $T\in\Sigma_0$ be a tangle and $D$ its diagram.

\begin{definition}\label{def:semiarc_probe}
A \emph{semiarc probe} is an unknotted embedding
\[
\gamma\colon(I;\, 0;\, \frac{1}{2};\, 1)\to (F\times I;\, F\times 0;\, T;\, F\times 1)
\]
such that $\gamma\cap N(T)$ is a simple curve in a meridian disk.

For a semiarc $s\in\mathcal R(D)$, choose a point $x\in s\subset F$. The curve $\gamma_s= x\times I$ (parametrised so that $\gamma_s(\frac 12)\in T$) is called a \emph{vertical semiarc probe} of the semiarc $s$.

\begin{figure}[h]
\centering
  \includegraphics[width=0.4\textwidth]{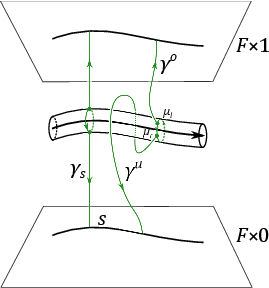}
  \caption{A vertical semiarc probe and a semiarc probe}\label{pic:semiarc_probe}
\end{figure}

For a semiarc probe $\gamma$, define the \emph{over-probe} $\gamma^o$ and the \emph{under-probe} $\gamma^u$ as the components of $\gamma\setminus N(T)$. We parameterize them as follows
\begin{gather*}
\gamma^o\colon(I; 0; 1)\to (M_T; \partial N(T); F\times 1),\\
\gamma^u\colon(I; 0; 1)\to (M_T; \partial N(T); F\times 0).
\end{gather*}
Then $\gamma^u(0)$ and $\gamma^o(0)$ are different points of the same meridian $\mu$. Denote the part of the meridian $\mu$ from $\gamma^o(0)$ to $\gamma^u(0)$ by $\mu_r$,  and the part from $\gamma^u(0)$ to $\gamma^o(0)$ by $\mu_l$. By definition of the semiarc probe, $(\gamma^u)^{-1}\mu_l\gamma^o$ and $(\gamma^u)^{-1}\mu_r^{-1}\gamma^o$ are unknotted in $F\times I$.
\end{definition}


We have a series of statements analogous to those for regions.

\begin{proposition}\label{prop:semiarc_affinity}
Let $s_1,s_2\in\mathcal{SA}(D)$. Then $s_1\approx_{\mathcal{SA}}s_2$ if and only if there is an embedded square $\Delta$ between the vertical semiarc probes $\gamma_{s_1}$ and $\gamma_{s_2}$:
\[
\Delta\colon(I\times I; 0\times I; 1\times I; I\times 0; I\times \frac 12;I\times 1)\hookrightarrow (F\times I; \gamma_{s_1}; \gamma_{s_2}; F\times 0; T; F\times 1),
\]
which intersects the tangle $T$ at one arc connecting $\gamma_{s_1}$ and $\gamma_{s_2}$.
\end{proposition}

We can assume that $\Delta\setminus N(T)=\Delta^o\cup\Delta^u$ consists of two disks which are parameterized as follows
{\small
\begin{gather*}
\Delta^o\colon(I\times I; 0\times I; 1\times I; I\times 0; I\times 1)\hookrightarrow (M_T; \gamma^o_{s_1}; \gamma^o_{s_2};  \partial N(T); F\times 1),\\
\Delta^u\colon(I\times I; 0\times I; 1\times I; I\times 0; I\times 1)\hookrightarrow (M_T; \gamma^u_{s_1}; \gamma^u_{s_2}; \partial N(T); F\times 0),
\end{gather*}
}
and for any $t\in[0,1]$  $\Delta^o(t,0)$ and $\Delta^u(t,0)$ belong to one meridian, and the restriction of the projection of $\partial N(T)\to T$ to $\Delta^o\cap\partial N(T)$ is injective.

For semiarc probes $\gamma, \gamma'$ such that there is an embedded square $\Delta$ between them, we denote $\gamma\approx^t_{\mathcal {SA}, T}\gamma'$ and say that  $\gamma$ and $\gamma'$ are \emph{topologically affined}.

Denote the transitive closure of $\approx^t_{\mathcal {SA}, T}$  by $\sim^{ts}_{\mathcal{SA}, T}$. The relation $\sim^{ts}_{\mathcal{SA}, T}$ in the set of semiarc probes is called  the \emph{topological strong equivalence}.

\begin{proposition}\label{prop:semiarc_top_affinity_closure}
    For any semiarc probes $\gamma, \gamma'$  $\gamma\sim^{ts}_{\mathcal{SA}, T}\gamma'$ if and only if $\gamma$ is isotopic to $\gamma'$ in $M_T$ (through semiarc probes).
\end{proposition}

\begin{proposition}\label{prop:semiarc_strong_transitivity}
    For any $s,s'\in\mathcal{SA}(D)$  $s\sim^s_{\mathcal{SA}, D} s'$ if and only if $\gamma_{s}\sim^{ts}_{\mathcal{SA}, T}\gamma_{s'}$.
\end{proposition}

\begin{remark}\label{rem:semiarc_strong_equivalence}
    1. The proof of Proposition~\ref{prop:semiarc_top_affinity_closure} goes along lines of the proof of Proposition~\ref{prop:region_top_affinity_closure} with the following modifications. An isotopy of semiarc probes $\phi=(\gamma_s)$ splits into isotopies of over- and under-probes $\phi^o$ and $\phi^u$ such that their restrictions to $\partial N(D)$ have support in two antipodal parallels $\lambda^o$ and $\lambda^u$ of $\partial N(D)$. Then for construction of the isotopy $\tilde\phi^o$ (and $\tilde\phi^u$), we choose the auxiliary probe $\hat\gamma^o$ ($\hat\gamma^o$) so that its end lies in $\lambda^o$ ($\lambda^u$).

    2. The proof of Proposition~\ref{prop:semiarc_strong_transitivity} is analogous the proof of Proposition~\ref{prop:region_strong_transitivity}, but we have to monitor that the corrections of the probes like in Fig.~\ref{pic:region_probe_correction} take place beyond the arc $\Delta_i\cap T$.

    After the curve $\tilde\gamma'$ and the square $\tilde\Delta'$ are constructed, we pull an arc of $T$ from $\gamma'_{n-1}$ to $\tilde\gamma'$ in $\tilde\Delta'$ and make $\tilde\gamma'$ a semiarc probe and $\tilde\Delta'$ an embedded square between two semiarc probes (Fig.~\ref{pic:semiarc_square_to_semiarc_square}).
\begin{figure}[ht]
\centering
  \includegraphics[width=0.6\textwidth]{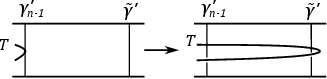}
  \caption{Turning an embedded square into a semiarc square}\label{pic:semiarc_square_to_semiarc_square}
\end{figure}

\end{remark}

\begin{definition}\label{def:semiarc_strong_tangle}
   The set of the isotopic classes of semiarc probes is called the \emph{(topological) set of semiarcs of the tangle $T$ (in the strong sense)} and denoted by ${\mathscr{SA}}^{ts}(T)$.
   Let $h^{ts}_{\mathcal{SA},D}$ be the projection $\mathcal{SA}(D)\to\mathscr{SA}^{ts}(T)$ determined by the formula $s\mapsto\gamma_s$.
\end{definition}

\begin{theorem}\label{thm:semiarc_universal_homotopy_coinvariant}
   The pair $(\mathscr{SA}^{ts}, h^{ts}_{\mathcal{SA}})$ is the universal $h$-coinvariant of the semiarc functor $\mathcal{SA}$.
\end{theorem}

\begin{proposition}\label{prop:semiarc_top_realizability}
 For any finite subset of $\{\gamma_i\}_{i=1}^n\subset\mathscr{SA}^{ts}(T)$ there is a morphism $f\colon T\to T'$ such that for any $i=1,\dots, n$, $f(\gamma_i)=\gamma_{s'_i}$ for some $s'_i\in\mathcal{SA}(D')$ where $D'=p(T')$.
\end{proposition}

\begin{definition}\label{def:semiarc_topological_weak_equivalence}
    A semiarc probe $\gamma$ of a tangle $T$ is \emph{topologically weak equivalent} to a semiarc probe $\gamma'$ of a tangle $T'$ ($\gamma\sim^{tw}_{\mathcal{SA}}\gamma'$) if  there is an isotopy $\phi$ of $F\times[0,1]$ such that $\phi(T)=T'$ and $\phi(\gamma)=\gamma'$. Denote the set of topologically weak equivalence classes of the semiarcs of the tangle $T$ by $\mathscr{SA}^{tw}(T)$.
\end{definition}

For a tangle $T$, the motion group $M(F\times I, T)$ acts on the set $\mathscr{SA}^{ts}(T)$ by composition: $f_t\times\gamma\mapsto f_1(\gamma)$, $f_t\in M(F\times I, T)$, $\gamma\in\mathscr{SA}^{ts}(T)$.
Then $\mathscr{SA}^{tw}(T)= \mathscr{SA}^{ts}(T)/M(F\times I, T)$.

\begin{theorem}\label{thm:semiarc_universal_homotopy_invariant}
The correspondence $T\mapsto\mathscr{SA}^{tw}(T)$ is the universal $h$-invariant of the semiarc functor $\mathcal{SA}$.
\end{theorem}

\begin{proposition}\label{prop:semiarc_spherical_trivial}
Let $K\subset S^2\times I$ be a knot in the thickened sphere. Then the set of semiarcs in the weak sense $\mathscr{SA}^{tw}(K)$ is trivial.
\end{proposition}

\begin{figure}[ht]
\centering
  \includegraphics[width=0.4\textwidth]{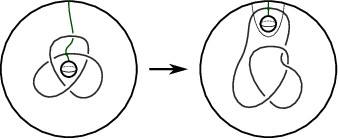}
  \caption{Isotopy of a semiarc probe of a spherical knot}\label{pic:semiarc_spherical_isotopy}
\end{figure}

\begin{proof}
Given semiarc probes $\gamma_i$ of isotopic knots $K_i$, $i=1,2$, in the thickened sphere, isotope them to a standard form as shown in Fig.~\ref{pic:semiarc_spherical_isotopy}. The long knots $K'_i$, $i=1,2$, which lie outside a standard neighbourhood $U$ of the probe ($U$ is marked by a dashed line in the figure), are equivalent. Then there is an isotopy between $K'_1$ and $K'_2$ with $U$ fixed. Then the composition of this isotopy with the standardizing isotopies maps $K_1$ to $K_2$ and  $\gamma_1$ to $\gamma_2$. Thus, $\gamma_1\sim_{\mathcal{SA}}^{tw}\gamma_2$.
\end{proof}

\begin{remark}
    The analogous result for knots in the thickened disk is not true. Like in Remark~\ref{rem:region_classical_knots}, one can construct an infinite series of non-equivalent semiarcs in $\mathscr{SA}^{tw}(K)$.
\end{remark}



\section{Crossing}\label{sect:crossing}

Let $T\in\Sigma_0$ be a tangle and $D$ its diagram.
\begin{definition}\label{def:crossing_probe}
A \emph{crossing probe} is an unknotted embedding
\[
\gamma\colon(I;\, 0;\, \{\frac 13, \frac 23\};\, 1)\to (F\times I;\, F\times 0;\, T;\, F\times 1)
\]
such that the intersection $\gamma\cap N(T)$ is two simple curves in two meridian disks, together with a framing of $\gamma$ in the interval $[\frac 13, \frac 23]$ which is collinear to $T$ in the points $\gamma(\frac 13), \gamma(\frac 13)\in\gamma\cap T$.

For a crossing $c\in\mathcal C(D)$,  its \emph{vertical crossing probe} is the curve $\gamma_c=c\times I$ with the framing that turns clockwise from the overcrossing to the undercrossing as shown in Fig.~\ref{pic:crossing_framing}.

\begin{figure}[h]
\centering
  \includegraphics[width=0.3\textwidth]{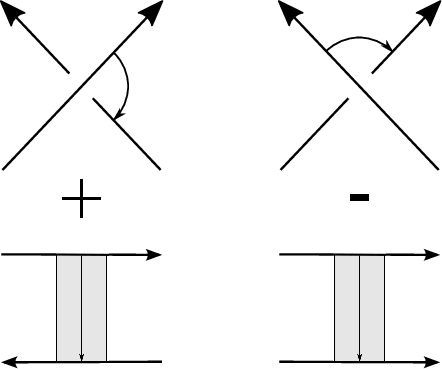}
  \caption{Crossing framing. The pictures below is the look at the crossing from the bottom (for the positive crossing) and from the right (for the negative crossing)}\label{pic:crossing_framing}
\end{figure}

For a crossing probe $\gamma$, the complement $\gamma\setminus N(T)$ splits into the \emph{over-probe} $\gamma^o$, \emph{mid-probe} $\gamma^m$ and \emph{under-probe} $\gamma^u$ (Fig.~\ref{pic:crossing_probe}). We parameterize them as embeddings
\begin{gather*}
\gamma^o\colon(I; 0; 1)\to (M_T; \partial N(T); F\times 1),\\
\gamma^m\colon(I; 0; 1)\to (M_T; \partial N(T); \partial N(T)),\\
\gamma^u\colon(I; 0; 1)\to (M_T;\partial N(T);  F\times 0)
\end{gather*}
such that
\begin{itemize}
\item $\gamma^o(0)$ and $\gamma^m(0)$ are different points of the same meridian $\mu^o$,
\item $\gamma^m(1)$ and $\gamma^u(0)$ are different points of the same meridian $\mu^u$,
\item $\gamma^m$ is framed by a transversal vector field which is collinear to the tangle at $\gamma^m(0)$ and $\gamma^m(1)$,
\item the curve $(\gamma^u)^{-1}\mu^u_l(\gamma^m)^{-1}\mu^o_l\gamma^o$ is unknotted in $F\times I$.
\end{itemize}
Here
\begin{itemize}
\item $\mu^o_r$ is the part of $\mu^o$ from $\gamma^o(0)$ to $\gamma^m(0)$,
\item $\mu^o_l$ is the part of $\mu^o$ from $\gamma^m(0)$ to $\gamma^o(0)$,
\item $\mu^u_r$ is the part of $\mu^u$ from $\gamma^m(1)$ to $\gamma^u(0)$,
\item $\mu^u_l$ is the part of $\mu^u$ from $\gamma^u(0)$ to $\gamma^m(1)$.
\end{itemize}
\begin{figure}[h]
\centering
  \includegraphics[width=0.5\textwidth]{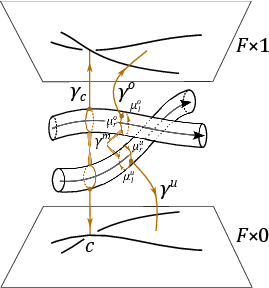}
  \caption{A crossing probe}\label{pic:crossing_probe}
\end{figure}
\end{definition}

\begin{proposition}\label{prop:crossing_affinity}
Let $c_1,c_2\in\mathcal{C}(D)$. Then $c_1\approx^h_{\mathcal{C}}c_2$ if and only if there is an embedded square $\Delta$ between the vertical crossing probes $\gamma_{c_1}$ and $\gamma_{c_2}$:
\[
\Delta\colon(I\times I; 0\times I; 1\times I; I\times 0; I\times \{\frac 13, \frac 23\};I\times 1)\hookrightarrow (F\times I; \gamma_{c_1}; \gamma_{c_2}; F\times 0; T; F\times 1),
\]
which intersects the tangle $T$ at two arc connecting $\gamma_{c_1}$ and $\gamma_{c_2}$, and is framed in the area  $\Delta([0,1]\times[\frac 13, \frac 23])$ by the horizontal framing $\partial_s D(s,t)$, so that the framing of the crossing probes is compatible with that of the embedded square as shown in Fig.~\ref{pic:crossing_square_framing_compatibility}.
\begin{figure}[h]
\centering
  \includegraphics[width=0.7\textwidth]{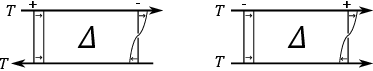}
  \caption{Framing compatibility condition for the embedded square (depending the tangle arcs and the signs of the crossings)  }\label{pic:crossing_square_framing_compatibility}
\end{figure}

In other words, the difference between the framings of a crossing probe $\gamma_{c_i}$ and that of the square is $\frac\pi 2(1-sgn(c_i))$ when the top and the bottom arcs of the tangle define the same orientation of the square (Fig.~\ref{pic:crossing_square_framing_compatibility} left), and the difference is $\frac\pi 2(1+sgn(c_i))$ when the top and the bottom arcs of the tangle define different orientations of the square (Fig.~\ref{pic:crossing_square_framing_compatibility} right).
\end{proposition}


\begin{proof}
If $c_1\approx_{\mathcal C}^h c_2$ then there exists a sequence of moves $f$ which transforms $c_1$ and $c_2$ to the configuration in Fig.~\ref{pic:homotopy_crossing_relation}. We can suppose that  $f$ does not move points $c_1$ and $c_2$. Let $\phi$ be the corresponding isotopy of the tangle $T$. Then $\phi$ extends to a spatial isotopy fixed on $\{c_1,c_2\}\times I$. For the configuration in Fig.~\ref{pic:homotopy_crossing_relation}, we can construct an embedded square $\Delta_0$ which satisfies the conditions of the proposition. Then take $\Delta=\phi^{-1}(\Delta_0)$.

Let $\Delta$ be a square between the vertical crossing probes $\gamma_{c_1}=c_1\times I$ and $\gamma_{c_2}=c_2\times I$. Isotope the square $\Delta$ rel $\gamma_{c_1}\cup\gamma_{c_2}$ to a square $\Delta_0= \alpha\times I$,  for some curve $\alpha\subset F$ connecting $c_1$ and $c_2$. Disturb $\Delta_0$ to get the configuration in Fig.~\ref{pic:homotopy_crossing_relation}. Extend the isotopy from $\Delta$ to the disturbed $\Delta_0$ to a spatial isotopy $\phi$ of $F\times I$. Then $\phi$ induces a sequence of Reidemeister moves on the tangle $T$, which places the crossings $c_1$ and $c_2$ in the configuration of Fig.~\ref{pic:homotopy_crossing_relation}. Thus, $c_1\approx_{\mathcal C}^h c_2$.
\end{proof}

For crossing probes $\gamma, \gamma'$ such that there is a framed embedded square $\Delta$ between them, we denote $\gamma\approx^t_{\mathcal {C}, T}\gamma'$ and say that  $\gamma$ and $\gamma'$ are \emph{topologically affined}.

Denote the transitive closure of $\approx^t_{\mathcal {C}, T}$  by $\sim^{ts}_{\mathcal{C}, T}$. The relation $\sim^{ts}_{\mathcal{C}, T}$ in the set of crossing probes is called  the \emph{topological strong equivalence}.

\begin{proposition}\label{prop:crossing_top_affinity_closure}
    For any crossing probes $\gamma, \gamma'$,  $\gamma\sim^{ts}_{\mathcal{C}, T}\gamma'$ if and only if $\gamma$ is isotopic to $\gamma'$ in $F\times I$ (through crossing probes).
\end{proposition}

\begin{lemma}\label{lem:crossing_top_affinity_closure}
 Let two crossing probes $\gamma_1$ and $\gamma_2$ be connected by an (unframed) embedded square
 \[
\Delta\colon(I\times I; 0\times I; 1\times I; I\times 0; I\times \{\frac 13, \frac 23\};I\times 1)\hookrightarrow (F\times I; \gamma_{1}; \gamma_{2}; F\times 0; T; F\times 1).
\]
Assume that the framings of $\gamma_1$ and $\gamma_2$ are compatible in the sense that the framings can be extended to a framing on $\Delta([0,1]\times[\frac 13,\frac 23])$ which is tangent to $T$ on $\Delta([0,1]\times\{\frac 13,\frac 23\})$ and is transversal to $\partial_t\Delta(s,t)$. Then $\gamma_1\sim_{\mathcal C,T}^{ts}\gamma_2$.
\end{lemma}
\begin{proof}
If the framings of $\gamma_1$ and $\gamma_2$ are compatible with the framing of the embedded square $\Delta$ then $\gamma_1\approx_{\mathcal C, T}^t\gamma_2$, hence, $\gamma_1\sim_{\mathcal C,T}^{ts}\gamma_2$. Assume that the framings of the probes and the square are not compatible. Assume that the sign of the probes is positive and the arcs of $T$ define the same orientation on $\Delta$ (Fig.~\ref{pic:crossing_uncompatible_square} left).

\begin{figure}[h]
\centering
  \includegraphics[width=0.5\textwidth]{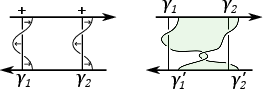}
  \caption{Equivalence of crossing probes connected by an unframed embedded square. The framed embedded square between the crossing probes $\gamma'_1$ and $\gamma'_2$ is shown by green color.}\label{pic:crossing_uncompatible_square}
\end{figure}

The framing of the probe $\gamma_1$ can be considered as a band connecting $\gamma_1$ with a probe $\gamma_1'$. The horizontal framing of the band is compatible with the framing of $\gamma_1$. Hence, $\gamma_1\approx_{\mathcal C, T}^t\gamma_1'$. Analogously, one defines a probe $\gamma_2'$ affined with $\gamma_2$. Isotope the probes $\gamma_1'$ and $\gamma_2'$ along $\Delta$ to make them linked (Fig.~\ref{pic:crossing_uncompatible_square} right). The linking coefficient is chosen so that the embedded square between $\gamma_1'$ and $\gamma_2'$ has compatible framing with the probes. Then we get the sequence of probes
\[
\gamma_1\approx_{\mathcal C, T}^t\gamma_1'\approx_{\mathcal C, T}^t\gamma_2'\approx_{\mathcal C, T}^t\gamma_2.
\]
Then $\gamma_1\sim_{\mathcal C,T}^{ts}\gamma_2$.
\end{proof}

\begin{proof}[Proof of Proposition~\ref{prop:crossing_top_affinity_closure}]
    The necessity is evident.

    Assume the crossing probe $\gamma$ is isotopic to the crossing probe $\gamma'$. Like in the proofs of Propositions~\ref{prop:region_top_affinity_closure} and ~\ref{prop:semiarc_top_affinity_closure}, we can construct a sequence of crossing probes between $\gamma$ and $\gamma'$ such that any two consequent probes are connected by an (unframed) embedded square. By Lemma~\ref{lem:crossing_top_affinity_closure}, we get a sequence of crossings
\[
\gamma = \gamma_0\sim_{\mathcal C, T}^{ts}\gamma_1\sim_{\mathcal C, T}^{ts}\cdots\sim_{\mathcal C, T}^{ts}\gamma_n=\gamma'.
\]
Then $\gamma\sim_{\mathcal C,T}^{ts}\gamma'$.
\end{proof}

We can formulate statements analogous to those for semiarcs or regions. The proofs are analogous.

\begin{proposition}\label{prop:crossing_strong_transitivity}
    For any $c,c'\in\mathcal{C}(D)$  $c\sim^s_{\mathcal{C}, D} c'$ if and only if $\gamma_{c}\sim^{ts}_{\mathcal{C}, T}\gamma_{c'}$.
\end{proposition}

\begin{definition}\label{def:crossing_strong_tangle}
   The set of the isotopic classes of crossing probes is called the \emph{(topological) set of crossings of the tangle $T$ (in the strong sense)} and denoted by ${\mathscr{C}}^{ts}(T)$.
   Let $h^{ts}_{\mathcal{C},D}$ be the projection $\mathcal{C}(D)\to\mathscr{C}^{ts}(T)$ determined by the formula $s\mapsto\gamma_s$.
\end{definition}

\begin{theorem}\label{thm:crossing_universal_homotopy_coinvariant}
   The pair $(\mathscr{C}^{ts}, h^{ts}_{\mathcal{C}})$ is the universal $h$-coinvariant of the crossing functor $\mathcal{C}$.
\end{theorem}

\begin{proposition}\label{prop:crossing_top_realizability}
 For any finite subset of $\{\gamma_i\}_{i=1}^n\subset\mathscr{C}^{ts}(T)$ there is a morphism $f\colon T\to T'$ such that for any $i=1,\dots, n$, $f(\gamma_i)=\gamma_{c'_i}$ for some $c'_i\in\mathcal{C}(D')$ where $D'=p(T')$.
\end{proposition}

\begin{definition}\label{def:crossing_topological_weak_equivalence}
    A crossing probe $\gamma$ of a tangle $T$ is \emph{topologically weak equivalent} to a crossing probe $\gamma'$ of a tangle $T'$ ($\gamma\sim^{tw}_{\mathcal{C}}\gamma'$) if  there is an isotopy $\phi$ of $F\times I$ such that $\phi(T)=T'$ and $\phi(\gamma)=\gamma'$. Denote the set of topologically weak equivalence classes of the crossings of the tangle $T$ by $\mathscr{C}^{tw}(T)$.
\end{definition}

For a tangle $T$, the motion group $M(F\times I, T)$ acts on the set $\mathscr{C}^{ts}(T)$ by composition: $f_t\times\gamma\mapsto f_1(\gamma)$, $f_t\in M(F\times I, T)$, $\gamma\in\mathscr{C}^{ts}(T)$.
Then $\mathscr{C}^{tw}(T)= \mathscr{C}^{ts}(T)/M(F\times I, T)$.

\begin{theorem}\label{thm:crossing_universal_homotopy_invariant}
The correspondence $T\mapsto\mathscr{C}^{tw}(T)$ is the universal $h$-invariant of the crossing functor $\mathcal{C}$.
\end{theorem}

\subsection{\texorpdfstring{$\Omega_2$}{R2}-equivalence of crossings and wrapping}\label{sect:crossing_R2_equivalence}

Theorem~\ref{thm:crossing_universal_homotopy_coinvariant} states that the universal $h$-coinvariant of crossings corresponds to the isotopy classes of crossing probes. The middle part of a crossing probe is framed, and we can consider an unframed crossing probe by forgetting the framing. What corresponds to the isotopy classes of unframed crossing probe? Below we describe the combinatorial counterparts of these classes.

\begin{definition}\label{def:crossing_probe_unframed}
An \emph{unframed crossing probe} is an unknotted embedding
\[
\gamma\colon(I;\, 0;\, \{\frac 13, \frac 23\};\, 1)\to (F\times I;\, F\times 0;\, T;\, F\times 1)
\]
such that the intersection $\gamma\cap N(T)$ is two simple curves in two meridian disks.

Denote the set of isotopy classes of unframed crossing probes of the tangle $T$ by $\mathscr C^{uf}(T)$, and let $h^{uf}(D)\colon\mathcal C(D)\to\mathscr C^{uf}(T)$ be the unframed vertical crossing probe map.
\end{definition}

\begin{definition}\label{def:crossing_R2_equivalence}
Consider the crossing functor $\mathcal C\colon\mathfrak D_s\to Rel$.

 The \emph{$\Omega_2$-strong equivalence relation} associated with the crossing functor  is the family of equivalence relations $\sim_{\mathcal C,D}^{s\Omega_2}$ on the sets $\mathcal C(D)$, $D\in Ob(\mathfrak D_s)$, generated by the rules:
\begin{itemize}
    \item $c_1\sim_{\mathcal C,D}^{s\Omega_2} c_2$ for any two crossing $c_1$, $c_2$ to which a second Reidemeister move can be applied;
    \item for any $f\in Mor_{\mathfrak D_s}(D,D')$, $c_1,c_2\in \mathcal C(D)$ and $c'_1\in \mathcal C(f)(c_1)$, $c'_2\in \mathcal C(f)(c_2)$  $c_1\sim_{\mathcal C,D}^{s\Omega_2} c_2$ implies $c'_1\sim_{\mathcal C,D'}^{s\Omega_2} c'_2$.
\end{itemize}

The \emph{$\Omega_2$-weak equivalence relation} associated with the crossing functor is the equivalence relation $\sim^{w\Omega_2}_{\mathcal C}$ on the set $\bigsqcup_{D\in Ob(\mathfrak D_s)}\mathcal C(D)$ generated by the rules
\begin{itemize}
    \item $c_1\sim_{\mathcal C}^{w\Omega_2} c_2$ for any two crossing $c_1$, $c_2$ to which a second Reidemeister move can be applied;
    \item for any morphism $f\in Mor_{\mathfrak D_s}(D,D')$,  $c\in \mathcal C(D)$ and  $c'\in \mathcal C(f)(c)$ one has $c\sim^{w\Omega_2}_{\mathcal C} c'$.
\end{itemize}
\end{definition}

Let us define another equivalence relation.

\begin{figure}[h]
\centering
  \includegraphics[width=0.6\textwidth]{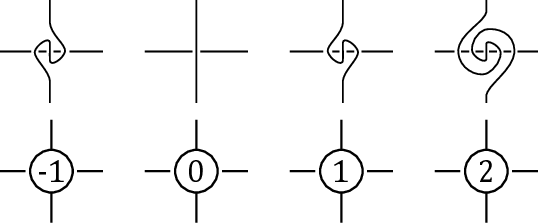}
  \caption{Crossing wrapping}\label{pic:crossing_wrappings}
\end{figure}

\begin{definition}\label{def:crossing_wrapping}
    Let $D$ be a tangle diagram and  $c\in\mathcal C(D)$. For an integer $n$, construct the \emph{$n$-wrapping} of $c$ by rotating the overcrossing by the angle $\pi n$ counterclockwise (see Fig.~\ref{pic:crossing_wrappings}). Denote the obtained diagram by $W_{c,n}(D)$ and the wrapped crossing in the new diagram by $W_n(c)$.
\end{definition}

Note that $W_{c,0}(D)=D$ and $W_{0}(c)=c$. The sign of the wrapped crossing is $sgn(W_n(c))=(-1)^n sgn(c)$.

\begin{remark}\label{rem:crossing_probe_wrapping}
    One can define the $n$-wrapping on the crossing probes as follows: for a crossing probe $\gamma\in\mathscr C^s(T)$, the wrapping $W_n(\gamma)$ is obtained from $\gamma$ by adding $n$ half-turns to the framing of the mid-probe of $\gamma$.
\end{remark}

\begin{definition}\label{def:crossing_W_equivalence}
 The \emph{$W$-strong (homotopical) equivalence relation} associated with the crossing functor is a family of equivalence relations $\sim_{\mathcal C,D}^{sW}$ on the sets $\mathcal C(D)$, $D\in Ob(\mathfrak D_s)$, generated by the rules:
\begin{itemize}
    \item $c_1\sim_{\mathcal C,D}^{sW} c_2$ for the configuration in Fig.~\ref{pic:homotopy_crossing_relation};
    \item for any $f\in Mor_{\mathfrak D_s}(D,D')$, $c_1,c_2\in \mathcal C(D)$ and $c'_1\in \mathcal C(f)(c_1)$, $c'_2\in \mathcal C(f)(c_2)$  $c_1\sim_{\mathcal C,D}^{sW} c_2$ implies $c'_1\sim_{\mathcal C,D'}^{sW} c'_2$;
    \item for any $n\in\mathbb Z$, $c_1\sim_{\mathcal C,D}^{sW} c_2$ implies $W_n(c_1)\sim_{\mathcal C,W_{c_1,n}(D)}^{sW} \mathcal C(W_{c_1,n})(c_2)$.
\end{itemize}

The \emph{$W$-weak equivalence relation} associated with the crossing functor is the equivalence relation $\sim^{wW}_{\mathcal C}$ on the set $\bigsqcup_{D\in Ob(\mathfrak D_s)}\mathcal C(D)$ generated by the rules
\begin{itemize}
    \item for any crossing $c$ and any $n\in\mathbb Z$, $c\sim_{\mathcal C}^{wW} W_n(c)$;
    \item for any morphism $f\in Mor_{\mathfrak D_s}(D,D')$,  $c\in \mathcal C(D)$ and  $c'\in \mathcal C(f)(c)$ one has $c\sim^{wW}_{\mathcal C} c'$.
\end{itemize}
\end{definition}

\begin{proposition}\label{prop:crossing_equivalence_R2_vs_W}
    1. The equivalence relations $\sim_{\mathcal C,D}^{sW}$ and $\sim_{\mathcal C,D}^{s\Omega_2}$ coincide for any diagram $D$;

    2. The equivalence relations $\sim_{\mathcal C}^{wW}$ and $\sim_{\mathcal C}^{w\Omega_2}$ coincide.
\end{proposition}

\begin{proof}
    Let us prove that $a\sim_{\mathcal C,D}^{sW}b$ implies $a\sim_{\mathcal C,D}^{s\Omega_2}b$. It is sufficient to consider a case when the third rule in the definition of the $W$-strong equivalence is applied. Let $a\sim_{\mathcal C,D}^{sW}b$, $b'''=W_n(b)$ and $a'=\mathcal C(W_{b,n})(a)$ (see Fig.~\ref{pic:crossing_R2_to_W}). By induction, we can assume that $a\sim_{\mathcal C,D}^{s\Omega_2}b$. The wrapping $W_{b,n}$ can be presented as a sequence of second Reidemeister moves, so that $\mathcal C(W_{b,n})(b)=b'$. Then $a'\sim_{\mathcal C,W_{b,n}(D)}^{s\Omega_2}b'$. On the other hand, we have $b'\sim_{\mathcal C,W_{b,n}(D)}^{s\Omega_2}b''$ and $b''\sim_{\mathcal C,W_{b,n}(D)}^{s\Omega_2}b'''$. Hence, $a'\sim_{\mathcal C,W_{b,n}(D)}^{s\Omega_2}b'''$.
\begin{figure}[h]
\centering
  \includegraphics[width=0.8\textwidth]{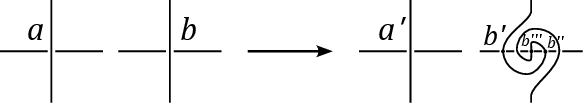}
  \caption{$W$-equivalence implies $\Omega_2$-equivalence}\label{pic:crossing_R2_to_W}
\end{figure}

Let us prove that $a\sim_{\mathcal C,D}^{s\Omega_2}b$ implies $a\sim_{\mathcal C,D}^{sW}b$. It is sufficient to consider the case when a second Reidemeister move can be applied to $a$ and $b$ (Fig.~\ref{pic:crossing_W_to_R2} left).

Wrap the crossing $b$. In the diagram $W_{b,1}(D)$ the crossings $a'$ and $b'$ are $W$-strong equivalent by the first rule of the definition. Apply the wrapping $W_{-1}$ to $b'$. Then the crossing $a''=\mathcal C(W_{b',-1})(a')$ and $b''=W_{-1}(b')$ are $W$-strong equivalent. By applying decreasing second Reidemeister moves, we find that $a\sim_{\mathcal C,D}^{sW}b$.
\begin{figure}[h]
\centering
  \includegraphics[width=0.6\textwidth]{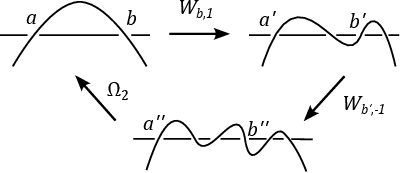}
  \caption{$\Omega_2$-equivalence implies $W$-equivalence}\label{pic:crossing_W_to_R2}
\end{figure}

The second statement of the proposition is proved analogously.
\end{proof}

\begin{theorem}\label{thm:crossing_universal_unframed_coinvariant}
For any $h$-coinvariant $(\mathcal G,h)$ of the crossing functor $\mathcal C$ such that for any $D\in\mathfrak Ob(D_s)$ and any $c_1\sim_{\mathcal C,D}^{s\Omega_2}c_2$ one has $h(D)(c_1)=h(D)(c_2)$, there is a unique single-valued natural transformation $\phi\colon\mathscr C^{uf}\Rightarrow\mathcal G$ such that $h=\phi\circ h^{uf}$.
\end{theorem}

\begin{proof}
By the universality of $\mathscr C$, there is a unique single-valued natural transformation $\tilde\phi\colon \mathscr C\Rightarrow\mathcal G$. The functor $\mathcal G$ is $\Omega_2$-invariant, hence, it is wrapping-invariant, i.e. the image $h\circ\tilde\phi(\gamma)$ of a crossing probe $\gamma$ does not depend on the framing. Thus, $\tilde\phi$ descends to a single-valued natural transformation $\phi\colon\mathscr C^{uf}\Rightarrow\mathcal G$.
\end{proof}

\begin{figure}[h]
\centering
  \includegraphics[width=0.6\textwidth]{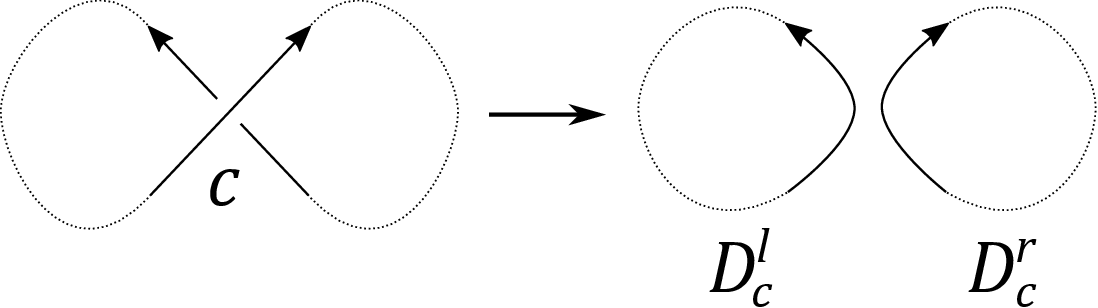}
  \caption{Diagram halves at a pure crossing}\label{pic:knot_halves}
\end{figure}

\begin{definition}\label{def:crossing_wrapping_index}
Let $c$ be a pure crossing in a tangle diagram $D$, i.e. a self-intersection of a component of the tangle. Smooth the diagram at the crossing $c$ according the orientation and consider the halves of the diagram (Fig.~\ref{pic:knot_halves}). The linking number $lk(D^l_c,D^r_c)$ is called the \emph{wrapping index} of the crossing $c$ and denoted by $wr(c)$.
\end{definition}

\begin{proposition}\label{prop:crossing_strong_equivalence_wrapping_index}
    1. Let $c_1, c_2\in\mathcal C(D)$ be crossings of a tangle diagram $D$. Then $c_1\sim_{\mathcal C,D}^{sh} c_2$ if and only if $c_1\sim_{\mathcal C,D}^{sW} c_2$, $sgn(c_1)=sgn(c_2)$ and $wr(c_1)=wr(c_2)$.

    2. For any crossings $c\in\mathcal C(D)$ and $c'\in\mathcal C(D')$, $c\sim_{\mathcal C}^w c'$ if and only if $c_1\sim_{\mathcal C}^{wW} c_2$, $sgn(c_1)=sgn(c_2)$ and $wr(c_1)=wr(c_2)$.
\end{proposition}
\begin{proof}
By induction, one can prove that  $c_1\sim_{\mathcal C,D}^{sh} c_2$ implies  $wr(c_1)=wr(c_2)$. Thus, the necessity is proved.

Let $c_1\sim_{\mathcal C,D}^{sW} c_2$, $sgn(c_1)=sgn(c_2)$ and $wr(c_1)=wr(c_2)$. By definition of the relation $\sim_{\mathcal C,D}^{sW}$, $c_1\sim_{\mathcal C,D}^{sW} c_2$ implies $c_1\sim_{\mathcal C,W_{c_2,n}(D)}^{s} W_n(c_2)$ for some integer $n$. Since $sgn(c_1)=sgn(c_2)$, $n=2k$ is even. Since $wr(c_1)=wr(W_{2k}(c_2))=wr(c_2)+k$,  we get $n=2k=0$. Then $c_1\sim_{\mathcal C,D}^{s} c_2$.

The second statement is proved analogously.
\end{proof}

\begin{remark}\label{rem:crossing_wrapping_monodromy}
The proposition above states that the wrapping of a pure crossing produces nonequivalent in the weak sense. The analogous result for mixed crossings is wrong as Fig.~\ref{pic:crossing_mixed_unwrapping} shows.
\begin{figure}[h]
\centering
  \includegraphics[width=0.7\textwidth]{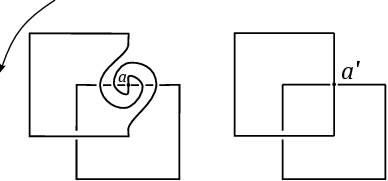}
  \caption{Unwrapping of a mixed crossing}\label{pic:crossing_mixed_unwrapping}
\end{figure}

For a crossing $c\in\mathcal C(D)$, we can define the \emph{wrapping monodromy} of the crossing $c$ as the minimal positive number $k$ such that $c\sim_{\mathcal C}^w W_{2k}(c)$. We suppose the monodromy to be $0$ if there is no such a number. Thus, pure crossings have zero wrapping monodromy. The wrapping monodromy of the crossing $a$ in Fig.~\ref{pic:crossing_mixed_unwrapping} is $1$. A question is: for an integer $n$ construct a crossing with the given wrapping monodromy.
\end{remark}

\section{Mid-crossing and tangle transformations}\label{sect:midcrossing}

By the definition of the crossing functor, the crossings which participate in a third Reidemeister move, have no relation to the corresponding crossings of the transformed tangle diagram. We can change this situation and allow a crossing to survive when an arc of the tangle passes over or under the crossing. Thus, we come to a new functor on tangle diagrams.

\begin{definition}\label{def:midcrossing_functor}
The \emph{midcrossing functor} is a functor $\mathcal{MC}\colon\mathfrak D_s\to Rel$ is a functor which extends the crossing functor $\mathcal C$ by relations $c\mapsto c'$ in third Reidemeister undermoves and overmoves.

\begin{figure}[h]
\centering
  \includegraphics[width=0.7\textwidth]{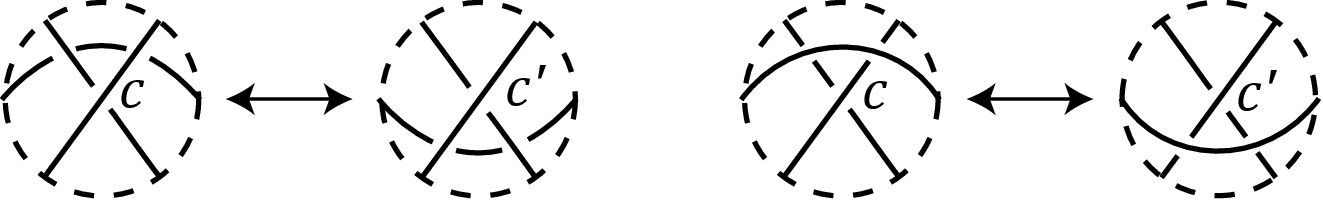}
  \caption{Third Reidemeister undermove and overmove}\label{pic:mdcrossing_moves}
\end{figure}
\end{definition}

Note that for any diagram $D$ $\mathcal{MC}(D)=\mathcal C(D)$ and for any Reidemeister move $f\colon D\to D'$, $\mathcal C(f)\subset\mathcal{MC}(f)$ as relations, i.e. subsets in $\mathcal{MC}(D)\times\mathcal{MC}(D')$.

\begin{definition}\label{def:midcrossing_probe}
A \emph{midcrossing probe} of a tangle $T$ in $F\times I$ is a framed embedded arc $\gamma\colon(I; 0; 1)\to (M_T; \partial N(T);  \partial N(T))$
such that the ends $\gamma(0)$ and $\gamma(1)$ belong to different meridians.

Let $T\in\Sigma_0$ and $D=p(T)$ its diagram.
For a crossing $c\in\mathcal{MC}(D)$, its \emph{vertical midcrossing probe} $\gamma_c$ is the middle component $\gamma^m$ of $c\times I\setminus N(T)$ oriented downwards and framed as shown in Fig.~\ref{pic:crossing_framing}.
\end{definition}

A midcrossing probe can be thought of as just a band glued to the tangle; one of the glued interval in the band boundary is the overcrossing and the other is the undercrossing.

For midcrossings, we can formulate statements analogous to those for arcs, regions, semiarcs and crossings.

\begin{definition}\label{def:midcrossing_affinity}
Two midcrossings $x_1,x_2\in \mathcal C(D')$, $D'\in Ob(\mathfrak D_s)$, are called \emph{homotopically affined} if there exists morphism $f\colon D'\to D$ and $c_i\in \mathcal{MC}(f)(x_i)$, $i=1,2$, such that the crossings $c_1$ and $c_2$ form the configuration shown in Fig.~\ref{pic:homotopy_crossing_relation}. In this case we write $x_1\approx_{\mathcal{MC},D}^h x_2$.
\end{definition}

\begin{proposition}\label{prop:midcrossing_affinity}
Let $c_1,c_2\in\mathcal{MC}(D)$. Then $c_1\approx^h_{\mathcal{MC},D}c_2$ if and only if there is an embedded square $\Delta$ between the vertical crossing probes $\gamma_{c_1}$ and $\gamma_{c_2}$:
\[
\Delta\colon(I\times I; 0\times I; 1\times I; I\times \{0,1\})\hookrightarrow (F\times I; \gamma_{c_1}; \gamma_{c_2}; T),
\]
which is framed by the horizontal framing $\partial_s D(s,t)$, so that the framing of the crossing probes is compatible with that of the embedded square as shown in Fig.~\ref{pic:crossing_square_framing_compatibility}.
\end{proposition}

\begin{proof}
    The necessity is proved like in Proposition~\ref{prop:crossing_affinity}.

    Let $\gamma_{c_1}$ and $\gamma_{c_2}$ are connected by a framed embedded square $\Delta$. We can consider $\Delta$ as a band or a framed arc with fixed ends. Ignoring the framing, this arc is homotopic (with fixed ends) to an arc whose projection to $F$ is a simple curve $\alpha$. Then we can realize this homotopy by a sequence of isotopies and crossing changes. The crossing changes can be resolved by isotopies as shown in Fig.~\ref{pic:midcrossing_probe_crossing_change}.

\begin{figure}[h]
\centering
  \includegraphics[width=0.7\textwidth]{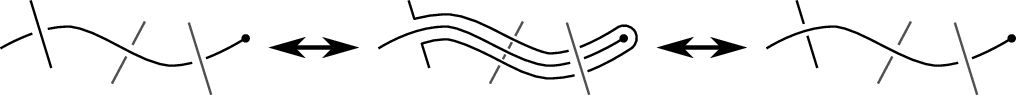}
  \caption{An isotopy which realizes a crossing change}\label{pic:midcrossing_probe_crossing_change}
\end{figure}

The constructed isotopy moves $\Delta$ to a band which can differ from the vertical band over $\alpha$ by some twists. We can remove these twists as shown in Fig.~\ref{pic:midcrossing_probe_unframing}.

\begin{figure}[h]
\centering
  \includegraphics[width=0.7\textwidth]{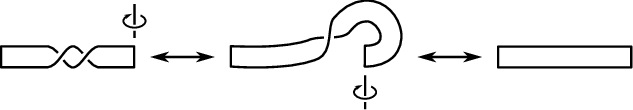}
  \caption{Removing a twist of the embedded square}\label{pic:midcrossing_probe_unframing}
\end{figure}

Then we get an isotopy $f$ which verticalizes the embedded square. By disturbing the isotopy, we can assume the crossings $c'_1=\mathcal{MC}(f)(c_1)$ and $c'_2=\mathcal{MC}(f)(c_2)$ form a configuration like in Fig.~\ref{pic:homotopy_crossing_relation}. Hence, $c_1\approx^h_{\mathcal{MC},D}c_2$.
\end{proof}

The following statement is proved analogously to Proposition~\ref{prop:crossing_strong_transitivity}.
\begin{proposition}\label{prop:midcrossing_strong_transitivity}
    For any $c,c'\in\mathcal{MC}(D)$  $c\sim^s_{\mathcal{MC}, D} c'$ if and only if $\gamma_{c}$ is isotopic to $\gamma_{c'}$ (though midcrossing probes).
\end{proposition}

\begin{definition}\label{def:midcrossing_set_of_tangle}
 Denote the isotopy classes of midcrossing probes of a tangle $T$ by $\mathscr{MC}^{s}(T)$, and call it the \emph{set of midcrossings in the strong sense} of the tangle $T$.  Denote the set of the orbits of $\mathscr{MC}^{s}(T)$ by the action of the motion group $M(F\times I, T)$ by $\mathscr{MC}^{w}(T)$  and call it the \emph{set of midcrossings in the weak sense} of the tangle. There are natural maps $h^s(T)\colon\mathcal{MC}(T)\to\mathscr{MC}^s(T)$ and $h^w(T)\colon\mathcal{MC}(T)\to\mathscr{MC}^w(T)$ which assign to a midcrossing the corresponding vertical probe.
\end{definition}

The following statement is analogous to Theorems~\ref{thm:crossing_universal_homotopy_coinvariant} and~\ref{thm:crossing_universal_homotopy_invariant}.

\begin{theorem}\label{thm:midcrossing_universal_invariant}
1. The pair $(\mathscr{MC}^{s}, h^{s}_{\mathcal{MC}})$ is the universal $h$-coinvariant of the midcrossing functor $\mathcal{MC}$.

2. The pair $(\mathscr{MC}^{w}, h^{w}_{\mathcal{MC}})$ is the universal $h$-invariant of the midcrossing functor $\mathcal{MC}$.
\end{theorem}

\subsection{Transformations}

\begin{figure}[h]
\centering
  \includegraphics[width=0.5\textwidth]{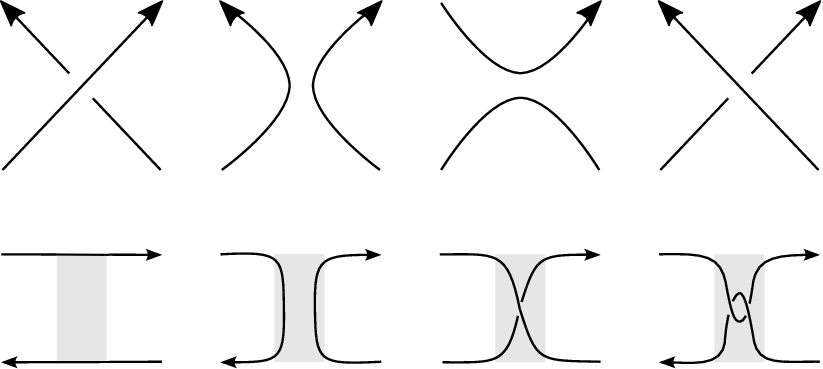}
  \caption{Transformation at a midcrossing: identity, oriented smoothing, non-oriented smoothing, crossing change}\label{pic:midcrossing_transformations}
\end{figure}

Given a midcrossing probe (in the weak sense) $\gamma$ of a tangle $T$, one can define several transformations which produce other tangles (Fig.~\ref{pic:midcrossing_transformations}).

Each transformation rule defines a map from $\mathscr{MC}^w(T)$ (hence, from $\mathscr{MC}^s(T)$) to the set of isotopy classes of tangles.

We can use transformation rule to introduce special types of midcrossings.

\begin{example}[Nugatory and cosmetic midcrossings] 

    An (unoriented) positive midcrossing probe $\gamma$ of a knot $K$ is called \emph{nugatory} if there is a sphere $S$ such that $\gamma\subset S$, $K\cap S=\partial \gamma$ and the framing is transversal to $S$ (Fig.~\ref{pic:midcrossing_nugatory}). In other words, the link obtained by the oriented smoothing at $\gamma$ is split.

    The set of nugatory midcrossing probes in $\mathscr{MC}^s(T)$ is finite.

\begin{figure}[h]
\centering
  \includegraphics[width=0.4\textwidth]{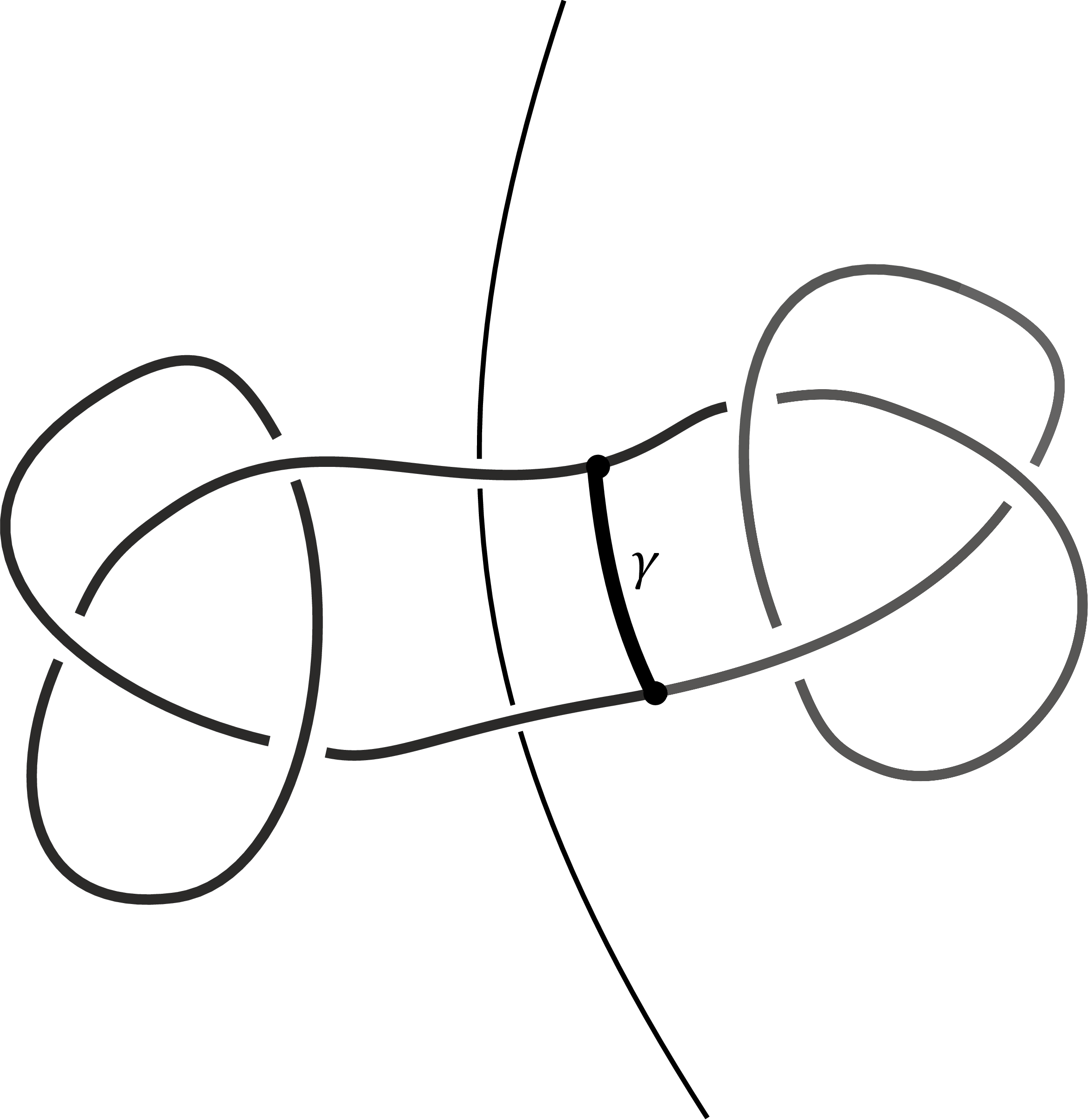}
  \caption{Nugatory midcrossing}\label{pic:midcrossing_nugatory}
\end{figure}

An (unoriented) midcrossing $\gamma$ of a knot $K$ is \emph{cosmetic} if the crossing change at $\gamma$ produces a knot isotopic to $K$.

The cosmetic conjecture~\cite[Problem 1.58]{Kirby} states that any cosmetic (mid)crossing is nugatory.
\end{example}

\begin{example}[Ribbon midcrossings]
A (classical) knot $K$ is called \emph{ribbon} if it spans a disc which has only ribbon self-intersections (Fig.~\ref{pic:ribbon_intersection}).

\begin{figure}[h]
\centering
  \includegraphics[width=0.35\textwidth]{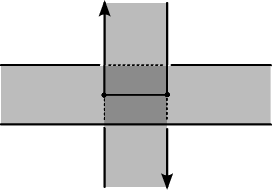}
  \caption{Ribbon intersection}\label{pic:ribbon_intersection}
\end{figure}
Note that any ribbon self-intersection produces an (unoriented) positive midcrossing.

A midcrossing tuple $(\gamma_1,\dots,\gamma_k)$ is called \emph{ribbon} if they are pairwise not interlaced and the oriented smoothing at all the midcrossings $\gamma_i$ produces a trivial link.

The proposition below follows from the definitions.
\begin{proposition}\label{prop:migcrossings_ribbon}
A knot is ribbon if and only if it has a ribbon tuple of midcrossings.
\end{proposition}
\begin{proof}
If the knot a ribbon then by cutting the ribbon disk along the ribbon midcrossings we get a set of non-intersecting disks. Then they boundary is an unlink.

Conversely, if we have a set of midcrossings that splits the knot into unlink then glue non-intersecting disks to the components of the unlink and move the folds of the disks beyond the midcrossing bands. Then the union of these disks and the midcrossing bands is a ribbon disk of the knot.
\end{proof}
\end{example}

\begin{remark}\label{rem:midcrossing_skein_module}
    The map from the set of midcrossings (in the weak sense) to the set of tangles which are induced by transformation rules, can be used in construction of skein modules~\cite{Turaev, Przytycki}. On can consider skein modules as topological description of polynomial knot invariants like Jones, Alexander, HOMFLY polynomials.
\end{remark}

{\color{red}
}

\section{Trait}\label{sect:trait}

Let us modify the crossing functor once again by allowing a crossing to survive when an arc passes between the overcrossing and undercrossing arcs.

\begin{definition}\label{def:trait_functor}
The \emph{trait functor} is a functor $\mathcal{T}\colon\mathfrak D_s\to Rel$ is a functor which extends the midcrossing functor $\mathcal {MC}$ by relations $c\mapsto c'$ in third Reidemeister midmoves.

\begin{figure}[h]
\centering
  \includegraphics[width=0.3\textwidth]{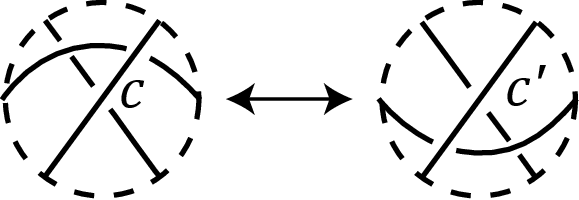}
  \caption{Third Reidemeister midmove}\label{pic:trait_move}
\end{figure}
\end{definition}

While a third Reidemeister midmove, an arc of the tangle passes through the vertical midcrossing probe of the crossing $c$. If we allow the probe to intersect the tangle during isotopy, then we can model a crossing change of the probe (see Fig.~\ref{pic:midcrossing_probe_crossing_change}) and change of the framing (see Fig.~\ref{pic:midcrossing_probe_unframing}). This leads us to the following definition.

\begin{definition}\label{def:trait_probe}
A \emph{trait probe} of a tangle $T$ in $F\times I$ is a pair $(\epsilon, \gamma)$ such that $\epsilon\in\{-1,1\}$ and
\[
\gamma\colon (I;0;1)\to (M_T, T,T)
\]
is a smooth map such that $\gamma(0)\ne\gamma(1)$.

Let $T\in\Sigma_0$ and $D=p(T)$ be its diagram. For a crossing $c\in\mathcal{T}(D)=\mathcal C(D)$, its \emph{vertical trait probe} is the pair $(sgn(c),\gamma_c)$ where $\gamma_c$ is the middle component of $(c\times I)\setminus T$ oriented downward.
\end{definition}

Similarly to the case of the midcrossing functor, we have the following statements.

\begin{proposition}\label{prop:trait_strong_transitivity}
    For any $c,c'\in\mathcal T(D)$, $c\sim^s_{\mathcal T,D}c'$ if and only if $sgn(c)=sgn(c')$ and $\gamma_c$ is homotopic to $\gamma_{c'}$ in the space of paths $\gamma\colon  I\to F\times I$ such that $\partial \gamma\in T$ and $\gamma(0)\ne\gamma(1)$.
\end{proposition}

\begin{definition}\label{def:trait_set_of_tangle}
 Denote the homotopy classes of trait probes of a tangle $T$ by $\mathscr{T}^{s}(T)$, and call it \emph{set of traits in the strong sense} of the tangle $T$.  Denote the set of orbits of $\mathscr{T}^{s}(T)$ by the action of the motion group $M(F\times I, T)$ by $\mathscr{T}^{w}(T)$  and call it the \emph{set of traits in the weak sense} of the tangle. There are natural maps $h^s(T)\colon\mathcal{T}(T)\to\mathscr{T}^s(T)$ and $h^w(T)\colon\mathcal{T}(T)\to\mathscr{T}^w(T)$ that assign to a crossing the corresponding vertical trait probe.
\end{definition}

The following statement is analogous to Theorems~\ref{thm:crossing_universal_homotopy_coinvariant} and~\ref{thm:crossing_universal_homotopy_invariant}.

\begin{theorem}\label{thm:trait_universal_invariant}
1. The pair $(\mathscr{T}^{s}, h^{s}_{\mathcal{T}})$ is the universal $h$-coinvariant of the trait functor $\mathcal{T}$.

2. The pair $(\mathscr{T}^{w}, h^{w}_{\mathcal{T}})$ is the universal $h$-invariant of the trait functor $\mathcal{T}$.
\end{theorem}

\begin{remark}
Similar to the case of crossings, we can define $W$- and $\Omega_2$-equivalence relations on the traits. Then the sets of $W$- and $\Omega_2$-strong equivalence classes of the trait probes coincide $\mathscr{T}^{sW}(T)=\mathscr{T}^{s\Omega_2}(T)$ and are equal to the set of homotopy classes of paths $\gamma\colon I\to F\times I$ such that $\partial\gamma\subset T$ and $\gamma(0)\ne\gamma(1)$.

The  sets of $W$- and $\Omega_2$-weak equivalence classes of the trait probes are equal to the set of orbits by the motion group:
\[
\mathscr{T}^{wW}(T)=\mathscr{T}^{w\Omega_2}(T)=\mathscr{T}^{sW}(T)/M(F\times I,T).
\]

Then we have $\mathscr{T}^{s}(T)=\mathbb Z_2\times\mathscr{T}^{sW}(T)$ and  $\mathscr{T}^{w}(T)=\mathbb Z_2\times\mathscr{T}^{wW}(T)$.
\end{remark}

Let us give an explicit description of the sets of traits of the tangle $T$. Let the tangle $T=T_1\sqcup\cdots\sqcup T_n$ has $n$ components. Then
\[
\mathscr{T}^{sW}(T)=\bigsqcup_{i,j=1}^n \mathscr{T}^{sW}_{ij}(T)
\]
where $\mathscr{T}^{sW}_{ij}(T)$ consists of paths $\gamma\colon I\to F\times I$ such that $\gamma(0)\in T_i$ and $\gamma(1)\in T_j$.

Choose arbitrary points $z_i\in T_i$, $1\le i\le n$, in the components. Denote by $z_i^+$ (resp. $z_i^-$) a point in $T_i$ that differs from $z_i$ by a small shift along (resp. against) the orientation of the component. Any path $\gamma$ is equivalent in $\mathscr{T}^{sW}_{ij}(T)$  to a path $\gamma'$ such that $\gamma'(0)=z_i$ and $\gamma'(1)=z_j$ if $i\ne j$, and to a path $\gamma'$ such that $\{\gamma'(0),\gamma'(1)\}=\{z^-_i,z^+_i\}$ if $i=j$.

If $i=j$ and the component $T_i$ is long (i.e. $\partial T_i\ne\emptyset$) then
\[
\mathscr{T}^{sW}_{ii}(T)=\mathscr{T}^{sW}_{ii,-}(T)\sqcup \mathscr{T}^{sW}_{ii,+}(T)
\]
where $\mathscr{T}^{sW}_{ii,\epsilon}(T)$, $\epsilon=\pm$, is formed paths $\gamma$ such that $\gamma(0)=z_i^\epsilon$ and $\gamma(1)=z_i^{-\epsilon}$. The traits in $\mathscr{T}^{sW}_{ii,-}(T)$ are called \emph{early overcrossings}, and those in $\mathscr{T}^{sW}_{ii,+}(T)$ are called \emph{early undercrossings}.

If $i=j$ and the component $T_i$ is closed (i.e. $\partial T_i=\emptyset$) then $\mathscr{T}^{sW}_{ii}(T)$ is presented by paths $\gamma$ such that $\gamma(0)=z_i^-$ and $\gamma(1)=z_i^+$.

By contracting the points $z_i^\pm$ to $z_i$, we can identify the set of the homotopy classes of paths $\gamma$ such that $\{\gamma(0),\gamma(1)\}=\{z_i^-,z_i^+\}$ with $\pi_1(F\times I, z_i)$.

Let $\tau_i\in\pi_1(F\times I, z_i)$ denotes the homotopy class of the component $T_i$ if the component is closed, and $\tau_i=1$ if the component is long. Given the ambiguity of moving the ends of paths to the base points $z_i^\pm$, we can write
\begin{gather*}
\mathscr{T}^{sW}_{ij}(T)=\tau_i\backslash\pi_1(F\times I,z_i,z_j)/\tau_j,\quad i\ne j,\\
\mathscr{T}^{sW}_{ii,\pm}(T)=\pi_1(F\times I,z_i),\quad T_i\mbox{ is long},\\
\mathscr{T}^{sW}_{ii,\pm}(T)=\pi_1(F\times I,z_i)/Ad(\tau_i),\quad T_i\mbox{ is closed},
\end{gather*}
where $\tau_i\backslash\pi_1(F\times I,z_i,z_j)/\tau_j$ is the set of orbits \[
\{\tau_i^n\gamma\tau_j^m\,\mid\, m,n\in\mathbb Z\},\quad  \gamma\in\pi_1(F\times I,z_i,z_j),
\]
and $\pi_1(F\times I,z_i)/Ad(\tau_i)$ is the set of orbits
\[
\{\tau_i^n\gamma\tau_i^{-n}\,\mid\, n\in\mathbb Z\},\quad  \gamma\in\pi_1(F\times I,z_i).
\]

Let $p\colon F\times I\to I$ be the projection map. Denote $x_i=p(z_i)$ and $\delta_i=p_*(\tau_i)\in\pi_1(F,x_i)$. Since $p$ is a homotopy equivalence, we come to the following statement.

\begin{theorem}\label{thm:trait_homotopic_description}
The set of traits of a tangle $T=T_1\sqcup\cdots\sqcup T_n$ is the union $\mathscr{T}^{s}(T)=\mathbb Z_2\times\mathscr{T}^{sW}(T)$ where $\mathscr{T}^{sW}(T)=\bigsqcup_{i,j=1}^n \mathscr{T}^{sW}_{ij}(T)$, and
\begin{gather*}
\mathscr{T}^{sW}_{ij}(T)=\delta_i\backslash\pi_1(F,x_i,x_j)/\delta_j,\quad\mbox{ if } i\ne j,\\
\mathscr{T}^{sW}_{ii,\pm}(T)=\pi_1(F,x_i),\quad\mbox{ if } T_i\mbox{ is long},\\
\mathscr{T}^{sW}_{ii,\pm}(T)=\pi_1(F,x_i)/Ad(\delta_i),\quad\mbox{ if } T_i\mbox{ is closed}.
\end{gather*}

For a trait $(\sigma,\gamma)\in\mathbb Z_2\times\mathscr{T}^{sW}(T)$, $\sigma$ is the sign of the trait and $\gamma$ is called the homotopy type of the trait. If $\gamma\in\mathscr{T}^{sW}_{ij}(T)$ then $(i,j)$ is the component type of the trait, and if $\gamma\in\mathscr{T}^{sW}_{ii,\epsilon}(T)$ for a long component $T_i$ then $\epsilon\in\mathbb Z_2$ is the order type of the trait.

Then we can say that the sign and the component, order, and homotopy types uniquely determine a trait of the tangle.
\end{theorem}

Theorem~\ref{thm:trait_homotopic_description} is a reformulation of~\cite[Theorem 1]{Ntribe}.


\section{Elements}\label{sect:elements}

\subsection{Incidence maps between diagram elements}\label{sect:element_incidence}

In a tangle diagram, there are incidence relations between diagram elements: an oriented arc has a region on its left and a region on its right; a crossing has four incident edges to it. We can give a topological interpretation of these relations.

\begin{definition}\label{def:element_incidence_maps}
Let $\gamma\in\mathscr{SA}^{ts}(T)$ be a semiarc probe, and $\gamma^u,\gamma^o,\mu$ the corresponding under-probe, over-probe, and meridian (see Definition~\ref{def:semiarc_probe}). Define \emph{left and right incidence maps} $R_r, R_l\colon\mathscr{SA}^{ts}(T)\to\mathscr{R}^{ts}(T)$ by the formulas
\[
R_r(\gamma)=(\gamma^u)^{-1}\mu_r^{-1}\gamma^o,\quad
R_l(\gamma)=(\gamma^u)^{-1}\mu_l\gamma^o.
\]

Let $\gamma\in\mathscr{C}^{ts}(T)$ be a crossing probe, and $\gamma^u,\gamma^m,\gamma^o,\mu^u,\mu^o$ the corresponding under-probe, mid-probe, over-probe and the meridians (see Definition~\ref{def:crossing_probe}). Define the \emph{incidence maps} $SA_{ur}, SA_{dr}, SA_{dl}, SA_{ul}\colon\mathscr{C}^{ts}(T)\to\mathscr{SA}^{ts}(T)$ for a positive crossing by the formulas
\begin{gather*}
SA_{ur}(\gamma)=(\gamma^m\mu^u_r\gamma^u,\gamma^o),\quad
SA_{dr}(\gamma)=(\gamma^u,(\gamma^m)^{-1}(\mu^o_r)^{-1}\gamma^o),\\
SA_{dl}(\gamma)=(\gamma^m\mu^u_l\gamma^u,\gamma^o),\quad
SA_{ul}(\gamma)=(\gamma^u,(\gamma^m)^{-1}(\mu^o_l)^{-1}\gamma^o),
\end{gather*}
where a semiarc probe is presented by its over-probe and under-probe. For a negative crossing, the incidence maps are defined by the formulas
\begin{gather*}
SA_{ur}(\gamma)=(\gamma^u,(\gamma^m)^{-1}(\mu^o_r)^{-1}\gamma^o),\quad
SA_{dr}(\gamma)=(\gamma^m\mu^u_r\gamma^u,\gamma^o),\\
SA_{dl}(\gamma)=(\gamma^u,(\gamma^m)^{-1}(\mu^o_l)^{-1}\gamma^o),\quad
SA_{ul}(\gamma)=(\gamma^m\mu^u_l\gamma^u,\gamma^o).
\end{gather*}

Compositions of the incidence maps from crossings to semiarcs and from semiarcs to regions produce the incidence maps $R_r, R_l, R_u, R_d\colon\mathscr{C}^{ts}(T)\to\mathscr{R}^{ts}(T)$:
\begin{gather*}
R_r=R_r\circ SA_{ur}=R_r\circ SA_{dr}, \quad R_l=R_l\circ SA_{ul}=R_l\circ SA_{dl},\\
R_u=R_r\circ SA_{ul}=R_l\circ SA_{ur}, \quad R_d=R_l\circ SA_{dl}=R_r\circ SA_{dr}.
\end{gather*}

\end{definition}

\begin{remark}
    The same formulas define incidence maps between tangle elements in the weak sense.
\end{remark}

On the other hand, there is a forgetting map from $\mathscr{SA}^{ts}(T)$ to $\mathscr{A}^{ts}(T)$ given by the formula $(\gamma^u,\gamma^o)\mapsto \gamma^o$.

\begin{proposition}
The incidence maps on the sets of diagram elements 
of a tangle $T$ form the following commutative diagram.
\begin{equation}\label{eq:element_incidence_maps_diagram}
\xymatrix{
\mathscr{A}(T) & & & & \mathscr{A}(T)\\
& \mathscr{SA}(T)\ar[ul] \ar[r]^{R_r}\ar[d]_{R_l} & \mathscr{R}(T) & \mathscr{SA}(T)\ar[ur]\ar[l]_{R_l}\ar[d]^{R_r} &\\
& \mathscr{R}(T) &
\mathscr{C}(T)\ar[l]|{R_l}\ar[ul]_{SA_{ul}}\ar[u]|{R_u}\ar[ur]^{SA_{ur}}\ar[r]|{R_r}\ar[dr]^{SA_{dr}}\ar[d]|{R_d}\ar[dl]_{SA_{dl}}
 & \mathscr{R}(T) & \\
& \mathscr{SA}(T)\ar[dl] \ar[r]_{R_r}\ar[u]^{R_l} & \mathscr{R}(T) & \mathscr{SA}(T)\ar[dr]\ar[l]^{R_l}\ar[u]_{R_r} &\\
\mathscr{A}(T) & & & & \mathscr{A}(T)
}
\end{equation}

The sets of diagram elements are either all considered in the strong sense or all considered in the weak sense.
\end{proposition}

\subsection{Compatibility of diagram elements}\label{sect:element_compatibility_relations}

\begin{definition}\label{def:region_adjacency}
    Two regions $\rho_0, \rho_1\in \mathscr R^s(T)$ of a tangle $T$ are \emph{adjacent} if there exists a morphism $f\colon T\to T'$ and two regions $r_0,r_1\in \mathcal R(D')$, $D'=p(T')$, such that $\gamma_{r_i}=\mathscr R^s(f)(\rho_i)$, $i=0,1$, and $r_0$ and $r_1$ are separated by an arc of the diagram $D'$ ($r_0$ is incident to the arc from the left, and $r_1$ is incident from the right). In this case, we write $\rho_0\uparrow\rho_1$.
\end{definition}

\begin{proposition}\label{prop:region_adjacency}
For a pair of regions $\rho_0,\rho_1\in\mathscr R^s(T)$, the following conditions are equivalent:
\begin{enumerate}
    \item $\rho_0\uparrow\rho_1$;
    \item there exists a semiarc $\sigma\in\mathscr{SA}^s(T)$ such that $\rho_0=R_l(\sigma)$ and $\rho_1=R_r(\sigma)$;
    \item there are region probes $\gamma_i$ representing $\rho_i$, $i=0,1$, such that there exists an embedded disk
\[
\Delta\colon (I\times I; 0\times I; 1\times I; I\times 0; I\times 1)\hookrightarrow (M_T; F\times 0; F\times 1; \gamma_0; \gamma_1)
\]
which intersects transversely $T$ in one point, and the intersection is positive.
\end{enumerate}
\end{proposition}

\begin{remark}\label{rem:region_semiarc_uniqueness}
    Two adjacent regions $\rho_0\uparrow\rho_1$ do not determine uniquely a semiarc $\sigma$ such that $\rho_0=R_l(\sigma)$ and $\rho_1=R_r(\sigma)$. For example, let $T$ be a composite knot, and let $\sigma_{in}$ and $\sigma_{out}$ be the incoming and outcoming semiarcs of a prime summand $K$ in $T$ (Fig.~\ref{pic:region_semiarc_nonuniqueness}). Then $\sigma_{in}\ne\sigma_{out}$ but they have the same incident regions.
\begin{figure}
    \centering
    \includegraphics[width=0.25\textwidth]{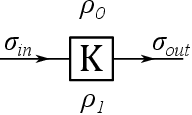}
    \caption{Semiarc with the same incident regions} \label{pic:region_semiarc_nonuniqueness}
\end{figure}

\end{remark}

\begin{proposition}\label{prop:semiarcs_adjacency}
Let $\epsilon\in\{-1,1\}$. For any semiarcs $\sigma_0,\sigma_1\in\mathscr{SA}^s(T)$ of a tangle $T$, the following conditions are equivalent:
\begin{enumerate}
    \item there exists a morphism $f\colon T\to T'$ and a crossing $c\in\mathcal C(D')$, $D'=p(T')$, such that $sgn(c)=\epsilon$, $SA_{dr}(\gamma_c)=\mathscr{SA}^s(f)(\sigma_0)$ and $SA_{ur}(\gamma_c)=\mathscr{SA}^s(f)(\sigma_1)$;
    \item there exists a crossing $\xi\in\mathscr C^s(T)$ such that $sgn(\xi)=\epsilon$, $SA_{dr}(\xi)=\sigma_0$ and $SA_{ur}(\xi)=\sigma_1$.
    \item $R_r(\sigma_0)=R_r(\sigma_1)$.
\end{enumerate}
\end{proposition}

There are statements analogous to those made for other cases of incident arcs to an oriented crossing.

\begin{remark}\label{rem:semiarc_adjacency}
    Two semiarcs do not determine a crossing. For example, the crossings $b'$ and $b'''$ in Fig.~\ref{pic:crossing_W_to_R2} have the same left semiarcs and the same top semiarcs (up to strong equivalence). But if $b'$ and $b'''$ are self-crossings of a component, they are equivalent by Proposition~\ref{prop:crossing_strong_equivalence_wrapping_index}.
\end{remark}

\begin{corollary}\label{cor:regions_form_crossing}
    For any regions $\rho_i\in\mathscr R^s(T)$, $i=1,2,3$, such that $\rho_1\uparrow\rho_2$ and $\rho_3\uparrow\rho_2$, and any $\epsilon\in\{-1,1\}$ there exist a morphism $f\colon T\to T'$ and a crossing $c\in\mathcal C(D')$, $D'=p(T')$, such that $sgn(c)=\epsilon$ and $R_d(\gamma_c)=\mathscr R^s(f)(\rho_1)$,  $R_r(\gamma_c)=\mathscr R^s(f)(\rho_2)$,  $R_u(\gamma_c)=\mathscr R^s(f)(\rho_3)$.
\end{corollary}
Informally speaking, three regions such that $\rho_1\uparrow\rho_2$ and $\rho_3\uparrow\rho_2$ form a crossing.

\begin{proof}
    Since $\rho_1\uparrow\rho_2$ and $\rho_3\uparrow\rho_2$, there exist semiarcs $\sigma_1,\sigma_2\in\mathscr{SA}^s(T)$ such that $R_l(\sigma_1)=\rho_1$, $R_r(\sigma_1)=\rho_2$, $R_l(\sigma_2)=\rho_3$ $R_r(\sigma_2)=\rho_2$. By Proposition~\ref{prop:semiarcs_adjacency}, there exists a morphism $f\colon T\to T'$ and a crossing $c\in\mathcal C(D')$, $D'=p(T')$, such that $sgn(c)=\epsilon$, $SA_{dr}(\gamma_c)=\mathscr{SA}^s(f)(\sigma_1)$ and $SA_{ur}(\gamma_c)=\mathscr{SA}^s(f)(\sigma_2)$. Then  $R_d(\gamma_c)=\mathscr R^s(f)(\rho_1)$,  $R_r(\gamma_c)=\mathscr R^s(f)(\rho_2)$,  $R_u(\gamma_c)=\mathscr R^s(f)(\rho_3)$.
\end{proof}

\begin{definition}\label{def:elements_reidemeister}
    1) Let $T$ be a tangle and $D$ its diagram to which a decreasing first Reidemeister move can be applied. Consider the disc of the move (see Fig.~\ref{pic:reidmove1elem} right). Then the region $r'_3$ is called the \emph{loop region} of the move, the semiarc $s'_2$ is the \emph{loop semiarc} of the move, and the (mid)crossing $c'_1$ is the \emph{loop (mid)crossing} of the move.

    2) Let $T$ be a tangle and $D$ its diagram to which a decreasing second Reidemeister move can be applied. Consider the disc of the move (see Fig.~\ref{pic:reidmove2elem} right). Then the region $r'_5$ is called the \emph{bigon region} of the move, the semiarcs $s'_3,s'_4$ are the \emph{bigon semiarcs} of the move, and the (mid)crossings $c'_1,c'_2$ are the \emph{bigon (mid)crossing} of the move.

    3) Let $T$ be a tangle, and $D$ its diagram to which a third Reidemeister move can be applied. Consider the disc of the move (see Fig.~\ref{pic:reidmove3elem} left). Then the unlabeled region is called the \emph{triangle region} of the move, the unlabeled semiarcs are the \emph{triangle semiarcs} of the move, and the (mid)crossings are the \emph{triangle (mid)crossing} of the move.
\end{definition}

\begin{proposition}\label{prop:midcrossing_loop}
    Let $T=K_1\cup\cdots\cup K_n$ be an $n$-component tangle. Then
\begin{enumerate}
    \item for any component $K_i$, $i=1,\dots,n$, there are four loop midcrossings (Fig.~\ref{pic:midcrossing_loops}) unless $K_i$ is an unlinked trivial component. In the latter case, there are two loop microssings $r_+=l_+$ and $r_-=l_-$ which differ by the sign;
    \item a midcrossing $\xi\in\mathscr{MC}^s(T)$ is a loop midcrossing if and only if the oriented smoothing $T'$ of the tangle $T$ at the midcrossing $\xi$ is isotopic to $T\sqcup\bigcirc$.
\end{enumerate}

\begin{figure}
    \centering
    \includegraphics[width=0.7\textwidth]{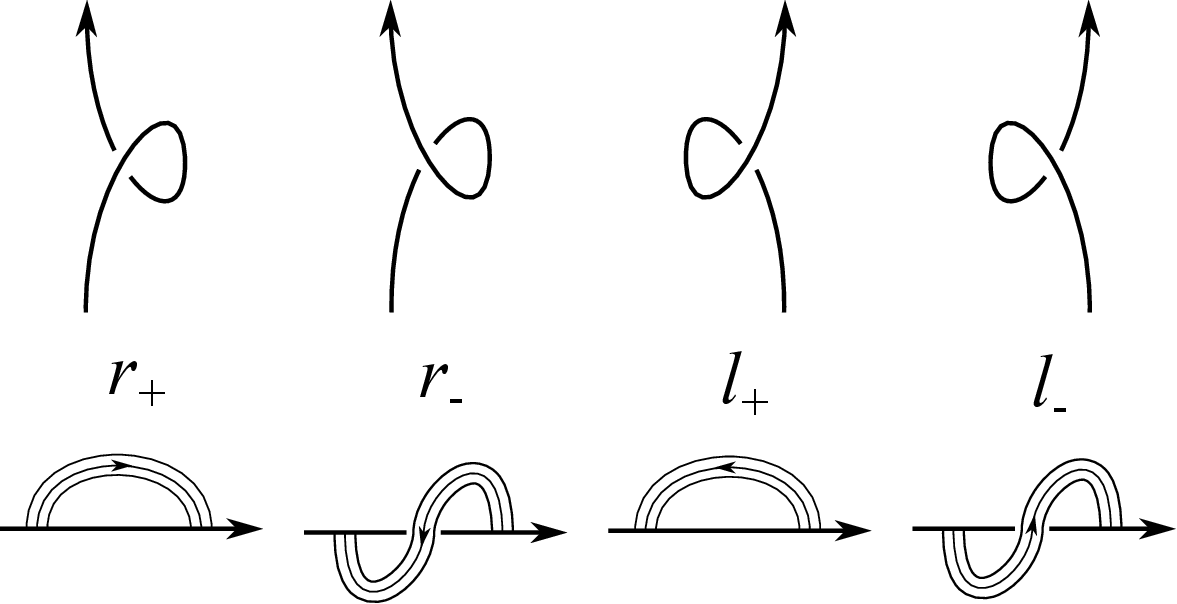}
    \caption{Loop midcrossings and their probes}
    \label{pic:midcrossing_loops}
\end{figure}
\end{proposition}

\begin{proof}
    It is easy to see that the small probes in Fig.~\ref{pic:midcrossing_loops} represents loop midcrossings.

    The midcrossings $r_+,l_+$ differ from $r_-,l_-$ by the sign. If $K_i$ is a nontrivial component, the midcrossings $r_\pm$ can be distinguished from $l_\pm$ by the oriented smoothing: at the midcrossing $r_\pm$, the left half of the tangle is $K_i$ and the right half is trivial, whereas at $l_\pm$ the left half is trival and the right half is $K_i$.

    If $K_i$ is trivial there is an isotopy between $r_\pm$ and $l_\pm$ in a neighborhood of the spanning disk of the component $K_i$.

    Let us show that any loop midcrossing $\xi\in\mathscr{MC}^s(T)$ is isotopic to $r_\pm$ or $l_\pm$. Let $f\colon T\to T'$ be a morphism and $c\in \mathcal C(D')$, $D'=p(T')$, a crossing such that $\gamma_c=\mathscr{MC}^s(f)(\xi)$. Then the probe $\gamma_c$ looks like probes in Fig.~\ref{pic:midcrossing_loops}. Contract $\gamma_c$ to a tiny probe $\gamma'$. Apply the inverse isotopy $f^{-1}$ to $\gamma'$. We can suppose that $f^{-1}$ is linear in a small neighborhood of $\gamma'$. Hence, the probe $f^{-1}(\gamma')$ looks like probes $r_\pm$ or $l_\pm$. Since $\mathscr{MC}^s(f)(\xi)=\gamma_c$ isotopic to $\gamma'$, the midcrossing $\xi$ is isotopic to $f^{-1}(\gamma')$, hence, to $r_\pm$ or $l_\pm$.

    It is clear that the oriented smoothing at the midcrossings $r_\pm$ and $l_\pm$ yields the tangle $T\sqcup\bigcirc$. On the other hand, if the oriented smoothing at a midcrossing $\xi$ is $T\sqcup\bigcirc$, after application of the contracting isotopy of the trivial component to the tangle $T$, one gets a tangle with a loop whose midcrossing is the image of $\xi$.
\end{proof}

\begin{proposition}\label{prop:elements_loop}
Let $T$ be a tangle.
    \begin{enumerate}
        \item Any region $\rho\in\mathscr{R}^s(T)$ is the loop region for some first Reidemeister move, i.e. there exists a morphism $f\colon T\to T'$ and a region $r\in\mathcal R(D')$, $D'=p(T')$, such that $\gamma_r=\mathscr{R}^s(f)(\rho)$ and $r$ is a loop region in the diagram $D'$;
        \item Any semiarc $\sigma\in\mathscr{SA}^s(T)$ is the loop semiarc for some first Reidemeister move, i.e. there exists a morphism $f\colon T\to T'$ and a semiarc $s\in\mathcal{SA}(D')$, $D'=p(T')$, such that $\gamma_s=\mathscr{SA}^s(f)(\sigma)$ and $s$ is a loop semiarc in the diagram $D'$;
        \item A crossing $\xi\in\mathscr C^s(T)$ is a loop crossing iff the midcrossing $\xi_m$ is a loop midcrossing;
        \item  A crossing $\xi\in\mathscr C^s(T)$ is a loop crossing iff the oriented smoothing $T'$ of the tangle $T$ at $\xi$ is $T\sqcup\bigcirc$.
    \end{enumerate}
\end{proposition}

\begin{proof}
    1. Let $\rho\in\mathscr{R}^s(T)$ and $f_0\colon T\to T'_0$ be a morphism realizing $\rho$, i.e. there exists $r\in\mathcal R(D'_0)$, $D'_0=p(T'_0)$, such that $\gamma_r=\mathscr{R}^s(f_0)(\rho)$ (Fig.~\ref{pic:loop_region_semiarc} left). Apply a second Reidemeister move to get a diagram $D'$ where the the region $r$ becomes a loop region (Fig.~\ref{pic:loop_region_semiarc} right). Let $f$ be the composition of $f_0$ and the move. Then $\mathscr{R}^s(f)(\rho)$ is the probe of a loop region.

    2. Analogously, one can make a loop semiarc from an explicit semiarc $s$ of a tangle diagram (Fig.~\ref{pic:loop_region_semiarc}).

    3. It is clear that the midprobe of a loop crossing probe is a loop midcrossing. On the other hand, if $\xi$ is a crossing such that its midprobe $\xi_m$ is a loop midcrossing, then consider an isotopy which contracts $\xi_m$ to a loop and verticalize the under- and the overprobes of $\xi$. Then after the isotopy $\xi$ will become a loop crossing.
\begin{figure}
    \centering
    \includegraphics[width=0.3\textwidth]{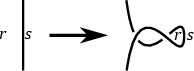}
    \caption{Loop region and semiarc}
    \label{pic:loop_region_semiarc}
\end{figure}

    4. The last statement of the proposition follows from Proposition~\ref{prop:midcrossing_loop}.
\end{proof}

\begin{proposition}\label{prop:elements_bigon}
    Let $T$ be a tangle. Then
\begin{enumerate}
    \item any region $\rho\in\mathscr{R}^s(T)$ is the bigon region of some second Reidmeister move;
    \item two semiarcs $\sigma_1,\sigma_2\in\mathscr{SA}^s(T)$ are the bigon semiarcs of some second Reidemeister move iff $R_\alpha(\sigma_1)=R_\beta(\sigma_2)$ for some $\alpha,\beta\in\{l,r\}$;
    \item two crossings $\xi_1,\xi_2\in\mathscr{C}^s(T)$ are the bigon crossings of some second Reidemeister move iff $\xi_2=W_{\pm 1}(\xi_1)$.
\end{enumerate}
\end{proposition}

\begin{proof}
    1. Let $f_0$ be a verticalizing isotopy of the region probe $\rho$. Apply two second Reidemeister moves to semiarcs bounding the verticalized region (Fig.~\ref{pic:bigon_region_semiarc}). Then the region becomes a bigon region in the new diagram.

\begin{figure}
    \centering
    \includegraphics[width=0.5\textwidth]{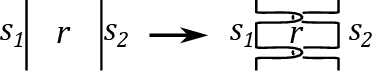}
    \caption{Bigon region and semiarcs}
    \label{pic:bigon_region_semiarc}
\end{figure}
    2. Let semiarcs $\sigma_1$ and $\sigma_2$ have one incident region $\rho$. Isotope the probes $\sigma_1$, $\sigma_2$ and $\rho$ to vertical probes of semiarcs $s_1$, $s_2$ and region $r$ of some diagram (Fig.~\ref{pic:bigon_region_semiarc} left). Apply two second Reidemeister moves to the diagram as shown in Fig.~\ref{pic:bigon_region_semiarc}. Then the semiarcs become bigon semiarcs in the new diagram.

    The inverse statement is clear.

    3. The third statement of the proposition follows from the proof of Proposition~\ref{prop:crossing_equivalence_R2_vs_W}.
\end{proof}

\begin{proposition}\label{prop:elements_trigon}
    Let $T$ be a tangle. Then
\begin{enumerate}
    \item any region $\rho\in\mathscr{R}^s(T)$ is the triangle region of some third Reidmeister move;
    \item three semiarcs $\sigma_1,\sigma_2,\sigma_3\in\mathscr{SA}^s(T)$ are the triangle semiarcs of some third Reidemeister move iff $R_\alpha(\sigma_1)=R_\beta(\sigma_2)=R_\gamma(\sigma_3)$ for some $\alpha,\beta,\gamma\in\{l,r\}$;
    \item two crossings $\xi_1,\xi_2\in\mathscr{C}^s(T)$ are the crossings of some third Reidemeister move iff $\xi_1$ and $\xi_2$ have a common incident semiarc.
\end{enumerate}
\end{proposition}

\begin{proof}
    1. Isotope the region probe $\rho$ to a vertical probe of some region $r$. Then one can use the incident semiarcs of $r$ to bound a triangle part $r'\subset r$ to which a third Ridemeister move can be applied. Then $\rho$ is the probe of the triangle region $r'$;

    2. Let $\sigma_1$, $\sigma_2$ and $\sigma_3$ be semiarcs which have a common incident region $\rho$. Verticalize the probe $\rho$. Then $\sigma_i$ is isotopic to the probe $\rho$ with a sprout attached to an arc of the tangle. After pulling the arcs along the sprouts to the probe $\rho$ (Fig.~\ref{pic:trigon_semiarcs}), one gets a triangle region whose semiarc probes are $\sigma_i$.
\begin{figure}
    \centering
    \includegraphics[width=0.7\textwidth]{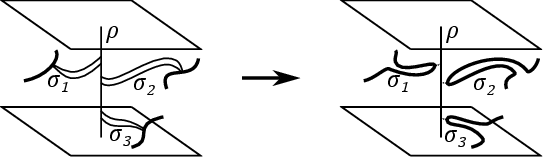}
    \caption{Pulling semiarcs to the region probe}
    \label{pic:trigon_semiarcs}
\end{figure}
    The inverse statement is clear.

    3. Let $\xi_1,\xi_2$ be crossing probes with a common semiarc probe $\sigma$. Verticalize the probe $\sigma$. Then $\xi_i$ is isotopic to the probe $\sigma$ with a sprout attached to the other incident arc of crossing. After pull these arcs to $\sigma$ along the sprouts and adjusting the framings of the midprobes of $\xi_1$, $\xi_2$, one gets a triangle region whose vertices are $\xi_1$, $\xi_2$.
\end{proof}

{\color{red}

}
\subsection{Representing spatial graphs of diagram elements}\label{subsect:element_representing_graph}

Description of universal homotopical invariants of diagram elements given in Theorems~\ref{thm:arc_universal_homotopy_invariant},\ref{thm:region_universal_homotopy_invariant},\ref{thm:semiarc_universal_homotopy_invariant},\ref{thm:crossing_universal_homotopy_invariant} and~\ref{thm:midcrossing_universal_invariant} can be reformulated for knots in the following way. For an element functor $\mathcal F=\mathcal A, \mathcal R, \mathcal{SA}, \mathcal C$ or $\mathcal{MC}$, consider the set of diagram elements (in the weak sense) of all knots
\[
\mathcal F_{knots}=\bigcup_{K\mbox{ \scriptsize is a knot}}\mathcal F^w(K).
\]

\begin{figure}[h]
\centering
  \includegraphics[width=0.8\textwidth]{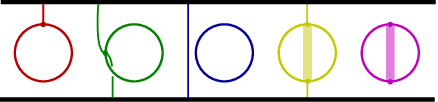}
  \caption{Representing spatial graphs of knot elements}\label{pic:element_spatial_graphs}
\end{figure}

On the other hand, consider the graphs (Fig.~\ref{pic:element_spatial_graphs})
\begin{itemize}
    \item $G_{\mathcal A}$ with vertices $V(G_{\mathcal A})=\{v,v_1\}$ and edges $E(G_{\mathcal A})=\{e,e_1\}$ such that $\partial e=v$, $\partial e_1=\{v,v_1\}$ and the ends of the edge $e$ are opposite at $v$;

    \item $G_{\mathcal{SA}}$ with vertices $V(G_{\mathcal{SA}})=\{v,v_0,v_1\}$ and edges $E(G_{\mathcal{SA}})=\{e,e_0,e_1\}$ such that $\partial e=v$ and $\partial e_i=\{v,v_i\}$, $i=0,1$. The ends of the edge $e$ are opposite in $v$ as well as the edges $e_0$ and $e_1$;

    \item $G_{\mathcal R}$ with vertices $V(G_{\mathcal R})=\{v_0,v_1\}$ and edges $E(G_{\mathcal R})=\{e,e_1\}$ such that $\partial e=\emptyset$ (i.e. $e$ is a circular edge), $\partial e_1=\{v_0,v_1\}$;

    \item $G_{\mathcal C}$ with vertices $V(G_{\mathcal C})=\{v_0,v_u,v_o,v_1\}$ and edges $E(G_{\mathcal C})=\{e_0,e,e',e_m,e_1\}$ such that $\partial e=\partial e'=\partial e_m=\{v_u,v_o\}$, $\partial e_0=\{v_0,v_u\}$ and $\partial e_1=\{v_o,v_1\}$. The ends of the edges $e$ and $e'$ are opposite in $v_u$ and $v_u$, the ends of $e_0$ and $e_m$ are opposite in $v_u$, and the ends of $e_m$ and $e_1$ are opposite at $v_o$;

    \item $G_{\mathcal{MC}}$ with vertices $V(G_{\mathcal{MC}})=\{v_u,v_o\}$ and edges $E(G_{\mathcal{MC}})=\{e,e',e_m\}$ such that $\partial e=\partial e'=\partial e_m=\{v_u,v_o\}$. The ends of the edges $e$ and $e'$ are opposite in $v_u$ and $v_u$.
\end{itemize}

For an element functor $\mathcal F=\mathcal A, \mathcal R, \mathcal{SA}, \mathcal C$ or $\mathcal{MC}$, denote by $Emb(G_{\mathcal F},F\times I)$ the set of isotopy classes of embeddings $f\colon G_{\mathcal F}\hookrightarrow F\times I$ such that $f(v_i)\in F\times i$, $i=0,1$, opposite edges of $G_{\mathcal F}$ go to opposite edges, and the edge $f(e_m)$ is framed. Then the following statements holds.

\begin{proposition}\label{prop:elemets_spatial_graph}
    For an element functor $\mathcal F=\mathcal A, \mathcal R, \mathcal{SA}, \mathcal C$ or $\mathcal{MC}$, there is a bijection between $\mathcal F_{knots}$ and $Emb(G_{\mathcal F},F\times I)$ established by the map
\[
x\mapsto K\cup\gamma_x,
\]
where $K$ is a knot with a diagram $D$, $x\in\mathcal F(D)$ is a diagram element, and $\gamma_x$ is the probe of $x$.
\end{proposition}

\begin{figure}[h]
\centering
  \includegraphics[width=0.5\textwidth]{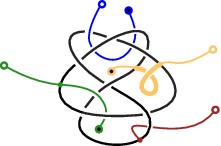}
  \caption{Diagrams of element probes: red arc, green semiarc, blue region, yellow crossing}\label{pic:diagram_probes}
\end{figure}

For a tangle $T\in\Sigma_0$ in $F\times I$, elements of the tangle in the weak or strong sense can be presented by probe diagrams (Fig.~\ref{pic:diagram_probes}).

\begin{definition}\label{def:element_probe_diagram}
Let $D$ be the diagram of $T$. Then
\begin{itemize}
    \item an \emph{arc probe diagram} is a diagram $D\cup\gamma$ where $\gamma\colon [0,1]\to F$ is a path such that $\gamma(0)$ is a non-crossing point of $D$, $\gamma(1)\in F\setminus D$ is marked as a float vertex, $\gamma$ has finite number of self-intersections and intersections with $D$ which are all double points with under-overcrossing structure, and $\gamma$ is transversal to $D$ in $\gamma(0)$;

    \item a \emph{region probe diagram} is a diagram $D\cup\gamma$ where $\gamma\colon [0,1]\to F$ is a path such that $\gamma(0)\in F\setminus D$ is marked as a sinker vertex, $\gamma(1)\in F\setminus D$ is marked as a float vertex, $\gamma$ has finite number of self-intersections and intersections with $D$ which are all transversal intersections with under-overcrossing structure;

    \item a \emph{semiarc probe diagram} is a diagram $D\cup\gamma$ where $\gamma\colon [0,1]\to F$ is a path such that $\gamma(0)\in F\setminus D$ is marked as a sinker vertex, $\gamma(1)\in F\setminus D$ is marked as a float vertex, $\gamma(\frac 12)$ is a non-crossing point of $D$, $\gamma(1)\in F\setminus D$ is marked as a float vertex, $\gamma$ has finite number of self-intersections and intersections with $D$ which are all transversal intersections with under-overcrossing structure (except $\gamma(\frac 12)$ which has no such structure);

    \item a \emph{crossing probe diagram} is a diagram $D\cup\gamma$ where $\gamma\colon [0,1]\to F$ is a path such that $\gamma(0)\in F\setminus D$ is marked as a sinker vertex, $\gamma(1)\in F\setminus D$ is marked as a float vertex, $\gamma(\frac 13)$ and $\gamma(\frac 23)$ are non-crossing points of $D$, $\gamma$ has finite number of self-intersections and intersections with $D$ which are all transversal intersections with under-overcrossing structure (except $\gamma(\frac 13), \gamma(\frac 23)$ which have no such structure). The part of the path between $\gamma(\frac 13)$ and $\gamma(\frac 23)$ is considered to be framed. We assume that the tangent vector to $\gamma$ and the tangent vector to $D$ form a negatively oriented frame in $\gamma(\frac 23)$ and the orientation of this frame in $\gamma(\frac 13)$ is equal to the sign of the crossing (Fig.~\ref{pic:graphoid_crossing_probes}). 
\begin{figure}[h]
\centering
  \includegraphics[width=0.8\textwidth]{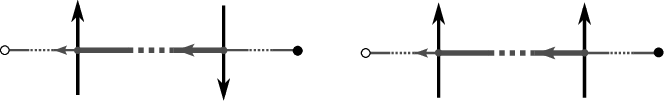}
  \caption{Probe diagrams of a positive crossing (left) and a negative crossing (right)}\label{pic:graphoid_crossing_probes}
\end{figure}
\end{itemize}
\end{definition}

Theorems~\ref{thm:arc_universal_homotopy_coinvariant}, \ref{thm:region_universal_homotopy_coinvariant}, \ref{thm:semiarc_universal_homotopy_coinvariant} and~\ref{thm:crossing_universal_homotopy_coinvariant} lead to the following statement.

\begin{figure}[h]
\centering
  \includegraphics[width=0.4\textwidth]{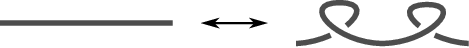}
  \caption{Framed first Reidemeister move}\label{pic:graphoid_framed_R1}
\end{figure}

\begin{figure}[h]
\centering
  \includegraphics[width=0.7\textwidth]{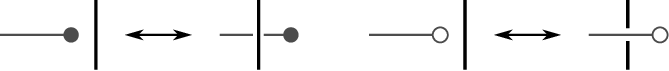}
  \caption{The sinker move (left) and the float move (right)}\label{pic:graphoid_float_sinker_moves}
\end{figure}

\begin{figure}[h]
\centering
  \includegraphics[width=0.6\textwidth]{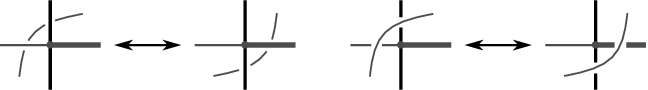}
  \caption{Third Reidemeister moves at vertex}\label{pic:graphoid_intersection_R3}
\end{figure}

\begin{figure}[h]
\centering
  \includegraphics[width=0.37\textwidth]{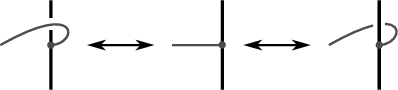}\qquad\qquad\includegraphics[width=0.3\textwidth]{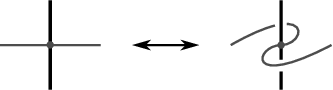}\\ \vspace{1em}
  \includegraphics[width=0.35\textwidth]{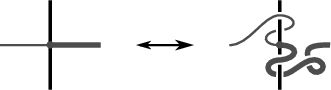}
  \caption{Vertex rotation moves}\label{pic:graphoid_vertex_rotation_moves}
\end{figure}

\begin{proposition}\label{prop:element_probe_diagram}
    The sets of arcs (regions, semiarcs, crossings) of the tangle $T$ in the strong sense can be identified with the classes of probe diagram modulo the following moves:
    \begin{itemize}
        \item isotopy of the probe;
        \item first, second and third Reidemeister moves in the unframed part of the probe;
        \item framed first (Fig.~\ref{pic:graphoid_framed_R1}), second and third Reidemeister moves in the framed part of the probe;
        \item float and sinker moves (Fig.~\ref{pic:graphoid_float_sinker_moves});
        \item third Reidemeister moves at vertices (Fig.\ref{pic:graphoid_intersection_R3});
        \item vertex rotation moves (Fig.~\ref{pic:graphoid_vertex_rotation_moves}).
    \end{itemize}

The sets of arcs (regions, semiarcs, crossings) of the tangle $T$ in the weak sense can be identified with the classes of probe diagram modulo the moves above and the moves:
\begin{itemize}
    \item isotopy of the diagram $D$;
    \item first, second and third Reidemeister moves on diagram $D$ which don't involve the probe.
\end{itemize}
\end{proposition}

{\color{red}




}

\subsection{Strong and homotopical strong equivalence}

Below we will prove Proposition~\ref{prop:h_strong_equivalence}. We start with the following lemma.

\begin{lemma}\label{lem:h-relation_is_isotopy}
    Let $T\in\Sigma_0$ be a tangle, $D=p(T)$ its diagram, and $\mathcal F=\mathcal A, \mathcal R, \mathcal{SA}$ or $\mathcal C$. For any two elements $x,y\in\mathcal F(D)$ the following conditions are equivalent:
\begin{enumerate}
    \item there exists a sequence of moves $f\in\mathrm{Mor}_{\mathfrak D_s}(D,D)$ such that $f$ is homotopic to the identity $id_D$ (i.e. $f=id_D\in\mathrm{Mor}_{\mathfrak D_h}(D,D)$) and  $y\in\mathcal F(f)(x)$;
    \item the probes $\gamma_x$ and $\gamma_y$ of the elements $x,y$ are isotopic.
\end{enumerate}
\end{lemma}

\begin{proof}
    1. Assume there is $f\sim id_D$ such that $y\in\mathcal F(f)(x)$. The morphism $f$ is realized by a tangle isotopy $\tau=(\tau_t)$, $t\in[0,1]$. Since $y\in\mathcal F(f)(x)$, there is a family $\gamma_t$, $t\in[0,1]$, of vertical probes on tangles $\tau_t$ such that $\gamma_0=\gamma_x$ and $\gamma_1=\gamma_y$.

    Since $f\sim id_D$ there is a homotopy $\Gamma=(\tau^s)$ such that $\tau^s_0=\tau^s_1=\tau^0_t=T$, $t,s\in[0,1]$, and $\tau^1=\tau$. By isotopy excision theorem, the homotopy $\Gamma$ extends to a map $H\colon I\times I\to\widetilde{Diff}(F\times I)$ such that $H_{t0}=H_{0s}=id$ and $\tau^s_t=H_{ts}(T)$ for all $t,s\in[0,1]$. Moreover, we can assume that $H_{t1}(\gamma_x)=\gamma_t$, $t\in[0,1]$. Then the family of probes $\gamma'_s=H_{s1}(\gamma_x)$, $s\in[0,1]$, is an isotopy of element probes $\gamma_x$ and $\gamma_y$.

    2. Now, assume that $\gamma_x$ and $\gamma_y$ are isotopic. This isotopy extends to an isotopy $h\colon I\to\widetilde{Diff}(F\times I)$ such that $h_0=id$, $h_t(T)=T$, $t\in[0,1]$ and $h_1(\gamma_x)=\gamma_y$. By Proposition~\ref{prop:probes_verticalizing}, there is a family of diffeomorphisms $H\colon I\times I\to\widetilde{Diff}(F\times I)$ such that $H_{t0}=h_t$, $H_{0s}=id$, $H_{1s}=h_1$ and $H_{t1}(\gamma_x)$ are vertical for all $t,s\in[0,1]$. We can suppose that the isotopy $H_{t1}(T)$ is $\Sigma_1$-transversal. Then it defines a morphism $f\in\mathrm{Mor}_{\mathfrak D_s}(D,D)$. Since $H_{11}(\gamma_x)=h_1(\gamma_x)=\gamma_y$,  we have $\mathcal F(f)(x)=y$. And by construction, $f\sim id_D$.
\end{proof}

We can reformulate the lemma above as follows: the additional relation on diagram elements introduced in Definition~\ref{def:h_strong_weak_equivalences} coincides with the topological strong equivalence on the corresponding probes.

\begin{proof}[Proof of Proposition~\ref{prop:h_strong_equivalence}]
    Let $T$ be a tangle, $D$ its diagram and $\mathcal F=\mathcal A, \mathcal R$ or $\mathcal {SA}$. We need to show that for $x,y\in\mathcal F(D)$, $x\sim^s_{\mathcal F,D}y$ is equivalent to $x\sim^{hs}_{\mathcal F,D}y$.
    The condition $x\sim^s_{\mathcal F,D}y$ implies $x\sim^{hs}_{\mathcal F,D}y$ by definition.

    Now, assume that $x\sim^{hs}_{\mathcal F,D}y$. We will show that $\gamma_x$ is isotopic to $\gamma_y$. The relation $\sim^{hs}$ is generated by the following rules:
\begin{itemize}
    \item (homotopy rule) for any $x,y\in\mathcal F(D)$ and any $f\in\mathrm{Mor}_{\mathfrak D_s}(D,D)$ such that $f= id_D$ in $\mathfrak D_h$ and $y\in\mathcal F(f)(x)$, one has $x\sim^{hs}_{\mathcal F,D}y$;
    \item (functorial rule) for any $f\in\mathrm{Mor}_{\mathfrak D_s}(D,D')$ $x_1,x_2\in\mathcal F(D)$ and $y_1\in\mathcal F(f)(x_1)$, $y_2\in\mathcal F(f)(x_2)$,  $x_1\sim^{hs}_{\mathcal F,D} x_2$ implies $y_1\sim^{hs}_{\mathcal F,D'} y_2$;
    \item (transitivity)  $x\sim^{hs}_{\mathcal F,D} y$ and $y\sim^{hs}_{\mathcal F,D} z$ implies $x\sim^{hs}_{\mathcal F,D} z$.
\end{itemize}
    It is enough to show that the topological strict equivalence
\[
x\sim^{ts} y \Leftrightarrow \gamma_x\mbox{ is isotopic to }\gamma_y
\]
satisfies the rules above. Homotopy rule holds by Lemma~\ref{lem:h-relation_is_isotopy}, and it is clear that the topological strict equivalence is transitive.

    Let us check the functorial rule. Let $x_1\sim^{ts}x_2$, $f\in\mathrm{Mor}_{\mathfrak D_s}(D,D')$ and $y_i\in\mathcal F(f)(x_i)$, $i=1,2$. Let $\phi\in\widetilde{Diff}(F\times I)$ be an isotopy realizing the morphism $f$. Since $y_i\in\mathcal F(f)(x_i)$, we can assume that $\phi(\gamma_{x_i})=\gamma_{y_i}$, $i=1,2$. Then the composition of the isotopy between the probes $\gamma_{x_1}$ and $\gamma_{x_2}$ and the map $\phi$ is an isotopy between $\gamma_{y_1}$ and $\gamma_{y_2}$.

    Thus, $x\sim^{hs}_{\mathcal F,D}y$ implies an isotopy between $\gamma_x$ and  $\gamma_y$. Then by Propositions~\ref{prop:arc_strong_transitivity},\ref{prop:region_strong_transitivity} and~\ref{prop:semiarc_strong_transitivity}, $x\sim^s_{\mathcal F,D}y$.

    The proof for the functor $\mathcal C$ is analogous.
\end{proof}

\section{Homotopy classes of diagram elements and colorings}\label{sect:colorings}

In the sections above, we described diagram elements of a tangle (up to natural equivalence) as \emph{isotopy classes} of arcs in the complement to the tangle. The sets of isotopy classes are complicated (there are surjections from these sets to the set of isotopy classes of tangles). Thus, it is impractical to use these sets directly as tangle invariants.

So, let us consider a rougher equivalence on diagram elements, and take \emph{homotopy classes} of elements' probes.

\begin{definition}\label{def:homotopy_classes_tangle_elements}
Let $F$ be a compact oriented surface and $T$ a tangle in the thickening $F\times I$. Denote a small tubular neighbourhood of $T$ by $N(T)$, and the complement manifold to $T$ by $M_T=\overline{F\times I\setminus N(T)}$.

Let us denote
\begin{itemize}
    \item $\mathdutchcal A(T)=[I,0,1; M_T,\partial N(T), F\times 1]$ the set of the homotopy classes of paths in $M_T$ from $\partial N(T)$ to $F\times 1$;
    \item $\mathdutchcal R(T)=[I,0,1; M_T,F\times 0, F\times 1]$ the set of the homotopy classes of paths in $M_T$ from $F\times 0$ to $F\times 1$;
    \item $\mathdutchcal{SA}(T)=[I,0,1; M_T,\partial N(T), F\times 1]\times_{\partial N(T)}[I,0,1; M_T,\partial N(T), F\times 0]$ --- the set of homotopy classes of pair of paths $(\gamma^o,\gamma^u)$ such that $\gamma^u(0)$ and $\gamma^o(0)$ are different points on the same meridian of $\partial N(T)$;
    \item  $\mathdutchcal{C}(T)$ --- the set of homotopy classes of triples $(\gamma^o,\gamma^m,\gamma^u)$ where
\begin{gather*}
\gamma^o\colon(I; 0; 1)\to (M_T; \partial N(T); F\times 1),\\
\gamma^m\colon(I; 0; 1)\to (M_T; \partial N(T); \partial N(T)),\\
\gamma^u\colon(I; 0; 1)\to (M_T;\partial N(T);  F\times 0)
\end{gather*}
such that
\begin{itemize}
\item $\gamma^o(0)$ and $\gamma^m(0)$ are different points of the same meridian $\mu^o$,
\item $\gamma^m(1)$ and $\gamma^u(0)$ are different points of the same meridian $\mu^u$,
\item the meridians $\mu^o$ and $\mu^u$ are distinct,
\item $\gamma^m$ is framed by a transversal vector field which is collinear to the tangle at $\gamma^m(0)$ and $\gamma^m(1)$.
\end{itemize}
\end{itemize}
We will call the sets $\mathdutchcal A(T), \mathdutchcal R(T), \mathdutchcal{SA}(T), \mathdutchcal C(T)$ the \emph{sets of homotopy classes of arcs, regions, semiarcs and crossings} of the tangle $T$.
\end{definition}

By forgetting the framing of the middle arc of the elements in $\mathdutchcal C(T)$, we get the set $\mathdutchcal C^{uf}(T)$ of unframed homotopy classes of crossings of the tangle.


Let us choose $x_0\in F$ such that $x_0\times I\cap N(T)=\emptyset$. Let $x^u=x_0\times 0$ and $x^o=x_0\times 1$.

Consider the set of homotopy classes
\[
\widetilde{\mathdutchcal A}(T)=[I,0,1; M_T,\partial N(T), x^o].
\]
Analogously, we define the \emph{sets of based homotopy classes of regions, semiarcs and crossings} $\widetilde{\mathdutchcal R}(T)$, $\widetilde{\mathdutchcal{SA}}(T)$, $\widetilde{\mathdutchcal C}(T)$
, imposing the condition that the paths end at the points $x^u$ and $x^o$.

There is a natural action of the group $\pi_1(F,x_0)$ on $\widetilde{\mathdutchcal A}(T)$:
\[
(\gamma,\alpha)\mapsto \gamma\cdot(\alpha\times 1),\quad \gamma\in\widetilde{\mathdutchcal A}(T), \alpha\in\pi_1(F,x_0),
\]
and there is an action of  $\pi_1(F,x_0)\times\pi_1(F,x_0)$ on $\widetilde{\mathdutchcal R}(T)$, $\widetilde{\mathdutchcal{SA}}(T)$, $\widetilde{\mathdutchcal C}(T)$ given by the formula
\[
(\gamma^u,\gamma^o;\alpha_0,\alpha_1)\mapsto (\gamma^u\cdot(\alpha_0\times 0), \gamma^o\cdot(\alpha_1\times 1)).
\]
Then $\mathdutchcal A(T)=\widetilde{\mathdutchcal A}(T)/\pi_1(F,x_0)$,
$\mathdutchcal R(T)=\widetilde{\mathdutchcal R}(T)/(\pi_1(F,x_0)\times\pi_1(F,x_0))$, analogous formulas hold for $\widetilde{\mathdutchcal{SA}}(T)$, $\widetilde{\mathdutchcal C}(T)$.

The projection $p\colon F\times I\to F$ induces maps $p_*$ from $\widetilde{\mathdutchcal R}(T)$, $\widetilde{\mathdutchcal{SA}}(T)$, $\widetilde{\mathdutchcal C}(T)$ to $\pi_1(F,x_0)$. Denote the preimages $p_*^{-1}(1)$ by $\widetilde{\mathdutchcal R}_0(T)$, $\widetilde{\mathdutchcal{SA}}_0(T)$, $\widetilde{\mathdutchcal C}_0(T)$ respectively. The group $\pi_1(F,x_0)$ acts on $\widetilde{\mathdutchcal R}_0(T)$, $\widetilde{\mathdutchcal{SA}}_0(T)$, $\widetilde{\mathdutchcal C}_0(T)$ by the formulas
\begin{gather*}
(\gamma,\alpha)\mapsto (\alpha^{-1}\times 0)\cdot\gamma\cdot(\alpha\times 1),\\
(\gamma^u,\gamma^o;\alpha)\mapsto (\gamma^u\cdot(\alpha\times 0), \gamma^o\cdot(\alpha\times 1)),\\
(\gamma^u,\gamma^m,\gamma^o;\alpha)\mapsto (\gamma^u\cdot(\alpha\times 0), \gamma^m, \gamma^o\cdot(\alpha\times 1)).
\end{gather*}

Then $\mathdutchcal R(T)=\widetilde{\mathdutchcal R}_0(T)/\pi_1(F,x_0)$, $\mathdutchcal{SA}(T)=\widetilde{\mathdutchcal {SA}}_0(T)/\pi_1(F,x_0)$, $\mathdutchcal C(T)=\widetilde{\mathdutchcal C}_0(T)/\pi_1(F,x_0)$.

In addition to the action of the fundamental group $\pi_1(F,x_0)$ on homotopy classes of diagram elements, we consider the following operations on them.

\begin{enumerate}
    \item Incidence relations. The formulas in Definition~\ref{def:element_incidence_maps} induce maps between $\mathdutchcal A(T)$, $\mathdutchcal R(T)$, $\mathdutchcal{SA}(T)$ and $\mathdutchcal C(T)$. These incidence maps form a diagram analogous to~\eqref{eq:element_incidence_maps_diagram}.

    \item Action of the tangle group. For a tangle $T$ in $F\times I$, denote $\pi^u(T)=\pi_1(M_T,x^u)$, $\pi^o(T)=\pi_1(M_T,x^o)$, and
\begin{gather*}
\pi^u_0(T)=\ker(p_*\colon\pi_1(M_T,x^u)\to\pi_1(F,x_0)),\\
\pi^o_0(T)=\ker(p_*\colon\pi_1(M_T,x^o)\to\pi_1(F,x_0)).
\end{gather*}
The tangle group $\pi^o(T)$ acts on the sets of homotopy classes $\widetilde{\mathdutchcal A}(T)$, $\widetilde{\mathdutchcal R}(T)$, $\widetilde{\mathdutchcal{SA}}(T)$ and $\widetilde{\mathdutchcal C}(T)$ according to the formulas
\begin{gather*}
(\gamma,\alpha)\mapsto \gamma\cdot\alpha,\quad \gamma\in\widetilde{\mathdutchcal A}(T)\mbox{ or }\widetilde{\mathdutchcal R}(T),\\
(\gamma^u,\gamma^o;\alpha)\mapsto (\gamma^u, \gamma^o\cdot\alpha),\quad (\gamma^u,\gamma^o)\in\widetilde{\mathdutchcal{SA}}(T),\\
(\gamma^u,\gamma^m,\gamma^o;\alpha)\mapsto (\gamma^u, \gamma^m, \gamma^o\cdot\alpha),\quad (\gamma^u,\gamma^m, \gamma^o)\in\widetilde{\mathdutchcal R}(T),
\end{gather*}
where $\alpha\in\pi^o(T)$. The group $\pi^u(T)$ acts on $\widetilde{\mathdutchcal R}(T)$, $\widetilde{\mathdutchcal{SA}}(T)$ and $\widetilde{\mathdutchcal C}(T)$ by the formulas
\begin{gather*}
(\gamma,\alpha)\mapsto \alpha^{-1}\cdot\gamma,\quad \gamma\in\widetilde{\mathdutchcal R}(T),\\
(\gamma^u,\gamma^o;\alpha)\mapsto (\gamma^u\cdot\alpha, \gamma^o),\quad (\gamma^u,\gamma^o)\in\widetilde{\mathdutchcal{SA}}(T)\\
(\gamma^u,\gamma^m,\gamma^o;\alpha)\mapsto (\gamma^u\cdot\alpha, \gamma^m, \gamma^o),\quad (\gamma^u,\gamma^m, \gamma^o)\in\widetilde{\mathdutchcal R}(T).
\end{gather*}
Using the same formulas, one defines an action of $\pi^o_0(T)$ on $\widetilde{\mathdutchcal R}_0(T)$, $\widetilde{\mathdutchcal{SA}}_0(T)$, $\widetilde{\mathdutchcal C}_0(T)$ and an action of $\pi^u_0(T)$ on $\widetilde{\mathdutchcal R}_0(T)$, $\widetilde{\mathdutchcal{SA}}_0(T)$ and $\widetilde{\mathdutchcal C}_0(T)$.

    \item Augmentation maps. Consider the maps $\epsilon^o$ from $\widetilde{\mathdutchcal A}(T)$, $\widetilde{\mathdutchcal{SA}}(T)$ and $\widetilde{\mathdutchcal C}(T)$ to the tangle group $\pi^o(T)$ given by the formulas
\begin{gather*}
\gamma\mapsto\gamma^{-1}\cdot\mu\cdot\gamma,\quad \gamma\in\widetilde{\mathdutchcal A}(T),\\
(\gamma^u,\gamma^o)\mapsto(\gamma^o)^{-1}\cdot\mu_r\cdot\mu_l\cdot\gamma^o,\quad (\gamma^u,\gamma^o)\in\widetilde{\mathdutchcal{SA}}(T)\\
(\gamma^u,\gamma^m,\gamma^o)\mapsto(\gamma^o)^{-1}\cdot\mu^o_r\cdot\mu^o_l\cdot\gamma^o,\quad (\gamma^u,\gamma^m, \gamma^o)\in\widetilde{\mathdutchcal R}(T).
\end{gather*}
where $\mu$ and $\mu^o$ are the meridians of the elements. Analogously on defines the augmentation map $\epsilon^u$ from $\widetilde{\mathdutchcal{SA}}(T)$ and $\widetilde{\mathdutchcal C}(T)$ to $\pi^u(T)$:
\begin{gather*}
(\gamma^u,\gamma^o)\mapsto(\gamma^u)^{-1}\cdot\mu_l\cdot\mu_r\cdot\gamma^u,\quad (\gamma^u,\gamma^o)\in\widetilde{\mathdutchcal{SA}}(T)\\
(\gamma^u,\gamma^m,\gamma^o)\mapsto(\gamma^u)^{-1}\cdot\mu^u_l\cdot\mu^u_r\cdot\gamma^u,\quad (\gamma^u,\gamma^m, \gamma^o)\in\widetilde{\mathdutchcal R}(T).
\end{gather*}
Note that the image of the maps $\epsilon^o$ and $\epsilon^u$ lies in $\pi^o_0(T)$ and $\pi^u_0(T)$, respectively.
\end{enumerate}


Let $T=T_1\sqcup\cdots\sqcup T_n$ be a tangle in $F\times I$. Denote a tubular neighborhood of the component $T_i$ by $N_i$.  Choose arbitrary points $z_i\in\partial N_i$ and paths $\delta_i$ in $M_T$ from $x^o$ to $z_i$, $i=1,\dots,n$. The paths $\delta_i$ 
allow us to identify the groups $\pi_1(\partial N_i,z_i)$ with subgroups in $\pi^o(T)$. 
Using an abuse of notation, we will denote these subgroups by $\pi_1(\partial N_i)$. Note that $\pi_1(\partial N_i)\simeq\Z$ when $T_i$ is a long or a trivial component of $T$, and $\pi_1(\partial N_i)\simeq\Z^2$ otherwise.

Given an arc probe $\gamma\in\widetilde{\mathdutchcal{A}}(T)$ such that $\gamma(0)\in\partial N_i$, $1\le i\le n$, we move the point $\gamma(0)$ in $\partial N_i$ to $z_i$ and get a homotopic path $\gamma'$ such that $\gamma'(0)=z_i$. Then $\delta_i\gamma'$ defines an element of $\pi^o(T)$.

Denote the vertical path from $x^o$ to $x^u$ by $\delta_0$.
Given a region probe $\gamma\in\widetilde{\mathdutchcal{R}}(T)$, the concatenation $\delta_0\gamma$ defines an element in $\pi^o(T)$.

If $(\gamma^u,\gamma^o)\in\widetilde{\mathdutchcal{SA}}(T)$ is a semiarc probe where $\gamma^u(0)\in\partial N_i$, then contract the part of meridian $\mu_r$ in $\partial N_i$ to a point and move this point to $z_i$. This homotopy transforms the paths $\gamma^u,\gamma^o$ into some paths $\gamma^u_1,\gamma^o_1$. Then we get a pair of elements $(\delta_i\gamma^u_1\delta_0^{-1},\delta_i\gamma^o_1)\in\pi^o(T)\times\pi^o(T)$.

Let $\gamma=(\gamma^u,\gamma^m,\gamma^o)\in\widetilde{\mathdutchcal{C}}(T)$ be a crossing probe such that $\gamma^m(0)\in\partial N_i$ and $\gamma^m(1)\in\partial N_j$.  If $i\ne j$, then contract the meridian parts $\mu^o_r\subset\partial N_i$ and $\mu^u_r\subset\partial N_j$ to the points and move these points to $z_i$ and $z_j$. Denote the paths obtained by this homotopy by $\gamma_1=(\gamma^u_1,\gamma^m_1,\gamma^o_1)$.

If $i=j$ then choose points $z_i^+$ (resp. $z_i^-$) in $T_i$ that differ from $z_i$ by a small shift along (resp. against) the orientation of the component. One can isotope the $\gamma$ in $\widetilde{\mathdutchcal{C}}(T)$ to a probe $\gamma_0=(\gamma^u_0,\gamma^m_0,\gamma^o_0)$ such that $\{\gamma^m_0(0),\gamma^m_0(1)\}=\{z_i^-,z_i^+\}$.

 If $T_i$ is a long component, one can distinguish the case of early overcrossing $\gamma^m_0(0)=z_i^-$ and the case of early undercrossing $\gamma^m_0(0)=z_i^+$. In the first case, we write $o(\gamma)=+$, and in the latter case, $o(\gamma)=-$. If $T_i$ is closed, we can suppose that $\gamma^m_0(0)=z_i^-$ and $\gamma^m_0(1)=z_i^+$. Then contract the meridian parts $\mu^u_r$ and $\mu^o_r$ of the probe $\gamma_0$ into the point $z_i$.  We get a path $\gamma_1=(\gamma^u_1,\gamma^m_1,\gamma^o_1)$ by the induced homotopy. We can assume that the linking number $lk(\gamma^m_1, T_i)$ vanishes.

 Then $\delta_i\gamma_1^o$ and $\delta_j\gamma_1^u\delta_0^{-1}$ can be viewed as elements of $\pi^o(T)$. Choosing a framing for each path $\delta_i$, we get a framed loop $\delta_i\gamma_1^m\delta_j^{-1}$ in the set $\pi^o_{fr}(T)$ of regular homotopy classes of framed loops when $i\ne j$, and in the set
 \[
 \pi^o_{fr,i}(T)=\ker\left(lk(\cdot,T_i)\colon \pi_1^{fr}(M_T,x^o)\to\Z\right)
 \]
 of regular homotopy classes of framed loops $\zeta$ such that $lk(\zeta,T_i)=0$ when $i=j$. By choosing a framing for each regular homotopy class, we get a bijection $\pi^u_{fr}(T)\simeq\pi^u(T)\times\Z_2\times Z_2$ where the first multiplier $\Z_2$ corresponds to the sign of the crossing, and the second $\Z_2$ corresponds to framing (we have $\Z_2$ instead of $\Z$ because we admit self-intersections of framed loops).

 Taking into account the ambiguity of reducing probes to loops in $\pi^o(T)$, we obtain the following statement.

\begin{proposition}\label{prop:element_homotopy_probes_loop_reduction}
Let $T=T_1\sqcup\cdots\sqcup T_n$ be a tangle in $F\times I$. The correspondences described above, establish the following bijections:
\begin{gather*}
    \widetilde{\mathdutchcal A}(T)\simeq \prod_{i=1}^n \pi_1(\partial N_i)\backslash\pi^o(T),\\
    \widetilde{\mathdutchcal R}(T)\simeq\pi^o(T),\\
    \widetilde{\mathdutchcal{SA}}(T)\simeq \prod_{i=1}^n \pi_1(\partial N_i)\backslash(\pi^o(T)\times\pi^o(T)),
\end{gather*}
and
\begin{multline*}
    \widetilde{\mathdutchcal C}(T)\simeq \prod_{i\ne j} (\pi_1(\partial N_i)\times\pi_1(\partial N_j))\backslash(\pi^o(T)\times\pi^o_{fr}(T)\times\pi^o(T)) \times\\ \prod_{i\colon T_i\mbox{ \scriptsize is long}}\left[\pi_1(\partial N_i)\backslash(\pi^o(T)\times\pi^o_{fr,i}(T)\times\pi^o(T))\right]^{\times 2}\times\\
    \prod_{i\colon T_i\mbox{ \scriptsize is closed}}\pi_1(\partial N_i)\backslash(\pi^o(T)\times\pi^o_{fr,i}(T)\times\pi^o(T)),
\end{multline*}
where the action of $\alpha\in\pi_1(\partial N_i)$ on $\zeta\in\pi^o(T)$ is given by the formula
\[
\alpha\cdot\zeta=\alpha\zeta;
\]
the action of $\alpha\in\pi_1(\partial N_i)$ on $(\zeta^u,\zeta^o)\in\pi^o(T)\times\pi^o(T)$  is given by the formula
\[
\alpha\cdot(\zeta^u,\zeta^o)=(\alpha\zeta^u,\alpha\zeta^o);
\]
the action of $(\alpha^o,\alpha^u)\in\pi_1(\partial N_i)\times\pi_1(\partial N_j)$ on $(\zeta^u,\zeta^m,\zeta^o)\in\pi^o(T)\times\pi^o_{fr}(T)\times\pi^o(T)$  is given by the formula
\[
(\alpha^o,\alpha^u)\cdot(\zeta^u,\zeta^m,\zeta^o)=(\alpha^u\zeta^u,\alpha^o\zeta^m(\alpha^u)^{-1},\alpha^o\zeta^o);
\]
and the action of $\alpha\in\pi_1(\partial N_i)$ on $(\zeta^u,\zeta^m,\zeta^o)\in\pi^o(T)\times\pi^o_{fr,i}(T)\times\pi^o(T)$  is given by the formula
\[
\alpha\cdot(\zeta^u,\zeta^m,\zeta^o)=(\alpha\zeta^u,\alpha\zeta^m\alpha^{-1}, \alpha\zeta^o).
\]
\end{proposition}

\begin{remark}\label{rem:element_homotopy_probe_loop_reduction}
    The sets $\mathdutchcal A(T), \mathdutchcal R(T), \mathdutchcal{SA}(T), \mathdutchcal C(T)$ can be described as sets of double cosets of the tangle group $\pi^u(T)$ and its products by the subgroups $\pi_1(\partial N_i)$ and $\pi_1(F,x_0)$.
\end{remark}


Given a coinvariant of diagram elements, we can look at it as a coloring of diagram elements with marks from some set $X$. Among all possible colorings, we take only those that satisfy some coloring propagation rule: the color of some diagram elements determines the color of other diagram elements.
In the subsequent sections, we will consider the following coloring propagation rules:
\begin{itemize}
    \item the colors of a pair of adjacent (semi)arc in a crossing uniquely determine the colors of the other (semi)arcs incident to the crossing;
    \item the colors of three regions incident to a crossing uniquely determine the color of the fourth region incident to the crossing;
    \item the color of a crossing incident to a region is uniquely determined by the colors of the other crossings of the region.
\end{itemize}
The propagation rules induce algebraic structures (such as a quandle or a group structure) in the set of marks $X$. It turns out that the set of homotopy classes of diagram elements defined above are universal objects among objects with these algebraic structures.

To illustrate the notion of propagation rule, consider the following examples.

\begin{example}[Adjacent region propagation rule]\label{exa:propagation_adjacent_region}
    For a coinvariant $\mathcal G$ of regions with values in a set $X$, consider the rule: the color of a region uniquely determines the color of an adjacent region. Formally, there is a map $\phi\colon X\to X$ such that for any regions $r,r'$ such that $r\uparrow r'$ with colors $c=\mathcal G(r)$ and $c'=\mathcal G(r')$, one has $c'=\phi(c)$ (Fig.~\ref{pic:propagation_rule_region_neighbour}).
\begin{figure}[h]
    \centering
    \includegraphics[width=0.1\textwidth]{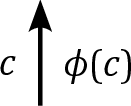}
    \caption{The color of a region determines the colors of the adjacent regions}
    \label{pic:propagation_rule_region_neighbour}
\end{figure}

    Let $T$ be a tangle in a connected compact oriented surface $F$ and $D=p(T)$ its diagram. Choose a point $x_0\in F\setminus D$ and denote the region containing $x_0$ by $r_0\in\mathcal R(D)$, and the color of $r_0$ by $c_0=\mathcal G(r_0)$. Then for any $r\in\mathcal R(D)$  $\mathcal G(r)=\phi^k(c_0)$, $k\in\Z$. The number $k$ can be calculated by the intersection formula $k=\gamma\cdot D$ where $\gamma$ is a path from $x_0$ to a point $x\in r$. The number $k$ is determined modulo
\[
m=gcd\{\gamma\cdot D\mid \gamma\in\pi_1(F,x_0)\}.
\]
Thus, the coinvariant $\mathcal G$ is determined by the coinvariant $N_{x_0}(r)=\gamma\cdot D\in\Z/m\Z$ that is the Alexander numbering of regions~\cite{Alexander}.
\end{example}

\begin{example}[Opposite semiarc propagation rule]\label{exa:propagation_rule_opposite_semiarc}
    For a coinvariant $\mathcal G$ of semiarcs with values in a set $X$, consider the rule: the color of a semiarc incident to a crossing uniquely determines the color of the opposite arc. Depending on the position of the incident semiarc and the sign of the crossing, there are four maps $\phi_{u+}$, $\phi_{u-}$, $\phi_{o+}$, $\phi_{o-}$ from $X$ to $X$ (Fig.~\ref{pic:propagation_rule_semiarc_opposite}).
\begin{figure}[h]
    \centering
    \includegraphics[width=0.5\textwidth]{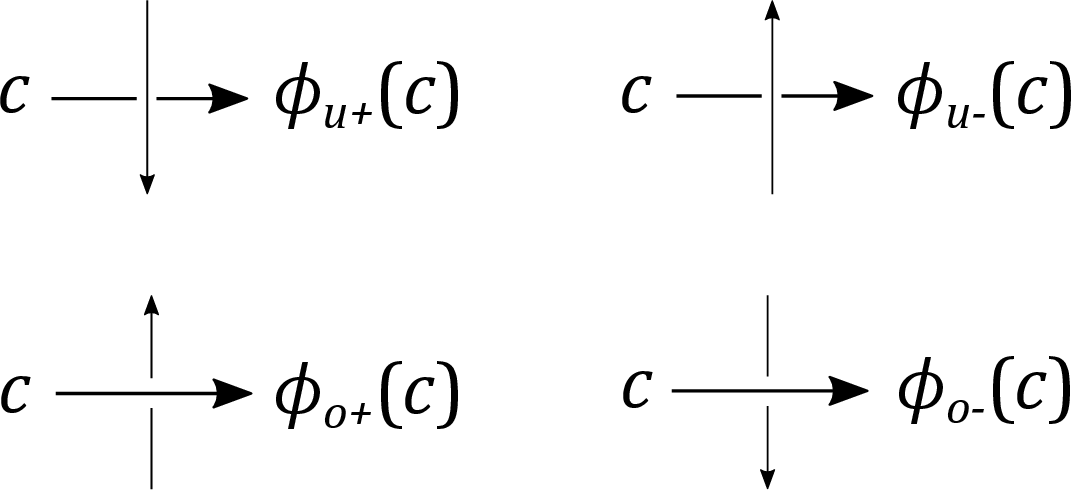}
    \caption{Maps of the propagation rule}
    \label{pic:propagation_rule_semiarc_opposite}
\end{figure}
    The invariance under third Reidemeister move implies that the operators $\phi_{u+}$, $\phi_{u-}$, $\phi_{o+}$, $\phi_{o-}$ commute. The second Reidemeister move implies $\phi_{\alpha+}\phi_{\alpha-}=id_X$ where $\alpha\in\{u,o\}$. The first Reidemeister move $\phi_{\alpha\epsilon}\phi_{\bar\alpha\bar\epsilon}=id_X$ where $\alpha\in\{u,o\}$, $\epsilon\in\{-,+\}$. Hence, $\phi_{u+}=\phi_{o+}=\phi_{u-}^{-1}=\phi_{o-}^{-1}$. Denote this map by $\phi$.

    Given a tangle diagram $D=D_1\cup\cdots\cup D_n$, choose a non-crossing point $x_i\in D_i$ in each component. Let $a_i\in\mathcal{SA}(D)$ be the semiarc containing $x_i$, and $c_i=\mathcal G(a_i)$ its color. Then for any semiarc $a\in\mathcal{SA}(D)$, $\mathcal G(a)=\phi^k(a_i)$ where $D_i$ is the component in which the semiarc $a$ lies. The number $k$ is calculated by the intersection formula $k=\gamma\cdot D\in\Z/m_i\Z$, where $\gamma\subset D_i$ is a path from $x_i$ to a point $x\in a$, and
\[
m_i=\left\{\begin{array}{cl}
    D_i\cdot D, & D_i\mbox{ is closed}, \\
     0, & D_i\mbox{ is long}.
\end{array}\right.
\]
Thus, the coinvariant value $\mathcal G(a)$ of an arc $a$ is determined by the component index $i$ and the Alexander numbering $N_{a_i}(a)=\gamma\cdot D$. The local behaviour of the Alexander numbering is shown in Fig.~\ref{pic:alexander_numbering}.
\begin{figure}[h]
    \centering
    \includegraphics[width=0.15\textwidth]{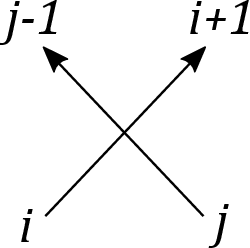}
    \caption{Alexander numbering rule~\cite{SW} (under-overcrossing structure does not matter)}
    \label{pic:alexander_numbering}
\end{figure}
\end{example}


\subsection{Colorings of arcs: Quandle}\label{subsect:quandle}

As we saw above, colorings of the arcs of a tangle diagram with elements from some set $X$ in a way that is compatible with Reidemeister moves, define coinvariants of the arc functor. Now, let us assume that the colors of two adjacent arcs at a crossing uniquely determine the color of the third arc. This assumption implies that an algebraic structure called a quandle exists in the color set $X$.

\begin{definition}[\cite{Joyce, Matveev}]\label{def:quandle}
A set $X$ with a binary operation $\ast\colon X\times X\to X$ is called a \emph{quandle} if it obeys the following conditions:
\begin{enumerate}
\item $x\ast x=x$ for any $x\in X$;
\item for any $y\in X$ the operator $\alpha_y\colon X\to X, x\mapsto x\ast y$ is invertible;
\item for any $x,y,z\in X$
\begin{gather*}
(x\ast y)\ast z=(x\ast z)\ast(y\ast z).
\end{gather*}
\end{enumerate}
\end{definition}

\begin{definition}\label{def:quandle_coloring}
Let $D$ be an oriented tangle diagram on an oriented connected compact surface $F$, and $(X,\ast)$ a quandle. Then a {\em coloring of the diagram $D$ with the quandle} is a map from the set of arcs of $D$ to $X$ such that the images of the arcs (colors) satisfy the coloring rule (Fig.~\ref{pic:quandle_coloring}).
\begin{figure}[h]
\centering\includegraphics[width=0.15\textwidth]{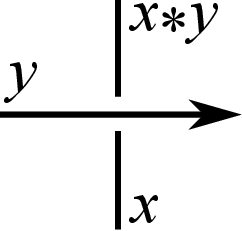}
\caption{The coloring rule}\label{pic:quandle_coloring}
\end{figure}
Let $Col_{X}(D)$ be the set of quandle colorings.
\end{definition}

\begin{theorem}[\cite{Joyce, Matveev}]\label{thm:quandle_coloring_coinvariance}
For any diagrams $D$ and $D'$ connected by a Reidemeister move, there is a bijection between $Col_{X}(D)$ and $Col_{X}(D')$.
\end{theorem}

\begin{example}[Alexander quandle]
Let $X$ be a module over the ring $\mathbb Z[t^{\pm 1}]$. Then the operation
\[
x\ast y = tx+(1-t)y
\]
defines a quandle operation called an \emph{Alexander quandle}.
\end{example}

\begin{example}[Fundamental quandle]
Let $D$ a diagram of a tangle $T$. The \emph{fundamental quandle} $FQ(D)$ of $T$ is determined by the universal rule: for any quandle $(X,\ast)$ and any coloring $c\in Col_X(D)$ there exists a unique quandle homomorphism $f\colon FQ(D)\to X$ such that $c=f(c_F)$ where $c_F\in Col_{FQ(D)}(D)$ is a fixed (fundamental) coloring. The fundamental quandle has a presentation of the form
\[
FQ(D)=\langle\mbox{arcs of $D$}\,\mid\,\mbox{coloring rule at the crossings}\rangle.
\]
\end{example}

\begin{example}[Topological quandle]\label{exa:topological_quandle}
Let $T$ be an oriented tangle in the thickening $F\times I$ of a connected oriented compact surface $F$.
Consider the set $\widetilde{\mathdutchcal A}(T)$ of based homotopy classes of arcs.
 Then the operation
\[
\gamma_1\ast\gamma_2=\gamma_1\cdot\epsilon^o(\gamma_2)
\]
defines a quandle structure on $\tilde{\mathscr A}(T)$.
\end{example}

We will prove that the topological quandle is fundamental in some sense. More precisely, the following theorem holds.

\begin{theorem}\label{thm:quandle_topological_fundamental}
Let $T\subset F\times I$ be a tangle, and $D$ its diagram. For any quandle $(X,\ast)$ there is a bijection between the set of colorings $Col_{X}(D)$ and the set of quandle homomorphisms $Hom(\widetilde{\mathdutchcal A}(T),X)^{\pi_1(F)}$ invariant under the action of $\pi_1(F,z)$ on $\widetilde{\mathdutchcal A}(T)$.
\end{theorem}

\begin{proof}
1. Let $\phi\colon \widetilde{\mathdutchcal A}(T)\to X$ be a $\pi_1(F,z)$-invariant quandle homomorphism. For an arc $a\in\mathcal A(D)$, choose a point $x\in a$ and the vertical arc probe $\gamma_a\subset x\times I$. Choose an arbitrary path $\delta_x\subset F$ from $x$ to $x_0$. Then $\gamma_a\delta_x\in\widetilde{\mathdutchcal A}(T)$. Define the color $c_\phi(a)$ of the arc $a$ by formula $c_\phi(a)=\phi(\gamma_a\delta_x)$. The element $c_\phi(a)$ does not depend on $x$ and $\delta_x$ by $\pi_1(F,z)$-invariance of $\phi$.

Let us check that the map $c_\phi\colon\mathcal A(D)\to X$ is a quandle coloring. Let  $x,y$ be arcs incident to a crossing of the diagram $D$ (Fig.~\ref{pic:quandle_coloring}). Denote the third arc by $z$. We need to show that $c_\phi(z)=c_\phi(x)\ast c_\phi(y)$. Choose points on the arcs close to the crossing, a path $\delta$ in $F$ from the crossing point to $x_0$. Using paths close to $\delta$, construct arc probes $\gamma_x,\gamma_y,\gamma_z\in\widetilde{\mathdutchcal A}(T)$. Then $\gamma_z=\gamma_x\ast\gamma_y$. Hence,
\[
c_\phi(z)=\phi(\gamma_z)=\phi(\gamma_x)\ast\phi(\gamma_y)=c_\phi(x)\ast c_\phi(y).
\]

2. Now, let $c\in Col_{X}(D)$ be a quandle coloring. We will construct a map $\phi_c\colon \mathdutchcal A(T)\to X$.

Let $D\cup\gamma$ be an arc probe diagram. Construct another arc probe diagram $D'\cup\gamma$ by pulling the arcs of $D$ overcrossing $\gamma$ to the initial point of $\gamma$ (Fig.~\ref{pic:quandle_probe_clearing}). The transformation of the diagram $D$ is a morphism $f\colon D\to D'$. This morphism induces a bijection $f_*\colon Col_X(D)\to Col_X(D')$. Then we set the values $\phi(\gamma)$ equal to the color $(f_*(c))(a')$ of the arc $a'$ in $D'$ where the probe $\gamma$ begins.
\begin{figure}[h]
    \centering
    \includegraphics[width=0.7\textwidth]{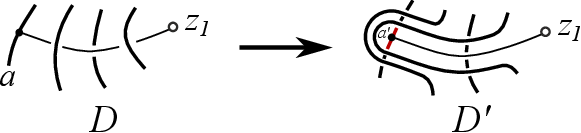}
    \caption{Modification of the diagram along the arc probe}
    \label{pic:quandle_probe_clearing}
\end{figure}

Let us check that the value $\phi(\gamma)$ does not change during isotopy of $\gamma$. If $g\colon D\to D_1$ is a second or a third Reidemeister move, then the corresponding transformed diagrams $D'$ and $D'_1$ differ by a sequence of Reidemeister moves of the same type that do not involve the initial point of $\gamma$ with a small neighborhood (Fig.~\ref{pic:quandle_R2_R3_enveloped}). Then the color of the arc $a'$ does not change during transformation from $D'$ to $D'_1$. If $g$ is a first Reidemeister move, then we can ignore the loop by Lemma~\ref{lem:quandle_reciprocial_crossing} below.

\begin{figure}[h]
    \centering
    \includegraphics[width=0.35\textwidth]{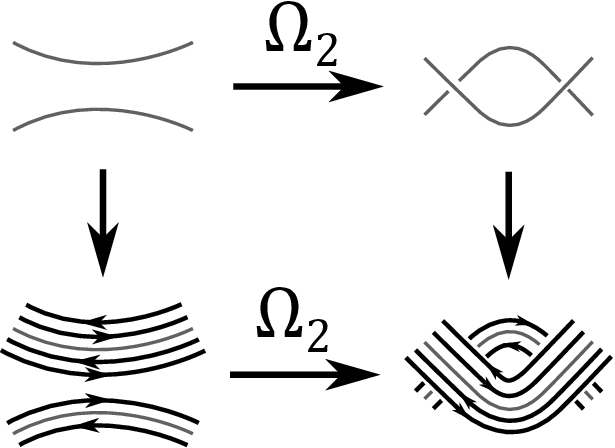}\qquad \includegraphics[width=0.45\textwidth]{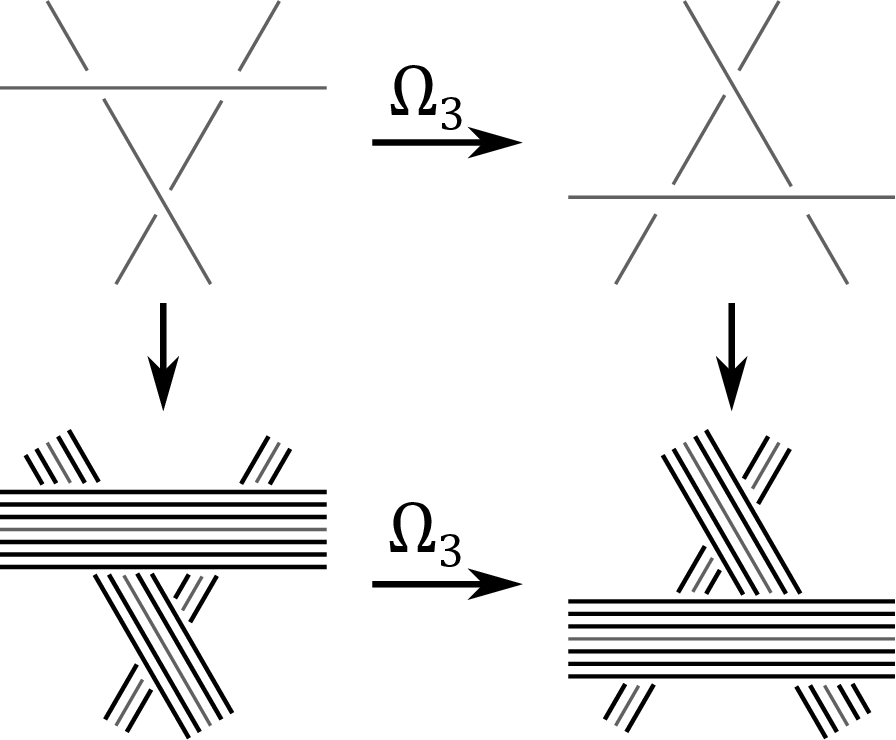}
    \caption{Second and third Reidemeister moves}
    \label{pic:quandle_R2_R3_enveloped}
\end{figure}

A second or a third Reidemeister move $f\colon D\to D_1$ including an arc of the diagram $D$, induces diagrams $D'$ and $D'_1$ connected by second and third Reidemeister moves (Fig.~\ref{pic:quandle_mixed_R2_R3_enveloped}).

\begin{figure}[h]
    \centering
    \includegraphics[width=0.45\textwidth]{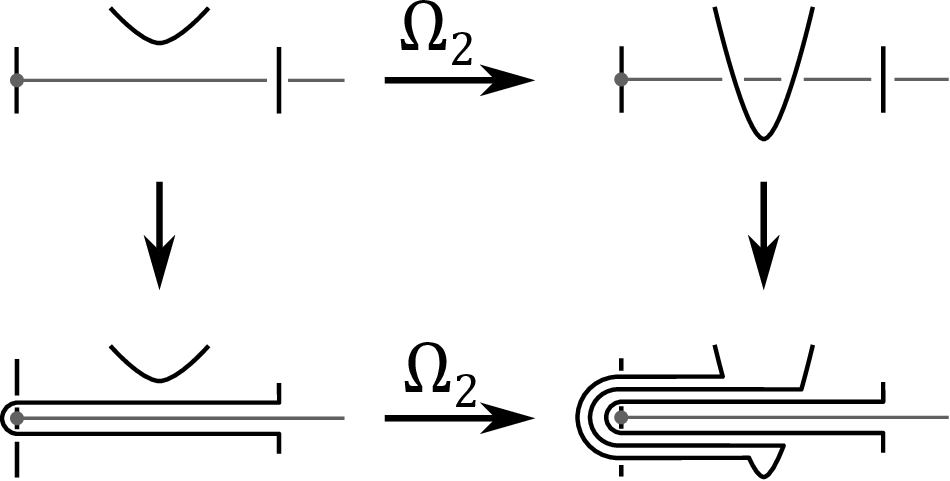}\quad \includegraphics[width=0.45\textwidth]{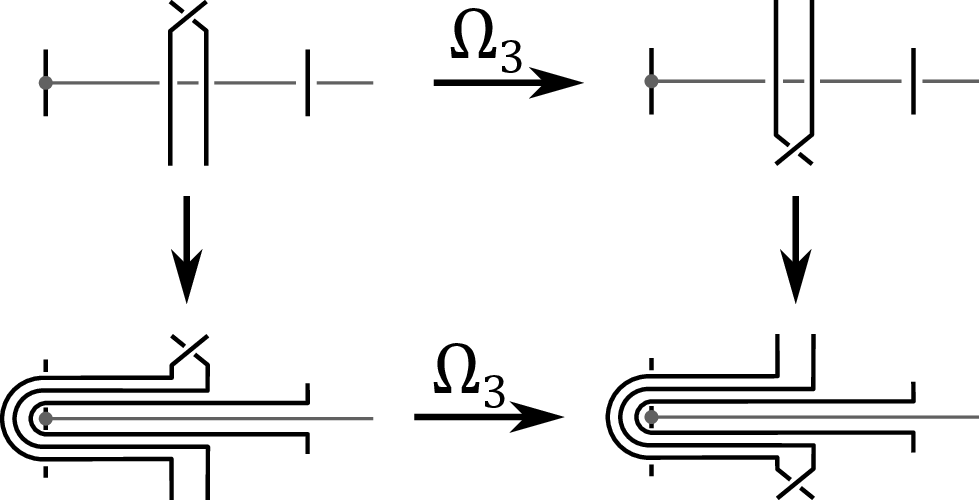}
    \caption{Second and third Reidemeister moves with an arc of the diagram}
    \label{pic:quandle_mixed_R2_R3_enveloped}
\end{figure}

If $f\colon D\to D_1$ is a float move (Fig.~\ref{pic:graphoid_float_sinker_moves}), then the corresponding transformed diagrams $D'$ and $D'_1$ differ by a float move which does not change the color of $a'$.

If $f\colon D\to D_1$ is a vertex rotation move (Fig.~\ref{pic:quandle_rotation_move} top line), then after transformation we get the diagram in Fig.~\ref{pic:quandle_rotation_move} bottom right. We see that a first Reidemeister move merges the old arc $a$ and the new arc $a'$. Hence, they have the same color.

\begin{figure}[h]
    \centering
    \includegraphics[width=0.3\textwidth]{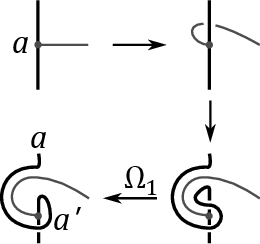}
    \caption{Rotation move}
    \label{pic:quandle_rotation_move}
\end{figure}

Thus, we get a map $\phi\colon\mathscr A^s(T)\to X$ from the isotopy classes of arc probes to the quandle.

Then let us show that self-intersections of arc probes do not change their colors. We need the following lemma.

\begin{lemma}\label{lem:quandle_reciprocial_crossing}
    Consider a part of the tangle that consists of $2n$ parallel arcs $a_1,\dots, a_{2n}$ and another arc $b$ such that the arcs $a_i$ and $a_{2n+1-i}$ have opposite orientation (Fig.~\ref{pic:quandle_reciprocial_crossing} left). Assume that the arcs are colored by the quandle $X$ so that the colors of the arcs $a_i$ and $a_{2n+1-i}$ coincide. Then after applying second Reidemiester moves, we get a colored part of the tangle (Fig.~\ref{pic:quandle_reciprocial_crossing} right) such that the arcs $b'$ and $b$ have the same color, as well as the arcs $a_i$ and $a'_i$, $i=1,\dots,2n$, as well as the arcs $b_i$ and $b'_i$, $i=1,\dots,2n-1$.
\begin{figure}[h]
    \centering
    \includegraphics[width=0.7\textwidth]{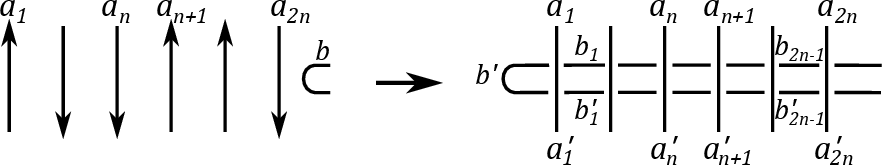}
    \caption{A colored tangle}
    \label{pic:quandle_reciprocial_crossing}
\end{figure}
\end{lemma}

\begin{proof}
The equality of arc colors for the arcs $a_i$ and $a'_i$, as well as $b_i$ and  $b'_i$ follows from the behavior of quandle colorings under Reidemeister moves. The proof for the arcs $b$ and $b'$ is given in Fig.~\ref{pic:quandle_reciprocial_crossing_proof}. Since the arcs $a_i$ and $a_{2n+1-i}$ have the same colors, we can connect them with an arc and then merge the arcs $b$ and $b'$ using second Reidemeister moves.
\begin{figure}[h]
    \centering
    \includegraphics[width=0.7\textwidth]{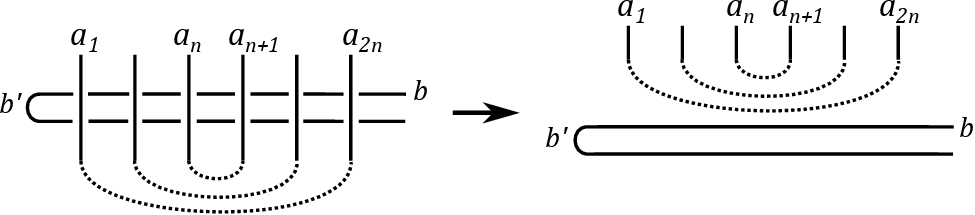}
    \caption{The arcs $b$ and $b'$ have the same color}
    \label{pic:quandle_reciprocial_crossing_proof}
\end{figure}
\end{proof}

The lemma implies that if an arc passes a self-intersection point of the probe during the pulling process, then its color does not change regardless of whether the arc passes above or below the crossing. Hence, if one switches the undercrossing and the overcrossing of a self-intersection of the probe, then the result does not change. Thus, the map $\phi$ induces a map from $\mathdutchcal A(T)$ to $X$, therefore, induces a $\pi_1(F)$-invariant map from $\phi\colon\widetilde{\mathdutchcal A}(T)\to X$.

The next claim is that the map $\phi$ is a natural transformation from the functor $\widetilde{\mathdutchcal A}$ to the constant functor $X$. It is enough to show that for any arc probe diagram $D\cup\gamma$ and for any Reidemeister move $f\colon D\to D_1$, the colors of the arcs $\gamma$ and $\mathdutchcal A(f)(\gamma)$ coincide.

Since the map $\phi$ gives the same value for isotopic arc probes, we can assume that $\gamma$ is distinct from the area where the move occurs. Then the transformed diagrams $D'$ and $D'_1$ are connected by the same move as $D$ and $D_1$, and this move does not involve the arc where the probe $\gamma$ starts. Hence, the color of this arc does not change with the move. Then the colors of the arcs $\gamma$ and $\mathdutchcal A(f)(\gamma)$ coincide.

Finally, let us show that the map $\phi\colon\widetilde{\mathdutchcal A}(T)\to X$ is a homomorphism of quandles. Let $\gamma_1,\gamma_2$ be arc probes of the tangle $T$. We need to prove that $\phi(\gamma_1\ast\gamma_2)=\phi(\gamma_1)\ast\phi(\gamma_2)$. Consider an isotopy $f\colon T\to T'$ which verticalize the probes and makes them form a crossing, i.e. $\mathdutchcal A(f)(\gamma_1)=\gamma_x$, $\mathdutchcal A(f)(\gamma_2)=\gamma_y$ where the arcs $x,y\in\mathcal A(D')$, $D'=p(T')$, incident to a crossing as in Fig.~\ref{pic:quandle_coloring}. Let $z$ be the third arc incident to the crossing. Then
\begin{multline*}
\phi(\gamma_1)\ast\phi(\gamma_2)=\phi_{D'}(\gamma_x)\ast\phi_{D'}(\gamma_y)=(f_*(c))(x)\ast (f_*(c))(y)=f_*(c)(z)=\\
\phi_{D'}(\gamma_x\ast\gamma_y)=\phi(\gamma_1\ast\gamma_2),
\end{multline*}
where we use the fact that $f_*(c)$ is a quandle coloring of $D'$, and that $\phi$ is a natural map.
The theorem is proved.
\end{proof}

Let $T$ be a tangle in the thickening $F\times I$ of a connected compact oriented surface $F$. Let $N(T)$ be a tubular neighborhood of $T$. Consider the space
\[
CF_T=(F\times I\setminus N(T))/F\times 1.
\]
There is a natural projection $M_T\to CF_T$. Let $x^o\in CF_T$ be the image of $F\times 1$ by this projection. Consider the set of homotopy classes of paths
\[
\bar{\mathdutchcal A}(T)=[I,0,1;CF_T,\partial N(T), x^o].
\]
The projection $M_T\to CF_T$ induces a surjection $\pi\colon\widetilde{\mathdutchcal A}(T)\to\bar{\mathdutchcal A}(T)$ that induces a quandle structure on $\bar{\mathdutchcal A}(T)$.

\begin{proposition}\label{prop:quandle_topological_reduction}
    For any quandle $X$ and a $\pi_1(F)$-invariant quandle homomorphism $\phi\colon\widetilde{\mathdutchcal A}(T)\to X$ there exists a unique quandle homomorphism $\bar\phi\colon\bar{\mathdutchcal A}(T)\to X$ such that $\phi=\bar\phi\circ\pi$.
\end{proposition}

\begin{proof}
    Since $\pi$ is a surjection, we define $\bar\phi=\phi\circ\pi^{-1}$. We need to show that this map is well defined, that is, for any $\gamma_1,\gamma_2\in\widetilde{\mathdutchcal A}(T)$ such that $\pi(\gamma_1)=\pi(y\gamma_2)$ one has $\phi(\gamma_1)=\phi(\gamma_2)$.

    One can homotope the probe $\gamma_2$ so that $\gamma_1(0)=\gamma_2(0)$ and the image of $\gamma_2^{-1}\gamma_1$ in $CF_T$ is contractible. The kernel of the map $\pi_1(M_T)\to\pi_1(CF_1)$ is generated by elements $\alpha\beta\alpha^{-1}$ where $\alpha\in\pi_1(M_T)$, and $\beta$ lies in $\pi_1(F\times 1,x^o)$. Hence, it is enough to check the case $\gamma_1=\gamma_2\alpha\beta\alpha^{-1}$.

    There is a decomposition $\alpha=\alpha_0\beta_0$ where $\alpha_0\in\pi^o_0(T)$ and $\beta_0\in\pi_1(F\times 1,x^u)$. Then
\[
\alpha\beta\alpha^{-1}=\alpha_0(\beta_0\beta\beta_0^{-1})\alpha_0^{-1}.
\]
    After replacing $\alpha$ with $\alpha_0$ and $\beta$ with $\beta_0\beta\beta_0^{-1}$, we can assume that $\alpha\in\pi^o_0(T)$. 

    Since $\pi_1(F\times I,x^o)=\pi_1(F,x_0)$,
\[
\pi^o_0(T)=\ker(p_*\colon \pi_1(M_T,x^o)\to\pi_1(F,x_0))=\ker(p_*\colon \pi_1(M_T,x^o)\to\pi_1(F\times I,x^o)).
\]
    Then group $\pi^o_0(T)$ is generated by elements $\epsilon^o(x)$, $x\in\widetilde{\mathdutchcal A}(T)$ (e.g., see the proof of~\cite[Proposition 3.2]{FR}).
    Then $\alpha^{-1}=\epsilon^o(x_1)\cdots\epsilon^o(x_k)$. Hence,
\begin{multline*}
\phi(\gamma_1)=\phi(\gamma_2\alpha\beta\alpha^{-1})=\phi(\gamma_2\alpha\beta\epsilon^o(x_1)\cdots\epsilon^o(x_k))=\\
\phi((\cdots(\gamma_2\alpha\beta)\ast x_1)\ast\cdots)\ast x_k)=    (\cdots(\phi(\gamma_2\alpha\beta)\ast \phi(x_1))\ast\cdots)\ast \phi(x_k)=\\
(\cdots(\phi(\gamma_2\alpha)\ast \phi(x_1))\ast\cdots)\ast \phi(x_k)=    \phi((\cdots(\gamma_2\alpha)\ast x_1)\ast\cdots)\ast x_k)=\\
\phi(\gamma_2\alpha\epsilon^o(x_1)\cdots\epsilon^o(x_k))=\phi(\gamma_2\alpha\alpha^{-1})=\phi(\gamma_2).
\end{multline*}
    The equality $\phi(\gamma_2\alpha\beta)=\phi(\gamma_2\alpha)$ follows from the $\pi_1(F)$-invariance of the map $\phi$.
\end{proof}

Theorem~\ref{thm:quandle_topological_fundamental} and Proposition~\ref{prop:quandle_topological_reduction} imply the following statement.

\begin{corollary}\label{cor:quandle_fundamental}
    The quandle $\bar{\mathdutchcal A}(T)$ is the fundamental quandle of the tangle $T$.
\end{corollary}

Recall the definition of stabilization map.

\begin{definition}\label{def:stabilization_map}
    Let $T\subset F\times I$ be a tangle in the thickening of a compact oriented surface $F$ and $D=p(T)$ the diagram of the tangle. Let $\delta$ be a closed $1$-submanifold such that $D\cap\delta=\emptyset$, and $D$ lies in a connected component $U$ of $F\setminus\delta$. Denote $F'=\overline{U}$. Then the pair $(F',T)$ is called a \emph{destabilization} of the pair $(F,T)$. The inverse operation is called stabilization.
\end{definition}

Recall that a \emph{virtual tangle} is an equivalence class of pairs $(F,T)$, where $F$ is a closed oriented surface and $T\subset F\times I$ is a tangle, modulo isotopies of the tangle, homeomorphisms, and (de)stabilizations.

\begin{proposition}\label{prop:quandle_stabilization}
Let $(F,T)\to(F',T)$ be a destabilization. Denote
\[
CF'(T)=(F'\times I\setminus N(T))/(F'\times 1)
\]
and
\[
\bar{\mathdutchcal A}'(T)=[I,0,1; CF'(T), \partial N(T),x^o].
\]
Then the map $i_*\colon\bar{\mathdutchcal A}'(T)\to\bar{\mathdutchcal A}(T)$ induced by the inclusion $i\colon CF'(T)\hookrightarrow CF(T)$ is an isomorphism of quandles.
\end{proposition}

\begin{proof}
    By definition of the quandle structures of $\bar{\mathdutchcal A}(T)$ and $\bar{\mathdutchcal A}'(T)$, the map $i_*$ is a quandle homomorphism.

    Let $\delta$ be the $1$-submanifold which yields the destabilization. Consider a collar neighborhood $U\subset F'$ of $\delta$ such that $D\cap U=\emptyset$ and a smooth function $h\colon F\to [0,1]$ such that $h|_{F'\setminus U}\equiv 0$ and $h|_{F\setminus F'}\equiv 1$. Consider the homotopy $H\colon (F\times I)\times I\to F\times I$ given by the formula
\[
H(x,t,s)=(x, t+(1-t)s\cdot h(x)),\quad x\in F, t,s\in I.
\]
Then $H$ defines a homotopy $\bar H$ in $CF$. Let $g=\bar H(\cdot,1)\colon CF\to CF$. Then $g(CF)=CF'$, $g\circ i=id_{CF}$, and $i\circ g\sim id_{CF}$. Hence, the induced map $g_*\colon\bar{\mathdutchcal A}(T)\to\bar{\mathdutchcal A}'(T)$ is inverse to $i_*$.
\end{proof}

\begin{corollary}\label{cor:quandle_virtual_invariant}
    The fundamental quandle $\bar{\mathdutchcal A}(T)$ (considered up to isomorphisms) is an invariant of virtual tangles.
\end{corollary}






\subsection{Colorings of regions: Partial tribracket}\label{subsect:partial_tribracket}

Consider coinvariants of regions that have the following property: the values of the invariant of three regions adjacent to a crossing determine the value of the invariant for the fourth region. This property leads to the following structure on the set of values of the invariant.

\begin{definition}[\cite{Niebrzydowski}]\label{def:tribracket}
A \emph{horizontal ternary quasigroup} is a pair of a set $X$ and a ternary operation $[]\colon X^3\to X$,  $(a, b, c)\mapsto [a, b, c]$ satisfying the following
property:
\begin{enumerate}
\item In the equation $[a,b,c]=d$, any three variables determine uniquely the fourth.
\item For any $a, b, c, d\in X$, it holds that
\[
[b, [a, b, c], [a, b, d]] = [c, [a, b, c], [a, c, d]] = [d, [a, b, d], [a, c, d]].
\]
\end{enumerate}
The operation $[ ]$ is called a \emph{horizontal tribracket}.
\end{definition}

For brevity, below we will use the word ``tribracket'' for horizontal ternary quasigroups.

\begin{definition}\label{def:tribracket_coloring}
Let $D$ be an oriented tangle diagram in an oriented connected compact surface $F$ and $(X,[])$ a tribracket. Then a {\em coloring of the diagram $D$ with the tribracket} is a map from the set of regions of $D$ to $X$ such that the images of the regions (colors) satisfy the coloring rule (Fig.~\ref{pic:tribracket_coloring}).
\begin{figure}[h]
\centering\includegraphics[width=0.4\textwidth]{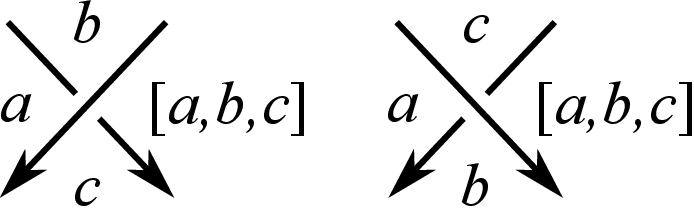}
\caption{The coloring rule}\label{pic:tribracket_coloring}
\end{figure}

Let $Col_{X}(D)$ be the set of tribracket colorings.
\end{definition}

\begin{figure}[h]
    \centering
    \includegraphics[width=0.8\textwidth]{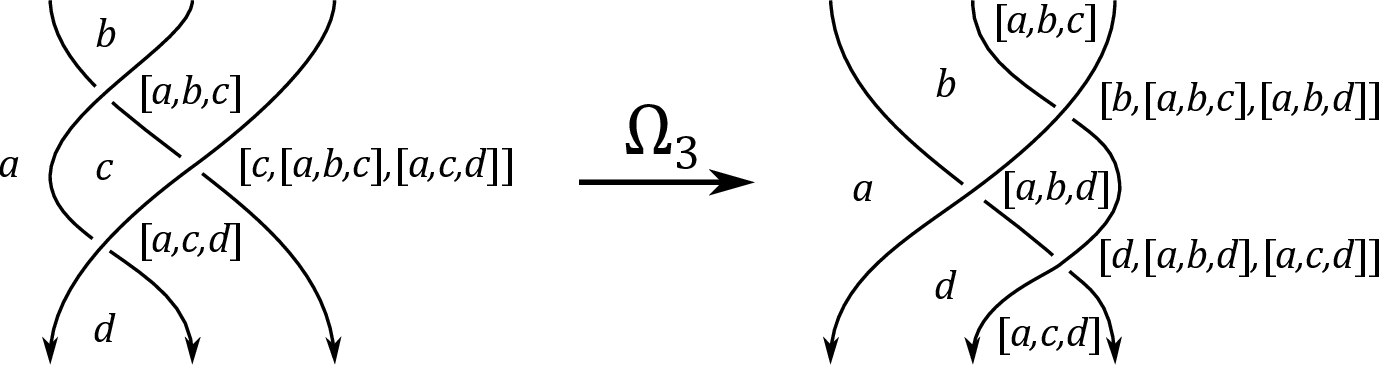}
    \caption{Third Reidemeister move}
    \label{pic:tribracket_R3}
\end{figure}

The second property in the definition of ternary quasigroups comes from the invariance under the third Reidemeister move (Fig.~\ref{pic:tribracket_R3}).

\begin{theorem}[\cite{Niebrzydowski}]\label{thm:tribracket_coloring_coinvariance}
For any diagrams $D$ and $D'$ connected by a Reidemeister move, there is a bijection between $Col_{X}(D)$ and $Col_{X}(D')$.
\end{theorem}

\begin{example}[Alexander tribracket]
Let $X$ be a module over the ring $\mathbb Z[t^{\pm 1},s^{\pm 1}]$. Then the operation
\[
[a, b, c] = tb + sc - tsa
\]
defines a tribracket operation called an \emph{Alexander tribracket}.
\end{example}

\begin{example}[Dehn tribracket]
Let $G$ be any group. Then the operation
\[
[a, b, c] = ca^{-1}b
\]
defines a tribracket structure called the \emph{Dehn tribracket} of $G$.
\end{example}

\begin{example}
Let $X=\{1,2,3\}$. It has a tribracket defined by the $3$-tensor

\[
\left[
\left[
\begin{array}{ccc}
1 & 3 & 2\\
2 & 1 & 3\\
3 & 2 & 1
\end{array}
\right],
\left[
\begin{array}{ccc}
2 & 1 & 3\\
3 & 2 & 1\\
1 & 3 & 2
\end{array}
\right],
\left[
\begin{array}{ccc}
3 & 2 & 1\\
1 & 3 & 2\\
2 & 1 & 3
\end{array}
\right]
\right].
\]

\end{example}

\begin{example}[Fundamental tribracket]
Let $D$ a diagram of a link $L$. The \emph{fundamental ternary quasigroup} $FT(D)$ of $L$ is determined by the universal rule: for any tribracket structure $(X,[])$ and any coloring $c\in Col_X(D)$ there exists a unique ternary quasigroup homomorphism $f\colon FT(D)\to X$ such that $c=f(c_F)$ where $c_F\in Col_{FT(D)}(D)$ is a fixed (fundamental) coloring. The fundamental ternary quasigroup has a presentation of the form
\[
FT(D)=\langle\mbox{regions of $D$}\,\mid\,\mbox{coloring rule at the crossings}\rangle.
\]
\end{example}

\begin{example}[Topological tribracket]\label{exa:topological_tribracket}
Let $F$ be a connected oriented compact surface. Let $T$ be an oriented tangle in the thickened surface $F\times I$, $M_T$ the complement to $T$, and $D=p(T)$ the diagram of $T$ in $F$. Choose $z\in F\setminus D$, denote $z_t=z\times t$, $t=0,1$. Consider the set $\widetilde{\mathdutchcal R}_0(T)$ of homotopy classes of paths
\[
\gamma\colon (I, 0, 1) \to (M_T, z_0, z_1)
\]
such that $p(\gamma)=1\in\pi_1(F,z)$. Then the operation
\[
[\alpha,\beta,\gamma]=\beta\alpha^{-1}\gamma
\]
defines a tribracket on $\widetilde{\mathdutchcal R}_0(T)$.

Note that the map $\gamma\mapsto\gamma^{-1}\gamma_0$ where $\gamma_0=z\times I$ defines an isomorphism from $\widetilde{\mathdutchcal R}_0(T)$ to the Dehn ternary quasigroup of the group 
\[
\pi^o(T)=ker(p_*\colon\pi_1(M_T,z_1)\to \pi_1(F,z)).
\]
\end{example}

\begin{remark}
The Dehn tribracket of any group $G$ possesses an involution-like property:
\[
[[a,b,c],c,b]=a\quad \forall a,b,c\in G.
\]
Consider a ternary quasigroup $(X,[])$ that contains elements $a,b,c$ which do not obey the condition above (we can take an Alexander tribracket). Let $D$ be the trivial plane diagram of the unlink $U_2$ with two components. Color the regions of the diagram $D$ with the colors $a,b,c$. Then there is no homomorphism from $\widetilde{\mathdutchcal R}_0(U_2)$ to $X$ that maps elements to $a,b,c$.

Thus, the topological tribracket is not fundamental.
\end{remark}

\begin{remark}
One explanation for why the topological tribracket is not fundamental is that the regions of diagrams on a surface cannot be arbitrarily combined, unlike the elements of the tribracket. For example, one cannot move the regions $a,b,c$ in Fig.~\ref{pic:tribracket_imcompatible_regions} into a position like in Fig.~\ref{pic:tribracket_coloring} (an obstruction is the checkerboard coloring of plane diagrams).

\begin{figure}
\centering\includegraphics[width=0.25\textwidth]{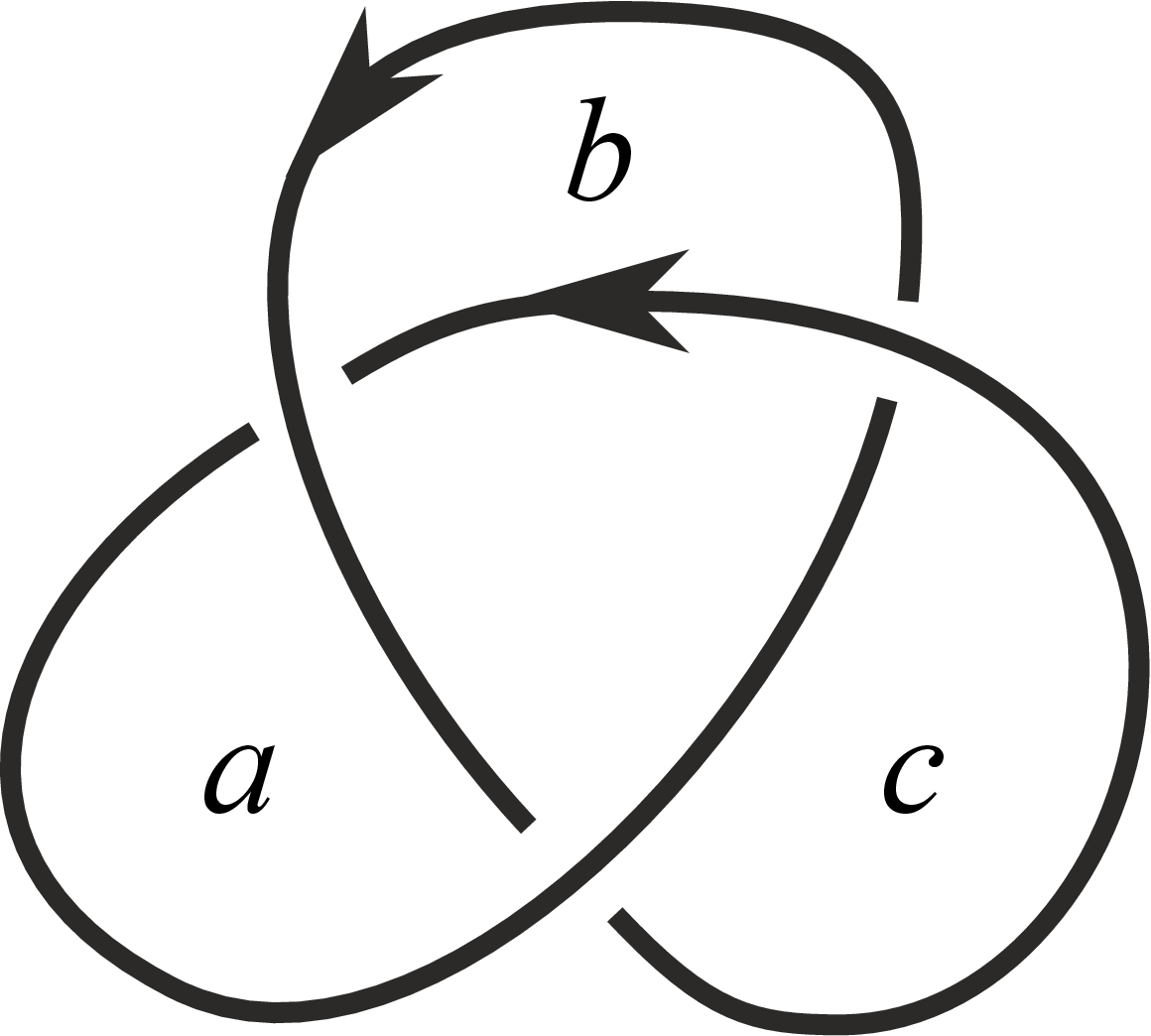}
\caption{The composition $[a,b,c]$ is not realizable on a classical diagram}\label{pic:tribracket_imcompatible_regions}
\end{figure}

In fact, regions $a,b,c$ can form a crossing only if $a$ and $b$ (and $a$ and $c$) are neighbor regions. We will use this compatibility condition to give a modified definition of tribrackets.
\end{remark}

\begin{definition}\label{def:tribracket_partial}
A \emph{partial ternary quasigroup} is a triple $(X,\uparrow, [])$ where $X$ is a set, $\uparrow$ is a relation on $X$ and $[]$ is a map from $X^3_\uparrow=\{(a,b,c)\in X^3 \mid b\uparrow a \mbox{ and } c\uparrow a\}$ to $X$ such that:
\begin{enumerate}
\item $[a,b,c]\uparrow b$ and $[a,b,c]\uparrow c$;
\item for any $a,b,d\in X$ such that $d\uparrow b\uparrow a$ there exists a unique $c$ such that $c\uparrow a$ and $[a,b,c]=d$;
\item for any $a,c,d\in X$ such that $d\uparrow c\uparrow a$ there exists a unique $b$ such that $b\uparrow a$ and $[a,b,c]=d$;
\item for any $b,c,d\in X$ such that $d\uparrow b$ and $d\uparrow c$ there exists a unique $a$ such that $b\uparrow a$, $c\uparrow a$ and $[a,b,c]=d$;
\item For any $a, b, c, d\in X$ such that $b\uparrow a$, $c\uparrow a$ and $d\uparrow a$, it holds that
\[
[b, [a, b, c], [a, b, d]] = [c, [a, b, c], [a, c, d]] = [d, [a, b, d], [a, c, d]].
\]
\end{enumerate}
\end{definition}

\begin{definition}\label{def:tribracket_partial_coloring}
Let $(X,\uparrow, [])$ be a partial ternary quasigroup. Then a {\em coloring of the diagram $D$ with $(B,R)$} is a map from the set of regions of $D$ to $X$ such that:
\begin{enumerate}
  \item for any two regions colored $a,b$ which are incident to one semiarc ($a$ from the left and $b$ from the right), $a\uparrow b$;
  \item the coloring rules (Fig.~\ref{pic:tribracket_coloring}) at the crossings hold.
\end{enumerate}
Denote the set of partial ternary quasigroup colorings by $Col_{(X,\uparrow)}(D)$.
\end{definition}

The following statement is proved analogously to Theorem~\ref{thm:tribracket_coloring_coinvariance}.

\begin{theorem}\label{thm:tribracket_partial_coloring}
For any diagrams $D$ and $D'$ connected by a Reidemeister move, there is a bijection between $Col_{(X,\uparrow)}(D)$ and $Col_{(X,\uparrow)}(D')$.
\end{theorem}

\begin{example}[Ternary quasigroup]
Let $(X,[])$ be a ternary quasigroup. Assume that $a\uparrow b$ for any $a,b\in X$. Then $(X,\uparrow, [])$ is a partial ternary quasigroup.
\end{example}

\begin{example}[Alexander numbering]
The triple $X=\mathbb Z$, $a\uparrow b$ iff $b=a-1$, and $[a,b,c]=a+2$ is a partial ternary quasigroup.
\end{example}

\begin{example}[Dehn partial tribracket]\label{exa:tribracket_partial_Dehn}
Let $G$ be a group, and $H\subset G$ be such that $g^{-1}Hg\subset H$ for any $g\in H$. Consider the relation $a\uparrow_H b$ iff $b^{-1}a\in H$. Then the Dehn tribracket operation $[a,b,c]=ba^{-1}c$ and the relation $\uparrow_H$ define a structure of partial ternary quasigroup on $G$.     
\end{example}

\begin{example}[Topological partial tribracket]\label{exa:tribracket_partial_topological}
Let $T$ be a tangle in $F\times I$, and $\widetilde{\mathdutchcal R}_0(T)$ the topological ternary quasigroup. For paths $\alpha,\beta\in\widetilde{\mathdutchcal R}_0$, set $\alpha\uparrow\beta$ iff there is a closed curve $\delta$ in $F$ such that the cycle 
$\beta^{-1}\alpha$ bounds a singular disk $\Delta\subset F\times I$ that intersects $T$ at one point which is a transversal positive intersection. Then $(\widetilde{\mathdutchcal R}_0(T),\uparrow, [])$ is a partial ternary quasigroup.
\end{example}

\begin{remark}\label{rem:tribracket_partial_topological}
 Let $H\subset\pi^o(T)$ be the set of homotopy classes of all meridians of the tangle $T$. Then the map 
 \[
 \widetilde{\mathdutchcal R}_0(T)\to\pi^o(T), \quad \gamma\mapsto\gamma^{-1}\gamma_0,
\]
defines an isomorphism from the topological partial ternary quasigroup to the Dehn partial ternary quasigroup of $(\pi^o(T),\uparrow_H)$.
\end{remark}

\begin{theorem}\label{thm:tribracket_partial_fundamental}
Let $T\subset F\times I$ be a tangle, and $D$ its diagram. For any partial ternary quasigroup $(X,\uparrow,[])$ there is a bijection between the set of colorings $Col_{(X,\uparrow)}(D)$ and the set of invariant homomorphisms $Hom(\widetilde{\mathdutchcal R}_0(T),(X,\uparrow))^{\pi_1(F)}$ under the action of $\pi_1(F,z)$ on $\widetilde{\mathdutchcal R}_0(T)$.
\end{theorem}

\begin{proof}
1. Let $\phi\colon \widetilde{\mathdutchcal R}_0(T)\to X$ be a $\pi_1(F,z)$-invariant homomorphism of partial ternary quasigroups. For a region $r\in\mathcal R(D)$, choose a point $x\in r$ and the vertical region probe $\gamma_r\subset x\times I$. Choose an arbitrary path $\delta_x\subset F$ from $x$ to $z$. Then $\gamma_x=(\delta_x\times 0)^{-1}\gamma_r(\delta_x\times 1)\in\widetilde{\mathdutchcal R}_0(T)$. Define the color $c_\phi(r)$ of the regon $r$ by formula $c_\phi(r)=\phi(\gamma_x)$. The element $c_\phi(r)$ does not depend on $x$ and $\delta_x$ by $\pi_1(F,z)$-invariance of $\phi$.

Let us check that the map $c_\phi\colon\mathcal R(D)\to X$ is a tribracket coloring. Let  $a,b,c$ be regions incident to a crossing of the diagram $D$ (Fig.~\ref{pic:tribracket_coloring}). Denote the fourth region by $d$. We need to show that $c_\phi(d)=[c_\phi(a),c_\phi(b),c_\phi(c)]$. Choose points in the regions close to the crossing, a path $\delta$ in $F$ from the crossing point to $z$. Using paths close to $\delta$, construct region probes $\gamma_a,\gamma_b,\gamma_c,\gamma_d\in\widetilde{\mathdutchcal R}_0(T)$. Then $\gamma_d=[\gamma_a,\gamma_b,\gamma_c]$ in $\widetilde{\mathdutchcal R}_0(T)$. Hence,
\[
c_\phi(d)=\phi(\gamma_d)=\phi([\gamma_a,\gamma_b,\gamma_c])=[\phi(\gamma_a),\phi(\gamma_b),\phi(\gamma_c)]
=[c_\phi(a),c_\phi(b),c_\phi(c)].
\]

2. Let $\xi\in Col_{(X,\uparrow)}(D)$. We construct a homomorphism $\phi=\phi_\xi$ from $\widetilde{\mathdutchcal R}_0(T)$ to $X$.

Let $D\cup\gamma$ a region probe diagram. Construct another region probe diagram $D'\cup\gamma$ by pulling the arcs of $D$ overcrossing $\gamma$ to the sinker point of $\gamma$, and pulling the arcs of $D$ undercrossing $\gamma$ to the float point of $\gamma$ (Fig.~\ref{pic:partial_tribracket_probe_clearing}).
The transformation of the diagram $D$ is a morphism $f\colon D\to D'$. This morphism induces a bijection $f_*\colon Col_{(X,\uparrow)}(D)\to Col_{(X,\uparrow)}(D')$. Then we set the values $\phi(\gamma)$ equal to the color $(f_*(\xi))(r')$ of the region $r'$ in $D'$ where the sinker point $z_0$ lies. (We can take the region of the float point $z_1$ for $r'$. The colors of these two regions coincide by Lemma~\ref{lem:tribracket_reciprocial_crossing}.)

\begin{figure}[h]
\centering\includegraphics[width=0.7\textwidth]{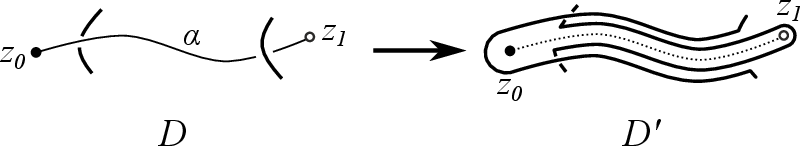}
\caption{Clearing a region probe}\label{pic:partial_tribracket_probe_clearing}
\end{figure}

Let us check that the value $\phi(\gamma)$ does not change during isotopy of $\gamma$. If $g\colon D\to D_1$ is a second or a third Reidemeister move, then the corresponding transformed diagrams $D'$ and $D'_1$ differ by a sequence of Reidemeister moves of the same type that do not involve the sinker point of $\gamma$. Then the color of the region $r'$ does not change during transformation from $D'$ to $D'_1$. If $g$ is a first Reidemeister move, then we can ignore the loop by Lemma~\ref{lem:tribracket_reciprocial_crossing} below.

A second or a third Reidemeister move $f\colon D\to D_1$ including an arc of the diagram $D$, induces diagrams $D'$ and $D'_1$ connected by an isotopy (and second and third Reidemeister moves if the probe has self-intersections), Fig.~\ref{pic:tribracket_mixed_R2_R3_enveloped}.

\begin{figure}[h]
    \centering
    \includegraphics[width=0.45\textwidth]{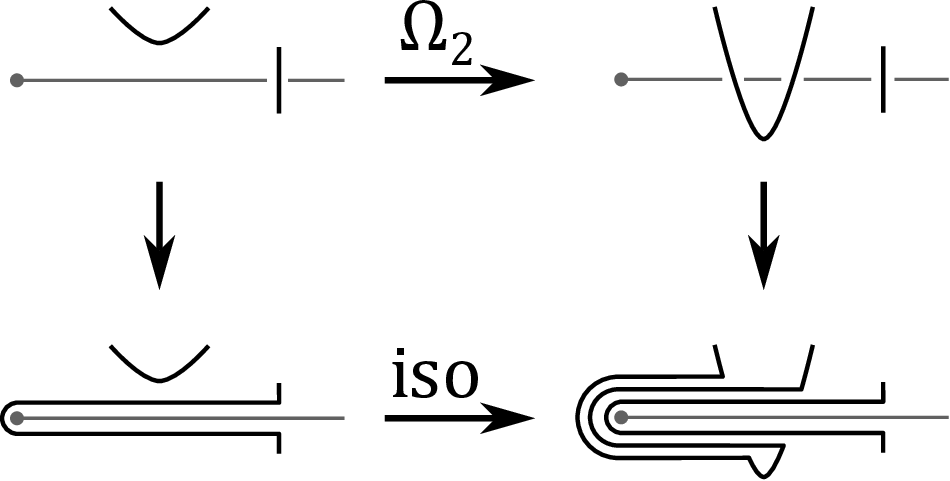}\quad \includegraphics[width=0.45\textwidth]{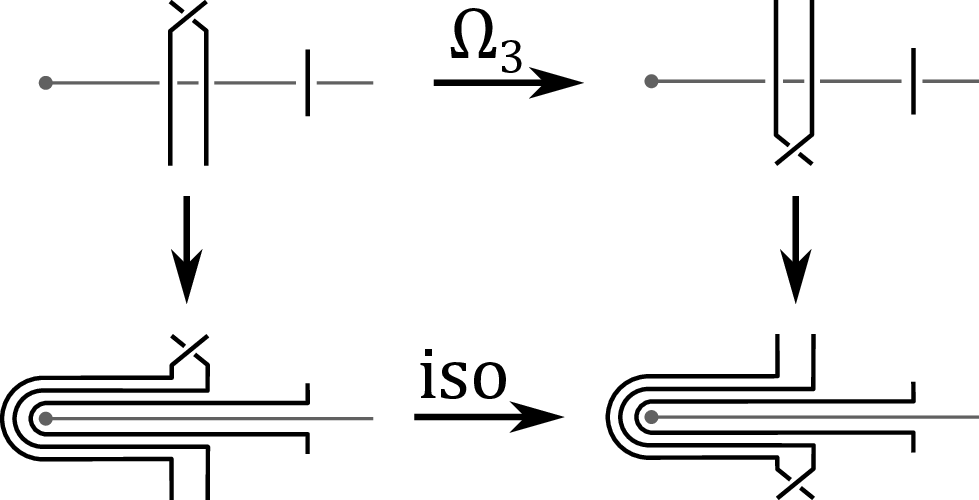}
    \caption{Second and third Reidemeister moves with an arc of the diagram}
    \label{pic:tribracket_mixed_R2_R3_enveloped}
\end{figure}

If $f\colon D\to D_1$ is a float or sinker move (Fig.~\ref{pic:graphoid_float_sinker_moves}), then the corresponding transformed diagrams $D'$ and $D'_1$ are isotopic.

Thus, we get a map $\phi\colon\mathscr R^s(T)\to X$ from the isotopy classes of region probes to the quandle.

Then let us show that self-intersections of region probes do not change their colors. We need the following lemma.

\begin{lemma}\label{lem:tribracket_reciprocial_crossing}
    Consider a configuration of regions colored by a partial quasigroup $(X,\uparrow)$ as shown in Fig.~\ref{pic:tribracket_reciprocial_crossing} left. Then after applying second Reidemeister moves, we get a colored part of the tangle (Fig.~\ref{pic:tribracket_reciprocial_crossing} right) such that the region colors $r$ and $r'$ coincide.
\begin{figure}[h]
    \centering
    \includegraphics[width=0.7\textwidth]{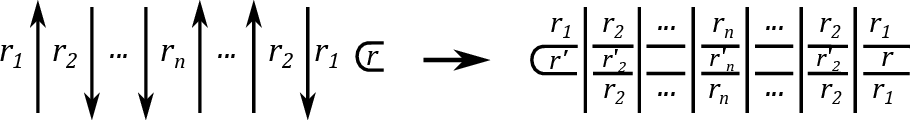}
    \caption{A colored tangle}
    \label{pic:tribracket_reciprocial_crossing}
\end{figure}
\end{lemma}
The proof is analogous to that of Lemma~\ref{lem:quandle_reciprocial_crossing}.

The lemma implies that if a region passes a self-intersection point of the probe during the pulling process, then its color does not change regardless of whether the region passes above or below the crossing. Hence, if one switches the undercrossing and the overcrossing of a self-intersection of the probe, then the result does not change. Thus, the map $\phi$ induces a map from $\mathdutchcal R(T)$ to $X$, therefore, induces a $\pi_1(F)$-invariant map from $\phi\colon\widetilde{\mathdutchcal R}_0(T)\to X$.

The next claim is that the map $\phi$ is a natural transformation from the functor $\widetilde{\mathdutchcal R}$ to the constant functor $X$. It is enough to show that for any region probe diagram $D\cup\gamma$ and for any Reidemeister move $f\colon D\to D_1$, the colors of the regions $\gamma$ and $\mathdutchcal R(f)(\gamma)$ coincide.

Since the map $\phi$ gives the same value for isotopic region probes, we can assume that $\gamma$ is distinct from the area where the move occurs. Then the transformed diagrams $D'$ and $D'_1$ are connected by the same move as $D$ and $D_1$, and this move does not involve the sinker region of the probe $\gamma$. Hence, the color of this region does not change with the move. Then the colors of the regions $\gamma$ and $\mathdutchcal A(f)(\gamma)$ coincide.

Finally, let us show that the map $\phi\colon\widetilde{\mathdutchcal R}_0(T)\to X$ is a homomorphism of partial quasigroups. Let $\gamma_1,\gamma_2,\gamma_3$ be region probes of the tangle $T$ such that $\gamma_2\uparrow\gamma_1$ and $\gamma_3\uparrow\gamma_1$. We need to prove that $\phi(\gamma_2)\uparrow\phi(\gamma_1)$, $\phi(\gamma_2)\uparrow\phi(\gamma_1)$, and $\phi([\gamma_1,\gamma_2,\gamma_3])=[\phi(\gamma_1),\phi(\gamma_2),\phi(\gamma_3)]$. Consider an isotopy $f\colon T\to T'$ which verticalize the probes and makes them form a crossing, i.e. $\mathdutchcal R(f)(\gamma_1)=\gamma_a$, $\mathdutchcal R(f)(\gamma_2)=\gamma_b$, $\mathdutchcal R(f)(\gamma_3)=\gamma_c$, where the regions $a,b,c\in\mathcal R(D')$, $D'=p(T')$, incident to a crossing as in Fig.~\ref{pic:tribracket_coloring}. Then
\[
\phi(\gamma_2)=\phi_{D'}(\gamma_b)=f_*(\xi)(b)\uparrow f_*(\xi)(a)=\phi_{D'}(\gamma_a)=\phi(\gamma_1),
\]
because $f_*(\xi)$ is a partial tribracket coloring.
Let $d$ be the fourth region incident to the crossing. Then
\begin{multline*}
[\phi(\gamma_1),\phi(\gamma_2),\phi(\gamma_3)]=[\phi_{D'}(\gamma_a),\phi_{D'}(\gamma_b),\phi_{D'}(\gamma_c)]=\\
[f_*(\xi)(a),f_*(\xi)(b),f_*(\xi)(c)]=f_*(\xi)(d)=
\phi_{D'}([\gamma_a,\gamma_b,\gamma_c])=\\
\phi_{D'}([\mathdutchcal R(f)(\gamma_1),\mathdutchcal R(f)(\gamma_2),\mathdutchcal R(f)(\gamma_3)])=\phi_{D'}(\mathdutchcal R(f)([\gamma_1,\gamma_2,\gamma_3]))=\phi([\gamma_1,\gamma_2,\gamma_3]),
\end{multline*}
where we use the fact that $f_*(\xi)$ is a partial tribracket coloring of $D'$, and that $\phi$ is a natural map.
The theorem is proved.
\end{proof}

The construction of homology of ternary quasigroups~\cite{Niebrzydowski_hom} can be extended to partial ternary quasigroups.

Let $(X,\uparrow, [])$ be a partial ternary quasigroup and $A$ an abelian group. Let $C_n(X,A)$ be the free $A$-module generated by
\[
X^{n+1}_\uparrow=\{(a_0,\dots,a_{n})\in X^{n+1}\mid a_{n}\uparrow a_{n-1}\uparrow\cdots\uparrow a_0\}.
\]
Define the differential by the formula
\begin{multline*}
  \partial_n(a_0,\dots,a_{n})=\sum_{i=1}^{n} (-1)^{i} [(a_0,\dots, a_{i-1},y_{(i,i)},\dots, y_{(i,n-1)})-\\(y_{(i,1)},y_{(i,2)}\dots,y_{(i,i-1)},a_{i},\dots,a_{n})],
\end{multline*}
where $[a_{j-1},a_j,y_{(i,j)}]=y_{(i,j+1)}$ for $j=1,\dots,i-2$, $[a_{i-2},a_{i-1},y_{(i,i-1)}]=a_{i}$, $[a_{i-1},a_{i},y_{(i,i)}]=a_{i+1}$, and
$[y_{(i,j-1)},a_{j},y_{(i,j)}]=a_{j+1}$ for $j=i+1,\dots,n-1$. For $n=1$, the differential is $\partial_1(a_0,a_1)=a_1-a_0$.
Then $(C_*(X,A),\partial_*)$ is a chain complex.
Let $D_n(X,A)$ the submodule generated by the elements
\[
\{(a_0,\dots,a_{n})\in C_n(X,A)\mid \exists i: [a_{i-1},a_i,a_i]=a_{i+1}\}.
\]



\begin{definition}\label{def:tribracket_partial_homology}
The homology $H^T(X,A)$ of the quotient complex $C_*(X,A)/D_*(X,A)$ is called \emph{partial tribracket homology}.
\end{definition}

\begin{definition}\label{def:tribracket_partial_cocycle_invariant}
Let $\theta\in Z^1(X,A)$ be a partial tribracket 1-cocycle. For a link diagram $D$, the \emph{cocycle invariant} of $D$ is the sum
\[
\Phi_\theta(D)=\sum_{c\in Col_{(X,R)}(D)}[\sum_{x\in\mathcal C(D)}\theta_x]\in\mathbb Z[A],
\]
where $\mathcal C(D)$ is the set of crossings and $\theta_x$ is the Boltzmann weight of the crossing $x$ defined as shown in Fig.~\ref{pic:tribracket_boltzmann_weight}.
\begin{figure}
\centering\includegraphics[width=0.3\textwidth]{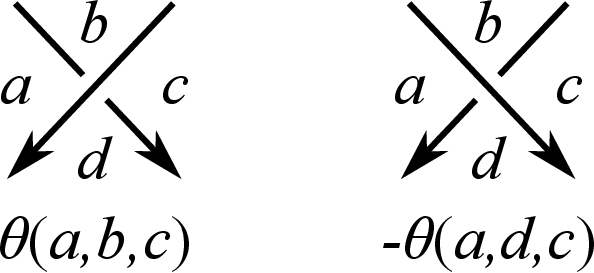}
\caption{The Boltzmann weights}\label{pic:tribracket_boltzmann_weight}
\end{figure}
\end{definition}

\begin{example}
Consider the ternary quasigroup given by the 3-tensor
\[
\left[
\left[\begin{array}{cc}
        1 & 2 \\
        2 & 1
      \end{array}
\right], \left[\begin{array}{cc}
        2 & 1 \\
        1 & 2
      \end{array}
\right]
\right]
\]
and the 1-cocycle $\theta$ such that $\theta(1,1,2)=\theta(2,2,1)$ and $\theta(a,b,c)=0$ otherwise.

For the Hopf link (Fig.~\ref{pic:hopf_link}), the cocycle invariant is equal to $\Phi_\theta(L)=4\cdot[0]+4\cdot[1]$.
\begin{figure}
\centering\includegraphics[width=0.25\textwidth]{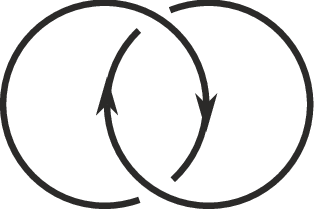}
\caption{The Hopf link}\label{pic:hopf_link}
\end{figure}
\end{example}

{\color{red}





}

Let $L\subset S^2\times I$ be a classical link. Let us see how strong the topological partial ternary quasigroup $\widetilde{\mathdutchcal R}_0(L)=\mathdutchcal R(L)$ is as a knot invariant. By Remark~\ref{rem:tribracket_partial_topological} the partial topological tribracket is determined by the knot group $\pi(L)$ and the set $H$ of (conjugacy classes of) meridians. We will show that the inverse is also true.

\begin{proposition}\label{prop:tribracket_partial_knot_group}
    Let $L$ and $L$ be classical links and $\mathdutchcal R(L)$, $\mathdutchcal R(L')$ their topological partial ternary quasigroups. If $\mathdutchcal R(L)$ and $\mathdutchcal R(L')$ are isomorphic then there is an isomorphism $\phi\colon\pi(L)\to\pi(L')$ of knot groups such that for any meridian $\mu\in\pi(L)$ of the link $L$ $\phi(\mu)$ is a meridian of the link $L'$. 
\end{proposition}

\begin{proof}
We will show that the link group $\pi(L)$ and the set $H\subset\pi(L)$ of meridians can be restored from $\mathdutchcal R(L)$.

Consider the set $\mathdutchcal R(L)^2_\uparrow=\{(a,b)\in\mathdutchcal R(L)\times\mathdutchcal R(L)\mid a\uparrow b\}$ and the equivalence relation on it generated by $(a,b)\sim ([b,c,a],c)$ for any $(a,b)\in \mathdutchcal R(L)^2_\uparrow$ and $c\uparrow b$. We claim that the equivalence relation has the properties:
\begin{enumerate}
    \item $(a,b)\sim(c,d) \Leftrightarrow b^{-1}\cdot a=d^{-1}\cdot c$;
    \item for any $(a,b)\in \mathdutchcal R(L)^2_\uparrow$ and $c\in\mathdutchcal R(L)$ there exists a unique $d$ such that $(a,b)\sim(c,d)$.
\end{enumerate}
The first property follows from
\[
c^{-1}\cdot[b,c,a]=c^{-1}\cdot c\cdot b^{-1}\cdot a=b^{-1}\cdot a.
\]
Let us prove the second property. Since the meridian set $H$ generates $\pi(L)$, there exists a sequence $a|^{\epsilon_1}a_1 |^{\epsilon_2}a_2 \cdots a_{n-1}|^{\epsilon_{n}}a_n=c$ where $\epsilon_i\in\{+,-\}$, $i=1,\dots,n$, and $|^+=\uparrow$, $|^-=\downarrow$. By induction, it is enough to consider the case $n=1$.

If $a\uparrow c$ then $(a,b)\sim(c,d)$ where $[d,b,c]=a$. If $c\uparrow a$ then $(a,b)\sim(c,d)$ where $[b,d,a]=c$.

Denote $\overline{\mathdutchcal R(L)}^2_\uparrow=\mathdutchcal R(L)^2_\uparrow/\sim$. Consider the group 
\[
G(L)=\langle (a,b)\in \overline{\mathdutchcal R(L)}^2_\uparrow\mid (b,a)\cdot(c,a)=(c,a)\cdot([a,b,c],c)\mbox{ for all } b\uparrow a, c\uparrow a\rangle.
\]
Since the map $(a,b)\mapsto b^{-1}\cdot a$ establishes a bijection between $\overline{\mathdutchcal R(L)}^2_\uparrow$ and $H$, the presentation of the group $G(L)$ is the presentation of the associated group of the link quandle which is the link group $\pi(L)$. Hence, the natural homomorphism $G(L)\to\pi(L)$ is a group isomorphism which maps the generator set $\overline{\mathdutchcal R(L)}^2_\uparrow$ onto $H$.

Thus, any partial ternary quasigroup isomorphism $\phi\colon\mathdutchcal R(L)\to\mathdutchcal R(L')$ induces a group isomorphism $\phi\colon\pi(L)=G(L)\to G(L')=\pi(L')$ which identifies the meridian sets $H$ and $H'$. 
\end{proof}

Proposition~\ref{prop:tribracket_partial_knot_group} and~\cite[Theorem 15.38]{BZ} imply the following statement.

\begin{corollary}\label{cor:tribracket_topological_isomorphism}
    Let $K,K'\subset S^2\times I$ be oriented classical knots, and $K=K_1\# \cdots \# K_n$ the prime decomposition of $K$. If the topological partial ternary quasigroups $\mathdutchcal R(K)$ and $\mathdutchcal R(K')$ are isomorphic, then $K'=K'_1\# \cdots \# K'_n$ where either $K'_i=K_i$ or $K'_i=-\bar K_i$ (the inverse mirror knot), $i=1,\dots,n$.
\end{corollary}

\subsection{Colorings of semiarcs: Biquandloid}\label{subsect:biquandloid}

Consider coinvariants of semiarcs that have the following property: the values of the invariant of two adjacent semiarcs at a crossing determine the value of the two opposite semiarcs. This property leads to the biquandle structure~\cite{EN} on the set of values of the invariant.

\begin{definition}\label{def:biquandle}
A set $B$ with two binary operations $\circ,\ast\colon B\times B\to B$ is called a \emph{biquandle} if it obeys the following conditions:
\begin{enumerate}
\item $x\circ x=x\ast x$ for any $x\in B$;
\item for any $y\in B$ the operators $\alpha_y\colon B\to B, x\mapsto x\circ y$, and $\beta_y\colon B\to B, x\mapsto x\ast y$, are invertible;
\item the map $S\colon B\times B\to B\times B$, $(x,y)\mapsto (x\circ y, y\ast x)$, is a bijection;
\item for any $x,y,z\in B$
\begin{gather*}
(x\circ z)\circ(y\circ z)=(x\circ y)\circ(z\ast y),\\
(x\circ z)\ast(y\circ z)=(x\ast y)\circ(z\ast y),\\
(x\ast z)\ast(y\ast z)=(x\ast y)\ast(z\circ y).
\end{gather*}
\end{enumerate}
\end{definition}

\begin{definition}\label{def:biquandle_coloring}
Let $B$ be a biquandle. Then a {\em coloring of the diagram $D$ with the biquandle $B$} is a map from the set of semi-arcs of $D$ (edges of the graph $D$) to $B$ such that the images of the arcs (colors) satisfy the coloring rule (Fig.~\ref{pic:biquandle_coloring}).
\begin{figure}[h]
\centering\includegraphics[width=0.4\textwidth]{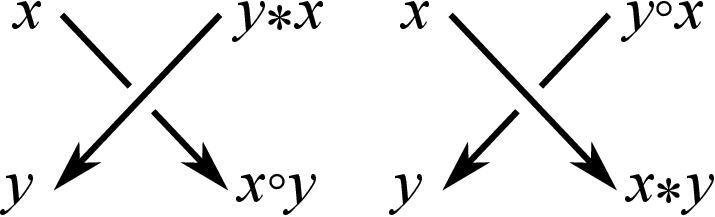}
\caption{The biquandle coloring rule}\label{pic:biquandle_coloring}
\end{figure}
Let $Col_{B}(D)$ denote the set of biquandle colorings.
\end{definition}

\begin{theorem}[\cite{EN}]\label{thm:biquandle_coloring_coinvariance}
For any diagrams $D$ and $D'$ connected by a Reidemeister move, there is a bijection between $Col_{B}(D)$ and $Col_{B}(D')$.
\end{theorem}

\begin{example}[Quandle]
If for any $x,y\in B$ one has $x\ast y=x$ then the biquandle $B$ is called a \emph{quandle}.
\end{example}

\begin{example}[Fundamental biquandle]
Let $D$ a diagram of a link $L$. The \emph{fundamental biquandle} $FB(D)$ of $L$ is determined by the universal rule: for any biquandle $B$ and any coloring $c\in Col_B(D)$ there exists a unique biquandle homomorphism $f\colon FB(D)\to B$ such that $c=f(c_F)$ where $c_F\in Col_{FB(D)}(D)$ is a fixed (fundamental) coloring. The fundamental biquandle has a presentation of the form
\[
FB(D)=\langle\mbox{semiarcs of $D$}\,\mid\,\mbox{coloring rule at the crossings}\rangle.
\]
\end{example}

\begin{example}[Topological biquandle~\cite{Horvat}]\label{exa:biquandle_topological}
Let $L\subset S^2\times I$ be a classical link, $N(L)$ a tubular neighborhood of $L$, and $E_L=\overline{(S^2\times I)\setminus N(L)}$. Choose points $z_0\in S^2\times 0$ and $z_1\in S^2\times 1$. Consider the set $B_L$ of the homotopy classes of pairs of paths $(a_0,a_1)$ such that
\[
a_0\colon (I,0,1)\to (E_L,\partial N(L),z_0),\qquad a_1\colon (I,0,1)\to (E_L,\partial N(L),z_1),
\]
and $a_0(0)=a_1(0)$. Define a biquandle structure by the formulas
\begin{gather*}
(a_0,a_1)\circ(b_0,b_1)=(a_0, a_1b_1^{-1}m_{b_1(0)}b_1),\\
 (a_0,a_1)\ast(b_0,b_1)=(a_0b_0^{-1}m_{b_0(0)}^{-1}b_0,a_1)
\end{gather*}
where $m_x$, $x\in\partial N(L)$, is a meridian of $\partial N(L)$ which passes over the point $x$.

The biquandle $B_L$ is called the \emph{topological biquandle} of the link $L$.
\end{example}

Unlike the topological quandle, the topological biquandle is not fundamental, as the following example shows.

\begin{example}[\cite{IT}]
Let $U$ be the unknot. Then $B_U$ is isomorphic to the biquandle $\mathbb Z$ with the operations
\[
x\circ y = x\ast y=x+1.
\]

On the other hand, $FB(U)$ is the free biquandle with 1 generator. Hence, there is a biquandle homomorphism of $FB(U)$ onto the biquandle $B=\mathbb Z^2$ with operations
\[
(x,a)\circ(y,b) = (x,a)\ast(y,b)=(x+1,a+y),
\]
which is generated by the element $(0,0)$.

But in $B$ the result of the operations depends on the second argument.
\end{example}

\begin{remark}
The topological biquandle is not fundamental because semiarcs of diagrams on a surface cannot be arbitrarily combined, unlike the elements of the biquandle. For example, the arcs $x$ and $y$ in Fig.~\ref{pic:tribracket_imcompatible_regions} cannot be positioned as shown in Fig.~\ref{pic:tribracket_coloring} (an obstruction is the checkerboard coloring of plane diagrams).
\begin{figure}[h]
\centering\includegraphics[width=0.25\textwidth]{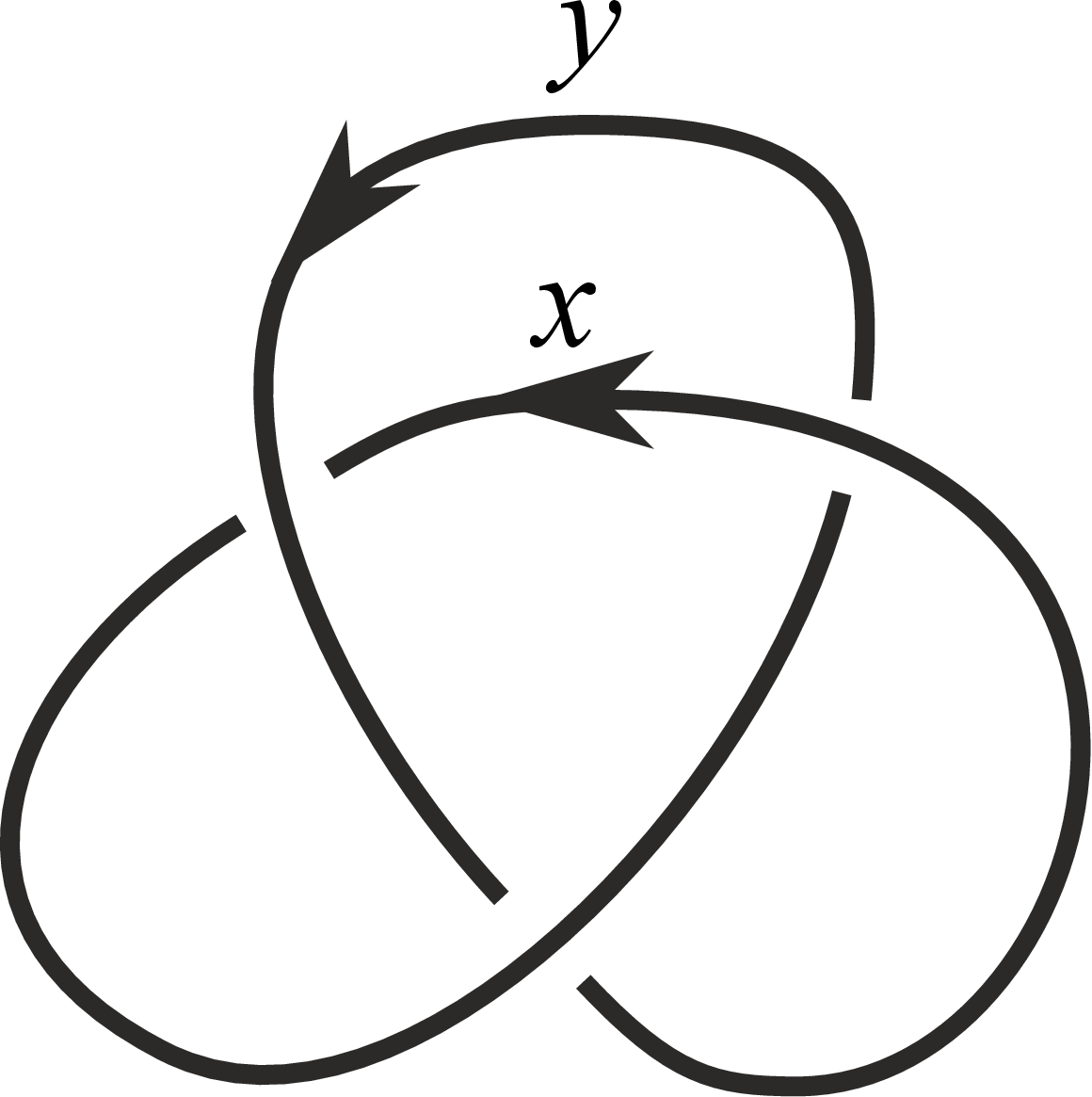}
\caption{The composition $x\ast y$ is not realizable in a classical diagram}\label{pic:trefoil_marked}
\end{figure}

Thus, the compatibility condition should be taken into account, which can be formulated as follows: two semiarcs can form a crossing if and only if they are incident (on the left or the right) to the same region.

Thereby, the notion of a biquandle is not well-suited for describing colorings of semiarcs in diagrams on a fixed surface, particularly for classical diagrams. However, biquandles are appropriate for colorings of virtual knots, as any two arcs can be paired in a virtual diagram.
\end{remark}

\begin{definition}\label{def:biquandloid}
A \emph{biquandloid} is a pair of sets $(B,R)$ with maps $\sigma_{l},\sigma_r\colon B\to R$ (called \emph{left} and \emph{right shadow maps}) and operations $\uast,\oast\colon B\times_R B\to B$ where $B\times_R B=\{(a,b)\in B\times B\mid\sigma_r(a)=\sigma_r(b)\}$ such that
\begin{enumerate}
\item for any $a,b\in B$ such that $\sigma_r(a)=\sigma_r(b)$ one has $\sigma_r(a\uast b)=\sigma_r(a\oast b)=\sigma_l(b)$ and $\sigma_l(a\uast b)=\sigma_l(b\oast a)$.
\item for any $b,c\in B$ such that $\sigma_l(b)=\sigma_r(c)$ there exist unique $a', a''\in B$ such that $\sigma_r(a')=\sigma_r(a'')=\sigma_r(b)$ and $c = a'\uast b = a''\oast b$.
\item for any $c,d\in B$ such that $\sigma_l(c)=\sigma_l(d)$ there exist unique $a, b\in B$ such that $\sigma_r(a)=\sigma_r(b)$ and $c = a\oast b$, $d = b\uast a$.

\item for any $a\in B$  $a\uast a=a\oast a$.
\item for any $a,b,c\in B$ such that $\sigma_r(a)=\sigma_r(b)=\sigma_r(c)$
\begin{gather*}
 (a\uast b)\uast(c\uast b)=(a\uast c)\uast(b\oast c),\\
 (a\uast b)\oast(c\uast b)=(a\oast c)\uast(b\uast c),\\
(a\oast b)\oast(c\oast b)=(a\oast c)\oast(b\uast c).
\end{gather*}
\end{enumerate}

\end{definition}

\begin{definition}\label{def:biquandloid_coloring}
Let $(B,R)$ be a biquandloid. Then a {\em coloring of the diagram $D$ with the biquandloid $(B,R)$} is a map from the set of semi-arcs of $D$ to $B$ such that:
\begin{enumerate}
  \item for any two arcs with colors $a,b$ incident to one region, $\sigma(a)=\sigma(b)$ where $\sigma(x)=\sigma_r(x)$, if the arc $x$ has positive orientation with respect to the region, and $\sigma(x)=\sigma_l(x)$ if the orientation is negative, $x=a,b$;
  \item the coloring rules (Fig.~\ref{pic:biquandle_coloring}) at the crossings hold.
\end{enumerate}
Denote the set of biquandloid colorings by $Col_{(B,R)}(D)$.
\end{definition}

\begin{remark}
For classical diagrams, the condition 2) implies 1).

\end{remark}

\begin{theorem}\label{thm:biquandloid_coloring_coinvariance}
For any diagrams $D$ and $D'$ connected by a Reidemeister move, there is a bijection between $Col_{(B,R)}(D)$ and $Col_{(B,R)}(D')$.
\end{theorem}

\begin{example}[Biquandle]
Let $R$ is a one-point set. Then all compatibility conditions are trivial, and the definition of biquandloid reduces to the definition of biquandle.
\end{example}

\begin{example}[Local biquandle~\cite{NOO}]

A {\it local biquandle} is defined as a triple $(X,  \{\uast_a\}_{a\in X}, \{\oast_a\}_{a\in X})$ of a set $X$ and two families of operations
\[
\uast_a, \oast_a: (\{a\} \times X)^2 \to X^2
\]
satisfying the following property:
\begin{itemize}
\item[] \hspace{-0.5cm}($\mathcal{L}$1) For any $a,b,c \in X$,
\begin{itemize}
\item[(i)] the first component of the result of $(a, b) \uast_a (a,c)$ is $c$,
\item[(ii)] the first component of the result of $(a, b) \oast_a (a,c)$ is $c$,
\item[(iii)] the second component of the result of $(a, b) \uast_a (a,c)$ coincides with that of the result of  $(a, c) \oast_a (a,b)$.
\end{itemize}
\item[] \hspace{-0.5cm}($\mathcal{L}$2)
\begin{itemize}
\item[(i)] For any $a, b\in  X$, the map $\uast_a(a, b) :\{a\}\times X \to \{b\}\times X$ sending $(a,c)$ to $(a,c)\uast_a(a,b)$ is bijective.
\item[(ii)]
For any $a, b\in X$, the map $\oast_a(a,b):\{a\}\times X \to \{b\}\times  X$ sending $(a,c)$ to $(a,c)\oast_a(a,b)$ is bijective.
\item[(iii)]
The map $S: \bigcup_{a\in X}(\{a\} \times X)^2\to \bigcup_{d\in X} (X\times \{d\})^2$ defined by $S\big((a,b),(a,c)\big)=\big((a,c)\oast_a(a,b),(a,b)\uast_a(a,c) \big)$ is bijective.
\end{itemize}
\item[] \hspace{-0.5cm}($\mathcal{L}$3)
For any $a,b,c\in X$, it holds that
\begin{itemize}
\item[(i)] $\big((a,b)\uast_a(a,c)\big)\uast_c\big((a,d)\oast_a(a,c)\big)=\big((a,b)\uast_a(a,d)\big)\uast_d\big((a,c)\uast_a(a,d)\big)$,
\item[(ii)]
$\big((a,b)\oast_a(a,c)\big)\uast_c\big((a,d)\oast_a(a,c)\big)=\big( (a,b)\uast_a(a,d)\big)\oast_d\big((a,c)\uast_a(a,d)\big)$,
\item[(iii)] $\big((a,b)\oast_a(a,c)\big)\oast_c\big((a,d)\uast_a(a,c)\big)=\big((a,b)\oast_a(a,d)\big)\oast_d\big((a,c)\oast_a(a,d)\big)$.
\end{itemize}
\end{itemize}

For simplicity, we often omit the subscript by $a$ as $\uast=\uast_a$, $\oast=\oast_a$,  $\{\uast\}=\{\uast_a\}_{a\in X}$ and $\{\oast\}=\{\oast_a\}_{a\in X}$ unless it causes confusion.

Given a local biquandle $(X,  \{\uast_a\}_{a\in X}, \{\oast_a\}_{a\in X})$, we set $B=X\times X$ and $R=X$ with the shadow maps defined by the formulas
\[
\sigma_r(a,b)=a,\quad \sigma_l(a,b)=b.
\]
\end{example}

\begin{remark}
There is a bijection between local biquandles and tribrackets of M. Niebrzydowski~\cite{Niebrzydowski}.
\end{remark}

\begin{example}\label{exa:biquandloid_tribracket}
Let $(X,[],\uparrow)$ be a partial ternary quasigroup. Consider the sets $B=\{(x,y)\in X^2\mid x\uparrow y\}$ and $R=X$. The maps
\begin{gather*}
    \sigma_l(x,y)=x,\quad \sigma_r(x,y)=y,\\
    (x,y)\uast(z,y)=([y,x,z],z),\quad (x,y)\oast(z,y)=([y,z,x],z)
\end{gather*}
establish a structure of a biquandloid on $(B,R)$. This biquandloid is called \emph{tribracket biquandloid}.   
\end{example}

\begin{example}[Biquandle double]
Let $(B,\circ,\ast)$ be a biquandle. Consider the sets $B_2=B\times\mathbb Z_2$ and $R=\mathbb Z_2$ with the operations
\begin{gather*}
  \sigma_r(x,\epsilon)=\epsilon,\quad \sigma_l(x,\epsilon)=\epsilon+1\\
  (x,\epsilon)\circ(y,\epsilon)=(x\circ y ,\epsilon+1),\quad (x,\epsilon)\ast(y,\epsilon)=(x\ast y ,\epsilon+1).
\end{gather*}
Then $(B_2,R)$ is a biquandloid which is not a local biquandle.

The set of colorings is $Col_{(B_2,R)}(L)=Col_B(L)\times\mathbb Z_2$ if the link $L$ has a checkerboard colored diagram (all classical links have one), and $Col_{(B_2,R)}(L)=\emptyset$ otherwise.
\end{example}

\begin{example}[Shadow biquandloid]
Let $(B,\circ,\ast)$ be a biquandle and $X$ is a set with an action $\triangleleft\colon X\times B\to X$ such that
\begin{enumerate}
  \item for any $a\in B$ the map $\triangleleft a\colon X\to X$ is bijective;
  \item for any $a,b\in B$ and $x\in X$ $(x\triangleleft a)\triangleleft(b\ast a)=(x\triangleleft b)\triangleleft(a\circ b)$.
\end{enumerate}

Consider the set $B_X=X\times B$ with the maps $\sigma_r,\sigma_l\colon B_X\to X$ and $\circ,\ast\colon B_X\times_X B_X\to B_X$ given by the formulas
\begin{gather*}
  \sigma_r(x,a)=x,\quad \sigma_l(x,a)=x\triangleleft a, \\
  (x,a)\circ(x,b)=(x\triangleleft b, a\circ b),\quad (x,a)\ast(x,b)=(x\triangleleft b, a\ast b).
\end{gather*}
Then $(B_X,X)$ is a biquandloid.
\end{example}

\begin{example}[Topological biquandloid]\label{exa:fundamental_biquandloid}
Let $T$ be a tangle in $F\times I$, $\widetilde{\mathdutchcal{SA}}_0(T)$ the set of based homotopy classes of semiarc probes, and $\widetilde{\mathdutchcal R}_0(T)$ the set of based homotopy classes of region probes. Consider the maps $\sigma_r, \sigma_l$ and $\ast,\circ$ defined by the formulas
\begin{gather*}
\sigma_r(\gamma^u,\gamma^o)=(\gamma^u)^{-1}\mu_r^{-1}\gamma^o,\quad \sigma_l(\gamma^u,\gamma^o)=(\gamma^u)^{-1}\mu_l\gamma^o,\\
(\gamma_1^u,\gamma_1^o)\circ(\gamma_2^u,\gamma_2^o)=(\gamma_1^u,\gamma_1^o(\gamma_2^o)^{-1}(\mu_2)_r(\mu_2)_l\gamma_2^o),\\
(\gamma_1^u,\gamma_1^o)\ast(\gamma_2^u,\gamma_2^o)=(\gamma_1^u(\gamma_2^u)^{-1}(\mu_2)_r^{-1}(\mu_2)_l^{-1}\gamma_2^u,\gamma_1^o).
\end{gather*}
and restrict the biquandle operations to $\widetilde{\mathdutchcal{SA}}_0(T)\times_{\widetilde{\mathdutchcal{R}}_0(T)}\widetilde{\mathdutchcal{SA}}_0(T)$. The result is called \emph{topological biquandloid} of the link $L$.
\end{example}

\begin{theorem}\label{thm:topological_biquandloid_universal}
Let $T\subset F\times I$ be a tangle, and $D$ its diagram. For any biquandloid $(B,R)$ there is a bijection between the set of colorings $Col_{(B,R)}(D)$ and the set of invariant homomorphisms $Hom((\widetilde{\mathdutchcal{SA}}_0(T),\widetilde{\mathdutchcal{R}}_0(T)),(B,R))^{\pi_1(F)}$ under the action of $\pi_1(F,x_0)$ on $\widetilde{\mathdutchcal{SA}}_0(T)$.
%
\end{theorem}

\begin{proof}
1. Let $\phi\colon \widetilde{\mathdutchcal{SA}}_0(T)\to X$ be a $\pi_1(F,x_0)$-invariant homomorphism of biquandloids. For a semiarc $a\in\mathcal{SA}(D)$, choose a point $x\in a$ and the vertical semiarc probe $\gamma_a=(\gamma_a^u,\gamma_a^o)\subset x\times I$. Choose an arbitrary path $\delta_x\subset F$ from $x$ to $x_0$. Then $\gamma_x=(\gamma_a^u(\delta_x\times 0),\gamma_a^o(\delta_x\times 1))\in\widetilde{\mathdutchcal{SA}}_0(T)$. Define the color $c_\phi(a)$ of the semiarc $a$ by formula $c_\phi(a)=\phi(\gamma_x)$. The element $c_\phi(a)$ does not depend on $x$ and $\delta_x$ by $\pi_1(F,x_0)$-invariance of $\phi$.

Let us check that the map $c_\phi\colon\mathcal{SA}(D)\to X$ is a biquandloid coloring. Let $x,y$ be semiarcs incident to a crossing of the diagram $D$ (Fig.~\ref{pic:biquandle_coloring}). Denote the opposite semiarcs to $x,y$ by $z$ and $w$. Denote $a=c_\phi(x)$, $b=c_\phi(y)$, $c=c_\phi(z)$, and $d=c_\phi(w)$. We need to show that $\sigma_l(a)=\sigma_r(d)$, $\sigma_l(b)=\sigma_r(c)$, $\sigma_l(c)=\sigma_l(d)$, $c=a\circ b$, and $d=b\ast a$. Choose points in the regions close to the crossing, a path $\delta$ in $F$ from the crossing point to $x_0$. Using paths close to $\delta$, construct region probes $\gamma_x,\gamma_y,\gamma_z,\gamma_w\in\widetilde{\mathdutchcal{SA}}_0(T)$. Then $\sigma_l(\gamma_x)=\sigma_r(\gamma_w)$, $\sigma_l(\gamma_y)=\sigma_r(\gamma_z)$, $\sigma_l(\gamma_z)=\sigma_l(\gamma_w)$, $\gamma_z=\gamma_x\circ \gamma_y$, and $\gamma_w=\gamma_y\ast\gamma_x$. in $\widetilde{\mathdutchcal{SA}}_0(T)$. Hence,
\begin{gather*}
\sigma_l(a)=\sigma_l(\phi(\gamma_x))=\phi(\sigma_l(\gamma_x))=\phi(\sigma_r(\gamma_w))=\sigma_r(\phi(\gamma_w))=\sigma_r(d),\\
c=\phi(\gamma_z)=\phi(\gamma_x\circ\gamma_y)=\phi(\gamma_x)\circ\phi(\gamma_y)=a\circ b.
\end{gather*}
The other equalities are checked analogously.


2. Let $\xi\in Col_{(B,R)}(D)$. We construct a biquandloid homomorphism $\phi=\phi_\xi$ from $\widetilde{\mathdutchcal{SA}}_0(T)$ to $X$.

Let $D\cup\gamma$ a semiarc probe diagram. Construct another semiarc probe diagram $D'\cup\gamma$ by pulling the arcs of $D$ overcrossing $\gamma$ to the sinker point of $\gamma$, and pulling the arcs of $D$ undercrossing $\gamma$ to the float point of $\gamma$ (Fig.~\ref{pic:biquandloid_probe_clearing}). 
\begin{figure}[h]
\centering\includegraphics[width=0.8\textwidth]{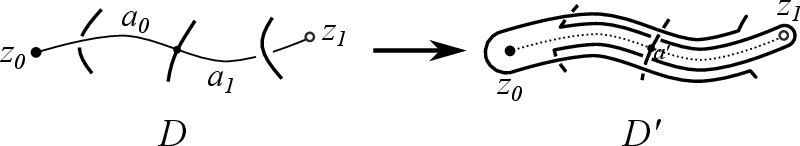}
\caption{Clearing a semiarc probe}\label{pic:biquandloid_probe_clearing}
\end{figure}

The transformation of the diagram $D$ is a morphism $f\colon D\to D'$. This morphism induces a bijection $f_*\colon Col_{(B,R)}(D)\to Col_{(B,R)}(D')$. Then we set the values $\phi(\gamma)$ equal to the color $(f_*(\xi))(a')$ of the semiarc $a'$ in $D'$ where the probe intersects the diagram $D'$.

Let us check that the value $\phi(\gamma)$ does not change during isotopy of $\gamma$. If $f\colon D\to D_1$ is a second or a third Reidemeister move, then the corresponding transformed diagrams $D'$ and $D'_1$ differ by a sequence of Reidemeister moves of the same type that do not involve the vertex of $\gamma$. Then the color of the semiarc $a'$ does not change during transformation from $D'$ to $D'_1$. If $f$ is a first Reidemeister move, then we can ignore the loop by Lemma~\ref{lem:biquandloid_reciprocial_crossing} below.

A second or a third Reidemeister move $f\colon D\to D_1$ including an arc of the diagram $D$, induces diagrams $D'$ and $D'_1$ connected by an isotopy (and second and third Reidemeister moves if the probe has self-intersections).

If $f\colon D\to D_1$ is a float or sinker move (Fig.~\ref{pic:graphoid_float_sinker_moves}), then the corresponding transformed diagrams $D'$ and $D'_1$ are isotopic.

If $f\colon D\to D_1$ is a vertex rotation move, then the corresponding transformed diagrams $D'$ and $D'_1$ are connected by isotopy, first and second Reidemeister moves which do not involve the vertex of the semiarc probe (Fig.~\ref{pic:biquandloid_rotation_enveloped}).
\begin{figure}[h]
    \centering
    \includegraphics[width=0.9\textwidth]{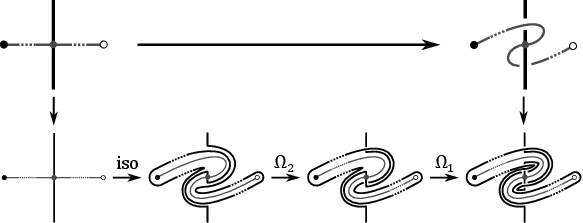}
    \caption{Rotation move}
    \label{pic:biquandloid_rotation_enveloped}
\end{figure}

Thus, we get a map $\phi\colon\mathscr{SA}^s(T)\to X$ from the isotopy classes of semiarc probes to the quandle.

Then let us show that self-intersections of semiarc probes do not change their colors. We need the following lemma.

\begin{lemma}\label{lem:biquandloid_reciprocial_crossing}
 Consider a part of the tangle that consists of $2n$ parallel arcs $a_1,\dots, a_{2n}$ and another arc $b$ such that the arcs $a_i$ and $a_{2n+1-i}$ have opposite orientation (Fig.~\ref{pic:biquandloid_reciprocial_crossing} left). Assume that the arcs are colored by the biquandloid $(B,R)$ so that the colors of the arcs $a_i$ and $a_{2n+1-i}$ coincide. Then after applying second Reidemiester moves, we get a colored part of the tangle (Fig.~\ref{pic:biquandloid_reciprocial_crossing} right) such that the arcs $b'$ and $b$ have the same color, as well as the arcs $a_i$ and $a'_i$, $i=1,\dots,2n$, as well as the arcs $b_i$ and $b'_i$, $i=1,\dots,2n-1$.
\begin{figure}[h]
    \centering
    \includegraphics[width=0.7\textwidth]{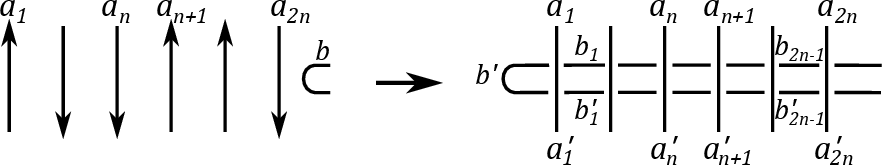}
    \caption{A colored tangle}
    \label{pic:biquandloid_reciprocial_crossing}
\end{figure}
\end{lemma}
The proof is analogous to that of Lemma~\ref{lem:quandle_reciprocial_crossing}.

The lemma implies that if a semiarc probe passes a self-intersection point during the pulling process, then its color does not change regardless of whether the semiarc passes above or below the crossing. Hence, if one switches the undercrossing and the overcrossing of a self-intersection of the probe, then the result does not change. Thus, the map $\phi$ induces a map from $\mathdutchcal{SA}(T)$ to $X$, therefore, induces a $\pi_1(F)$-invariant map from $\phi\colon\widetilde{\mathdutchcal{SA}}_0(T)\to X$.

The next claim is that the map $\phi$ is a natural transformation from the functor $\widetilde{\mathdutchcal{SA}}$ to the constant functor $B$. It is enough to show that for any semiarc probe diagram $D\cup\gamma$ and for any Reidemeister move $f\colon D\to D_1$, the colors of the semiarcs $\gamma$ and $\mathdutchcal{SA}(f)(\gamma)$ coincide.

Since the map $\phi$ gives the same value for isotopic semiarc probes, we can assume that $\gamma$ is distinct from the area where the move occurs. Then the transformed diagrams $D'$ and $D'_1$ are connected by the same move as $D$ and $D_1$, and this move does not involve the vertex of the probe $\gamma$. Hence, the color of this semiarc does not change with the move. Then the colors of the semiarcs $\gamma$ and $\mathdutchcal{SA}(f)(\gamma)$ coincide.

Next, we define a map $\phi\colon\widetilde{\mathdutchcal{R}}_0(T)\to R$. Given a region probe $\gamma\in\widetilde{\mathdutchcal{R}}_0(T)$, consider a semiarc probe $\gamma_a\in\widetilde{\mathdutchcal{SA}}_0(T)$ such that $\gamma=\sigma_r(\gamma_a)$, and set $\phi(\gamma)=\sigma_r(\phi(\gamma_a))$. 

Let us check that the image does not depend on the choice of the semiarc. Let $\gamma_b$ be another semiarc such that $\gamma=\sigma_r(\gamma_b)$. Consider an isotopy $f\colon T\to T'$ which verticalize the probes $\gamma$, $\gamma_a$ and $\gamma_b$. Let $\xi'=f_*(\xi)$ be the corresponding coloring of the diagram $D'=p(T')$ by the biquandloid $(B,R)$. Denote $\gamma'=\mathdutchcal{SA}(f)(\gamma)$, $\gamma_{a'}=\mathdutchcal{SA}(f)(\gamma_a)$, and $\gamma_{b'}=\mathdutchcal{SA}(f)(\gamma_b)$. Since the semiarcs $a'$ and $b'$ are incident to the same region, and $\xi'$ is a biquandloid coloring, we get the equality $\sigma_r(\xi'(a'))=\sigma_r(\xi'(b'))$. Then
\[
\sigma_r(\phi(\gamma_a))=\sigma_r(\phi(\gamma_{a'}))=\sigma_r(\xi'(a'))=\sigma_r(\xi'(b'))=\sigma_r(\phi(\gamma_{b'}))=\sigma_r(\phi(\gamma_b)).
\]

Finally, let us show that the map $\phi\colon\widetilde{\mathdutchcal{SA}}_0(T)\to X$ is a homomorphism of biquandloids. Let $\gamma_1,\gamma_2$ be semiarc probes of the tangle $T$ such that $\sigma_r(\gamma_2)=\sigma_r(\gamma_1)$. Denote $\gamma_3=\gamma_1\circ\gamma_2$, $\gamma_4=\gamma_2\ast\gamma_1$. We need to prove that $\sigma_l(\phi(\gamma_1))=\sigma_r(\phi(\gamma_4))$, $\sigma_l(\phi(\gamma_2))=\sigma_r(\phi(\gamma_3))$, $\sigma_l(\phi(\gamma_3))=\sigma_l(\phi(\gamma_4))$, and $\phi(\gamma_3)=\phi(\gamma_1)\circ\phi(\gamma_2)$, $\phi(\gamma_4)=\phi(\gamma_2)\ast\phi(\gamma_1)$.

Since $\sigma_r(\gamma_1)=\sigma_r(\gamma_2)$, there is a positive crossing probe $\gamma_c$ such that $\gamma_1$ is the homotopy class of the probe $SA_{dr}(\gamma_c)$, and $\gamma_2$ is the homotopy class of $SA_{ur}(\gamma_c)$. Then the homotopy classes of $SA_{dl}(\gamma_c)$ and $SA_{ul}(\gamma_c)$ are $\gamma_4$ and $\gamma_3$. 

Consider an isotopy $f\colon T\to T'$ which verticalize the probe $\gamma_c$. Denote $\gamma'_i=\mathdutchcal{SA}(f)(\gamma_i)$, $i=1,2,3,4$. Then $\gamma'_i$ are the vertical probes of the semiarc $a'_i\in\mathcal{SA}(D')$ incident to a crossing of the diagram $D'=p(T')$. Let $\xi'=f_*(\xi)$ be the biquandloid coloring of $D'$ corresponding to $\xi$. Then
\[
\sigma_l(\phi(\gamma_1))=\sigma_l(\phi(\gamma'_1))=\sigma_l(\xi'(a'_1))=\sigma_r(\xi'(a'_4))=\sigma_r(\phi(\gamma'_4))=\sigma_r(\phi(\gamma_4)),
\]
and
\[
\phi(\gamma_1)\circ\phi(\gamma_2)=\phi(\gamma'_1)\circ\phi(\gamma'_2)=\xi'(a'_1)\circ\xi(a'_2)=\xi'(a'_3)=\phi(\gamma'_3)=\phi(\gamma_3).
\]
In the third equality we used the fact that $\xi'$ is a biquandle coloring. The other relations are proved analogously.
The theorem is proved.
\end{proof}

The construction of the homology of biquandles~\cite{CES} can be extended to biquandloids.

\begin{definition}\label{def:biquandloid_homology}
Let $(B,R)$ be a biquandloid and $A$ an abelian group. Let $C_n(B,A)$ be the free $A$-module generated by
\[
B^{\times_R n}=\{(a_1,\dots,a_n)\in B^{\times n}\mid \sigma_r(a_1)=\sigma_r(a_2)=\cdots=\sigma_r(a_n)\},
\]
and $C_0(B,A)=A[R]$. Define the differential by the formula
\begin{multline*}
  \partial_n(a_1,\dots,a_{n})=\sum_{i=1}^n(-1)^i [(a_1,\dots,\widehat{a_i},\dots,a_n)- \\
  (a_1\circ a_i,\dots,a_{i-1}\circ a_i,a_{i+1}\ast a_i, \dots, a_n\ast a_i)]
\end{multline*}
when $n>1$, and $\partial_1(a_1)=\sigma_r(a_1)-\sigma_l(a_1)$.
Then $(C_*(B,A),\partial_*)$ is a chain complex.

Let $D_n(B,A)$, $n\ge 2$, be the submodule generated by the elements
\[
\{(a_1,\dots,a_n)\in C_n(B,A)\mid \exists i: a_i=a_{i+1}\}.
\]
The homology $H_*(B,A)$ of the quotient complex $C_*(B,A)/D_*(B,A)$ is called \emph{biquandloid homology}.
\end{definition}

\begin{definition}\label{def:biquandloid_cocycle_invariant}
Let $\theta\in Z^2(B,A)$ be a biquandloid 2-cocycle. For a link diagram $D$, the \emph{cocycle invariant} of $D$ is the sum
\[
\Phi_\theta(D)=\sum_{c\in Col_{(B,R)}(D)}[\sum_{c\in\mathcal C(D)}\pm\theta(x,y)]\in\mathbb Z[A],
\]
where $\mathcal C(D)$ is the set of crossings and $\pm\theta(x,y)$ is the Boltzmann weight of the crossing defined in Fig.~\ref{pic:biquandle_bolzmann_weight}.

\begin{figure}
\centering\includegraphics[width=0.3\textwidth]{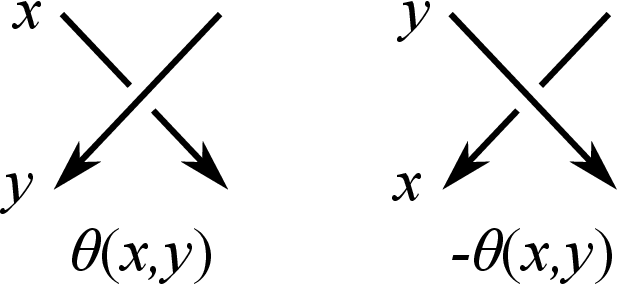}
\caption{The Boltzmann weights}\label{pic:biquandle_bolzmann_weight}
\end{figure}
\end{definition}


\begin{example}[\cite{CG}]
Consider the quandle $S_4=\{0,1,2,3\}$ with the $\ast$-operation given by the table
\[
\left(\begin{array}{cccc}
        0 & 2 & 3 & 1 \\
        3 & 1 & 0 & 2 \\
        1 & 3 & 2 & 0 \\
        2 & 0 & 1 & 3
      \end{array}
\right).
\]
Consider the biquandle double $B_2$ of $S_4$ and the 2-cocycle $\theta\in Z^2(B_2,\mathbb Z_2)$ given by the matrices
\[
\Theta^0=\Theta^1=\left(\begin{array}{cccc}
        0 & 1 & 1 & 0 \\
        1 & 0 & 1 & 0 \\
        1 & 1 & 0 & 0 \\
        0 & 0 & 0 & 0
      \end{array}
\right).
\]
The element $\theta^\epsilon_{ij}$ of $\Theta^\epsilon$ is equal to $\theta((i-1,\epsilon),(j-1,\epsilon))$.

For the trefoil, the cocycle invariant is equal to $\Phi_\theta(3_1)=8\cdot[0]+24\cdot[1]$.
\end{example}

For classical links, biquandles generate invariants that are no stronger than those of quandles~\cite{IT}. In particular, there is a bijection between the set of coloring of a diagram by a biquandle and the set of coloring of it by the associated quandle. This result is generalized to the case of biquandloids.

\begin{definition}\label{def:biquandloid_associated_quandle}
    Let $(B,R)$ be a biquandloid and $z\in R$ an arbitrary element. The set $B_z=\sigma_r^{-1}(z)$ with the binary operation $\triangleleft\colon B_z\times B_z\to B_z$, $x\triangleleft y=(x\ast y)\circ^{-1}y$, $x,y\in B_z$, is called an \emph{associated quandle} of the biquandloid $(B,R)$. 
\end{definition}
A direct check shows that $B_z$ is indeed a quandle, that is, the subset $B_z\subset B$ is closed with respect to the operation $\triangleleft$, and the operation $\triangleleft$ satisfies the properties of the quandle.

\begin{proposition}\label{prop:biquandloid_vs_quandle_colorings}
    Let $D\subset \R^2$ be a diagram of a classical oriented link and $(B,R)$ a biquandloid. Then for any element $r\in R$ there is a bijection between the set of quandle colorings $Col_{B_r}(D)$ and the set $Col_{(B,R)}(D,r)$ of biquandloid colorings such that the color of the unbounded region is $r$. 
\end{proposition}

\begin{proof}
    Following~\cite[Theorem 3.1]{IT}, we construct maps $\Phi\colon Col_{(B,R)}(D,r)\to Col_{B_r}(D)$ and $\Psi\colon Col_{B_r}(D)\to Col_{(B,R)}(D,r)$.

    The map $\Phi$ is defined as follows. Let $c\in  Col_{(B,R)}(D,r)$. Given an arc $a\in\mathcal A(D)$, consider a sequence of second Reidemeister moves $f\colon D\to D'$ that pulls the arc to the unbounded region (Fig.~\ref{pic:biquandloid_to_quandle}). Let $c'=f_*(c)$ be the biquandloid coloring of $D'$ that corresponds to the coloring $c$. Then the color $\Phi(c)(a)$ of the arc $a$ is the color $c'(a')$.
\begin{figure}[h]
    \centering
    \includegraphics[width=0.6\textwidth]{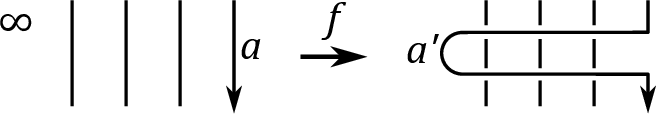}
    \caption{Pulling an arc}
    \label{pic:biquandloid_to_quandle}
\end{figure}

    The inverse map $\Psi$ is defined as follows. Let $c\in Col_{B_r}(D)$. Consider the diagram $W(D)=D\cup -\bar D$ which is the union of $D$ and the inverse mirror diagram as shown in Fig.~\ref{pic:diagram_doubling} middle. The coloring $c$ induces a $(B,R)$-coloring $\tilde c$ of the diagram $W(D)$ (Fig.~\ref{pic:crossing_doubling}). Pull away the mirror link $-\bar D$ as shown in Fig.~\ref{pic:diagram_doubling} right. Then the image $\Psi(c)$ is the restriction to the component $D$ of the corresponding coloring $\tilde c'$ of the obtained diagram $W'(D)$.
\begin{figure}[h]
    \centering
    \includegraphics[width=0.7\textwidth]{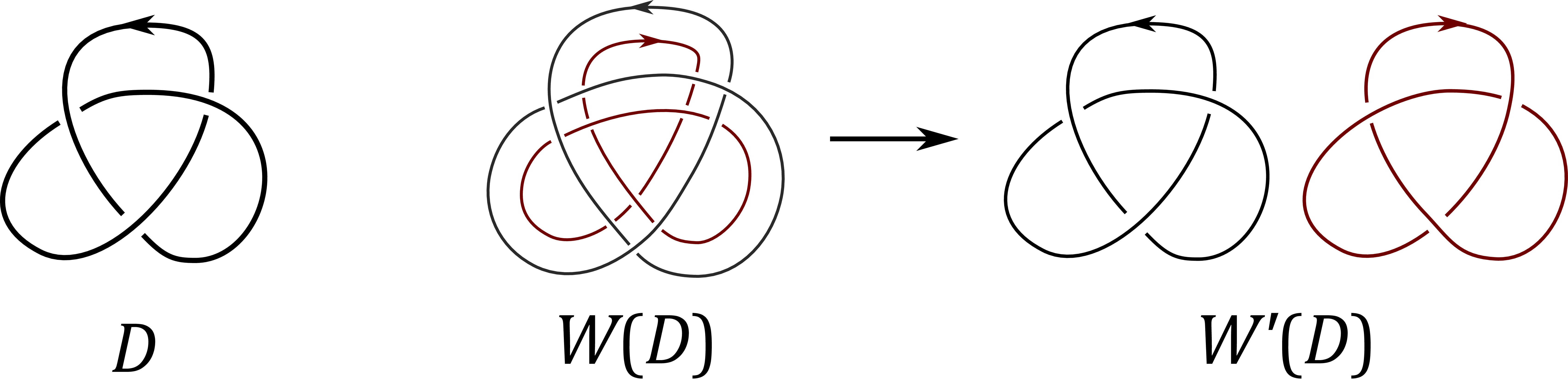}
    \caption{Diagrams $W(D)=D\cup -\bar D$ and $W'(D)=D\sqcup -\bar D$}
    \label{pic:diagram_doubling}
\end{figure}
\begin{figure}[h]
    \centering
    \includegraphics[width=0.6\textwidth]{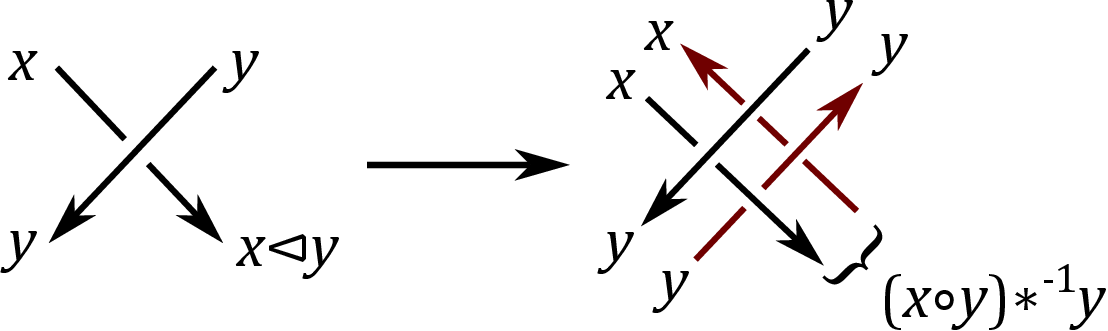}
    \caption{The coloring of the diagram $W(D)$}
    \label{pic:crossing_doubling}
\end{figure}

    The proof that $\Phi$ and $\Psi$ are inverse to each other is analogous to the proof in~\cite[Theorem 3.1]{IT}. 
\end{proof}


{\color{red}





}

\subsection{Coloring of crossings: Crossoid}\label{subsect:crossoid}

Let us introduce an algebraic structure which one can impose on the homotopy classes of crossings. For simplicity, we will consider the set $\mathdutchcal C^{uf}(T)$ of unframed classes that immolate the bit of information responsible for framing.

\begin{figure}[h]
    \centering
    \includegraphics[width=0.2\textwidth]{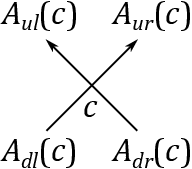}
    \caption{Incidence relations between a crossing and the semiarcs}
    \label{pic:crossoid_incidence}
\end{figure}

An outline of the structure is as follows. Given a tangle diagram $D$, we will color its crossings by elements of some set $X$, and color the semiarcs of $D$ by elements of some set $A$. The colors of a crossing and the four arcs incident to it (Fig.~\ref{pic:crossoid_incidence}) are connected by some incidence maps from $X$ to $A$ .

\begin{figure}[h]
    \centering
    \includegraphics[width=0.25\textwidth]{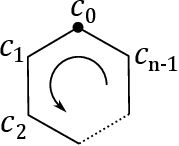}
    \caption{A polygon}
    \label{pic:crossoid_polygon}
\end{figure}

The core of the structure is the \emph{polygon map}: given an $n$-gonal region of the diagram (Fig.~\ref{pic:crossoid_polygon}), the color $c_0$ of one crossing of the region is uniquely determined by the colors $(c_1,\dots, c_{n-1})$ of the other crossings. Herewith, the crossing colors $c_0,\dots,c_{n-1}$ must be compatible, i.e. the adjacent crossings must be incident to the same semiarc.

Now, let us give a formal definition of the structure.

\begin{definition}\label{def:crossoid}
A \emph{crossoid} is a pair of sets $(C,A)$ ($C$ is called the \emph{crossing set} and $A$ is the \emph{arc set}) and the following maps:
\begin{enumerate}
    \item a sign map $sgn\colon C\to\{-1,1\}$. We denote $C_\pm=sgn^{-1}(\pm 1)$;
    \item four incidence maps $A_{\alpha\beta}\colon C\to A$, $\alpha\in\{u,d\}$, $\beta\in\{l,r\}$.

        Using the incidence map, we define the incoming and the outcoming incidence maps $A_{in}, A_{out}\colon C\times\{\leftarrow,\rightarrow\}^2\to A$ according to the table
\begin{equation}\label{eq:table_incoming_outcoming_arc}
\begin{array}{|c|c|c|}
\hline
o_1o_2 & A_{in}(c,o_1,o_2) & A_{out}(c,o_1,o_2)\\
\hline
\rightarrow\rightarrow &  A_{dl}(c) & A_{ul}(c)\\
\leftarrow\rightarrow &  A_{ul}(c) & A_{ur}(c)\\
\rightarrow\leftarrow &  A_{dr}(c) & A_{dl}(c)\\
\leftarrow\leftarrow &  A_{ur}(c) & A_{dr}(c)\\
\hline
\end{array}
\end{equation}
    \item polygonal maps $P_{\bm o}^{\bm\epsilon}\colon C(\hat{\bm\epsilon},\bm o)\to C_{\epsilon_0}$, where $n\ge 2$, $\bm\epsilon=(\epsilon_0,\epsilon_1,\dots,\epsilon_{n-1})$,
    $\hat{\bm\epsilon}=(\epsilon_1,\dots,\epsilon_{n-1})$, $\epsilon_i\in\{+,-\}$,  $\bm o =(o_0,\dots,o_{n-1})$, $o_i\in\{\leftarrow,\rightarrow\}$, and
\begin{multline*}
C(\hat{\bm\epsilon},\bm o)=\{(c_1,\dots,c_{n-1})\in C_{\epsilon_1}\times\cdots\times C_{\epsilon_{n-1}}\mid \\
A_{out}(c_i,o_{i-1},o_i)=A_{in}(c_{i+1},o_i,o_{i+1}), i=1,\dots,n-2\},
\end{multline*}
and $P^\epsilon_o\colon A\to C_\epsilon$, $\epsilon\in\{+,-\}$, $o\in\{\leftarrow,\rightarrow\}$.
\end{enumerate}

The maps must satisfy the following relations:
\begin{enumerate}

\item polygon incidence relations
\begin{eqnarray*}
A_{in}(P_{\bm o}^{\bm\epsilon}(\bm c),o_{n-1},o_0)&=&A_{out}(c_{n-1},o_{n-2},o_{n-1}),\\
A_{out}(P_{\bm o}^{\bm\epsilon}(\bm c),o_{n-1},o_0)&=&A_{in}(c_{1},o_0,o_1),
\end{eqnarray*}
where $\bm c=(c_1,\dots,c_{n-1})\in C_{n-1}(\bm\epsilon,\bm o)$, and
\[
A_{in}(P^\epsilon_o(a),o,o)=A_{out}(P^\epsilon_o(a),o,o)=a,
\]
where $\epsilon\in\{+,-\}$, $o\in\{\leftarrow,\rightarrow\}$, $a\in A$;

\item polygon rotation relations
\begin{equation}\label{eq:crossoid_rotation_equation}
P^{(\epsilon_1,\dots,\epsilon_{n-1},\epsilon_0)}_{(o_1,\dots,o_{n-1},o_0)}\left(c_2,\dots,c_{n-1},P^{(\epsilon_0,\dots,\epsilon_{n-1})}_{(o_0,\dots,o_{n-1})}(c_1,\dots,c_{n-1})\right)=c_1;
\end{equation}

\item incidence $\Omega_1$-relations:
\[
A_{ur}\circ P^\epsilon_\rightarrow = A_{dr}\circ P^\epsilon_\rightarrow,\quad A_{ul}\circ P^\epsilon_\leftarrow = A_{dl}\circ P^\epsilon_\leftarrow,\quad \epsilon\in\{+,-\},
\]
     and the maps $A_{ur}\circ P^\epsilon_\rightarrow$, $A_{ul}\circ P^\epsilon_\leftarrow$ are bijections.

Define the \emph{loop maps} by the formulas
\begin{equation}\label{eq:crossoid_loop_maps}
L^\epsilon_\rightarrow=P^\epsilon_\rightarrow\circ(A_{ur}\circ P^\epsilon_\rightarrow)^{-1},\quad L^\epsilon_\leftarrow=P^\epsilon_\leftarrow\circ(A_{ul}\circ P^\epsilon_\leftarrow)^{-1},\quad \epsilon\in\{+,-\}.
\end{equation}

\item inner $\Omega_1$-relations
\begin{multline}\label{eq:crossoid_R1in}
P^{(\epsilon_0,\epsilon,\epsilon,\epsilon_1,\dots,\epsilon_{n-1})}_{(o_0,\overline{o_0},o_0,o_1,\dots, o_{n-1})}\left(L^\epsilon_{o_0}(A_{in}(c_1,o_0,o_1)),L^\epsilon_{o_0}(A_{in}(c_1,o_0,o_1)),c_1,\dots,c_{n-1}\right)=\\
P^{(\epsilon_0,\dots,\epsilon_{n-1})}_{(o_0,\dots,o_{n-1})}(c_1,\dots,c_{n-1}),
\end{multline}
where $(c_1,\dots,c_{n-1})\in C((\epsilon_1,\dots,\epsilon_{n-1}),(o_0,\dots,o_{n-1}))$, $\epsilon\in\{+,-\}$, and
\[
P^{(\epsilon,\epsilon',\epsilon')}_{(o,\bar o, o)}(L^{\epsilon'}_o(a),L^{\epsilon'}_o(a))=P^\epsilon_o(a);
\]

\item outer $\Omega_1$-relations
\begin{multline}\label{eq:crossoid_R1out}
P^{(\epsilon_0,\epsilon,\epsilon_1,\dots,\epsilon_{n-1})}_{(o_0,o_0,o_1,\dots, o_{n-1})}\left(L^\epsilon_{\overline{o_0}}(A_{in}(c_1,o_0,o_1)),c_1,\dots,c_{n-1}\right)=\\
P^{(\epsilon_0,\dots,\epsilon_{n-1})}_{(o_0,\dots,o_{n-1})}(c_1,\dots,c_{n-1}),
\end{multline}
and
\[
P^{(\epsilon,\epsilon')}_{(o,o)}(L^{\epsilon'}_{\bar o}(a))=P^\epsilon_o(a);
\]

\item incidence $\Omega_2$-relations:
\begin{gather*}
A_{in}\left(P^{(\epsilon,-\epsilon)}_{(o_0,o_1)}(c),\overline{o_1},\overline{o_0}\right)=A_{out}(c,\overline{o_0},\overline{o_1}),\\
A_{out}\left(P^{(\epsilon,-\epsilon)}_{(o_0,o_1)}(c),\overline{o_1},\overline{o_0}\right)=A_{in}(c,\overline{o_0},\overline{o_1});
\end{gather*}

\item inner $\Omega_2$-relations
\begin{multline}\label{eq:crossoid_R2in}
P^{(\epsilon_0,-\epsilon,\epsilon_{i+1},\dots,\epsilon_{n-1})}_{(o_0,o_i,\dots, o_{n-1})}\left(P^{(-\epsilon,\epsilon)}_{(\overline{o_i},\overline{o_0})}\left(P^{(\epsilon,\epsilon_1,\dots,\epsilon_i)}_{(o_0,o_1,\cdots,o_i)}(c_1,\dots,c_i)\right),c_{i+1},\dots,c_{n-1}\right)=\\
P^{(\epsilon_0,\epsilon_1,\dots,\epsilon_{n-1})}_{(o_0,o_1,\dots,o_{n-1})}\left(c_1,\dots,c_{n-1}\right);
\end{multline}

\item outer $\Omega_2$-relations
\begin{multline}\label{eq:crossoid_R2out}
P^{(\epsilon_0,-\epsilon,\epsilon,\epsilon_1,\dots,\epsilon_{n-1})}_{(o_0,o,o_0,o_1,\dots,o_{n-1})}\left(P^{(-\epsilon,\epsilon)}_{(o_0,\overline o)}(c),c,c_1,\dots,c_{n-1}\right)= \\
P^{(\epsilon_0,\epsilon_1,\dots,\epsilon_{n-1})}_{(o_0,o_1,\dots,o_{n-1})}\left(c_1,\dots,c_{n-1}\right),
\end{multline}
where $(c,c_1,\dots,c_{n-1})\in C((\epsilon,\epsilon_1,\dots,\epsilon_{n-1}),(o,o_0,o_1,\dots,o_{n-1}))$, and
\[
P^{(\epsilon,-\epsilon',\epsilon')}_{(o,o',o)}\left(P^{(-\epsilon',\epsilon')}_{(o,\overline{o'})}(c),c\right)=P^\epsilon_o(A_{out}(c,o',o)),
\]
where $c\in C_{\epsilon'}$, $\epsilon,\epsilon'\in\{+,-\}$, $o,o'\in\{\leftarrow,\rightarrow\}$;

\item $\Omega_3$-relations

For any $y,z\in C_+$ and $o\in\{\leftarrow,\rightarrow\}$ such that $A_{out}(y,o,o)=A_{in}(z,o,o)$ denote $x=P_{(o,o,o)}^{(-,+,+)}(y,z)$ and
\begin{gather}\label{eq:crossoid_R3_new_crossing}
    x'=P_{(o,o)}^{(-,+)}\circ P_{(o,\bar o,o)}^{(+,+,+)}(z,y),\quad
    y'=P_{(o,o)}^{(+,-)}\circ P_{(o,\bar o,o)}^{(-,-,+)}(x,z),\nonumber\\  
    z'=P_{(o,o)}^{(+,-)}\circ P_{(o,\bar o,o)}^{(-,+,-)}(y,x). 
\end{gather}
Then $A_{out}(y',\bar o,\bar o)=A_{in}(z',\bar o,\bar o)$ and $x'=P_{(\bar o,\bar o,\bar o)}^{(-,+,+)}(y',z')$.


\end{enumerate}
\end{definition}

\begin{definition}\label{def:crossoid_coloring}
    Let $(C,A)$ be a crossoid and $D$ a filling tangle diagram in $F$ (i.e. $F\setminus D$ is a union of cells). A \emph{coloring of the diagram $D$ with the crossoid $(C,A)$} is a pair of maps $\chi_C\colon \mathcal C(D)\to C$, $\chi_A\colon \mathcal{SA}(D)\to A$ such that:
\begin{enumerate}
    \item for any crossing $c$ the colors of the incident semiarcs satisfy the incidence relation in Fig.~\ref{pic:crossoid_incidence};
    \item for any polygonal region (Fig.~\ref{pic:crossoid_polygon}) the colors $c_0,\dots,c_{n-1}$ satisfy the polygon relation
\[
c_0 = P^{(\epsilon_0,\dots,\epsilon_{n-1})}_{(o_0,\dots,o_{n-1})}(c_1,\dots,c_{n-1}),
\]
where $\epsilon_i$ is the sign of the $i$-th crossing, and $o_i=\rightarrow$ if the edge from the $i$-th to the $(i+1)$-th crossing is oriented counterclockwise, otherwise $o_i=\leftarrow$.
\end{enumerate}
Denote the set of colorings by $Col_{(C,A)}(D)$.
\end{definition}

\begin{remark}
Now, having the definition of crossoid coloring, we can give a geometric interpretation of the relations in the definition of a crossoid.

The incoming and the outcoming arcs~\eqref{eq:table_incoming_outcoming_arc} of a crossings in a polygonal face are the edges of the polygon (oriented counterclockwise) incident to the crossing (Fig.~\ref{pic:crossoid_polygon_incidence}). The polygon rotation relation~\eqref{eq:crossoid_rotation_equation} corresponds to the change of base crossing in a polygonal region (Fig.~\ref{pic:crossoid_rotation}). The inner and outer first Reidemeister move relations~\eqref{eq:crossoid_R1in},\eqref{eq:crossoid_R1out} correspond the polygon transformation in Fig.~\ref{pic:crossoid_R1_move}. The inner and outer second Reidemeister move relations describe the transformations in Fig.~\ref{pic:crossoid_R2_move}.

According to M. Polyak~\cite{Polyak}, the oriented version of the third Reidemeister move in Fig.~\ref{pic:reidemeister_R3a} using second Reidemeister moves generates the other versions of oriented $\Omega_3$-moves. Thus, when we look at how a third Reidemeister move impacts a polygonal region (Fig.~\ref{pic:crossoid_R3_move}), we can consider only the move $\Omega_{3a}$. To calculate the colors of the new crossings $x', y', z'$, we use an auxiliary crossing as shown in Fig.~\ref{pic:crossoid_R3_auxiliary}. then $\Omega_3$-relation means that the new crossings form a correct triangle, i.e. satisfy the polygon relation.

\begin{figure}[h]
    \centering
    \includegraphics[width=0.3\textwidth]{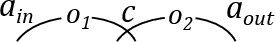}
    \caption{The incoming and outcoming incident arcs}
    \label{pic:crossoid_polygon_incidence}
\end{figure}

\begin{figure}[h]
    \centering
    \includegraphics[width=0.5\textwidth]{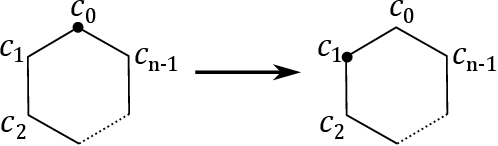}
    \caption{Change of the base crossing in a polygon}
    \label{pic:crossoid_rotation}
\end{figure}

\begin{figure}[h]
    \centering
    \includegraphics[width=0.45\textwidth]{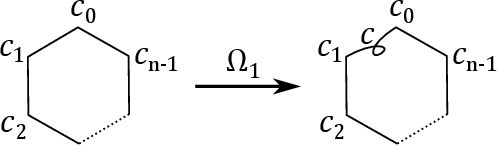}\quad
    \includegraphics[width=0.45\textwidth]{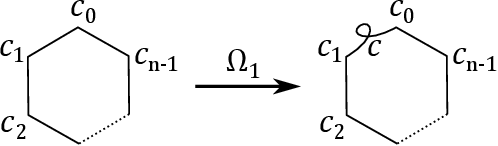}
    \caption{Inner and outer $\Omega_1$-moves}
    \label{pic:crossoid_R1_move}
\end{figure}

\begin{figure}[h]
    \centering
    \includegraphics[width=0.45\textwidth]{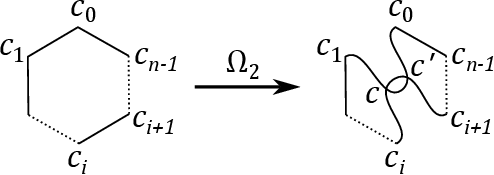}\quad
    \includegraphics[width=0.45\textwidth]{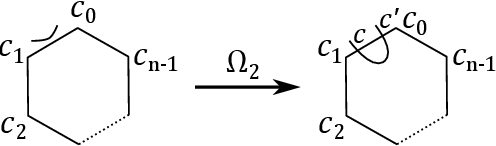}
    \caption{Inner and outer $\Omega_2$-moves}
    \label{pic:crossoid_R2_move}
\end{figure}

\begin{figure}[h]
    \centering
    \includegraphics[width=0.3\textwidth]{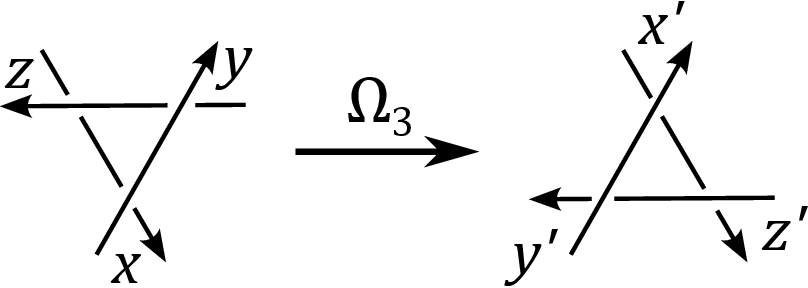}
    \caption{The oriented move $\Omega_{3a}$}
    \label{pic:reidemeister_R3a}
\end{figure}

\begin{figure}[h]
    \centering
    \includegraphics[width=0.5\textwidth]{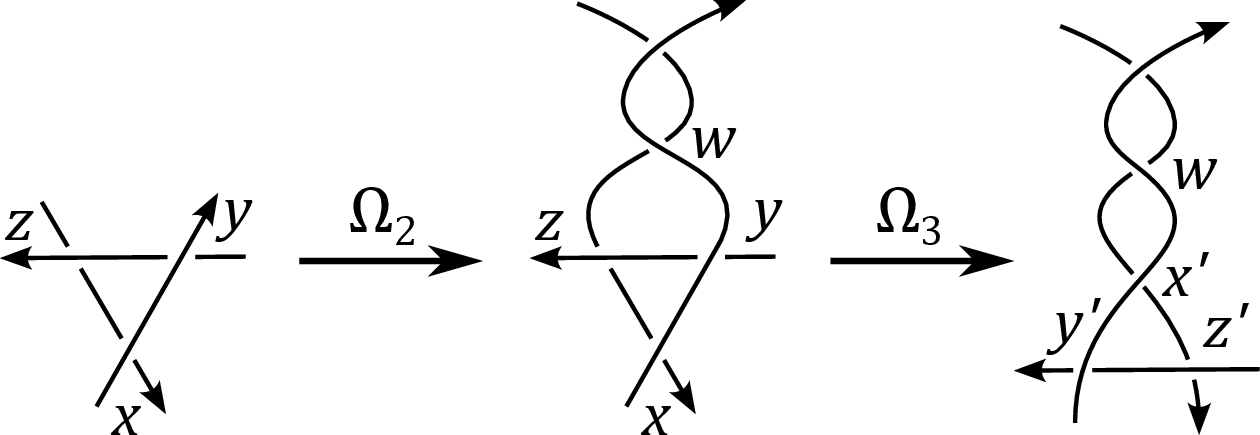}
    \caption{Scheme of finding the color of a new crossing in the move $\Omega_{3a}$}
    \label{pic:crossoid_R3_auxiliary}
\end{figure}

\end{remark}

\begin{remark}\label{rem:crossoid_inner_R2_to_h_relation}
The inner $\Omega_2$-relations imply the $h$-relation
\begin{equation}\label{eq:crossoid_h}
    P_{(o,o')}^{(\epsilon,-\epsilon)}=P_{(\bar o,\bar o')}^{(\epsilon,-\epsilon)},\quad \epsilon\in\{+,-\}, o,o'\in\{\leftarrow,\rightarrow\},
\end{equation}
which corresponds to the homotopical equivalence of crossings (Fig.~\ref{pic:homotopy_crossing_relation}). Indeed, for any $x\in C_\epsilon$, consider $y=P^{(-\epsilon,\epsilon)}_{(o',o)}$. Then
\[
x=P^{(\epsilon,-\epsilon)}_{(o,o')}(y)=P^{(\epsilon,-\epsilon)}_{(o,o')}\circ P^{(-\epsilon,\epsilon)}_{(\bar o',\bar o)}\circ P^{(\epsilon,-\epsilon)}_{(o,o')}(y)=P^{(\epsilon,-\epsilon)}_{(o,o')}\circ P^{(-\epsilon,\epsilon)}_{(\bar o',\bar o)}(x),
\]
hence, $P_{(o,o')}^{(\epsilon,-\epsilon)}=\left(P^{(-\epsilon,\epsilon)}_{(\bar o',\bar o)}\right)^{-1}=P^{(\epsilon,-\epsilon)}_{(\bar o,\bar o')}$.
\end{remark}

\begin{remark}\label{rem:crossoid_R3_polygon_relation}
    The formula~\eqref{eq:crossoid_R3_new_crossing} for the new crossings in a third Reidemeister move implies the following relation (see Fig.~\ref{pic:crossoid_R3_move}):
\begin{multline}\label{eq:crossoid_R3}
P^{(\epsilon_0,\epsilon_1,\epsilon_2,\epsilon_3\dots,\epsilon_{n-1})}_{(o,\overline o,o,o_3,\dots,o_{n-1})}(c_1,c_2,c_3,\dots,c_{n-1})=\\
P^{(\epsilon_0,\epsilon,\epsilon_3\dots,\epsilon_{n-1})}_{(o,o,o_3,\dots,o_{n-1})} \left(P^{(\epsilon,-\epsilon)}_{(o,o)} \left(P^{(-\epsilon,\epsilon_1,\epsilon_2)}_{(o,\overline o,o)}(c_1,c_2)\right), c_3,\dots,c_{n-1}\right)
\end{multline}
where $\{\epsilon,\epsilon_1,\epsilon_2\}=\{+,+,-\}$ as multisets, $o\in\{\leftarrow,\rightarrow\}$.
\begin{figure}[h]
    \centering
    \includegraphics[width=0.45\textwidth]{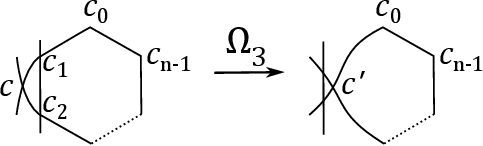}
    \caption{$\Omega_3$-move}
    \label{pic:crossoid_R3_move}
\end{figure}
\end{remark}

From the geometric interpretation of the crossoid relations, we immediately get the following statement.

\begin{theorem}\label{thm:crossoid_coinvariance}
    Let $(C,A)$ be a crossoid and $f\colon D\to D'$ a Reidemeister move. Then $f$ induces a bijection between the coloring sets $Col_{(C,A)}(D)$ and $Col_{(C,A)}(D')$.
\end{theorem}

Moreover, the correspondence $D\mapsto Col_{(C,A)}(D)$ defines a coinvariant 
\[
Col_{(C,A)}\colon \mathfrak D_s\to Rel,
\]
$(Map(Col_{(C,A)}(D),C),ev_C)$ is an $h$-coinvariant of the crossing functor $\mathcal C$, and $(Map(Col_{(C,A)}(D),A),ev_A)$ is an $h$-coinvariant of the semiarc functor $\mathcal{SA}$. The single-valued natural transformations $ev_A$, $ev_C$ are defined by the formulas
\[
ev_C(c)(\chi)=\chi_C(c),\quad ev_A(a)(\chi)=\chi_A(a), 
\]
where $c\in\mathcal C(D)$, $a\in\mathcal{SA}(D)$, $\chi=(\chi_C,\chi_A)\in  Col_{(C,A)}(D)$.

\begin{remark}\label{rem:crossoid_unframed}
    In Section~\ref{sect:crossing_R2_equivalence} we considered $\Omega_2$-equivalence on crossings. For crossoids, $\Omega_2$-equivalence takes the form of an unframing relation
\begin{equation}\label{eq:crossoid_unframing}
    P^{(\epsilon,-\epsilon)}_{(\rightarrow,\rightarrow)}=P^{(\epsilon,-\epsilon)}_{(\leftarrow,\rightarrow)}=P^{(\epsilon,-\epsilon)}_{(\rightarrow,\leftarrow)}=P^{(\epsilon,-\epsilon)}_{(\leftarrow,\leftarrow)},\ \epsilon\in\{+,-\}.
\end{equation}
 Due to the $h$-relation, the second and third equalities in~\eqref{eq:crossoid_unframing} follow from the first. 

We say that a crossoid $(C,A)$ is \emph{unframed} is it satisfies the unframing relation~\eqref{eq:crossoid_unframing}.
\end{remark}

Using internal $\Omega_2$-relations, one can split polygonal maps into composition of monogonal, bigonal and trigonal maps. Then we can reformulate the definition of an (unframed) crossoid in the following form.

\begin{definition}\label{def:crossoid_reduced}
An \emph{crossoid} is a pair of sets $(C,A)$ ($C$ is called the \emph{crossing set} and $A$ is the \emph{arc set}) and the following maps:
\begin{enumerate}
    \item a sign map $sgn\colon C\to\{-1,1\}$. We denote $C_\pm=sgn^{-1}(\pm 1)$;
    \item four incidence maps $A_{\alpha\beta}\colon C\to A$, $\alpha\in\{u,d\}$, $\beta\in\{l,r\}$ (Fig.~\ref{pic:crossoid_incidence});
    
    \item four loop maps $L^\epsilon_{\beta}\colon A\to C_\epsilon$, $\beta\in\{l,r\}$, $\epsilon\in\{+,-\}$ (Fig.~\ref{pic:crossoid_loop_maps});
\begin{figure}[h]
    \centering
    \includegraphics[width=0.5\textwidth]{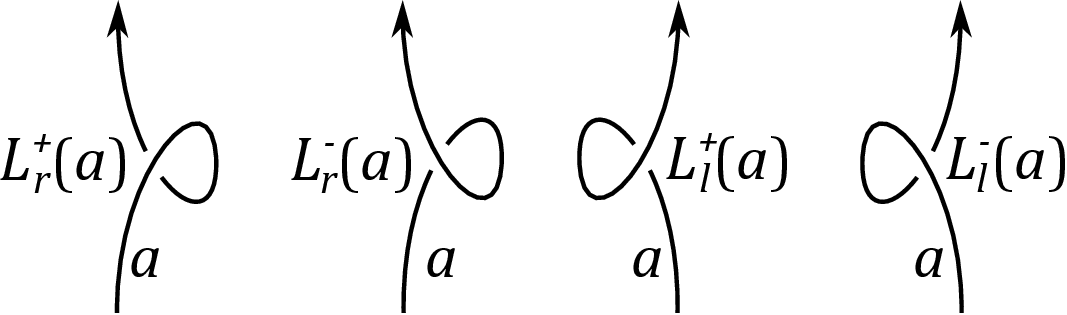}
    \caption{Loop maps}
    \label{pic:crossoid_loop_maps}
\end{figure}

    \item two non-alternating bigon map $i_s, i_w\colon C\to C$ (Fig.~\ref{pic:crossoid_nonalternating_bigons});
\begin{figure}[h]
    \centering
    \includegraphics[width=0.35\textwidth]{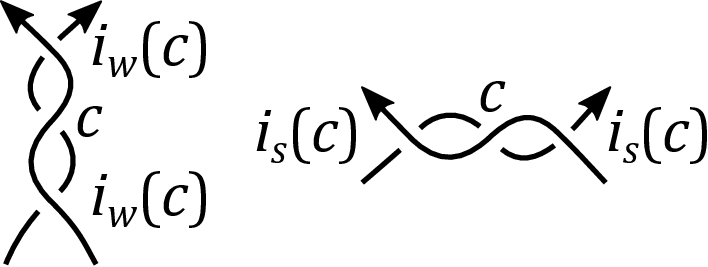}
    \caption{Non-alternating bigon maps}
    \label{pic:crossoid_nonalternating_bigons}
\end{figure}

    \item four alternating bigon maps $B_\alpha\colon C\to C$, $\alpha\in\{u,d,l,r\}$ (Fig.~\ref{pic:crossoid_alternating_bigons});
\begin{figure}[h]
    \centering
    \includegraphics[width=0.6\textwidth]{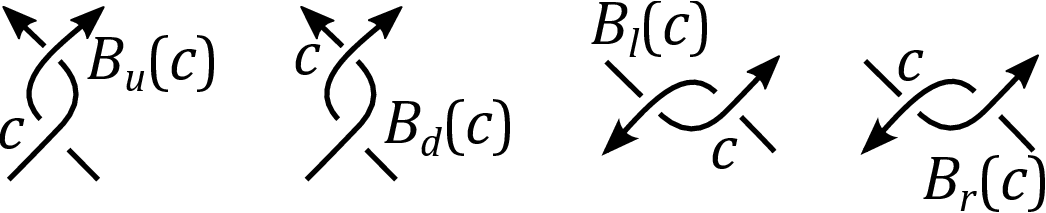}
    \caption{Alternating bigon maps}
    \label{pic:crossoid_alternating_bigons}
\end{figure}

    \item eight triangle maps (Fig.~\ref{pic:crossoid_trigons})
\[
T_{o_0o_1o_2}\colon \{(c_1,c_2)\in C_+^2\mid A_{out}(c_1,o_0,o_1)=A_{in}(c_2,o_1,o_2)\}\to C_-,
\]
$o_0,o_1,o_2\in\{\leftarrow,\rightarrow\}$, and the maps $A_{out}$, $A_{in}$ are defined by the table~\eqref{eq:table_incoming_outcoming_arc}.
\begin{figure}[h]
    \centering
    \includegraphics[width=0.8\textwidth]{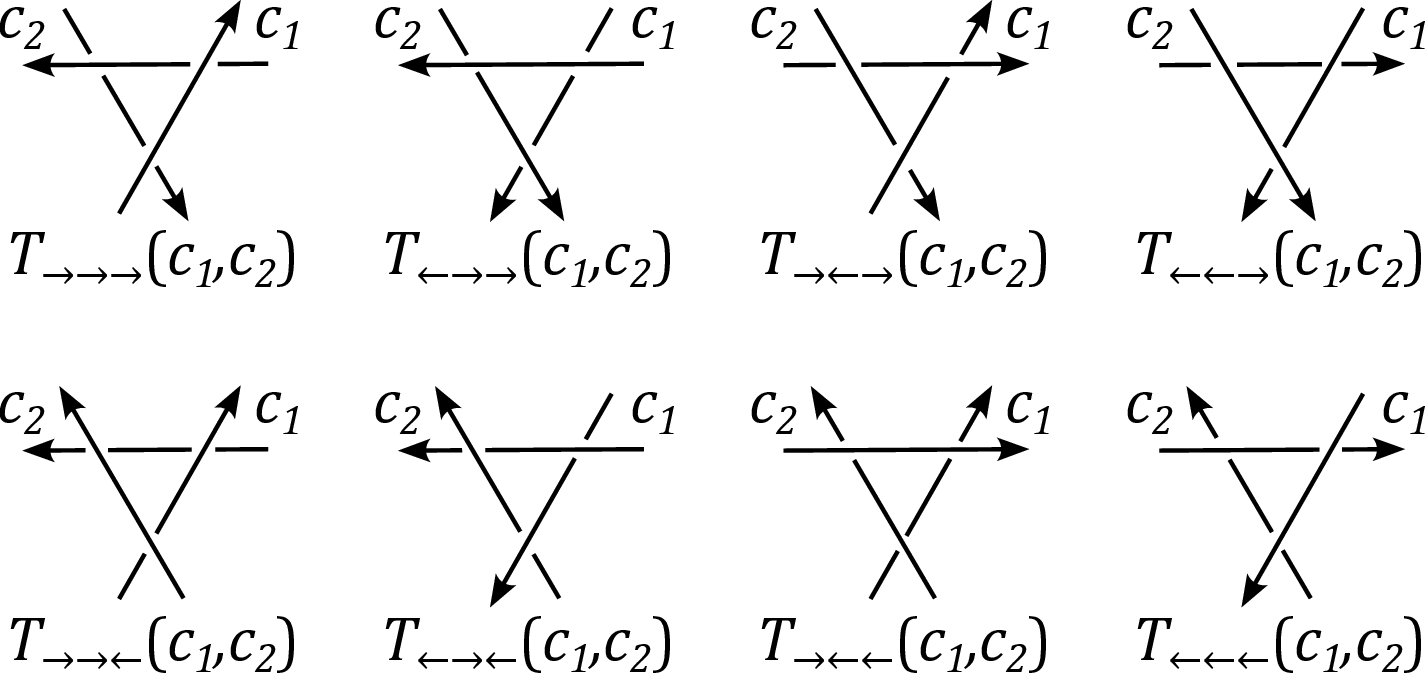}
    \caption{Triangle maps}
    \label{pic:crossoid_trigons}
\end{figure}
\end{enumerate}

The maps satisfy the following set of relations:
\begin{enumerate}
    \item relations with the sign map
\[
sgn\circ i_s=sgn\circ i_w=-sgn,\quad sgn\circ B_\alpha=sgn,\quad \alpha\in\{u,d,l,r\};
\]
    \item relations with the incidence maps:
    \begin{itemize}
        \item with the loop maps
\[
A_{ul}\circ L_l^\epsilon=A_{dl}\circ L_l^\epsilon,\quad A_{ur}\circ L_r^\epsilon=A_{dr}\circ L_r^\epsilon,\quad\epsilon=\pm;
\]
        \item with the bigon maps
\[
A_{dr}\circ i_\beta= A_{ur},\ A_{dl}\circ i_\beta= A_{ul},\ A_{ur}\circ i_\beta= A_{dr},\ A_{ul}\circ i_\beta= A_{dl},\ \beta=s,w,
\]
\begin{gather*}
A_{dr}\circ B_u= A_{ur},\ A_{dl}\circ B_u= A_{ul},\quad A_{ur}\circ B_d= A_{dr},\ A_{ul}\circ B_d= A_{dl},\\
A_{ul}\circ B_l= A_{dl},\ A_{dl}\circ B_l= A_{ul},\quad A_{ur}\circ B_r= A_{dr},\ A_{dr}\circ B_r= A_{ur};
\end{gather*}
        \item with the triangle maps
\begin{gather*}
A_{out}(T_{o_0o_1o_2}(c_1,c_2),o_2,o_0)=A_{in}(c_1,o_0,o_1),\\
A_{in}(T_{o_0o_1o_2}(c_1,c_2),o_2,o_0)=A_{out}(c_2,o_1,o_2);
\end{gather*}
    \end{itemize}
    \item $\Omega_1$-relations:
    \begin{itemize}
        \item with the incidence maps
\[
A_{ur}\circ L_l^\epsilon=A_{dr}\circ L_l^\epsilon=A_{ul}\circ L_r^\epsilon=A_{dl}\circ L_r^\epsilon=\mathrm{id}_A,\ \epsilon=\pm;
\]
        \item with the bigon maps (Fig.~\ref{pic:crossoid_R1_outer_12-gon} left)
\begin{gather*}
i_s\circ L_l^\epsilon=L^{-\epsilon}_r\circ A_{ul}\circ L_l^\epsilon,\quad i_s\circ L_r^\epsilon=L^{-\epsilon}_l\circ A_{ur}\circ L_r^\epsilon,\\
B_l\circ L_l^\epsilon=L^{\epsilon}_r\circ A_{ul}\circ L_l^\epsilon, \quad
B_r\circ L_r^\epsilon=L^{\epsilon}_l\circ A_{ur}\circ L_r^{\epsilon},
\end{gather*}

        \item with the triangle maps (Fig.~\ref{pic:crossoid_R1_outer_12-gon} right)
\begin{gather*}
T_{\rightarrow\rightarrow\rightarrow}(c,B_l(c))=L^-_r(A_{dl}(c)),\quad
T_{\rightarrow\leftarrow\rightarrow}(c,B_d(c))=L^-_r(A_{dr}(c)),\\
T_{\leftarrow\rightarrow\leftarrow}(c,B_u(c))=L^-_l(A_{ul}(c)),\quad
T_{\leftarrow\leftarrow\leftarrow}(c,B_r(c))=L^-_l(A_{ur}(c)),
\end{gather*}
for any $c\in C_+$.
\begin{figure}[h]
    \centering
    \includegraphics[width=0.3\textwidth]{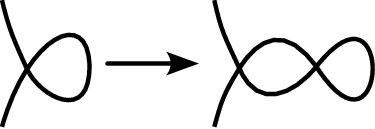}\qquad
    \includegraphics[width=0.43\textwidth]{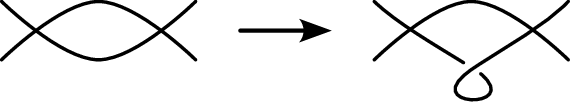}
    \caption{Outer $\Omega_1$-move in a loop and a bigon}
    \label{pic:crossoid_R1_outer_12-gon}
\end{figure}

        \item degenerated triangle relations (Fig.~\ref{pic:crossoid_R1_inner_loop})
\begin{gather*}
T_{\rightarrow\leftarrow\rightarrow}(L^+_l(a), L^+_l(a))=L^-_l\circ A_{ur}\circ L^+_r(a),\\
T_{\leftarrow\rightarrow\leftarrow}(L^+_r(a), L^+_r(a))=L^-_r\circ A_{ul}\circ L^+_l(a),\\
T_{\rightarrow\leftarrow\rightarrow}(i_w\circ B_d\circ L^-_l(a), i_w\circ B_u\circ L^-_l(a))=L^-_l\circ A_{ur}\circ L^+_r(a),\\
T_{\leftarrow\rightarrow\leftarrow}(i_w\circ B_u\circ L^-_r(a), i_w\circ B_d\circ L^-_r(a))=L^-_r\circ A_{ul}\circ L^+_l(a);
\end{gather*}
\begin{figure}[h]
    \centering
    \includegraphics[width=0.65\textwidth]{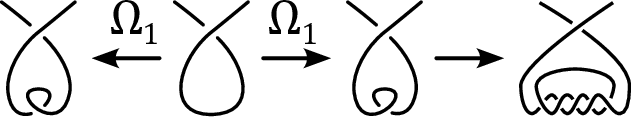}
    \caption{Inner $\Omega_1$-move in a loop}
    \label{pic:crossoid_R1_inner_loop}
\end{figure}

    \end{itemize}
    \item $\Omega_2$-relations
\begin{itemize}
    \item with the loop maps (Fig.~\ref{pic:crossoid_R2_outer_loop})
\begin{gather*}
    T_{\rightarrow\rightarrow\rightarrow}(c, i_s\circ B_l\circ i_w(c))=L^-_l\circ A_{ur}\circ L^+_r\circ A_{ul}(c),\\
    T_{\rightarrow\rightarrow\rightarrow}(i_s\circ B_l\circ i_w(c), c)=L^-_l\circ A_{ur}\circ L^+_r\circ A_{dl}(c),\\
    T_{\rightarrow\leftarrow\rightarrow}(c, i_w\circ B_u\circ i_s(c))=L^-_l\circ A_{ur}\circ L^+_r\circ A_{dr}(c),\\
    T_{\rightarrow\leftarrow\rightarrow}(i_w\circ B_d\circ i_s(c), c)=L^-_l\circ A_{ur}\circ L^+_r\circ A_{ur}(c),\\
    T_{\leftarrow\rightarrow\leftarrow}(c, i_w\circ B_d\circ i_s(c))=L^-_r\circ A_{ul}\circ L^+_l\circ A_{ul}(c),\\
    T_{\leftarrow\rightarrow\leftarrow}(i_w\circ B_u\circ i_s(c), c)=L^-_r\circ A_{ul}\circ L^+_l\circ A_{dl}(c),\\
    T_{\leftarrow\leftarrow\leftarrow}(c, i_s\circ B_r\circ i_w(c))=L^-_r\circ A_{ul}\circ L^+_l\circ A_{ur}(c),\\
    T_{\leftarrow\leftarrow\leftarrow}(i_s\circ B_r\circ i_w(c), c)=L^-_r\circ A_{ul}\circ L^+_l\circ A_{dr}(c)
\end{gather*}
for any $c\in C_+$.
\begin{figure}[h]
    \centering
    \includegraphics[width=0.6\textwidth]{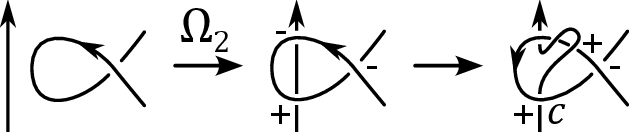}
    \caption{Outer $\Omega_2$-move in a loop}
    \label{pic:crossoid_R2_outer_loop}
\end{figure}

    \item with the bigon maps (Fig.~\ref{pic:crossoid_R2_inner_bigon})
\begin{eqnarray*}
    i_s=B_l\circ i_s\circ B_l,\quad i_s=B_r\circ i_s\circ B_r,\\
    i_w=B_u\circ i_w\circ B_u,\quad i_w=B_d\circ i_w\circ B_d.
\end{eqnarray*}
\begin{figure}[h]
    \centering
    \includegraphics[width=0.45\textwidth]{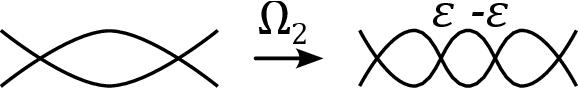}
    \caption{Inner $\Omega_2$-move in a bigon}
    \label{pic:crossoid_R2_inner_bigon}
\end{figure}

\end{itemize}

    \item rotation relations
\begin{gather*}
i_s^2=i_w^2=B_uB_d=B_dB_u=B_r^2=B_l^2=\mathrm{id}_C,\\
S_{o_1o_2}\circ T_{o_1o_2o_0}(c_2,S_{o_2o_0}\circ T_{o_0o_1o_2}(c_1,c_2))=c_1
\end{gather*}
where $S_{oo'}$ is the sign change map defined by the formulas
\begin{equation}\label{eq:crossoid_sign_change_map}
S_{\rightarrow\rightarrow}=i_s\circ B_l,\ S_{\leftarrow\rightarrow}=i_w\circ B_u,\ S_{\rightarrow\leftarrow}=i_w\circ B_d,\ S_{\leftarrow\leftarrow}=i_s\circ B_r.
\end{equation}

    \item $\Omega_3$-relation
    
For any $y,z\in C_+$ such that $A_{ul}(y)=A_{dl}(z)$ denote $x=T_{\rightarrow\rightarrow\rightarrow}(y,z)$ and
\begin{gather*}
    x'=B_l\circ T_{\rightarrow\leftarrow\rightarrow}(z,y),\quad
    y'=i_s\circ  T_{\rightarrow\leftarrow\rightarrow}(i_w\circ B_d(x),z),\\  
    z'=i_s\circ T_{\rightarrow\leftarrow\rightarrow}(y,i_w\circ B_u(x)). 
\end{gather*}
Then $A_{dr}(y')=A_{ur}(z')$ and $x'=T_{\leftarrow\leftarrow\leftarrow}(y',z')$.

    \item flip relation (Fig.~\ref{pic:crossoid_flip_move})
\[
T_{o_0o_2o_3}(i_{o_1o_3}\circ T_{o_0o_1o_2}(c_1,c_2),c_3)=T_{o_0o_1o_3}(c_1,i_{o_0o_2}\circ T_{o_1o_2o_3}(c_2,c_3)),
\]
where 
\begin{equation}\label{eq:crossoid_nonalternating_bigon_orientation}
i_{oo'}=\left\{\begin{array}{cl}
    i_s, & o=o', \\
    i_w, & o\ne o'.
\end{array}\right.
\end{equation}
\begin{figure}[h]
    \centering
    \includegraphics[width=0.45\textwidth]{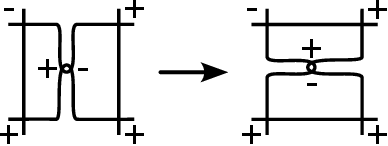}
    \caption{Flip move}
    \label{pic:crossoid_flip_move}
\end{figure}

\end{enumerate}

The crossoid is called \emph{unframed} if $i_w=i_s$.
\end{definition}

\begin{proposition}
    Definitions~\ref{def:crossoid} and~\ref{def:crossoid_reduced} are equivalent.
\end{proposition}

\begin{proof}
    1. Let $(C,A)$ be a crossoid in the sense of Definition~\ref{def:crossoid}. Consider the maps
\begin{gather*}
    L^\epsilon_l=L^\epsilon_\rightarrow,\quad L^\epsilon_r=L^\epsilon_\leftarrow,\quad i_s=P^{(-\epsilon,\epsilon)}_{(o,o)},\quad i_s=P^{(-\epsilon,\epsilon)}_{(o,\bar o)},\\
    B_l=P^{(\epsilon,\epsilon)}_{(\rightarrow,\rightarrow)},\quad B_r=P^{(\epsilon,\epsilon)}_{(\leftarrow,\leftarrow)},\quad 
    B_u=P^{(\epsilon,\epsilon)}_{(\leftarrow,\rightarrow)},\quad B_d=P^{(\epsilon,\epsilon)}_{(\rightarrow,\leftarrow)},\\
    T_{o_0o_1o_2}=P^{(-,+,+)}_{(o_0,o_1,o_2)},
\end{gather*}
$\epsilon\in\{+,-\}$, $o,o_0,o_1,o_2\in\{\leftarrow,\rightarrow\}$. A direct check shows that the relations in Definition~\ref{def:crossoid_reduced} are special cases of the relations in Definition~\ref{def:crossoid}.

2. Let $(C,A)$ be a crossoid in the sense of Definition~\ref{def:crossoid_reduced}. Using the loop and bigon maps, one can define the polygonal maps for loops and bigons:
\begin{gather*}
    P^\epsilon_\rightarrow=L^\epsilon_l\circ A_{ur}\circ L^{-\epsilon}_r,\quad P^\epsilon_\leftarrow=L^\epsilon_r\circ A_{ul}\circ L^{-\epsilon}_l,\quad
    P^{(-\epsilon,\epsilon)}_{(o,o)}=i_s,\quad P^{(-\epsilon,\epsilon)}_{(o,\bar o)}=i_w,\\
    P^{(\epsilon,\epsilon)}_{(\rightarrow,\rightarrow)}=B_l,\quad P^{(\epsilon,\epsilon)}_{(\leftarrow,\leftarrow)}=B_r,\quad 
    P^{(\epsilon,\epsilon)}_{(\leftarrow,\rightarrow)}=B_u,\quad P^{(\epsilon,\epsilon)}_{(\rightarrow,\leftarrow)}=B_d,
\end{gather*}
$\epsilon\in\{+,-\}$, $o\in\{\leftarrow,\rightarrow\}$.

Let us define the polygonal map for $n$-gons, $n\ge 3$. Using the sign change maps~\eqref{eq:crossoid_sign_change_map}, we reduce a general polygonal map to that with special signs:
\begin{multline*}
P^{(\epsilon_0,\dots,\epsilon_{n-1})}_{(o_0,\dots,o_{n-1})}(c_1,\dots,c_{n-1})=\\
S^{\frac{1+\epsilon_0}2}_{o_{n-1}o_0}\circ P^{(-,+,\dots,+)}_{(o_0,\dots,o_{n-1})}\left(S^{\frac{1-\epsilon_1}2}_{o_0o_1}(c_1),\dots,S^{\frac{1-\epsilon_{n-1}}2}_{o_{n-2}o_{n-1}}(c_{n-1})\right).
\end{multline*}

The maps $P^{(-,+,\dots,+)}_{(o_0,\dots,o_{n-1})}$ are defined by the induction on $n$:
\[
P^{(-,+,+)}_{(o_0,o_1,o_2)}=T_{o_0o_1o_2},
\]
and
\begin{multline*}
P^{(-,+,\dots,+)}_{(o_0,\dots,o_{n-1})}(c_1,\dots,c_n)=\\
P^{(-,+,\dots,+)}_{(o_0,\dots,o_i,o_j,\dots,o_{n-1})}\left(c_1,\dots,c_i,i_{o_io_j}\circ P^{(-,+,\dots,+)}_{(o_i,\dots,o_j)}(c_{i+1},\dots,c_j),c_{j+1},\dots,c_{n-1}\right)
\end{multline*}
for any $0\le i<j\le n-1$ such that $j-i\ne 1,n-1$. 

We need to show that the result of reduction to triangle maps does not depend on the choice of indices $i,j$ in the formula above. Consider a polygon (Fig.~\ref{pic:crossoid_triangulation} top left) and the dual graph (Fig.~\ref{pic:crossoid_triangulation} top right) whose vertices corresponds to the edges of the polygon. Any reduction of a polygon map to triangle maps corresponds to a triangulation of the dual graph (Fig.~\ref{pic:crossoid_triangulation} bottom). Any two triangulations are connected by a sequence of flips. Thus, the flip relation implies that the polygon map is uniquely defined.
\begin{figure}[h]
    \centering
    \includegraphics[width=0.4\textwidth]{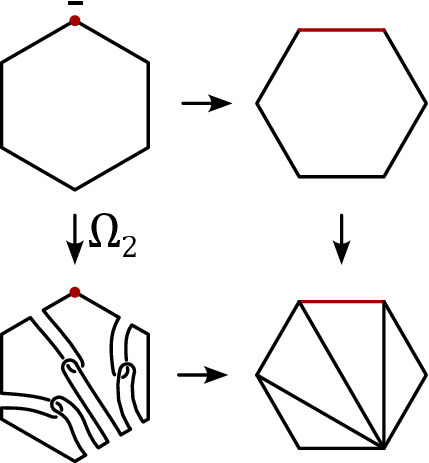}
    \caption{Reduction of a polygon map to triangle maps corresponds to a triangulation}
    \label{pic:crossoid_triangulation}
\end{figure}

Next, we must check that the polygonal maps obey the relations induced by Reidemeister moves. These relations for loops and bigons follow are the relations in Definition~\ref{def:crossoid_reduced}. The relations for $n$-gons, $n\ge 3$, can be reduced to those for bigons as shown in Fig.~\ref{pic:crossoid_relation_reduction}. Outer $\Omega_2$-relation is a part of definition of polygonal maps by reduction to triangle maps.
\begin{figure}[h]
    \centering
    \includegraphics[width=0.2\textwidth]{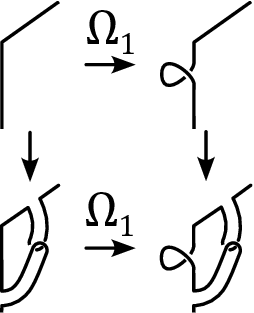}\qquad 
    \includegraphics[width=0.2\textwidth]{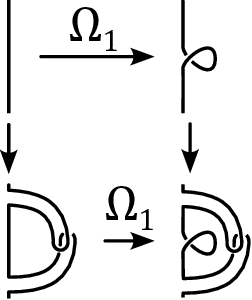}\qquad
    \includegraphics[width=0.2\textwidth]{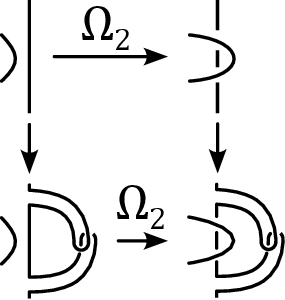}\qquad 
    \includegraphics[width=0.2\textwidth]{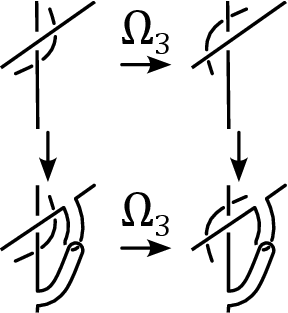}\qquad
    \caption{Reduction of $\Omega_1$-relations, outer $\Omega_2$-relation and $\Omega_3$-relation}
    \label{pic:crossoid_relation_reduction}
\end{figure}
\end{proof}

\begin{example}[Parity crossoid]\label{exa:crossoid_parity}
 Let $G$ be a group. Set $C=G\times\{-1,1\}$, and $A=\{1\}$. Then the incidence maps $A_{\alpha\beta}\colon C\to A$ are trivial. We set $sgn(g,\epsilon)=\epsilon$. Define the polygonal maps by the formulas $P^\epsilon_o=(1,\epsilon)$ and
\begin{equation}\label{eq:parity_crossoid_polygonal_map}
P^{(\epsilon_0,\dots,\epsilon_{n-1})}_{(o_0,\dots,o_{n-1})}((g_1,\epsilon_1),\dots,(g_{n-1},\epsilon_{n-1}))=\left(\left(\prod_{i=1}^{n-1}g_i^{e(o_{i-1},o_i)}\right)^{-e(o_{n-1},o_0)},\epsilon_0\right),
\end{equation}
where $e(o,o')$ is the incidence index defined by the formula
\[
e(o,o')=\left\{\begin{array}{cl}
    1, & o=o', \\
    -1, & o\ne o'.
\end{array}\right.
\]
A coloring of a diagram by the crossoid $(G\times\{-1,+1\},\{1\})$ is called an \emph{oriented parity with values in the group $G$} (cf.~\cite[Definition 1.10]{N2016}). 
\end{example}


\begin{example}[Biquandloid crossoid]\label{exa:crossoid_biquandloid}
    Let $(B,R)$ be a biquandloid. Consider the crossing set $C=(B\times_R B)\times\{-1,+1\}$ and the arc set $A=B$. Define the sign map by the formula $sgn(a,b,s)=s$ and the incidence maps as follows
\begin{gather*}
  A_{dr}(a,b,+)=a,\ A_{ur}(a,b,+)=b,\ A_{ul}(a,b,+)= a\uast b,\ A_{dl}(a,b,+)=b\oast a,\\
  A_{dr}(a,b,-)=a,\ A_{ur}(a,b,-)=b,\ A_{ul}(a,b,-)= a\oast b,\ A_{dl}(a,b,-)=b\uast a.
\end{gather*}
The polynomial maps are defined by the formulas $P^\epsilon_\rightarrow(a)=(k(a),k(a),\epsilon)$ where $k(a)\uast k(a)=a$, $P^\epsilon_\leftarrow(a)=(a,a,\epsilon)$. Let us define the map $P^{\bm\epsilon}_{\bm o}(c_1,\dots,c_{n-1})$, $\bm\epsilon=(\epsilon_0,\dots,\epsilon_{n-1})$, $\bm o=(o_0,\dots,o_{n-1})$. Denote $a=A_{in}(c_1,o_0,o_1)$ and $b=A_{out}(c_{n-1})$. Then the polynomial map is given by the table
\[
\begin{array}{|c|c|c|l|}
 \hline
\epsilon_0 & o_0 & o_{n-1} & P^{\bm\epsilon}_{\bm o}(c_1,\dots,c_{n-1})   \\
\hline
 + & \rightarrow & \rightarrow & (c,d,+)\mbox{ where }a=c\uast d, b=d\oast c  \\
 + & \leftarrow & \rightarrow & (b,c,+)\mbox{ where }a=c\oast b\\
 + & \rightarrow & \leftarrow & (d,a,+)\mbox{ where }b=d\uast a \\
 + & \leftarrow & \leftarrow & (a,b,+) \\
 - & \rightarrow & \rightarrow & (c,d,-)\mbox{ where }a=c\oast d, b=d\uast c \\
 - & \leftarrow & \rightarrow & (b,c,-)\mbox{ where }a=c\uast b\\
 - & \rightarrow & \leftarrow & (d,a,-)\mbox{ where }b=d\oast a \\
 - & \leftarrow & \leftarrow & (a,b,-) \\
\hline
\end{array}
\]
For a tangle diagram $D$, the map $(\chi_C,\chi_A)\mapsto \chi_A$ establishes a bijection between the set $Col_{(C,A)}(D)$ of crossoid colorings and the set $Col_{(B,R)}(D)$ of biquandloid colorings of the diagram.
\end{example}


\begin{example}[Topological crossoid of a tangle]\label{exa:crossoid_topological}
   Let $T$ be a tangle in $F\times I$. Consider the crossing set $C=\widetilde{\mathdutchcal C}_0(T)$ and the arc set $A=\widetilde{\mathdutchcal{SA}}_0(T)$. 

   The sign of a crossing $(\gamma^u,\gamma^m,\gamma^o)$ is negative if the translation of the tangent vector to the tangle along the framing of $\gamma^m$ is co-directed with the tangent to the tangle; otherwise the sign is positive (Fig.~\ref{pic:crossing_framing}).
   
   The incidence maps are defined by the formulas
\begin{gather*}
    A_{dl}(\gamma^u,\gamma^m,\gamma^o)=(\gamma^m
    (\mu^u_l)^{-1}\gamma^u,\gamma^o),\quad  A_{ul}(\gamma^u,\gamma^m,\gamma^o)=(\gamma^u,(\gamma^m)^{-1}
    \mu^o_l\gamma^o),\\
    A_{dr}(\gamma^u,\gamma^m,\gamma^o)=(\gamma^u,(\gamma^m)^{-1}
    (\mu^o_r)^{-1}\gamma^o),\quad  A_{ur}(\gamma^u,\gamma^m,\gamma^o)=(\gamma^m
    \mu^u_r\gamma^u,\gamma^o) 
\end{gather*}
for a positive crossing $(\gamma^u,\gamma^m,\gamma^o)$, and 
\begin{gather*}
    A_{dl}(\gamma^u,\gamma^m,\gamma^o)=(\gamma^u,(\gamma^m)^{-1}
    \mu^o_l\gamma^o),\quad  A_{ul}(\gamma^u,\gamma^m,\gamma^o)=(\gamma^m
    (\mu^u_l)^{-1}\gamma^u,\gamma^o),\\
    A_{dr}(\gamma^u,\gamma^m,\gamma^o)=(\gamma^m
    \mu^u_r\gamma^u,\gamma^o),\quad  A_{ur}(\gamma^u,\gamma^m,\gamma^o)=(\gamma^u,(\gamma^m)^{-1}
    (\mu^o_r)^{-1}\gamma^o) 
\end{gather*}
for a negative crossing $(\gamma^u,\gamma^m,\gamma^o)$.

The loop maps are given by the formulas (Fig.~\ref{pic:crossoid_top_loop_maps})
\begin{gather*}
    L^+_r(\gamma^u,\gamma^o)=(\gamma^u,(\mu_l)_{+1},\gamma^o)_{eo},\quad
    L^-_r(\gamma^u,\gamma^o)=(\gamma^u,(\mu_r)^{-1},\gamma^o)_{eu},\\
    L^+_l(\gamma^u,\gamma^o)=(\gamma^u,(\mu_r)_{+1}^{-1},\gamma^o)_{eu},\quad
    L^-_l(\gamma^u,\gamma^o)=(\gamma^u,\mu_l,\gamma^o)_{eo},
\end{gather*}
where $(\gamma)_{\epsilon}$ means changing the framing of the path $\gamma$ by $\epsilon$ half-twists; the subscript $eo$ (resp., $eu$) means that the overcrossing meridian is slightly shifted backward (resp., forward) along the orientation of the tangle compared to the undercrossing meridian. The terms $eo$ and $eu$ are abbreviations for early ovecrossing and early undercrossing.
\begin{figure}[h]
    \centering
    \includegraphics[width=0.8\textwidth]{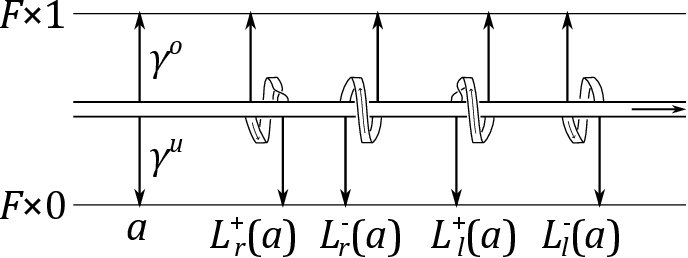}
    \caption{Loop maps, side view}
    \label{pic:crossoid_top_loop_maps}
\end{figure}

The formulas for non-alternating bigons are
\[
i_w(\gamma)=(\gamma^u,(\gamma^m)_{\sgn(\gamma)},\gamma^o),\quad i_s(\gamma)=(\gamma^u,(\gamma^m)_{-\sgn(\gamma)},\gamma^o),
\]
where $\gamma=(\gamma^u,\gamma^m,\gamma^o)\in\widetilde{\mathdutchcal C}_0(T)$ (Fig.~\ref{pic:bigon_nonalternating_top_view},\ref{pic:bigon_nonalternating_side_view}). 
\begin{figure}[h]
    \centering
    \includegraphics[width=0.7\textwidth]{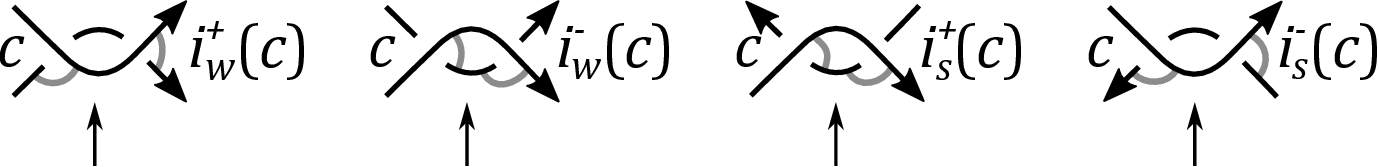}
    \caption{Non-alternating bigons, top view}
    \label{pic:bigon_nonalternating_top_view}
\end{figure}
\begin{figure}[h]
    \centering
    \includegraphics[width=0.7\textwidth]{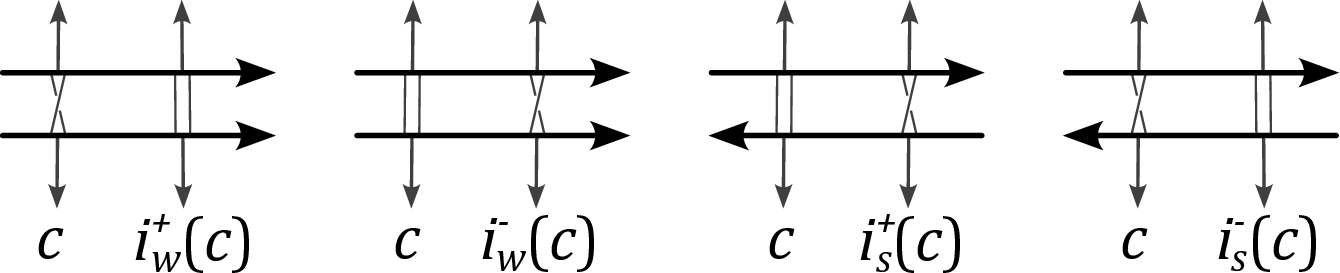}
    \caption{Non-alternating bigons, side view}
    \label{pic:bigon_nonalternating_side_view}
\end{figure}

The alternating bigons maps are described by the formulas
\begin{gather*}
    B_l(\gamma)=(\gamma^m(\mu^u_l)^{-1}\gamma^u,\mu^u_l(\gamma^m)^{-1}\mu_l^o, (\gamma^m)^{-1}\mu_l^o\gamma^o),\\
    B_r(\gamma)=(\gamma^m\mu^u_r\gamma^u,(\mu^u_r)^{-1}(\gamma^m)^{-1}(\mu_r^o)^{-1}, (\gamma^m)^{-1}(\mu_r^o)^{-1}\gamma^o),\\
    B_u^+(\gamma)=B_d^-(\gamma)=
    (\gamma^m\mu^u_r\gamma^u,(\mu^u_r)^{-1}(\gamma^m)^{-1}\mu_l^o, (\gamma^m)^{-1}\mu_l^o\gamma^o),\\
    B_u^-(\gamma)=B_d^+(\gamma)=(\gamma^m(\mu^u_l)^{-1}\gamma^u,\mu^u_l(\gamma^m)^{-1}(\mu_r^o)^{-1}, (\gamma^m)^{-1}(\mu_r^o)^{-1}\gamma^o),
\end{gather*}
    where $B^\epsilon_\alpha=B_\alpha|_{C_\epsilon}$, $\alpha\in\{u,d\}$, $\epsilon\in\{+,-\}$, is the restriction of the bigon map to the crossings of the given sign (Fig.~\ref{pic:bigon_alternating_top_view},\ref{pic:bigon_alternating_side_view}). 
\begin{figure}[h]
    \centering
    \includegraphics[width=0.7\textwidth]{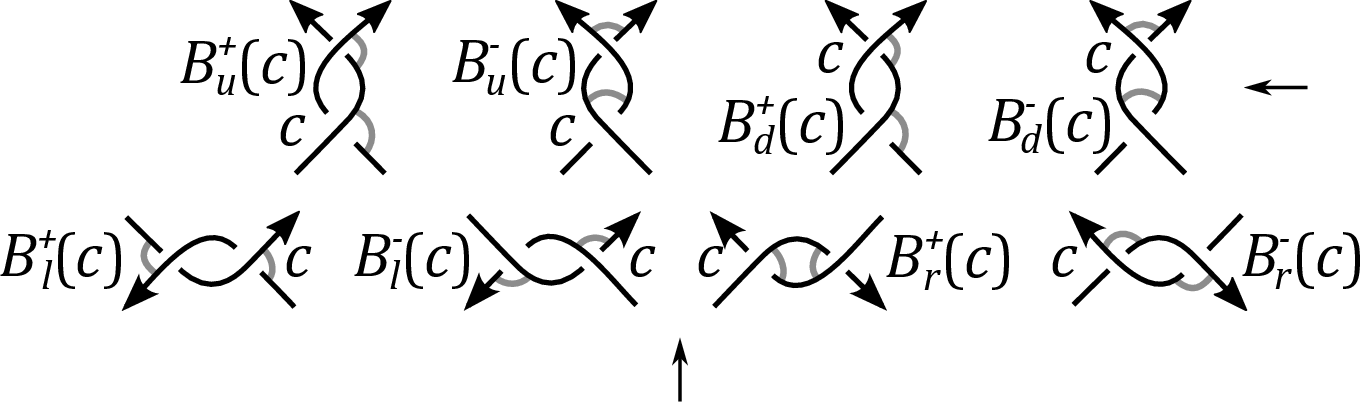}
    \caption{Alternating bigons, top view}
    \label{pic:bigon_alternating_top_view}
\end{figure}
\begin{figure}[h]
    \centering
    \includegraphics[width=0.7\textwidth]{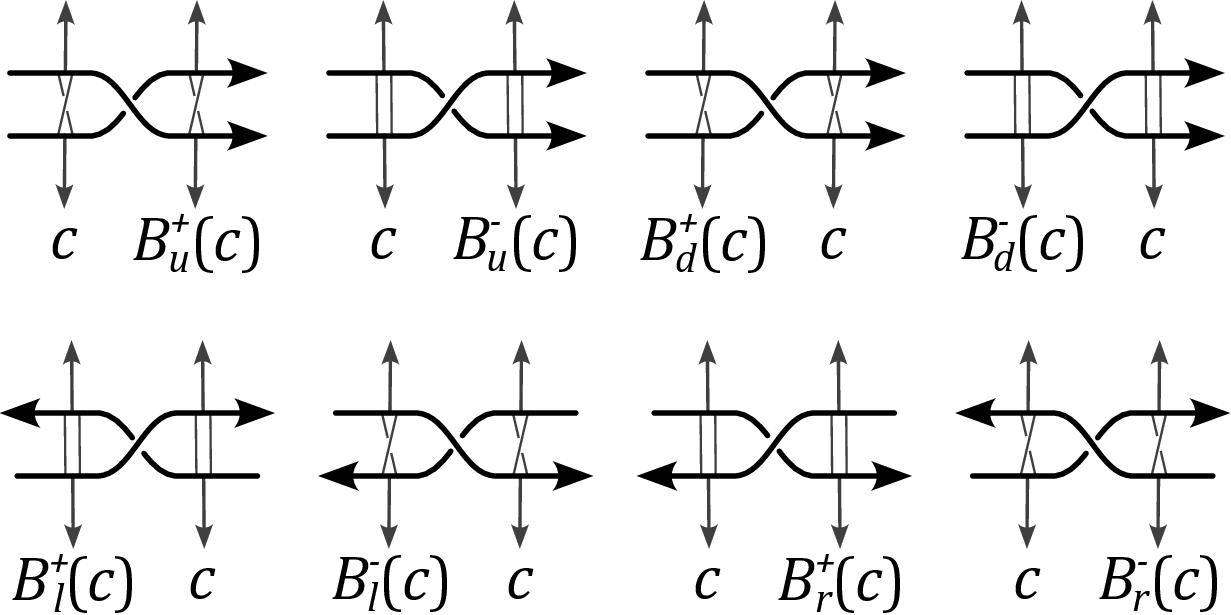}
    \caption{Alternating bigons, side view}
    \label{pic:bigon_alternating_side_view}
\end{figure}

In order to define the triangle maps, consider two positive crossings $\gamma_1=(\gamma^u_1,\gamma^m_1,\gamma^o_1)$ and $\gamma_2=(\gamma^u_2,\gamma^m_2,\gamma^o_2)$ in $\widetilde{\mathdutchcal C}_0(T)$ that are incident to an arc $a\in\widetilde{\mathdutchcal {SA}}_0(T)$, i.e. $A_{\alpha_1}(\gamma_1)=A_{\alpha_2}(\gamma_2)=a$ for some $\alpha_1,\alpha_2\in\{dl,ul,dr,ur\}$. We homotop $\gamma_1$ and $\gamma_2$ in order to match the meridians of these crossings that correspond to the arc $a$. We require that the underprobes $A_{\alpha_1}(\gamma_1)^u$ and $A_{\alpha_2}(\gamma_2)^u$ (as well as the overprobes $A_{\alpha_1}(\gamma_1)^o$ and $A_{\alpha_2}(\gamma_2)^o$) are homotopic as curves with fixed ends. Then the triangle maps are defined as follows (Fig.~\ref{pic:trigon_side_view})
\begin{gather*}
    T_{\rightarrow\rightarrow\rightarrow}(\gamma_1,\gamma_2)=(\gamma_2^u,\gamma_1^m(\mu^u_{1l})^{-1}\gamma^m_2,\gamma_1^o),\\
    T_{\leftarrow\rightarrow\rightarrow}(\gamma_1,\gamma_2)=(\gamma_1^u,\mu^u_{2l}(\gamma_2^m)^{-1}\cdot_{-}\gamma_1^m,(\gamma_2^m)^{-1}\mu^o_{2l}\gamma_2^o),\\
    T_{\rightarrow\leftarrow\rightarrow}(\gamma_1,\gamma_2)=(\gamma_2^m\mu^u_{2r}\gamma_2^u,\mu^u_{1l}(\gamma_1^m)^{-1}(\mu^o_{1r})^{-1}(\gamma_2^m)^{-1}(\mu^o_{2r})^{-1},(\gamma^m_1)^{-1}(\mu^o_{1r})^{-1}\gamma_1^o),\\
    T_{\rightarrow\rightarrow\leftarrow}(\gamma_1,\gamma_2)=(\gamma^m_1(\mu^u_{1l})^{_1}\gamma_1^u,\gamma^m_2\cdot_-(\gamma^m_1)^{-1}\mu^o_{1l},\gamma_2^o),\\
    T_{\rightarrow\leftarrow\leftarrow}(\gamma_1,\gamma_2)=(\gamma_1^u,(\mu^u_{2r})^{-1}(\gamma_2^m)^{-1}\cdot_{+}\gamma_1^m,(\gamma_2^m)^{-1}(\mu^o_{2r})^{-1}\gamma_2^o),\\
    T_{\leftarrow\rightarrow\leftarrow}(\gamma_1,\gamma_2)=(\gamma_2^m(\mu^u_{2l})^{-1}\gamma_2^u,(\mu^u_{1r})^{-1}(\gamma_1^m)^{-1}\mu^o_{1l}(\gamma_2^m)^{-1}(\mu^o_{2r})^{-1},(\gamma^m_1)^{-1}\mu^o_{1l}\gamma_1^o),\\
    T_{\leftarrow\leftarrow\rightarrow}(\gamma_1,\gamma_2)=(\gamma^m_1\mu^u_{1r}\gamma_1^u,\gamma^m_2\cdot_+(\gamma^m_1)^{-1}(\mu^o_{1r})^{-1},\gamma_2^o),\\
    T_{\leftarrow\leftarrow\leftarrow}(\gamma_1,\gamma_2)=(\gamma_2^u,\gamma_1^m\mu^u_{1r}\gamma^m_2,\gamma_1^o),
\end{gather*}
where $\cdot_-$ means smoothing against the orientation of the meridian, and $\cdot_+$ means smoothing along the orientation of the meridian (Fig.~\ref{pic:framing_smoothing}).
\begin{figure}[h]
    \centering
    \includegraphics[width=0.95\textwidth]{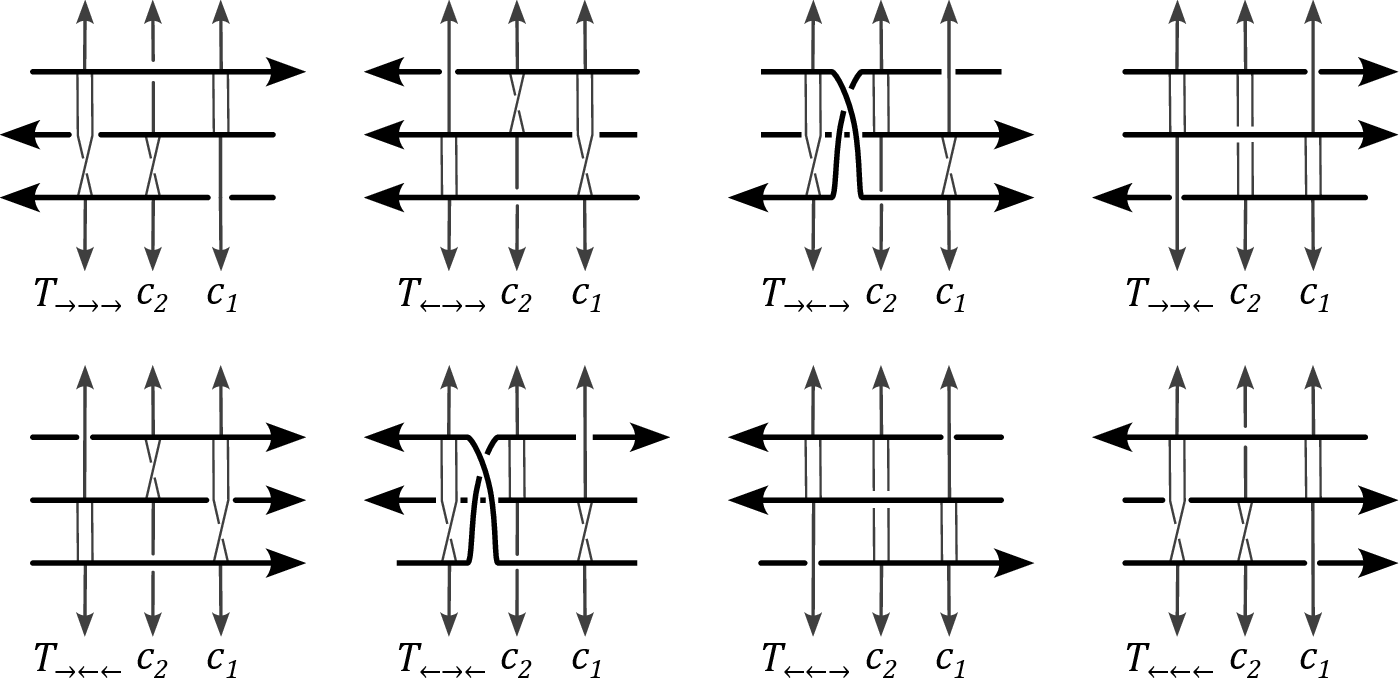}
    \caption{Triangle maps, side view}
    \label{pic:trigon_side_view}
\end{figure}
\begin{figure}[h]
    \centering
    \includegraphics[width=0.4\textwidth]{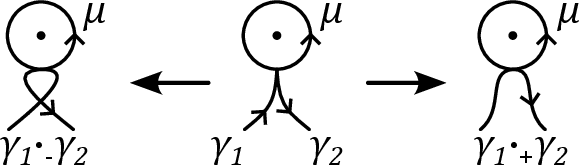}
    \caption{Smoothing of a concatenation of paths. The tangle is oriented towards the reader.}
    \label{pic:framing_smoothing}
\end{figure}

\end{example}

\begin{example}[Unframed topological crossoid of a tangle]\label{exa:crossoid_topological_unframed}
   Let $T$ be a tangle in $F\times I$. Consider the crossing set $C=\widetilde{\mathdutchcal C}^{uf}_0(T)\times\{-1,1\}$ and the arc set $A=\widetilde{\mathdutchcal{SA}}_0(T)$. 

   The sign map is defined by the formula $sgn(\gamma,\epsilon)=\epsilon$. The formulas for the incidence maps, loop maps, bigon maps, and triangle maps are obtained from those in Example~\ref{exa:crossoid_topological} by forgetting the framing.
\end{example}

The topological crossoid of a tangle is fundamental in the following sense.

\begin{theorem}\label{thm:crossoid_topological_universal}
Let $T\subset F\times I$ be a tangle, and $D$ its diagram. For any crossoid $(C,A)$ there is a bijection between the set of colorings $Col_{(C,A)}(D)$ and the set of invariant crossoid homomorphisms $Hom((\widetilde{\mathdutchcal{C}}_0(T),\widetilde{\mathdutchcal{SA}}_0(T)),(C,A))^{\pi_1(F)}$ under the action of $\pi_1(F,x_0)$ on $\widetilde{\mathdutchcal{C}}_0(T)$.
\end{theorem}

\begin{proof}
1. Let $\phi\colon (\widetilde{\mathdutchcal{C}}_0(T), \widetilde{\mathdutchcal{SA}}_0(T)))\to (C,A)$ be a $\pi_1(F,x_0)$-invariant homomorphism of crossoids. For a crossing $c\in\mathcal{C}(D)$, consider the vertical crossing probe $\gamma_c=(\gamma_c^u,\gamma_c^m,\gamma_c^o)$. Choose an arbitrary path $\delta_x\subset F$ from $x=p(c)$ to $x_0$. Then $\gamma_x=(\gamma_c^u(\delta_x\times 0),\gamma_c^m(\delta_x\times 0),\gamma_c^o(\delta_x\times 1))\in\widetilde{\mathdutchcal{SA}}_0(T)$. Define the color $\chi_{\phi,C}(c)$ of the crossing $c$ by the formula $\chi_{\phi,C}(c)=\phi(\gamma_x)$. The element $\chi_{\phi,C}(c)$ does not depend on $\delta_x$ by $\pi_1(F,x_0)$-invariance of $\phi$. Analogously, one defines the coloring $\chi_{\phi,A}$ of the semiarcs by elements of $A$.

Let us check that the map $\chi_\phi=(\chi_{\phi,C},\chi_{\phi,A})$ is a crossoid coloring. Let  $c_0,c_1,c_2\in\mathcal C(D)$ form a triangle in the diagram $D$. Let $\gamma_{c_0}$, $\gamma_{c_1}$, $\gamma_{c_2}$ be their vertical probes. Choose a path $\delta$ in $F$ from a point in the triangle to $x_0$. Using paths close to $\delta$, construct crossing probes $\gamma_0,\gamma_1,\gamma_2\in\widetilde{\mathdutchcal{C}}_0(T)$. Then $\gamma_0=T_{o_0o_1o_2}(\gamma_1,\gamma_2)$ for the appropriate orientations $o_0,o_1,o_2$. Hence,
\begin{multline*}
\chi_{\phi,C}(c_0)=\phi(\gamma_0)=\phi(T_{o_0o_1o_2}(\gamma_1,\gamma_2))=T_{o_0o_1o_2}(\phi(\gamma_1),\phi(\gamma_2))=\\
T_{o_0o_1o_2}(\chi_{\phi,C}(c_1),\chi_{\phi,C}(c_2)).
\end{multline*}
The other equalities are checked analogously.

2. Let $\chi\in Col_{(B,R)}(D)$. We construct a crossoid homomorphism $\phi=\phi_\chi$ from $\widetilde{\mathdutchcal{C}}_0(T)$ to $C$.

Let $D\cup\gamma$ a crossing probe diagram. Construct another crossing probe diagram $D'\cup\gamma$ by pulling the arcs of $D$ overcrossing $\gamma$ to the sinker point of $\gamma$, and pulling the arcs of $D$ undercrossing $\gamma$ to the float point of $\gamma$, and pulling the undercrossing arc of the probe to the overcrossing arc so that they form a pair of crossings (Fig.~\ref{pic:crossoid_probe_clearing}). Mark one of those two crossings that has the sign of the probe.
\begin{figure}[h]
\centering\includegraphics[width=0.8\textwidth]{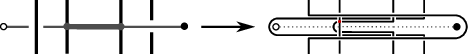}
\caption{Clearing a crossing probe}\label{pic:crossoid_probe_clearing}
\end{figure}

The transformation of the diagram $D$ is a morphism $f\colon D\to D'$. This morphism induces a bijection $f_*\colon Col_{(C,A)}(D)\to Col_{(C,A)}(D')$. Then we set the values $\phi(\gamma)$ equal to the color $(f_*(\chi))(c')$ of the marked crossing $c'$ in the diagram $D'$.

Let us check that the value $\phi(\gamma)$ does not change during isotopy of $\gamma$. If $f\colon D\to D_1$ is a second or a third Reidemeister move, then the corresponding transformed diagrams $D'$ and $D'_1$ differ by a sequence of Reidemeister moves of the same type that do not involve the marked crossing of $\gamma$. Then the color of the marked crossing does not change during the transformation from $D'$ to $D'_1$. 

A second or a third Reidemeister move $f\colon D\to D_1$ including an arc of the diagram $D$, induces diagrams $D'$ and $D'_1$ connected by an isotopy (and second and third Reidemeister moves if the probe has self-intersections).

If $f$ is a first Reidemeister move on $\gamma^u$ or $\gamma^o$, then we can move the loop to the sinker or to the float point of $\gamma$ using second and third Reidemeister moves. For a first Reidemeister move near the sinker or float point, the corresponding diagrams $D'$ and $D'_1$ are connected by second Reidemeister moves (Fig.~\ref{pic:crossoid_R1_enveloped}). 
\begin{figure}[h]
    \centering
\includegraphics[width=0.45\textwidth]{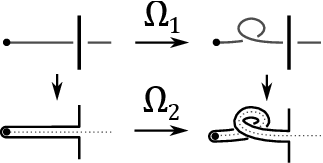}
    \caption{First Reidemeister moves near the sinker point}
    \label{pic:crossoid_R1_enveloped}
\end{figure}


If $f\colon D\to D_1$ is a float or sinker move (Fig.~\ref{pic:graphoid_float_sinker_moves}), then the corresponding transformed diagrams $D'$ and $D'_1$ are isotopic.

Let $f\colon D\to D_1$ be a vertex rotation move (Fig.~\ref{pic:crossoid_framed_intersection_rotation} top line). Consider the corresponding transformed diagrams $D'$ and $D'_1$ (Fig.~\ref{pic:crossoid_framed_intersection_rotation} middle line). Using Reidemeister moves outside the marked crossing, one can simplify the diagram $D'_1$ to the diagram in Fig.~\ref{pic:crossoid_framed_intersection_rotation} bottom right. Using Lemma~\ref{lem:crossoid_reciprocial_crossing} below, one can move a pair of counter-directed arcs below the marked crossing as shown in Fig.~\ref{pic:crossoid_framed_intersection_rotation} bottom line. The diagram obtained is isotopic to $D'$.
\begin{figure}[h]
    \centering
    \includegraphics[width=0.9\textwidth]{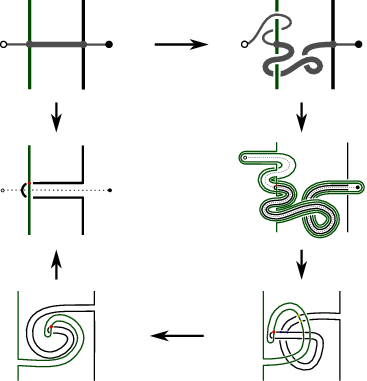}
    \caption{Framed rotation move}
    \label{pic:crossoid_framed_intersection_rotation}
\end{figure}

Thus, we get a map $\phi\colon\mathscr{SA}^s(T)\to X$ from the isotopy classes of semiarc probes to the quandle.

Then let us show that self-intersections of semiarc probes do not change their colors. We need the following lemma.

\begin{lemma}\label{lem:crossoid_reciprocial_crossing}
1. Let crossings $c_0, c_1, c_2, c_3\in C$ are the crossoid colors of a quadrilateral in a tangle (Fig.~\ref{pic:crossoid_tetragon}) that satisfies the conditions $\epsilon_1=-\epsilon_2$, $\epsilon_0=-\epsilon_3$ on the signs of the crossings, and the condition $o_0=o_2$ on the orientaions. Assume that $c_2=i_\alpha(c_1)$, $\alpha\in\{w,s\}$. Then $c_0=i_\beta(c_3)$, $\beta\in\{w,s\}$, where $\beta=\alpha$ if $o_1\ne o_3$, and $\beta\ne\alpha$ if $o_1=o_3$.
\begin{figure}[h]
    \centering
    \includegraphics[width=0.15\textwidth]{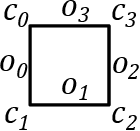}
    \caption{A square tangle}
    \label{pic:crossoid_tetragon}
\end{figure}

2. Consider a part of the tangle that consists of $2n$ vertical arcs $a_1,\dots, a_{2n}$ and two intersecting horizontal arcs such that the arcs $a_i$ and $a_{2n+1-i}$ have opposite orientation (Fig.~\ref{pic:crossoid_crossing_transfer} left). Assume that the arcs are colored by the crossoid $(C,A)$ so that the $c_{1,2n+1-i}=i_{\alpha_i}(c_{1i})$, $\alpha_i\in\{w,s\}$. Then after applying third Reidemeister moves, we get a colored part of the tangle (Fig.~\ref{pic:crossoid_crossing_transfer} right) such that $c'=c$ and $c'_{2,2n+1-i}=i_{\alpha_i}(c'_{2i})$, $i=1,\dots,2n-1$. 
\begin{figure}[h]
    \centering
    \includegraphics[width=0.7\textwidth]{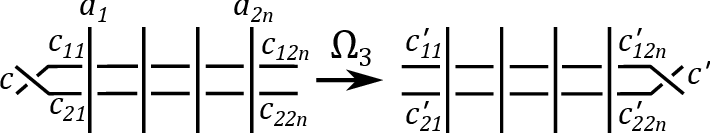}
    \caption{A colored tangle}
    \label{pic:crossoid_crossing_transfer}
\end{figure}
Analogous statement holds in the case of overcrossing horizontal arcs.
\end{lemma}

\begin{proof}
    1. Assume that $o_0=o_1=o_2=o_3=\rightarrow$. If $c_2=i_w(c_1)=P_{(\leftarrow,\rightarrow)}^{(\epsilon_2,\epsilon_1)}(c_1)$. Then by outer $\Omega_2$-relation
\[
c_0=P_{(\rightarrow,\rightarrow,\rightarrow,\rightarrow)}^{(\epsilon_0,\epsilon_1,\epsilon_2,\epsilon_3)}(c_1,c_2,c_3)=P_{(\rightarrow,\rightarrow)}^{(\epsilon_0,\epsilon_3)}(c_3)=i_s(c_3).
\]
If $c_2=i_s(c_1)=$ then by outer $\Omega_2$-relation
\begin{multline*}
c_2=P_{(\rightarrow,\rightarrow)}^{(\epsilon_2,\epsilon_1)}(c_1)=P_{(\rightarrow,\rightarrow,\rightarrow,\rightarrow)}^{(\epsilon_2,\epsilon_3,-\epsilon_3,\epsilon_1)}(c_3,P_{(\leftarrow,\rightarrow)}^{(-\epsilon_3,\epsilon_3)}(c_3),c_1)=\\
P_{(\rightarrow,\rightarrow,\rightarrow,\rightarrow)}^{(\epsilon_2,\epsilon_3,\epsilon_0,\epsilon_1)}(c_3, i_w(c_3),c_1).
\end{multline*}
Then by rotation relation
\[
i_w(c_3)=P_{(\rightarrow,\rightarrow,\rightarrow,\rightarrow)}^{(\epsilon_0,\epsilon_1,\epsilon_2,\epsilon_3)}(c_1,c_2,c_3)=c_0.
\]
Other cases of orientations of the quadrilateral sides are proved analogously.

2. Let us use induction on $n$. Let $n=1$. Assume that arcs are oriented as shown in Fig.~\ref{pic:crossoid_crossing_transfer_proof} left. Assume that $c_{12}=i_w(c_{11})$. Then by the first statement $c_{22}=i_s(c_{21})$. Then $A_{ul}(c_{21})=A_{dl}(c_{22})$ and $A_{ur}(c_{21})=A_{dr}(c_{22})$. Then a ``crossing'' $x=P_{(\rightarrow,\leftarrow,\rightarrow,\leftarrow)}^{(+,+,-,-)}(c_{22},c_{21},c)$ is defined. By outer $\Omega_2$-relation, $x=i_w(c)$. 
\begin{figure}[h]
    \centering
    \includegraphics[width=0.9\textwidth]{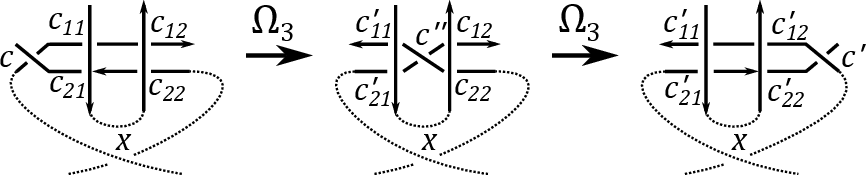}
    \caption{Case $n=1$}
    \label{pic:crossoid_crossing_transfer_proof}
\end{figure}

Since  $c_{22}=i_s(c_{21})$, the crossings with colors $c_{21}$ and $c_{22}$ form a ``bigon''. By applying twice the relation~\eqref{eq:crossoid_R3} to this ``bigon'', we get a ``bigon'' with colors $c'_{21}$ and $c'_{22}$. Hence, $c'_{22}=P_{(\leftarrow,\rightarrow)}^{(-,+)}(c'_{21})=i_w(c'_{21})$. 

Analogously, by applying~\eqref{eq:crossoid_R3} to the ``$4$-gon'' $x,c_{22},c_{21},c$, one gets ``$4$-gon'' $x,c',c'_{22},c'_{21}$. Then by outer $\Omega_2$-relation
$c'=i_w(x)=i_w^2(c)=c$.

For $n>1$, the induction step is performed by removing the crossings $c_{1,n}$, $c_{1,n+1}$, $c_{2,n}$, $c_{2,n+1}$ using outer $\Omega_2$-relations.
\end{proof}

The lemma implies that if a crosing probe passes a self-intersection point during the pulling process, then the crossing color does not change regardless of whether it passes above or below the intersection. Hence, if one switches the undercrossing and the overcrossing of a self-intersection of the probe, then the result does not change. Thus, the map $\phi$ induces a map from $\mathdutchcal{C}(T)$ to $X$, therefore, induces a $\pi_1(F)$-invariant map from $\phi\colon\widetilde{\mathdutchcal{C}}_0(T)\to X$.

The next claim is that the map $\phi$ is a natural transformation from the functor $\widetilde{\mathdutchcal{C}}$ to the constant functor $C$. It is enough to show that for any semiarc probe diagram $D\cup\gamma$ and for any Reidemeister move $f\colon D\to D_1$, the colors of the crossings $\gamma$ and $\mathdutchcal{C}(f)(\gamma)$ coincide.

Since the map $\phi$ gives the same value for isotopic crossing probes, we can assume that $\gamma$ is distinct from the area where the move occurs. Then the transformed diagrams $D'$ and $D'_1$ are connected by the same move as $D$ and $D_1$, and this move does not involve the marked crossing of the probe $\gamma$. Hence, the color of this crossing does not change with the move. Then the colors of the crossings $\gamma$ and $\mathdutchcal{C}(f)(\gamma)$ coincide.

Next, we define a map $\phi\colon\widetilde{\mathdutchcal{SA}}_0(T)\to A$. Given a semiarc probe $\gamma\in\widetilde{\mathdutchcal{SA}}_0(T)$, consider a crossing probe $\gamma_c\in\widetilde{\mathdutchcal{C}}_0(T)$ such that $\gamma=SA_{\alpha}(\gamma_c)$, $\alpha\in\{ul,dl,ur,dr\}$, and set $\phi(\gamma)=A_\alpha(\phi(\gamma_c))$. 

Let us check that the image does not depend on the choice of the crossing. Let $\gamma_{c_1}$ be another semiarc such that $\gamma=SA_\beta(\gamma_{c_1})$. Consider an isotopy $f\colon T\to T'$ which verticalize the probes $\gamma$, $\gamma_c$ and $\gamma_{c_1}$. Let $\chi'=f_*(\chi)$ be the corresponding coloring of the diagram $D'=p(T')$ by the crossoid $(C,A)$. Denote $\gamma'=\mathdutchcal{SA}(f)(\gamma)$, $\gamma_{c'}=\mathdutchcal{C}(f)(\gamma_c)$, and $\gamma_{c_1'}=\mathdutchcal{C}(f)(\gamma_{c'_1})$. Since the crossings $c'$ and $c_1'$ are incident to the same semiarc of $D'$, and $\chi'$ is a crossoid coloring, we get the equality $A_\alpha(\chi'(c'))=A_\beta(\chi'(c'_1))$. Then
\[
A_\alpha(\phi(\gamma_c))=A_\alpha(\phi(\gamma_{c'}))=A_\alpha(\chi'(c'))=A_\beta(\chi'(c'_1))=A_\beta(\phi(\gamma_{c_1'}))=A_\beta(\phi(\gamma_{c_1})).
\]

Finally, let us show that the map $\phi\colon\widetilde{\mathdutchcal{C}}_0(T)\to C$ is a homomorphism of crossoids. Let $\gamma_1,\gamma_2$ be crossing probes of the tangle $T$ such that $sgn(\gamma_1)=sgn(\gamma_2)=+1$ and $A_{ul}(\gamma_1)=A_{dl}(\gamma_2)$. Denote $\gamma_0=T_{\rightarrow\rightarrow\rightarrow}(\gamma_1,\gamma_2)$. We need to prove that  $\phi(\gamma_0)=T_{\rightarrow\rightarrow\rightarrow}(\phi(\gamma_1),\phi(\gamma_2))$.

Consider an isotopy $f\colon T\to T'$ which verticalize appropriate representatives of the probes $\gamma_1$, $\gamma_2$. Denote $\gamma'_i=\mathdutchcal{C}(f)(\gamma_i)$, $i=1,2$. Then $\gamma'_i$ are the vertical probes of the crossings $c'_i\in\mathcal{C}(D')$ of the diagram $D'=p(T')$. We can assume that the upper left semiarc of $c'_1$ is the lower left semiarc of $c'_2$. Assume also that the down-left semiarc of $c'_1$ and the uppel-left semiarc of $c'_2$ intersect in a negative crossing $c'_0$ (if not, then apply the second Reidemeister move). Denote the corresponding vertical probe of $c'_0$ by $\gamma'_0$. Then $\gamma'_0=T_{\rightarrow\rightarrow\rightarrow}(\gamma'_1,\gamma'_2)$, hence, $\gamma'_0=\mathdutchcal{C}(f)(\gamma_0)$.

Let $\chi'=f_*(\chi)$ be the crossoid coloring of $D'$ corresponding to $\chi$. Then
\begin{multline*}
T_{\rightarrow\rightarrow\rightarrow}(\phi(\gamma_1),\phi(\gamma_2))=T_{\rightarrow\rightarrow\rightarrow}(\phi(\gamma'_1),\phi(\gamma'_2))=T_{\rightarrow\rightarrow\rightarrow}(\chi'(c'_1),\chi(c'_2))=\\
\chi'(c'_0)=\phi(\gamma'_0)=\phi(\gamma_0).
\end{multline*}
In the third equality we used the fact that $\chi'$ is a crossoid coloring. 

The other relations are proved analogously.
The theorem is proved.
\end{proof}

\begin{corollary}\label{thm:crossoid_topological_unframed_universal}
Let $T\subset F\times I$ be a tangle, and $D$ its diagram. For any unframed crossoid $(C,A)$ there is a bijection between the set of colorings $Col_{(C,A)}(D)$ and the set of invariant crossoid homomorphisms 
\[
Hom((\widetilde{\mathdutchcal{C}}^{uf}_0(T)\times\{-1,+1\},\widetilde{\mathdutchcal{SA}}_0(T)),(C,A))^{\pi_1(F)}
\]
under the action of $\pi_1(F,x_0)$ on $\widetilde{\mathdutchcal{C}}^{uf}_0(T)$.  
\end{corollary}

\begin{proof}
    Since the maps $i_w$ and $i_s$ on a crossing probe $\gamma=(\gamma^u,\gamma^m,\gamma^o)\in\widetilde{\mathdutchcal{C}}_0(T)$ act by shifting framing of the midprobe:
\[
i_w(\gamma)=(\gamma^u,(\gamma^m)_{\sgn(\gamma)},\gamma^o),\quad i_s(\gamma)=(\gamma^u,(\gamma^m)_{-\sgn(\gamma)},\gamma^o),
\]
then the relation $i_w=i_s$ is equivalent to the equality 
\[
(\gamma^u,\gamma^m,\gamma^o)=(\gamma^u,(\gamma^m)_{+2k},\gamma^o),\ k\in\Z.
\]
Hence, the space $\widetilde{\mathdutchcal{C}}^{uf}_0(T)\times\{-1,+1\}$ is the quotient of the space $\widetilde{\mathdutchcal{C}}_0(T)$ by the relation $i_w=i_s$. Thus, any homomorphism $\phi$ from  $\widetilde{\mathdutchcal{C}}_0(T)$ to an unframed crossoid $(C,A)$ has a unique decomposition $\phi=\bar\phi\circ\pi$ where $\pi\colon\widetilde{\mathdutchcal{C}}_0(T)\to \widetilde{\mathdutchcal{C}}^{uf}_0(T)\times\{-1,+1\}$ is the natural projection. Then
\begin{multline*}
Col_{(C,A)}(D)\simeq Hom((\widetilde{\mathdutchcal{C}}_0(T),\widetilde{\mathdutchcal{SA}}_0(T)),(C,A))^{\pi_1(F)}\simeq\\ Hom((\widetilde{\mathdutchcal{C}}^{uf}_0(T)\times\{-1,+1\},\widetilde{\mathdutchcal{SA}}_0(T)),(C,A))^{\pi_1(F)}.
\end{multline*}
\end{proof}

\begin{definition}\label{def:crossoid_cocycle}
    A \emph{cocycle} on a crossoid $(C,A)$ is a map $\phi\colon C\to H$ to an abelian group $H$ such that:
\begin{enumerate}
    \item $\phi(P^\epsilon_o(a))=0$ for any $a\in A$, $\epsilon\in\{+,-\}$, $o\in\{\leftarrow,\rightarrow\}$;
    \item $\phi(P^{(-\epsilon,\epsilon)}_{(o,o')}(c))=\phi(c)$ for any $\epsilon\in\{+,-\}$, $o,o'\in\{\leftarrow,\rightarrow\}$ and $c\in C_\epsilon$;
    \item for any $y,z\in C_+$ and $o\in\{\leftarrow,\rightarrow\}$ such that $A_{out}(y,o,o)=A_{in}(z,o,o)$ denote $x=P_{(o,o,o)}^{(-,+,+)}(y,z)$ and
\begin{gather*}
    x'=P_{(o,o)}^{(-,+)}\circ P_{(o,\bar o,o)}^{(+,+,+)}(z,y),\quad
    y'=P_{(o,o)}^{(+,-)}\circ P_{(o,\bar o,o)}^{(-,-,+)}(x,z),\nonumber\\  
    z'=P_{(o,o)}^{(+,-)}\circ P_{(o,\bar o,o)}^{(-,+,-)}(y,x). 
\end{gather*}
Then 
\[
-\phi(x)+\phi(y)+\phi(z)=-\phi(x')+\phi(y')+\phi(z').
\]
\end{enumerate}
\end{definition}

\begin{example}\label{exa:crossoid_cocycle}
    1. Let $C=\Z_2\times\{-1,+1\}$ be the parity crossoid with values in $\Z_2$. Then the map $\phi\colon C\to\Z_2$, $\phi(x,\epsilon)=x$, is a crossoid cocycle.

    2. Let $(B,R)$ be a biquandloid and $\theta\in Z^2(B,H)$ a $2$-cocycle with coefficients in an abelian group $H$. Consider the biquandloid crossoid $C=(B\times_R B)\times\{-1,1\}$ from Example~\ref{exa:crossoid_biquandloid}. Then the map $\phi\colon C\to H$, 
\[
\phi(x,y,\epsilon)=\left\{\begin{array}{cl}
    \theta(x,y), & \epsilon=+1, \\
     -\theta (y,x), & \epsilon=-1,
\end{array}\right.
\]
is a crossoid cocycle.
\end{example}

\begin{proposition}\label{prop:crossoid_cocycle_invariant}
    Let $(C,A)$ be a crossoid and $\phi\colon C\to H$ a crossoid cocycle. Let $D$ be a tangle diagram in a surface $F$. Then the formula 
\[
I_\phi(D)=\sum_{\chi\in Col_{(C,A)}(D)}\left[\sum_{c\in\mathcal C(D)}sgn(c)\phi(\chi(c))\right]\in\Z[H]
\]
defines a tangle invariant.
\end{proposition}
\begin{proof}
    We need to check that the element $I_\phi(D)$ does not change under Reidemeister moves. But the invariance follows directly from the definitions of crossoid coloring and crossoid cocycle.
\end{proof}

\begin{definition}\label{def:crossoid_homology}
Let $(C,A)$ be an unframed crossoid and $H$ an abelian group. Consider the free $H$-module $C_n(C,H)$, $n\ge 1$, generated by the set
\[
\{(c_1,\dots,c_{n-1})\in C_+^{n-1}\mid A_{ur}(c_i)=A_{dr}(c_{i+1}),\ i=1,\dots,n-1\}
\]
when $n\ge 2$, and by $A$ when $n=1$. For $n\ge 3$, define a differential $d$ by the formula
\begin{multline*}
    d(c_1,\dots,c_{n-1})=(c_2,\dots,c_{n-1})+(-1)^{n-1}(c_1,\dots,c_{n-2})+\\
    \sum_{i=2}^{n-1}(-1)^{i-1}(c_1,\dots,c_{i-2},\iota P^{-++}_{\leftarrow\leftarrow\leftarrow}(c_i,c_{i-1}),c_{i+1},\dots,c_{n-1})+\\
\sum_{i=1}^n (-1)^i(c_{i,1},\dots, c_{i,i-2}, P^{+++}_{\leftarrow\rightarrow\leftarrow}(c_{i-1},c_i),c_{i,i+1},\dots,c_{i,n-1}),
\end{multline*}
where 
\[
c_{i,k}=\left\{\begin{array}{cl}
   \iota P^{-++}_{\leftarrow\rightarrow\rightarrow}(P^{+++}_{\leftarrow\rightarrow\leftarrow}(c_k,x_{i,k}),x_{i,k}), & k\le i-2, \\
    \iota P^{-++}_{\rightarrow\rightarrow\leftarrow}(x_{i,k},P^{+++}_{\leftarrow\rightarrow\leftarrow}(x_{i,k},c_k)), & k\ge i+1, 
\end{array}\right.
\]
and the elements $x_{i,k}$ are defined inductively
\[
x_{i,k}=\left\{\begin{array}{cl}
    \iota P^{-++}_{\leftarrow\leftarrow\leftarrow}(x_{i,k+1},c_{k+1}), & k<i-2, \\
    c_{i-1}, & k=i-2,\\
    c_i, & k=i+1,\\
    \iota P^{-++}_{\leftarrow\leftarrow\leftarrow}(c_{k-1},x_{i,k-1}), & k>i+1.
\end{array}\right.
\]
For $n=2$ the differential is defined by the formula
\[
d(c_1)=A_{ul}(c_1)+A_{ur}(c_1)-A_{dl}(c_1)-A_{dr}(c_1).
\]
Then $(C_*(C,H),d)$ is a chain commplex. It contains a subcomplex $D_n(C,H)$, $n\ge 2$, which is the free $H$-module generated by the set
\[
\{(c_1,\dots,c_{n-1})\mid \exists i\ \exists a\in A\ \exists o\in\{\leftarrow,\rightarrow\}: c_i=P^+_o(a)\}.
\]

The homology $H_*(C,H)$ of the quotient complex $C_*(C,H)/D_*(C,H)$ is called the \emph{crossoid homology} with the coefficients in $H$.
Analogously, one defines the \emph{crossoid cohomology} $H^*(C,H)$.
\end{definition}

The formulas for the differential of the crossoid complex are explained in the next subsection.

By definition, the set of crossoid cocycles is identified with the set of $2$-cocycles $Z^2(C,H)$ of the complex $C^*(C,H)$. 

\begin{remark}
1.  If $(C,A)$ is the biquandloid crossoid of a biquandloid $(B,R)$ then the maps
\[
((a_1,a_2,+),(a_2,a_3,+),\dots,(a_{n-1},a_n,+))\mapsto (a_1,a_2,\dots,a_n),\quad a_i\in B,
\]
establish an isomorphism between the complexes $C_*(C,H)$ and $C_*(B,H)$. Hence, the crossoid (co)homology and the biquandloid (co)homology are isomorphic.

2. If $(C,A)$ is a parity crossoid, then the differential $d$ in the complex $C_*(C,H)$ vanishes. Hence,
\[
H_*(C,H)=C_*(C,H)=H[(G\setminus 1)^{\times n-1}].
\]
\end{remark}

\subsection{Homotopical multicrossing complex}\label{subsect:multicrossing_complex}

Let us define a homology construction that is behind all those tribracket, biquandle and crossoid homologies.

Let $T$ be a tangle in the thickened surface $F\times I$. Choose a point $x_0\in F$ and denote $x^u=x_0\times 0$ and $x^o=x_0\times 1$.

For $n\ge 0$, consider the space $\widetilde{\mathdutchcal C}(T,n)$ that consists of homotopy classes of tuples $\gamma=(\gamma_0,\dots,\gamma_n)$ of paths
\begin{gather*}
\gamma_0\colon (I,0,1)\to (M_T, x^u,\partial N(T)),\\  
\gamma_i\colon (I,0,1)\to (M_T, \partial N(T),\partial N(T)),\ i=1,\dots,n-1,\\
\gamma_n\colon (I,0,1)\to (M_T, \partial N(T), x^o);
\end{gather*}
such that
\begin{itemize}
\item $\gamma_{i-1}(1)$ and $\gamma_{i}(0)$ are different points of the same meridian $\mu_i$, $1\le i\le n$. 

Denote the part of the meridian $\mu_i$ from $\gamma_{i-1}(1)$ to $\gamma_{i}(0)$ by $\mu_{i,l}$, and the part from $\gamma_{i}(0)$ to $\gamma_{i-1}(1)$ by $\mu_{i,r}$.
\item the meridians $\mu_i$, $i=1,\dots,n$, are all distinct;
\item the loop $p(\gamma_0\mu_{1,l}\gamma_1\cdots\gamma_{n-1}\mu_{n,l}\gamma_n)$ is trivial in $\pi_1(F,x)$

\end{itemize}

The sequence $(\widetilde{\mathdutchcal C}(T,n))_{n\ge 0}$ with the maps $\partial^n_{i,l},\partial^n_{i,l}\colon \widetilde{\mathdutchcal C}(T,n)\to\widetilde{\mathdutchcal C}(T,n-1)$, $i=1,\dots,n$, given by the formulas
\begin{gather*}
    \partial^n_{i,l}(\gamma_0,\dots,\gamma_n)=(\gamma_0,\dots,\gamma_{i-2},\gamma_{i-1}\mu_{i,l}\gamma_i,\gamma_{i+1},\dots,\gamma_n),\\
    \partial^n_{i,r}(\gamma_0,\dots,\gamma_n)=(\gamma_0,\dots,\gamma_{i-2},\gamma_{i-1}(\mu_{i,r})^{-1}\gamma_i,\gamma_{i+1},\dots,\gamma_n),
\end{gather*}
is a semi-cubical set.

\begin{remark}\label{rem:homotopical_crossing_space_012}
The space $\widetilde{\mathdutchcal C}(T,0)$ is identified with $\widetilde{\mathdutchcal R}_0(T)$, $\widetilde{\mathdutchcal C}(T,1)$ is identified with $\widetilde{\mathdutchcal{SA}}_0(T)$. and $\widetilde{\mathdutchcal C}(T,2)$ is identified with $\widetilde{\mathdutchcal C}^{uf}_0(T)$. Using this identification, the map $\partial^1_l$ is the incidence map $R_l$, and $\partial^1_r$ is $R_r$. For a positive crossing $c\in\widetilde{\mathdutchcal C}_0(T)$, the map $\partial^2_{1,l}(c)$ is $A_{dl}(c)$, $\partial^2_{1,r}(c)=A_{ur}(c)$, $\partial^2_{2,l}(c)=A_{ul}(c)$, $\partial^2_{2,r}(c)=A_{dr}(c)$. For a negative crossing, $\partial^2_{1,l}(c)=A_{ul}(c)$, $\partial^2_{1,r}(c)=A_{dr}(c)$, $\partial^2_{2,l}(c)=A_{dl}(c)$, $\partial^2_{2,r}(c)=A_{ur}(c)$.
\end{remark}

There is an action of $\pi_1(F,x_0)$ on $\widetilde{\mathdutchcal C}(T,n)$ given by the formula
\[
\alpha\cdot(\gamma_0,\dots,\gamma_n)=((\alpha\times 0)\gamma_0,\gamma_1\dots,\gamma_{n-1},\gamma_n(\alpha^{-1}\times 1)),\quad \alpha\in\pi_1(F,x_0).
\]
The maps $\partial^n_{i,l}, \partial^n_{i,r}$ commutes with this action.

For $n\ge 2$, consider a subset $\widetilde{\mathdutchcal D}(T,n)=\bigcup_{i=1}^{n-1}\widetilde{\mathdutchcal D}_i(T,n) \subset \widetilde{\mathdutchcal C}(T,n)$ where $\widetilde{\mathdutchcal D}_i(T,n)$ consists of elements $\gamma$ such that 
\[
\gamma_{\alpha_1\dots\alpha_{i-1}}^{\alpha_i\dots\alpha_{n-2}}=\partial^3_{3,\alpha_{n-2}}\circ\cdots\circ\partial^{n-i+1}_{3,\alpha_{i}}\circ\partial^{n-i+2}_{1,\alpha_{i-1}}\circ\cdots\partial^n_{1,\alpha_{1}}(\gamma),
\]
is a loop crossing for some $\alpha_j\in\{l,r\}$, $j=1,\dots,n-2$. Note that in this case $\gamma_{\alpha_1\dots\alpha_{i-1}}^{\alpha_i\dots\alpha_{n-2}}$ is a loop crossing for any $\alpha_1,\dots,\alpha_{n-2}\in\{l,r\}$. The sequence $(\widetilde{\mathdutchcal D}(T,n))_{n\ge 2}$ is a semicubical subset in $(\widetilde{\mathdutchcal D}(T,n))$.

For an abelian group $M$, consider the complex 
\[
    C_n(T,M)=\left(\widetilde{\mathdutchcal C}(T,n)\otimes M/\widetilde{\mathdutchcal D}(T,n)\otimes M\right)^{\pi_1(F)}
\]
with the differential 
\[
d(\gamma)=\sum_{i=1}^n(-1)^i(\partial^n_{i,l}(\gamma)-\partial^n_{i,r}(\gamma)),
\]
and consider the complex
\begin{multline*}
    C^n(T,M)=\{\phi\in Hom(\widetilde{\mathdutchcal C}(T,n),M)\mid \forall \gamma\in\widetilde{\mathdutchcal C}(T,n) \forall \gamma_1\in \widetilde{\mathdutchcal D}(T,n) \forall \alpha\in\pi_1(F)\\  \phi(\alpha\cdot\gamma)=\phi(\gamma)\ \&\ \phi(\gamma_1)=0\}
\end{multline*}
with the differential $(d\phi)(\gamma)=\phi(d\gamma)$.

\begin{definition}\label{def:homotopical_crossing_homology}
    The homology $H_*(T,M)=H(C_*(T,M))$ and $H^*(T,M)=H(C^*(T,M))$ are called the \emph{homotopical multicrossing homology} and \emph{homotopical multicrossing cohomology} of the tangle $T$ with coefficients in $M$.
\end{definition}

Let $D=p(T)$ be the diagram of the tangle $T$. Consider the chain
\[
z(D)=\sum_{c\in\mathcal C(D)}sgn(c)\gamma_c\in C_2(T,\Z).
\]

\begin{proposition}\label{prop:homotopy_crossing_class}
\begin{enumerate}
    \item the chain $z(D)$ is a cycle $d(z(D))=0$;
    \item the homology class $cr(D)=[z(D)]\in H_1(T,\Z)$ is invariant, i.e. for any isotopy $f\colon T\to T'$ we have $f_*(cr(D))=cr(D')$ where $D=p(T')$ is the diagram of the tangle $T$ and $f_*\colon H_*(T,\Z)\to H_*(T',\Z)$ is the homomorphism induced by the isotopy $f$.
\end{enumerate} 
\end{proposition}
\begin{proof}
   By Remark~\ref{rem:homotopical_crossing_space_012} for a crossing $c\in\mathcal C(D)$, 
\[
sgn(c)d(\gamma_c)=A_{ul}(\gamma_c)+A_{ur}(\gamma_c)-A_{dl}(\gamma_c)-A_{dr}(\gamma_c)
\]
Since any semiarc $a\in\mathcal{SA}(D)$ starts at one crossing and ends at another, the total coefficient of this arc in the sum $d(z(D))$ is zero. Hence, $d(z(D))=0$.

Let $f\colon T\to T'$ be an isotopy that represents a Reidemeister move. If $f$ is an increasing first Reidemeister move and $\gamma$ is the probe of the appearing crossing, then
$z(D')=f_*(z(D))\pm\gamma=f_*(z(D))$ because $\gamma$ is a loop crossing and is equal to zero in $C_n(T',\Z)$.

If $f$ is an increasing second Reidemeister move and $\gamma_{1}$, $\gamma_{2}$ are the probes of the new crossings $c_1$, $c_2$, then $\gamma_1=\gamma_2\in\widetilde{\mathdutchcal C}(T,2)$ and $sgn(c_1)=-sgn(c_2)$. Hence,
\[
z(D')=f_*(z(D))+sgn(c_1)\gamma_1+sgn(c_2)\gamma_2=f_*(z(D)).
\]

If $f$ is a third Reidemeister move, then $z(D')=f_*(z(D))\pm d\omega$, where $\omega\in\widetilde{\mathdutchcal C}(T,3)$ is the probe of the triple crossing that appears during the move. Hence, $cr(D')=f_*(cr(D))$.
\end{proof}

\begin{definition}\label{def:homotopical_crossing_class}
    The element $cr(D)\in H_1(T,\Z)$ is called the \emph{homotopical crossing class} of the tangle $T$.
\end{definition}

\begin{remark}\label{rem:isotopical_crossing_class}
    We can consider the sets of isotopy classes of path sequences $\gamma=(\gamma_0,\dots,\gamma_n)$ instead of the homotopy classes. This will lead us to definitions of isotopical multicrossing homology and the isotopical crossing class of the tangle.
\end{remark}

\begin{remark}
    The cycle $z(D)$ that represents the crossing cycle, has a geometric description as the set of singular values of the projection $p\colon T\to D$.
\end{remark}

\begin{corollary}
    Let $\phi$ be a crossing cocycle of a crossing $(C,A)$ valued in an abelian group $H$. Then
\begin{enumerate}
    \item any coloring $\chi\in Col_{(C,A)}(D)$ induces a cocycle in $\phi_\chi\in C^2(T,H)$;
    \item the cocycle invariant $I_\phi$ is equal to
\[
I_\phi(D)=\sum_{\chi\in Col_{(C,A)}(D)}[\phi_\chi(cr(D))]\in\Z[H].
\]
\end{enumerate}
\end{corollary}
\begin{proof}
    The cocycle $\phi_\chi$ is the composition $\phi\circ\psi_\chi$ of the crossoid cycle and the crossoid homomorphism $\psi_\chi\colon \widetilde{\mathdutchcal C}_0(T)\to C$ induced by the coloring $\chi$.

    The second statement follows from the definition of $I_\phi$.
\end{proof}

Now, let us demonstrate how multicrossing homology relates to tribracket, biquandloid and crossoid homology.

\begin{figure}
    \centering
    \includegraphics[width=0.4\textwidth]{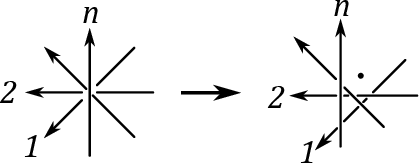}
    \caption{Splitting a multicrossing}
    \label{pic:multicrossing_resolution}
\end{figure}

Consider a vertical $n$-crossing probe such that the arcs of the probe rotate clockwise as we move from the bottom to the top (Fig.~\ref{pic:multicrossing_resolution} left). Shift the arcs to left to split the multicrossing (Fig.~\ref{pic:multicrossing_resolution} right).

1. Let $(X,[],\uparrow)$ be a partial ternary quasigroup.
Any coloring $\chi\in Col_{(X,\uparrow)}(D)$ of the diagram $D$ induces a homomorphism $\phi_\chi\colon\widetilde{\mathdutchcal R}_0(T)\to X$ from the topological partial ternary quasigroup. Then the maps $\phi\colon C_n(T,H)\to C_n(X,H)$,
\[
\gamma\mapsto (r_0,r_1,\dots,r_n),
\]
where $r_i=\phi_\chi\circ(\partial_{1,r})^{n-i}(\partial_{1,r})^i(\gamma)$,  $i=0,\dots,n$, define a chain map between the complexes.

\begin{figure}
    \centering
    \includegraphics[width=0.9\textwidth]{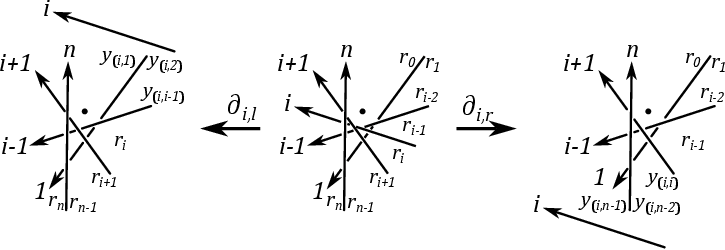}
    \caption{Geometric interpretation of the differential for a ternary quasigroup}
    \label{pic:multicrossing_tribracket}
\end{figure}

Let us give a geometric interpretation of this homomorphism. Assume that the multicrossing in (Fig.~\ref{pic:multicrossing_resolution}) is described by the $n+1$ region colors $r_n\uparrow r_{n-1}\uparrow\cdots\uparrow r_0$ as in Fig.~\ref{pic:multicrossing_tribracket} middle. The left and right cubic maps produce the tangles (Fig.~\ref{pic:multicrossing_tribracket} left and right) whose region colors correspond to the terms of the differential of the partial tribracket complex (see Definition~\ref{def:tribracket_partial_homology}).

2. Let $(B,R)$ be a biquandloid.
A coloring $\chi\in Col_{(B,R)}(D)$ of the diagram $D$ by the biquandle induces a homomorphism $\phi_\chi\colon\widetilde{\mathdutchcal{SA}}_0(T)\to B$ from the topological biquandloid. Then the maps $\phi\colon C_n(T,H)\to C_n(B,H)$,
\[
\gamma\mapsto (a_1,\dots,a_n),
\]
where $a_i=\phi_\chi\circ(\partial_{2,r})^{n-i}(\partial_{1,r})^{i-1}(\gamma)$,  $i=1,\dots,n$, define a chain map between the complexes.

\begin{figure}
    \centering
    \includegraphics[width=0.9\textwidth]{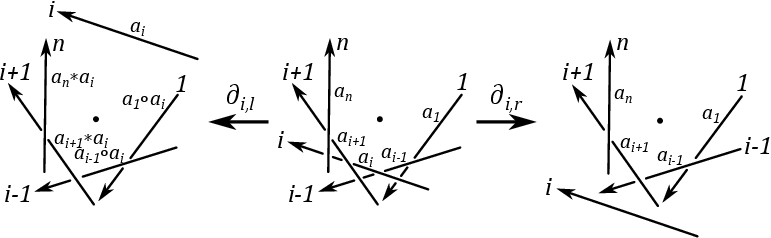}
    \caption{Geometric interpretation of the differential for a biquandloid}
    \label{pic:multicrossing_biquandloid}
\end{figure}

The geometric interpretation of this homomorphism is as follows. Assume that the multicrossing in (Fig.~\ref{pic:multicrossing_resolution}) is described by the $n$ semiarc colors $a_1,\dots, a_n$ such that $\sigma_r(a_1)=\cdots=\sigma_r(a_n)$, see Fig.~\ref{pic:multicrossing_biquandloid} middle. The left and the right cubic maps produce the tangles (Fig.~\ref{pic:multicrossing_biquandloid} left and right), whose semiarc colors correspond to the terms of the differential in the biquandloid complex (Definition~\ref{def:biquandloid_homology}).

3. Let $(C,A)$ be an unframed crossoid.
Let $\phi_\chi\colon\widetilde{\mathdutchcal{SA}}_0(T)\to B$ be the homomorphism from the topological crossoid that corresponds to a coloring $\chi\in Col_{(C,A)}(D)$ of the diagram $D$. The maps $\phi\colon C_n(T,H)\to C_n(C,H)$,
\[
\gamma\mapsto (c_1,\dots,c_{n-1}),
\]
where $a_i=\phi_\chi\circ(\partial_{3,r})^{n-i-1}(\partial_{1,r})^{i-1}(\gamma)$,  $i=1,\dots,n-1$, define a chain map between the complexes.

\begin{figure}
    \centering
    \includegraphics[width=0.9\textwidth]{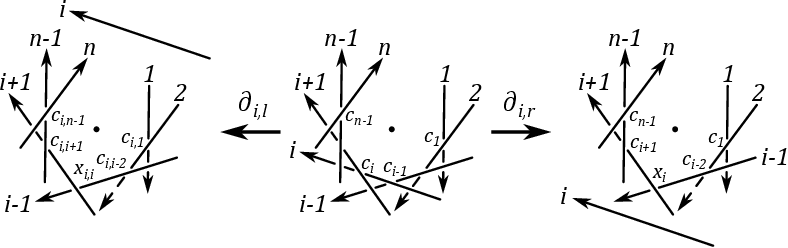}
    \caption{Geometric interpretation of the differential for a crossoid}
    \label{pic:multicrossing_crossoid}
\end{figure}

We assign to the multicrossing in (Fig.~\ref{pic:multicrossing_resolution}) a sequence of $n-1$ positive crossing colors $c_1,\dots, c_n\in C_+$ such that $A_{ur}(a_i)=A_{dr}(a_{i+1})$, $1,\dots,n-2$ as shown in Fig.~\ref{pic:multicrossing_crossoid} middle.
The geometric interpretation of this homomorphism above comes from the colorings of the tangle which appear after applying the left and right cubic maps (Fig.~\ref{pic:multicrossing_crossoid} left and right). The obtained sequences of crossing colors correspond to the terms of the differential in the crossoid complex (Definition~\ref{def:crossoid_homology}).

{\color{red}
}


\section{Invariants}\label{sect:invariants}

After obtaining a topological interpretation of the diagram elements, the next step is to "topologize" the known combinatorial invariants of knots. This task will be left for future work. In this section, we will briefly mention some invariants that use elements of knot diagrams in their construction.


\subsection{Diagram sum invariants}

\begin{definition}\label{def:diagram_chain}
Let $\mathcal F\colon \mathfrak D_s\to Rel$ be a functor and $(\mathcal G, h)$ its invariant. For a diagram $D\in Ob(\mathfrak D_s)$ the multiset $h_D(\mathcal F(D))\in\Z[\mathcal G(D)]$ is called the \emph{diagram chain} of the functor $\mathcal F$ with values in $\mathcal G$. If $\mathcal G(D)$ is an abelian group, the image of $h_D(\mathcal F(D))$ under the homomorphism $\Z[\mathcal G(D)]\to \mathcal G(D)$ is the \emph{diagram-sum} of $\mathcal F$ with values in $\mathcal G$.

If the diagram sum does not depend on the diagram, it is the \emph{diagram-sum invariant}.
\end{definition}

\begin{example}\label{ex:diagram_sum_invariant}
1. For the crossing functor, diagram-sum invariants come from the crossing class $cr(D)$. Examples are the linking coefficient, index polynomials and quandle cycle invariants.

2. For the arc and region functors, all diagram-sum invariants are trivial. Indeed, by applying a first Reidemeister move, we can duplicate any arc (any region). Invariance of the diagram-sum implies that the invariant value of this arc (region) is zero. Hence, the invariant is zero.

3. For the semiarc functor, a nontrivial example of invariant diagram-sum is following. Let $L$ be a link and $K\subset L$ one of its components. Consider the invariant with constant coefficients $\mathcal G=\Z_2$ that assigns $1$ to the semiarcs which belong to the component $K$, and $0$ to the rest of the semiarcs. In this case, the diagram-sum is invariant and equal to 1 when the component $K$ is odd (passes through odd number of classical crossings), and equal to 0 when $K$ is even.

\end{example}

\subsection{Global sum invariants}

Let $\mathcal F\colon\mathfrak D_s\to Rel$ be a functor and $(\mathcal G,h)$ its invariant. Let $\mu_{\mathcal G}$ be a measure on the set $\mathcal G(D)$. (For example, one can take a measure $\mu$ on the set of isotopy classes of tangle diagrams and consider the measure on $\mathcal G(D)$ induced by the diagram chain.) Let $f\colon \mathcal G\to\R$ be a function, then $I(D)=\int_{\mathcal G(D)}f(x)d\mu_{\mathcal G}$ is called a \emph{global sum invariant}.

\begin{example}\label{ex:global_sub_invariant}
Consider a classical knot $K$. Let $\mathcal G=\mathscr{MC}^w(K)$ be the universal midcrossing invariant and $\mu$ the discrete uniform measure on it. Consider the function $f\colon \mathscr{MC}^w(K)\to\Z$ such that $f(\gamma)=1$ if $\gamma$ is a nugatory midcrossing and $0$ otherwise. The corresponding invariant $I(K)=\int_{\mathscr{MC}^w(K)}f(\gamma)d\mu$ is equal to the number of nugatory midcrossings up to weak equivalence. If $K=\#_{i=1}^p n_iK_i$ is the prime decomposition of the knot $K$, then $I(K)=\prod_{i=1}^p(n_i+1)$, since the equivalence class of a nugatory midcrossing is determined by the isotopy classes of the components of the oriented smoothing at the midcrossing.
\end{example}

The example shows that using the global sums allows us to take into account the hidden elements of the tangle.

\subsection{Skein invariants}

\begin{definition}\label{def:skein_invariant}
Let $\mathcal F\colon\mathfrak D_s\to Rel$ be a functor and $(\mathcal G,h)$ its invariant. Assume that there is a skein-relation map $sk\colon \mathcal G(D)\to R[\mathfrak T]$ where $R$ is a ring. Denote $\mathcal S= R[\mathfrak T]/\mathrm{Im} sk$. The natural projection $\mathfrak T\to \mathcal S$ is called a \emph{skein invariant}~\cite{Nskein}.
\end{definition}

Examples of skein invariants are polynomial invariants (Alexander, Jones, HOMFLY).

\begin{remark}
    Using the (semi)arc functor, one can construct skein-invariants of order $1$ (by replacing a semiarc with a linear combination of $1$-tangles). Analogously, the (mid)crossing functor yields skein-invariants of order $2$ (including the polynomial invariants). To get skein-invariants of higher orders, one needs to use multicrossing functors on the tangle diagram category.
\end{remark}










\section{Further directions}\label{sect:open_question}

Let us outline possible directions for the future development of the framework presented in this paper.


\begin{itemize}
\item Extension to other knot theories
\begin{itemize}
\item Higher dimension knots. The extension of the approach described in the article to $2$-knots looks straightforward. In a broken surface diagram of a $2$-knot, one can distinguish such elements as sheets, double lines, triple points and cusps. Then one defines the elements of the $2$-knot as the isotopy classes of the corresponding probes.
\item Knots in $3$-manifolds. Given a knot in a $3$-manifold, its arcs, regions, and crossings can be defined as isotopy classes of probes. The ends of the probes should lie on a substrate whose role is naturally assigned to a spine of the manifold.
\item Virtual knots. The initial motivation for this study was to describe the parities of virtual knots. According to one definition, virtual knots are knots in thickened surfaces considered up to stabilizations. Therefore, it is necessary to combine the topological description of diagram elements with the stabilization moves. An example of the biquandle of the unknot (which is the free biquandle with one generator) shows that the structure of diagram elements (for example, semiarcs) of a virtual knot can be quite complex.
\end{itemize}

\item Topologization of combinatorial invariants. The topological description of diagram elements of knots allows us to reformulate combinatorial invariants such as the quandle cycle invariant and Khovanov homology in topological terms and thus gain a new perspective on these invariants and the possibilities of their application. 

\item Monodromy groups. They are responsible for the difference between the tangle invariants and the tangle coinvariants. The task is to describe monodromy at three different levels: 
\begin{enumerate}
    \item the motion groups of tangles;
    \item the monodromy groups of the diagram elements (arcs, regions, crossings);
    \item the monodromy of colorings (by tribrackets, biquandles and crossoids).
\end{enumerate}
\item Finite type invariants of diagram elements. In this paper, we consider homotopy classes of diagram elements and relate them to colorings. These homotopy classes can be seen as the lowest level of a hierarchy that leads to isotopic classes. Each step in this hierarchy corresponds to a finite-order invariant in the sense of Vasiliev. The next step is to define first-order invariants, such as first-order quandles.
\item Cobordisms. The term "functoriality" first appeared in the context of Khovanov homology, where a TQFT acted as a category, with its morphisms being cobordisms. One can check how the functors of diagram elements and their coinvariants behave when cobordisms are added to the diagram category.

\item Non-Reidemeister knot theories. In the Reidemeister approach, diagrams and their elements arise as a result of stratification of the projection of a knot onto the plane. An interesting task is to adapt the framework discussed in the article to non-Reidemeister approaches to knot theory, such as the one-parameter approach of T. Fiedler and V. Kurlin~\cite{FiedKur} or the theory of groups $G^k_n$ by V.O.Manturov~\cite{MN}.

\item Coloring propagation rules. When considering colorings of knot diagrams, the choice of propagation rule determines the algebraic structure of the set of colors. Naive distribution rules, as discussed in Section~\ref{sect:colorings}, lead to Alexander numbering. Other rules lead to known algebraic structures, such as quandles, biquandles, and ternary quasigroups. The question arises: what other propagation rules are there and what algebraic structures do they correspond to?

\item Tuples of diagram elements. In this article, we consider the topological description of a single diagram element. However, there are times when it is necessary to simultaneously consider several diagram elements. Examples are ribbon singularities of a ribbon knot which are described by a tuple of midcrossings, and Gauss diagrams of knots which are set by a tuple of traits which can intersect but the order of the ends of the traits in the knot is fixed.

\item Universal cycle formulas in the multicrossing complex. Are there universal formulas for cycles in the multicrossing complex, other than the crossing cycle?


\end{itemize}

{\color{red}
%
%
%
%
%
%
%
%
}

\end{document}